\documentclass[11pt]{amsart}

\renewcommand{\familydefault}{\sfdefault}
\newcommand{\dateversion}{Ismar27July2012}% change also the name of
\newcommand{\version}{\dateversion.
This is the last private version. Set \texttt{setboolean\{privateversion\}}\
to \texttt{false} to get the public version.\\
Continued by Ismar from Pascal's revisions of 24 July 2012.\\
pl in the margin means a change done by Pascal after Ismar's
23July2012 revision.
There are still a very few footnotes to resolve.
The \texttt{$\backslash$p}, \texttt{$\backslash$pfn}  and
\texttt{$\backslash$q} became invisible.
Instead the few new changes use \texttt{$\backslash$pp} and
\texttt{$\backslash$ppfn}.\\
Ismar's comments and footnote of July 23  are now introduced by  $\backslash$ivfn
}

\usepackage{ifthen}\newboolean{privateversion}
\setboolean{privateversion}{false}% to be set to "false" to publish a PUBLIC version
\newcommand{\previousfn}[1]{}%{\private{\marginpar{\dag}\footnote{{#1}}}}
\newcommand{\ivfn}[1]{\private{\footnote{\iv {#1}}}}
\newcommand{\pp}{\private{\textbf{\copyright}\marginpar{pl}{{}}}}
\newcommand{\ppfn}[1]{\private{\p\footnote{\pl 24july'12:{#1}}}}
\newcommand{\p}{}%{\private{\textbf{\copyright}\marginpar{pl}{{}}}}
\newcommand{\pfn}[1]{}%{\private{\p\footnote{\pl:{#1}}}}
\newcommand{\q}{}%{\private{\marginpar{iv}{{}}}}
\newcommand{\chaug}{\previousoldfn{\pl {\large\textbf{the following was changed in August 2010}}}}
\newcommand{\D}{}%{\private{{\textbf{\Large{\copyright}}}   }}
\newcommand{\Z}{}%{\private{{\textbf{\Large{\copyright}}}   }}%\Z INDIQUE
\newcommand{\previousoldfn}[1]{}%{\private{\marginpar{\dag}\endnote{{#1}}}} 
\newcommand{\optproof}[2][this]{\private{(a private detailed proof
is at the following endnote\endnote{\textbf{\Large Detailed proof of {#1}:}\\
{#2}
\textbf{\Large End of detailed proof of {#1}}})}}
\newcommand{\pl}{{\bf Pascal:} }
\newcommand{\iv}{{\bf Ismar: }}

\newcommand{\todo}{{\Large\bf TO DO! }}

\usepackage{amssymb, amsmath, amscd, amsthm, stmaryrd, color, epsfig, pdflscape, longtable}
\usepackage{xspace,url}
\usepackage[all]{xy}          % for xy-pic pictures
\xyoption{dvips}              % for xy-pic pictures
\usepackage{graphicx}
\usepackage{accents}% for \undertilde
\usepackage{esint} % pour la commande \fint pour fiber-integral
\usepackage{endnotes}%see "Latex, apprentissage, guide er reference" 2eme edition page 88-90

\reversemarginpar
\ifthenelse{\boolean{privateversion}}{% 
        \usepackage[notref,notcite]{showkeys}% shows keys of labels and bibliographies 
        \oddsidemargin=0.5in%
        \evensidemargin=0.5in% leave enough space to show the keys 
        \newcommand{\private}[1]{{#1}}     % text appearing onlyin private versions
}{     % FOR A PUBLIC VERSION 
        \newcommand{\private}[1]{}
        \oddsidemargin=0in
        \evensidemargin=0in
}
\newcommand{\dontcompile}[1]{}

\makeatletter
\renewcommand\l@subsection{\@tocline{2}{0pt}{2pc}{5pc}{}}
\renewcommand\l@subsubsection{\@tocline{3}{0pt}{4pc}{10pc}{}}
\makeatother
\newcommand{\BN}{{\mathbb N}}  % natural  integers
\newcommand{\BZ}{{\mathbb Z}}  % ring of integers
\newcommand{\BR}{{\mathbb R}}  % field of real numbers
\newcommand{\BQ}{{\mathbb Q}}  % field of rational numbers
\newcommand{\BC}{{\mathbb C}}  %  field of complex numbers
\newcommand{\BK}{{\mathbb K}}  % a genreic field
\newcommand{\unit}{{\mathbf{1}}}
\newcommand{\calB}{{\mathcal{B}}}
\newcommand{\calT}{{\mathcal{T}}}

\newcommand{\calA}{{\mathcal{A}}}
\newcommand{\calV}{\mathcal{V}}
\newcommand{\calC}{\mathcal{C}}
\newcommand{\BF}{\mathcal{BF}}
\newcommand{\calU}{\mathcal{U}}
\newcommand{\calO}{\mathcal{O}}

\newcommand{\GD}{\widehat{\mathcal{D}}}% (generalized) diagrams
\newcommand{\NAI}{{\mathcal{N}}}% ideal of non admissible diagrams
\newcommand{\AD}{{\mathcal{D}}}%admissible diagrams
\newcommand{\cAD}{{\undertilde{\mathcal{D}}}}%connected admissible diagrams
\newcommand{\IK}{{\operatorname{I}}}%Kontsevich integral
\newcommand{\bIK}{{\overline{\operatorname{I}}}}%Kontsevich integral in homology
\newcommand{\GIK}{\widehat{\operatorname{I}}}%generalized Kontsevich integral
\newcommand{\dvol}{\operatorname{vol}}

\newcommand{\ExtVertbar}{\overline{\ExtVert}}
\newcommand{\SemiAlg}{\operatorname{SemiAlg}}
\newcommand{\CompSemiAlg}{\operatorname{CompactSemiAlg}}
\newcommand{\Ebar}{\overline{E}}\newcommand{\CSA}{\operatorname{C}}
\newcommand{\Ibar}{\overline{I}}
\newcommand{\thetabar}{\overline{\theta}}
\newcommand{\mubar}{\overline{\mu}}
\newcommand{\pibar}{\overline{\pi}}
\newcommand{\tbar}{\overline{t}}
\newcommand{\sbar}{\overline{s}}
\newcommand{\ExtVert}{A}
\newcommand{\Gbar}{\overline{\Gamma}}
\newcommand{\oprd}{}%{(\bullet)} %notation for an operad: too heavy
\newcommand{\omin}{\Omega_{\min}}
\newcommand{\ompa}{\Omega_{PA}}\newcommand{\omDR}{\Omega_{DR}}

\newcommand{\bbr}[1]{\llbracket{#1}\rrbracket}%{{\mathbf{[}}{#1}{\mathbf{]}}}
\newcommand{\IsoFin}{\operatorname{IsoFin}}
\newcommand{\Top}{\operatorname{Top}}  % category of topological spaces
\newcommand{\Tree}{\operatorname{Tree}}
\newcommand{\CDGA}{\operatorname{CDGA}}
\newcommand{\htree}{\operatorname{height}}
\newcommand{\pos}{{\operatorname{pos}}}  % position in a linear order
\newcommand{\ChZ}{{\operatorname{Ch}_\BZ}} % categorie des complexes de chaines sur Z
\newcommand{\ChK}{{\operatorname{Ch}_\BK}} % categorie des complexes de chaines sur K
\newcommand{\Econtr}{E^{\textrm{\tiny{contr}}}}
\newcommand{\Ho}{{\operatorname{H}}}  % homology
\newcommand{\Sing}{{\operatorname{S}}}  %singular chains
\newcommand{\Hom}{\operatorname{Hom}}
\newcommand{\Apl}{A_{PL}}
\newcommand{\Ball}{\operatorname{B}}
\newcommand{\proj}{\mathrm{proj}}

\newcommand{\id}{\mathrm{id}}\newcommand{\ev}{\mathrm{ev}}

\newcommand{\setn}[1]{\underline{{#1}}}%\setn{k} is the notation for the set {1,...,k}
\newcommand{\Cond}{{\operatorname{Cond}}}
\newcommand{\AdmCond}{{\operatorname{AdmCond}}}
\newcommand{\EssCond}{{\operatorname{EssCond}}}
\newcommand{\rel}{\operatorname{rel}}
\newcommand{\im}{\operatorname{im}}
\newcommand{\fiber}{\operatorname{fiber}}
\newcommand{\sign}{\operatorname{sign}}
\newcommand{\Inj}{\operatorname{Inj}}
\newcommand{\theroot}{\operatorname{root}}
\newcommand{\predec}{\operatorname{pred}}
\newcommand{\outputset}{\operatorname{output}}
\newcommand{\Emb}{\operatorname{Emb}}
\newcommand{\Conf}{\operatorname{C}}%espace de configuration
\newcommand{\Confsing}{\operatorname{C}^{\operatorname{sing}}}%espace de configuration
\newcommand{\barEmb}{{\overline{\Emb}}}

\newcommand{\Perm}{\operatorname{Perm}}
\newcommand{\ordsum}{\olessthan}
\newcommand{\bary}{\operatorname{barycenter}}
\newcommand{\radius}{\operatorname{radius}}

\newcommand{\quism}{\stackrel{\simeq}{\longrightarrow}}
\newcommand{\iso}{\stackrel{\cong}{\longrightarrow}}

\newcommand{\current}{chain}% pour remplacer SA currents by SA chain
\newcommand{\citePADC}{\cite[Definition 3.1]{HLTV:RHTSAS}}
\newcommand{\citePASMinimalForms}{\cite[Section 5.2]{HLTV:RHTSAS}}
\newcommand{\citePASSPAForms}{\cite[Section 5.4]{HLTV:RHTSAS}}
\newcommand{\citePADPAForm}{\cite[Definition 5.20]{HLTV:RHTSAS}}
\newcommand{\citePASSAbundle}{\cite[Section 8]{HLTV:RHTSAS}}
\newcommand{\citePASpropintfiber}{\cite[Section 8.2]{HLTV:RHTSAS}}
\newcommand{\citePADSAbundle}{\cite[Definition 8.1]{HLTV:RHTSAS}}
\newcommand{\citePADintegrfiber}{\cite[Definition 8.3]{HLTV:RHTSAS}}
\newcommand{\citePAPnatpush}{\cite[Proposition 8.10]{HLTV:RHTSAS}}
\newcommand{\citePAPaddpush}{\cite[Proposition 8.11]{HLTV:RHTSAS}}
\newcommand{\citePAPPAdegdimz}{\cite[Proposition 5.24]{HLTV:RHTSAS}}
\newcommand{\citePAPpushfactcodim}{\cite[Proposition 8.14]{HLTV:RHTSAS}}
\newcommand{\citePAPdpush}{\cite[Proposition  8.12]{HLTV:RHTSAS}}%fiberwise stokes
\newcommand{\citePAPdoublepush}{\cite[Proposition 8.13]{HLTV:RHTSAS}}
\newcommand{\citePAPmonompaapl}{\cite[Theorem 7.1]{HLTV:RHTSAS}}
\newcommand{\citePAPcompositebdl}{\cite[Proposition 8.5]{HLTV:RHTSAS}}
\newcommand{\citePAPevalmostmon}{\cite[Proposition 7.3]{HLTV:RHTSAS}}
\newcommand{\citePAPprodpushPB}{\cite[Proposition 8.15]{HLTV:RHTSAS}}
\newcommand{\citePAPmonCHSing}{\cite[Proposition 7.2]{HLTV:RHTSAS}}
\newcommand{\citePAPpbSAbdl}{\cite[Proposition 8.4]{HLTV:RHTSAS}}
\newcommand{\citePAPpushpullback}{\cite[Proposition 8.9]{HLTV:RHTSAS}}
\theoremstyle{plain}
\newtheorem{thm}{Theorem}[section]
\newtheorem{prop}[thm]{Proposition}
\newtheorem{lemma}[thm]{Lemma}
\newtheorem{cor}[thm]{Corollary}

\newtheorem{setting}[thm]{Setting}

\theoremstyle{definition}
\newtheorem{defin}[thm]{Definition}
\newtheorem{eg}[thm]{Example}
\newtheorem{remark}[thm]{Remark}
\theoremstyle{remark}

\newcommand{\refS}[1]{Section~\ref{sec:#1}}
\newcommand{\refT}[1]{Theorem~\ref{T:#1}}
\newcommand{\refC}[1]{Corollary~\ref{C:#1}}
\newcommand{\refP}[1]{Proposition~\ref{P:#1}}
\newcommand{\refR}[1]{Remark~\ref{R:#1}}
\newcommand{\refD}[1]{Definition~\ref{D:#1}}
\newcommand{\refL}[1]{Lemma~\ref{L:#1}}
\newcommand{\refX}[1]{Example~\ref{X:#1}}
\newcommand{\refF}[1]{Figure~\ref{fig:#1}}
\newcommand{\refE}[1]{Equation~$(\ref{eq:#1})$}
\newcommand{\refN}[1]{$(\ref{#1})$}
\begin{document}
%\dontcompile{
%% =========== TITLE ============
\ifthenelse{\boolean{privateversion}}{
%  FOR A PRIVATE VERSION
\title[\dateversion]{Formality of the little $N$-disks operad}
}{%
% FOR A PUBLIC VERSION
\title{Formality of the little $N$-disks operad}
}
%%
%%=========== AUTHORS ============
\author{Pascal Lambrechts}
\address{Universit\'e catholique de Louvain, IRMP, 2 Chemin du Cyclotron, B-1348 Louvain-la-Neuve, Belgium}
\email{pascal.lambrechts@uclouvain.be}
\urladdr{http://milnor.math.ucl.ac.be/plwiki}
\author{Ismar Voli\'c}
\address{Department of Mathematics, Wellesley College, Wellesley, MA 02482}
\email{ivolic@wellesley.edu}
\urladdr{http://palmer.wellesley.edu/~ivolic}

% -----------------------------------------------------------------
%%
%%============ CLASSIFICATION AND KEYWORDS


\thanks{The first author is Ma\^\i tre de Recherches au F.R.S-FNRS.  
The second author was supported in part by the National Science Foundation grants DMS 0504390 and DMS 1205786.}

%% =========== ABSTRACT ============
\private{\version}

\begin{abstract}
%\private{\Huge CONFIDENTIAL\\}
%\private{See the boolean "privateversion" at the beginning of the latex file.\\}
%\private{\privateabstract}
We develop the details of Kontsevich's proof of the formality  of
little $N$-disks operad over the field of real numbers. 
Formality holds in the category of  operads of chain complexes
and also in some sense in the category of commutative differential graded algebras,
which is the  category encoding rational homotopy theory (or 
 ``real'' homotopy  theory in this case).  We also prove a relative version of the formality for the inclusion  of the little $m$-disks operad in the little $N$-disks operad when $N\geq2m+1$. 
\end{abstract}

%%
%% ================== BODY OF THE PAPER =======================
%\dontcompile{
%}%dontcompile

\maketitle

\pfn{in the title I suggest \texttt{Rational formality} better than \texttt{The
    formality}.  \iv I added ``The".} 
\ppfn{I had second thoughts about the title. Our preprint is known on
  arXiv  as 
``Formality of the little N-disks operads `` since 2008 (!)''. Moreover
many paper are entitled ``Formality ...'' without an article ``The'': see the long
list below. Therefore, despite referee comment (1), I suggest to keep
the original title without ``The'' and ``Rational''. Among the
advantage there is the fact that robots making citations, like on
MathReview, will be able to connect the arXiv preprint with the
published memoirs.
\\Here is the list from MathSciNet:
\begin{itemize}
\item
MR2946936 Prelim Giansiracusa, Jeffrey; Salvatore, Paolo; Cyclic operad formality for compactified moduli spaces of genus zero surfaces
\item
MR2836399 Pending Willwacher, Thomas Formality of cyclic chains.
\item
MR2821234 Reviewed Dolgushev, Vasily Formality theorem for Hochschild
cochains via transfer. 
\item
MR2823091 Pending Miller, Matthew; Wakefield, Max Formality of Pascal
arrangements
\item
MR2728693 Reviewed Zuber, Hugues Non-formality of Milnor fibres of
line arrangements.
\item
MR2661521 Reviewed Severa, Pavol Formality of the chain operad of
framed little disks.
\item 
MR2648708 Reviewed Giansiracusa, Jeffrey; Salvatore, Paolo Formality
of the framed little 2-discs operad and semidirect products.
\item
MR2668829 Reviewed Calaque, Damien; Van den Bergh, Michel Global
formality at the G-level.
\item etc...
\end{itemize}
\iv Totally fine with me.
}
\ifthenelse{\boolean{privateversion}}%
{\textbf{TABLE OF CONTENTS} does not work in the private version but
  will appear here in the public version. \iv If we don't want
  subsubsections to be displayed in the table of contents, the
  commands about secdepth and tocdepth should be commented
  out. \pl: commented out because it is better not to have these subsubsection in  TOC
since I turned these subsubtitle into textbf }%
{\tableofcontents}

\parskip=5pt
\parindent=0cm

% \cm{\pl 08july09. There is some confusion btween the maps $\Psi$ and $\Phi$ in the text. $\Psi$ should be the operad structure  maps for $\Conf[\bullet]$
% and $\Phi$ is the cooperad structure maps for $\GD$. This has to be corrected}

% \cm{\pl 08july09. Greg suggested to add up the equivalence between FultonMacPherson operad and little disks}
%\newpage

\section{Introduction}
\label{sec:intro}
\previousoldfn{\iv We should have table of contents.\pl There is a technical pbm to insert the tableofcontents. Anyway there will be a roadmap at the end of intro
\pl Now there is a roadmap at the end of the intro; I guess that it is better than a toc. \iv I think having toc is crucial, in addition to a roadmap.  Such a long paper is really hard to navigate.  I have a hard time finding things easily, let alone a reader. \pl OK the toc is added; the problem was only in the private version}
\previousfn{\iv Maybe we should make a ``clean" tex file for arXiv and for submission where all this public/private crap has been taken out?\pl OK I'll do that before submitting and arXiving. \iv Or maybe we can wait until we have to, if we have to.}

In this paper we give a detailed proof of Kontsevich's theorem on the formality of the
little $N$-disks operad.  The theorem, whose proof was sketched in \cite[Theorem 2]{Kon:OMDQ},  asserts that
the  singular chains on the little $N$-disks operads is
weakly \D equivalent to its homology in the category of operads of\p\ chain complexes. We  also improve that result in
three directions:
\begin{enumerate}
\item  Formality is in the category of CDGA (commutative 
  differential graded \p algebras) which, following
  Sullivan and Quillen, models rational homotopy theory; 
\item  For us, the little disks operad has an operation in arity $0$ while Kontsevich discards that nullary operation;
\item We establish a \emph{relative} formality result, namely formality of the
  inclusion of the little $m$-disks operad into the little $N$-disks
  operad for $N\geq 2m+1$.
\end{enumerate}

Our motivation  for proving these results comes from  applications 
to the study of the rational homology of the space $\Emb(M,\BR^N)$ of smooth embeddings
of a  compact manifold $M$ into $\BR^N$. Goodwillie-Weiss manifold calculus \cite{Wei:EI1,GoWe:EI2}\previousoldfn{\iv There should also be a reference to Weiss' Part I here.\pl done}\previousoldfn{\iv  Some other bibliography notes are in the last footnote. \pl Any to add?}  approximates
 this embedding space 
by homotopical constructions based on a category $\calO_\infty$
of open subsets of $M$ diffeomorphic to finitely many open balls with inclusions as morphisms.  This category is closely related to the little 
balls operad.
On the other hand, formality theorems can 
often lead to collapse  results for spectral
sequences.  Combining manifold calculus with formality, the authors, along with Greg Arone, were thus able to prove in \cite{ALV:HQE} the collapse of a spectral sequence
computing $\Ho_*(\barEmb(M,\BR^N);\BQ)$, where ${\barEmb(M,\BR^N)}$ is a slight variation
of ${\Emb(M,\BR^N)}$.
 A special case of this approach also led the authors, jointly with Victor Turchin, to the proof in \cite{LTV:HQLK} of the
collapse of the Vassiliev spectral sequence computing the rational homology of
the space of long knots in $\BR^N$ for $N\geq4$.

To explain the formality results that we prove here,  fix an integer $N\geq1$ and 
recall the  classical  little $N$-disks operad $\calB_N\oprd=\{\calB_N(n)\}_{n\geq0}$, where 
$\calB_N(n)$ is the space of configurations of $n$ closed $N$-disks with disjoint interiors
contained in the unit disk  of $\BR^N$  \cite{BoVo:hom}\previousoldfn{precise the page in this reference, because it is a long book. \iv Actually, the first reference for little cubes is Boardman, J. M.; Vogt, R. M. Homotopy-everything $H$-spaces. 
Bull. Amer. Math. Soc. 74 1968 1117--1122.  The definition is not numbered, but it's in Section 2 of that paper.\pl I suggest to leave it like this. ok?}.
The integer $N$ will usually be understood so we will just denote this operad by $\calB$ and often simply say ``little balls\previousoldfn{\pl balls or disk? should remains disks in the title since it was title so in arXiv and the paper is cited} operad". 
This operad is homotopy equivalent to many other operads, such as the little $N$-cubes operad,
or the Fulton-MacPherson operad $\Conf[\bullet]=\{\Conf[n]\}_{n\geq0}$ of compactified configurations of points in $\BR^N$.
The latter will be important in our proofs and  we will say more about it in \refS{FMoperad}.

Fix a unital commutative ring $\BK$.
The functor 
\[\Sing_*(-;\BK)\colon\Top\longrightarrow\ChK
\]
of singular chains with coefficients in $\BK$ is symmetric monoidal.  
Therefore $\Sing_*(\calB;\BK)$ is an operad of chain complexes.  In addition, its homology
$\Ho_*(\calB;\BK)$ can be viewed as an operad of chain complexes with differential $0$.
One of the main results that we will prove in detail is 
\begin{thm}[Kontsevich \cite{Kon:OMDQ}; Tamarkin for $N=2$ \cite{Tam:for}]\label{T:stableformality}
The   little $N$-disks operad is stably formal over the real numbers, that is, there exists a chain of weak equivalences of operads of chain complexes
\[
\Sing_*(\calB_N;\BR)\stackrel{\simeq}\longleftarrow\cdots\quism\Ho_*(\calB_N;\BR).
\]
\end{thm}
The proof of this theorem was sketched in \cite[Section 3.3]{Kon:OMDQ} but we felt 
that it would be useful to develop it in full detail.  
%This\previousfn{\pl this paragraph should be rephrased once $N=2$ case is clear}
%proof seems to
%break down\previousfn{\pl That is a strong assertion! We should better be sure that Kontseich approach fails for $N=2$; otherwise we
%should just say that Tamarkin has another proof in $N=2$}
% for $N=2$, but in that dimension the formality has been proved by Tamarkin \cite{Tam:for} using a different approach.\previousfn{The following was here: 
%``Actually it seems that minor corrections'' Waat is this?}
In this paper, $\calB(0)$ is the one-point space, contrary to \cite{Kon:OMDQ} where it is the empty set.\previousoldfn{what about tamrkin?.}
This fact makes our \D proof more delicate, but in the application we have in mind it will be important to
have $\calB(0)=*$ (operad composition with this corresponds to the operation of forgetting a ball from a configuration of little balls).
\previousoldfn{citer travail Roig and co on their proof of formality}

 Morally, \D singular chains with coefficients in $\BQ$
encode the rational stable  homotopy type of spaces or topological operads,
and with coefficients in $\BR$ we get the ``real stable homotopy
type''.
This is why in \refT{stableformality} we \D talk about \emph{stable} formality.
The unstable  real (or more correctly, rational) homotopy type of spaces is encoded by
commutative differential graded algebras (CDGAs for short), as was discovered by Sullivan
using the functor $\Apl$ of polynomial forms (see \refS{CDGAmodel}).
One then has the important notion of
a \emph{CDGA model} for a space $X$, which by definition 
is a CDGA weakly equivalent to $\Apl(X)$. Any CDGA model (over the
field $\BQ$) for a simply-connected space with finite Betti numbers\previousoldfn{\pl Global
  comment: sometimes we say ``model of'', sometimes ``model
  for''. Shoult it be uniformized?}
contains all the information about its rational homotopy type.
We can define an analogous notion of a CDGA model for a topological operad, 
although the definition
is a little bit more intricate (see \refD{CDGAmodeloperad}).
We then have the following unstable version of \refT{stableformality}.
\begin{thm}\label{T:unstableformality}
For \D $N\not=2$, a CDGA model over $\BR$ of the little $N$-disks operad
is given by its cohomology algebra, that is, it is formal \D over $\BR$ (in the sense of \refD{CDGAmodeloperad}). \q
\end{thm}
%We do not know whether the result is true for $N=2$. There is only one place in the proof where the hypothesis
%$N\not=2$ is needed (see Remark \ref{R:N=2break}).\previousfn{\pl to restate once  $N=2$ settled}

As explained in \refS{CDGAmodel}, one reason for which our definition of a CDGA model for an operad is not
as direct as one might wish is that  $\Apl(\calB)$ is not a cooperad.  This is because 
the contravariant functor $\Apl$ is not  comonoidal.
It might be better to consider the coalgebra of singular chains $\Sing_*(\calB;\BR)$,
 which is indeed an operad of differential coalgebras. However, we do not know how to prove
that this operad is weakly equivalent to its homology in the category of differential coalgebras. 
Moreover, that category is 
 not  very suitable  
for doing real homotopy theory because of the lack of strict cocommutativity.

In \refT{unstableformality}, \D we assumed $N\not=2$. 
Our proof in the case $N=2$ fails because some of
our CDGAs become $\BZ$-graded instead \D of non-negatively graded as
required in rational homotopy theory.
 We still however obtain some results in the case $N=2$ \p \q and we
 believe that our proof can be adapted to include that case as well; see \refS{proofform}.\previousoldfn{\pl check if done in
  that section.\pl ok}

\previousoldfn{\iv Should we address the formality of configuration spaces,
  namely say that we prove they're formal in the course of our proofs?
  Paolo said to me last month in Tunisia how configuration spaces
  aren't formal, but I didn't get a chance to ask him what he
  meant. \pl See a further note}

We now state a relative version of the above theorems. 
Let $1\leq m\leq N$ be integers and suppose given a linear isometry
\[\epsilon\colon\BR^m\longrightarrow\BR^N.
\]
Define  the map
\[
\calB_\epsilon(n)\colon\calB_m(n)\longrightarrow\calB_N(n)
\]
that sends a configuration of $n$ $m$-disks to the configuration of $n$ $N$-disks 
where the center of each $N$-disk is the image under $\epsilon$ 
of the center of the corresponding $m$-disk and has the same radius.
This clearly defines a morphism of operads.

\begin{defin}\label{D:StableFormality}
A morphism of  topological operads
\[\alpha\oprd\colon\calA\oprd\longrightarrow\calA'\oprd
\]
is \emph{stably formal} over \D $\BK$  if there exists a zig-zag of
 quasi-isomorphisms of operads in $\ChK$ connecting the singular chains
$\Sing_*(\alpha\oprd;\BK)$ to its homology
$\Ho_*(\alpha\oprd;\BK)$
as in the following diagram:
\[\xymatrix{
\Sing_*(\calA\oprd;\BK)
\ar[d]_{\Sing_*(\alpha)}
&
\ar[l]_-{\simeq}
\calC_1\oprd
\ar[r]^-{\simeq}
\ar[d]
&
\cdots
&
\ar[l]_-{\simeq}
\calC_k\oprd
\ar[r]^-{\simeq}
\ar[d]
&
\Ho_*(\calA\oprd;\BK)
\ar[d]^{\Ho_*(\alpha)}
\\
\Sing_*(\calA'\oprd;\BK)
&
\ar[l]_-{\simeq}
\calC'_1\oprd
\ar[r]^-{\simeq}
&
\cdots
&
\ar[l]_-{\simeq}
\calC'_k\oprd
\ar[r]^-{\simeq}
&
\Ho_*(\calA'\oprd;\BK)
}
\]

When $\BK$ is a field of characteristic $0$,
we say that $\alpha\oprd$ is \emph{formal} over $\BK$
if the morphism of CDGA cooperads $\Ho^*(\alpha;\BK)$ is a model for
$\alpha$ (see \refS{CDGAmodel} for the precise definition of a model for CDGA cooperads).

\end{defin}

\begin{thm}\label{T:relfor}
Assume that $m\geq1$ and $N\geq 2m+1$.
Then %\previousfn{\pl what about $m=2$? to be settled} 
the morphism  of operads
\[
\calB_\epsilon\colon\calB_m\longrightarrow\calB_N
\]
is stably formal over $\BR$. If $m\not=2$, it is \D also formal over $\BR$.
\end{thm}
\p \q 
% It is possible that his theorem does not hold when $N\leq2m$.

There is also a notion of \emph{coformality} which is Eckman-Hilton dual to that of
(unstable) formality \cite{NeMi:for}. \previousoldfn{\pl add reference?  \iv The one we used before was J. Neisendorfer and T. Miller, Formal and coformal spaces, Illinois J. Math. 22 (1978), no. 4, 565--580.  I didn't add this to the bibliography since you should probably do it in your .bib file directly.  I also checked that our citations from the Santos, Navarro, Pascual, and Roig paper were correct with respect to the version that appeared in Duke.\pl Neisendorfer-Miller added}
 Roughly speaking, coformality of a space $X$ means that its rational homotopy type
 is determined
by its rational homotopy Lie algebra $\pi_*(\Omega X)\otimes\BQ$ (instead of its 
rational cohomology algebra in the case of formality). 
In some sense, the operad of little $N$-disks also seems to be coformal, although there is difficulty
in making this idea precise because of the lack of a basepoint for the operad.
We refer the reader to \cite{ALTV:cof} for a discussion of coformality of the little
$N$-disks operad.

All of the above formality results are over the field of real numbers. It would be more
convenient to have rational formality because localization over $\BQ$ is topologically meaningful, contrary to
localization over $\BR$. This descent of fields for stable formality of operads is always possible when one considers
operads\previousoldfn{\pl ``reduced'' operad usually means that in arity $0$ we
  have a single operation, which is not what we mean here.So I
  replaced it in the text by paraphrases like ``with no nullary operations'' etc...}  in which the zeroth term (corresponding to $0$-ary operations)
is empty,  as proved in \cite[Theorem 6.2.1]{GNPR:mod}.
In particular, we can consider the operad $\widetilde{\calB}$ \D defined
by $\widetilde{\calB}(0)=\emptyset$ and $\widetilde{\calB}(n)=\calB(n)$ for $n\geq1$.
Our formality results for $\calB$ are clearly also true for
$\widetilde{\calB}$; the latter was the operad considered by Kontsevich in \cite{Kon:OMDQ}.
Moreover, since this operad has no nullary \D operations, stable
formality for $\widetilde{\calB}$ over $\BR$ descends to $\BQ$.

\q
For \pfn{I rephrased the next paragraph; check please; it correponds
  to referee comment (6). \iv Checked and slightly modified.} our applications to embedding spaces \cite{ALV:HQE, LTV:HQLK},
however, it is important to take the usual little balls operad, \p$\calB$,
which is only formal over $\BR$. In
those applications, this weaker  formality  is sufficient essentially because the main results there are about 
collapse of spectral sequences, and these collapse results do not depend on which field of
characteristic $0$ is used.
%Our results about embeddings spaces are over $\BR$ since in this case
%we do not know how to prove formality  over $\BQ$; 
The proof of descent of formality  in \cite[Section 6]{GNPR:mod}
 does not generalize easily to the case with nullary operations because of the lack
of minimal models when these degeneracy operations occur\p.

The formality of the operad $\widetilde\calB$ implies the formality over
$\BQ$ of each
space  $\calB(n)$, in the sense that the CDGA $\Apl(\calB(n))$ is
weakly equivalent to its cohomology algebra, $\Ho^*(\calB(n);\BQ)$. Paolo Salvatore has 
recently proved using a computer that, for $n=4$ and $N=2$, the space $\calB_2(4)$  is
not formal over the ring $\BZ/2$, i.e.~its cohomology
algebra, $\Ho^*(\calB_2(4);\BZ/2)$, and its algebra of singular cochains,
$\Sing^*(\calB_2(4);\BZ/2)$, are not quasi-isomorphic.  We do not know
whether the ({non-symmetric}) little disks  operad is stably formal
over some field of positive characteristic.\pp\ppfn{I added this paragraph. I just exchanged email
  with Paolo for that and remember his talk in Lille about it. \iv  Good, this makes sense.  I just changed a couple of minor things.}

As a final comment, the Tamarkin's and
Kontsevich's  proofs of formality for $N=2$ have been compared in \cite{SeWi:EFL} where it is
proved that the weak equivalences obtained in those two proofs are homotopic.\previousoldfn{ \iv Should we cite the Severa-Willwacher paper that proves that
  Kontsevich's and Tamarkin's formalities are homotopic for $N=2$?
  \textbf{pascal} Here it is}

We end this introduction by explaining the general idea of Kontsevich's proof of formality that we develop in this paper.
The main ingredient is a combinatorial CDGA
cooperad $\AD\oprd=\{\AD(n)\}_{n\geq0}$ of  \emph{admissible diagrams} and an
  explicit CDGA map 
\begin{equation}\label{eq:Iintro}
\IK\colon\AD(n)\longrightarrow\ompa (\Conf[n])
\end{equation}
which we will call the \emph{Kontsevich configuration space
  integral}\pfn{\pl is the ``'s'' ok after kontsevich? otherwise global
  change on ``Kontsevich's configuration space
  integral'' \iv I changed it to ``Kontsevich".  We also needed ``the"
  in front, which I added everywhere where it was missing.}. \q   Here
$\Conf[n]$ are compact manifolds homotopy equivalent to $\calB(n)$, and $\ompa $
is a semi-algebraic analog of the deRham CDGA of differential forms $\omDR $. 
A combinatorial argument will show that the cooperad $\AD$ is quasi-isomorphic 
to the cohomology of the little balls operad.
We will also show that $\IK$ is a quasi-isomorphism and, since $\IK$ also respects the cooperad structures,
 the desired result will follow.

Let us elaborate on  $\AD(n)$  and $\IK$ a bit further.
We will work with the Fulton-MacPherson operad $\Conf[\bullet]=\{\Conf[n]\}_{n\geq0}$ which is homotopy
equivalent to the little balls operad. The space $\Conf[n]$ is a compact manifold with corners obtained by adding a boundary to the open manifold $F_n(\BR^N)$, the space of configurations of $n$ points in
$\BR^N$, that is, 
\[F_n(\BR^N):=\{(z_1,\dots,z_n)\in(\BR^N)^n:z_i\not=z_j\text{ for }i\not=j\}\]
(after \D normalizing by modding out by translations and positive dilations).
 Arnold \cite{Arn:coh}  computed the cohomology algebra of $F_n(\BR^2)=F_n(\BC)$ and in fact
 proved that these spaces are formal over $\BC$.
His argument is as follows:

Consider the complex smooth differential one-forms\p \q
\begin{equation}\label{eq:omegaij}
\omega_{ij}:=\frac{d(z_j-z_i)}{z_j-z_i}=d\log(z_j-z_i)\in\omDR^1(F_n(\BC);\BC)
\end{equation}
which are cocycles and can easily be  shown to be
cohomologically independent for $1\leq i<j\leq n$.
%These forms are  the pullback of the volume form $\dvol$ on the circle
%through the map
% \begin{eqnarray*}
% \theta_{ij}\colon& F_n(\BC)&\to S^1\\
% &(z_1,\dots,z_n)&\mapsto\frac{z_j-z_i}{|z_j-z_i|}.
%\end{eqnarray*}
A direct computation shows that these forms satisfy the  \emph{$3$-term relation}
\begin{equation}\label{eq:arnold}
\omega_{ij}\wedge\omega_{jk}+\omega_{jk}\wedge\omega_{ki}+\omega_{ki}\wedge\omega_{ij}=0.
\end{equation}
It is  convenient to represent this relation by the diagram
pictured
\pfn{the order of the 3 terms of the diagram does not fit with
  the above arnold relation :( \iv Fixed.} 
  \q
\previousoldfn{\pl the labels of tex external vertices are wrong: it
  should vbe $1,i,j,l,n$ instead of $1,i,j,n,n$. Also the horizontal
  line is too short. \textbf{\Large{CAUTION}} Actually in january
 2011 I changed all the arguments of operads in the intro from $k$ to  
 $n$, to be uniformized with rest of the text. Therefore we should
 \textbf{keep} the old picture for 3 terms relation with labels of 
 middle points $i,j,k$ and last point $n$. Sorry for the unneeded extra
 work.}
in \refF{3terms}.

%\begin{figure}
%  \centering
%  \fbox{\includegraphics[width=80mm,height=30mm]{3terms.pdf}}
%  \caption{Diagrammatic description of the $3$-term relation}
%  \label{fig:3terms}
%\end{figure}

\begin{figure}[h]
\input{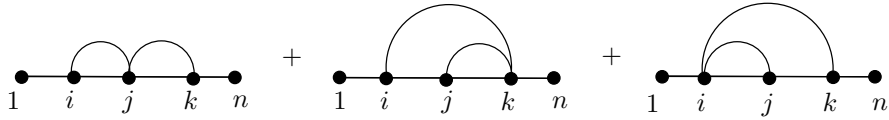}
\caption{Diagrammatic description of the $3$-term relation.}
\label{fig:3terms}
\end{figure}

In this figure, the vertices on the line correspond to the labels of the points $z_1,\dots,z_n$ of 
a configuration and each edge $(u,v)$ between two vertices represents a differential form 
$\omega_{uv}$.

The subalgebra of $\omDR (F_n(\BC);\BC)$ generated by the $\omega_{ij}$ is 
\[\frac{\wedge(\omega_{ij}:1\leq i<j\leq n)}
{\left(\omega_{ij}\wedge\omega_{jk}+\omega_{jk}\wedge\omega_{ki}+\omega_{ki}\wedge\omega_{ij}\right)}.
\]
This algebra has a trivial differential and it maps to the cohomology
algebra $\Ho^*(F_n(\BC);\BC)$.
 A Serre spectral sequence argument shows that this map is
actually an isomorphism.
 In other words, the cohomology embeds in the deRham algebra of forms, and hence $F_n(\BC)$ is formal.

Arnold's argument for $N=2$ can be generalized to all $N$ as follows. Consider the differential forms $\omega_{ij}=\theta^*_{ij}(\dvol)$
where
\begin{eqnarray*}
\theta_{ij}\colon &F_n(\BR^N)&\longrightarrow S^{N-1}\\
&(z_1,\dots,z_n)&\longmapsto\frac{z_j-z_i}{\|z_j-z_i\|},
\end{eqnarray*}
and $\dvol\in\omDR^{N-1}(S^{N-1})$ \p is the symmetric volume form on the sphere $S^{N-1}$ that integrates to $1$.
For $N=2$, these are  analogous to \D \refN{eq:omegaij}.
It is well known by work of F. Cohen that these forms generate the cohomology algebra of $F_n(\BR^N)$
and that the $3$-term relation holds  in cohomology. However, the relation is not always true at the level of forms.
One only knows that, for each $i$, $j$, and $k$, there exists some differential form  $\beta$
such that
\begin{equation}\label{eq:dbeta}
d\beta\,=\,\omega_{ij}\wedge\omega_{jk}+\omega_{jk}\wedge\omega_{ki}+\omega_{ki}\wedge\omega_{ij}.
\end{equation}
The key idea now is to describe an algorithm which constructs  \emph{in a natural way} such a cobounding
form  $\beta$.
To explain this, suppose that $n=3$ and $(i,j,k)=(1,2,3)$.
Consider the projection
\begin{equation}\label{eq:piF4F3}
\pi\colon F_4(\BR^N)\longrightarrow F_3(\BR^N)
\end{equation}
that forgets the fourth point of the configuration. It is a fibration\previousoldfn{\iv Give reference to Fadell and Husseini\pl No} with fiber  $$F=\BR^N\setminus\{z_1,z_2,z_3\}.$$

We will obtain $\beta$ by integration along the fiber of $\pi$ of some
suitable differential form $\alpha$ on $F_4(\BR^N)$.
To ensure convergence of the integral, we replace the spaces in the fibration \refN{eq:piF4F3} by their Fulton-MacPherson compactifications $\Conf[4]$ and $\Conf[3]$ 
so  that  the fiber becomes 
diffeomorphic to a closed disk in $\BR^N$ with three small open disks removed. We will denote this fiber by $\overline{F}$. 
Intuitively, each of the three inner boundary spheres of $\overline{F}$ corresponds 
to points $z_4$  becoming infinitesimaly close to $z_1$, $z_2$, or
$z_3$, \p(which we denote by $z_4\simeq z_i$), and the outer
boundary sphere of $\overline{F}$ corresponds to the point $z_4$ going
to infinity
 (which we denote by \p$z_4\simeq\infty)$).\q

Now consider the map
\begin{equation}\label{eq:theta}
\theta:=(\theta_{14},\theta_{24},\theta_{34})\colon \Conf[4]\longrightarrow S^{N-1}\times S^{N-1}\times S^{N-1}.
\end{equation}
The pullback form 
\[\theta^*(\dvol\times\dvol\times\dvol)\]
is a cocycle in $\omDR^{3N-3}(\Conf[4])$ 
and is exactly
\[\omega_{14}\wedge\omega_{24}\wedge\omega_{34}.
\]

{Integration along the fiber} of $\pi$ is a linear map
\begin{eqnarray*}
\pi_*=\fint_{\overline{F}}\colon&\omDR^{3N-3}(\Conf[4])&\longrightarrow\omDR^{2N-3}(\Conf[3])\\
&\alpha&\longmapsto \fint_{\overline{F}}\alpha.
\end{eqnarray*}
\p \q The integration takes place along the variable $z_4$ in the fiber
$\overline{F}$ which corresponds to the fourth component of a
configuration $z\in \Conf[4]$\p.  The map $\pi_*$ satisfies a fiberwise Stokes formula \q
\begin{equation}\label{eq:Stokesintro}
d( \fint_{\overline{F}}\alpha)=\fint_{\overline{F}}d(\alpha)\pm\fint_{\partial \overline{F}}\alpha.
\end{equation}
When $\alpha=\omega_{14}\wedge\omega_{24}\wedge\omega_{34}$, the first term on the right side of 
\refN{eq:Stokesintro} vanishes because $\alpha$ is a cocyle.
We study  its second term. One of the boundary  components  of 
$\overline{F}$ corresponds to \q\previousoldfn{\pl 17/01/11: there was no need to introduce two notations
    for that part of the boundary so I cancelled the notation
    $\partial_{14}F$ ... Good idea or not? \iv Looks ok to me.}
% $\partial_{14}\overline{F}$ or
 $\{z_4\simeq z_1\}\subset\partial\overline{F}$, and $\theta_{14}$ restricts
 to \Z  a \D diffeomorphism
$$
\theta_{14}\colon \{z_4\simeq z_1\}\iso S^{N-1}.
$$
We then \D  have
\Z 
\[\fint\limits_{\{z_4\simeq z_1\}}\omega_{14}\wedge\omega_{24}\wedge\omega_{34}=
\fint\limits_{\{z_4\simeq z_1\}}\omega_{14}\wedge\omega_{21}\wedge\omega_{31}=
\left(\ \int\limits_{S^{N-1}}\dvol\right)\cdot\omega_{21}\wedge\omega_{31}
=\omega_{21}\wedge\omega_{31}.
\]
Similarly the components corresponding to $z_4\simeq z_2$ and $z_4\simeq z_2$ give the two other
summands of the $3$-term relation \refN{eq:arnold}. Another argument shows that the integral along the
outer  boundary
corresponding to $z_4\simeq\infty$ vanishes. Thus
\[\beta:=\fint_{\overline{F}}\alpha\]
satisfies \refE{dbeta} and is naturally defined.

This algorithm for constructing $\beta$ can be encoded by a diagram $\Gamma$ as pictured in \refF{killer3terms}.
In this diagram, \D vertices $1,2,3$ (pictured on a line segment) are called \emph{external} and  vertex $4$ is called \emph{internal}. 

%\begin{figure}
%  \centering
%  \fbox{\includegraphics[width=50mm,height=30mm]{killer3terms.pdf}}
%  \caption{The diagram $\Gamma$ that cancels the $3$-term relation from \refF{3terms}.}
%  \label{fig:killer3terms}
%\end{figure}

\begin{figure}[h]
\input{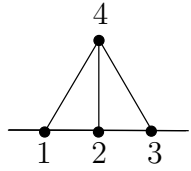}
\caption{The diagram $\Gamma$ that cancels the $3$-term relation from \refF{3terms}.}
\label{fig:killer3terms}
\end{figure}

The three edges $(1,4)$, $(2,4)$, and $(3,4)$ correspond to the
three components of the map $\theta$ from \refN{eq:theta}. To such a diagram we have associated the differential
form
\begin{equation}\label{eq:IGamma}
\IK(\Gamma):=\fint\limits_{\fiber}\theta_{14}^*(\dvol)\wedge\theta_{24}^*(\dvol)\wedge\theta_{34}^*(\dvol)=\pi_*(\theta^*(\underset{3}\times\dvol))
\end{equation}
where the points of the fiber are those labeled by internal vertices in the diagram $\Gamma$ (that is, not on the
\D horizontal line, which is $z_4$ in this case).

We define the coboundary of such a diagram $\Gamma$ by taking the sum over all possible contractions of an edge
with not all endpoints \D  on the line. In particular, for $\Gamma$ as in \refF{killer3terms},
its coboundary is exactly the diagrams of \refF{3terms} corresponding to the $3$-term relation specialized \D to
$n=3$ and $(i,j,k)=(1,2,3)$. Applying  $\IK$,
defined similarly as in \refN{eq:IGamma},   to  \D the diagrams  of  \refF{3terms} gives the right hand side of
\refN{eq:dbeta}, which we have shown to be $d(\IK(\Gamma))$.
In other words, $\IK$ commutes with the differential in this example.

The vector space of \D all such ``admissible'' diagrams will be
denoted by $\AD$ and will be endowed with the structure of a cooperad
in CDGA.
The generalization of Formula \refN{eq:IGamma} will define \D the
Kontsevich configuration space integral $\IK$ from \refN{eq:Iintro}.
An algebraic computation will show that $\AD(n)$, where $n$ is the number of external vertices
(the ones drawn on the
horizontal line segment),
is quasi-isomorphic to
$\Ho^*(\Conf[n])$, from which we will deduce that $\IK$ in \refN{eq:Iintro} is a quasi-isomorphism and hence that
$\Conf[n]$ is formal. Since these quasi-isomorphisms respect the
cooperadic structure, this will prove the formality of the operad $\Conf[\bullet]$ which is equivalent to the little disks operad.

There is  one last technical issue. The operad structure 
\previousoldfn{\pl dec2010: in the previous version there was an erroneous
  global change of "operad structure" into "fiber". I made a reverse changing.
 To easily spot the changes there are new commands \texttt{newcommand operadicstructure} and alike}
 on $\Conf[n]$ corresponds  to the
inclusions of various faces of the boundary of  $\Conf[n]$. Therefore, in order for $\IK$ to be a map of cooperads,
it is essential that the forms $\IK(\Gamma)$ are well-defined on this boundary.
However, the projection 
\[\pi\colon\Conf[n+l]\longrightarrow\Conf[n]\]
is unfortunately not a smooth submersion on the boundary $\partial\Conf[n]$
(see \refX{notsubm}),
 and hence
$\IK(\Gamma)$ need not be a smooth form on this boundary.
To fix this problem we will replace $\omDR$ by the CDGA $\ompa$  of PA
forms as defined in \cite[Appendix]{KoSo:def}.  These were studied in \D great detail in
\cite{HLTV:RHTSAS} and are reviewed in \refS{RHTSAS}.

%vefore

\subsection{Plan of the paper}
\ 

\p \q For a faster run through this paper, the reader could, after reading the Introduction,  jump directly to the beginning of \refS{KCSI}
to get a better idea of the construction of the quasi-isomorphism of operads
\[\IK\colon\AD(\bullet)\longrightarrow\ompa(\Conf[\bullet])\]
which is central to our proofs.  Along the way, 
 a quick look at Sections \ref{sec:compconf} and
\ref{sec:diagrams1}--\ref{sec:modulediagrams} will supply a better sense
of the Fulton-MacPherson operad $\Conf[\bullet]$ and the CDGA cooperad of
admissible diagrams $\AD(\bullet)$, respectively.\pfn{this paragraph answer referee
  comment (82)}

\p The plan of the paper is as follows (see also the Table of Contents at
the beginning of the paper).

\textbullet\ \ 
In \refS{notation} we fix some notation, and in particular establish some terminology relating to linear orders and weak ordered partitions which will be useful in describing the operad structure maps.

\textbullet\ \ 
In \refS{CDGAmodel} we define in detail what we mean by formality for operads. This is not as straighforward as one might wish 
 because the Sullivan-deRham functor $\Apl$ (or its semi-algebraic
 analog $\ompa$) does not turn operads into genuine
 cooperads  of \D CDGAs.  Our definition, however, is practical enough for applications. 

\textbullet\ \ 
In \refS{RHTSAS} we review the functor $\ompa$ of PA forms.  This is the analog for semi-algebraic spaces 
of the deRham functor $\omDR$ of differential forms for smooth manifolds. We review the main results we will
need from this theory, such as the notion of semi-algebraic {\current}s
$\CSA_*(X)$ on a{\D} semi-algebraic set $X$, which are weakly equivalent
to singular chains;  
 the fact that $\ompa $ encodes (monoidaly) the real homotopy type of compact semi-algebraic 
sets; and the important notion of integration along the fiber, or pushforward, of a ``minimal'' PA form along
a semi-algebraic bundle.

\textbullet\ \ 
In \refS{FMoperad} we define and study in detail the Fulton-MacPherson operad $\Conf[\bullet]$ and prove the results about this operad that are necessary for establishing certain properties of the Kontsevich configuration space integral. 
We also review \q the fact that the Fulton-MacPherson operad is equivalent to the little balls operad.

\textbullet\ \ 
In \refS{diagrams} we construct the combinatorial CDGA $\AD(n)$ of \emph{admissible diagrams}
 (on $n$ external vertices), built from a{\D} larger companion CDGA $\GD(n)$ of \emph{diagrams}. The 
CDGA $\AD(n)$  will later be shown to be quasi-isomorphic to both $\ompa (\Conf[n])$ and its cohomology.

\textbullet\ \ 
In \refS{diagcoop} we endow first $\GD$ and then $\AD$ with the structure of a cooperad.
The cooperad structure is obtained by considering
\emph{condensations}\pfn{check that AdmLoc, Loc, EssLoc, etc... were
  changed cosinsistently into Cond etc... all over the paper. \iv I did a search inside the tex document and also a pdf search in the pdf file and didn't find any issues.}, which will have already appeared in the study of the
Fulton-MacPherson operad in \refS{FMoperad}.

\textbullet\ \ 
In \refS{HoAD}, we prove that the cooperad $\AD$ is quasi-isomorphic to the cohomology of the Fulton-MacPherson operad.  

\textbullet\ \ 
In \refS{KCSI} we construct the Kontsevich configuration space integrals, which are CDGA maps
\[\GIK\colon\GD(n)\longrightarrow\ompa(\Conf[n])\quad\textrm{and}\quad\IK\colon\AD(n)\longrightarrow\ompa(\Conf[n]).\]
We prove that they are (almost) morphisms of cooperads. The arguments use many properties of the 
Fulton-MacPherson operad developed in  \refS{FMoperad}.

\textbullet\ \ 
In \refS{proofform} we collect the  results of the previous two sections 
to deduce our main formality results.  In particular,
we prove that $\IK$ is a quasi-isomorphism.

\textbullet\ \ 
\p \q Lastly, for the convenience of the reader we have included an index of
notation in the Appendix.

%In \refS{proofpibdl} of the Appendix we give a  proof of the technical fact that the canonical projection between Fulton-MacPherson compactifications of configuration spaces are semi-algebraic b%undles (this result is stated in \refS{FMoperad}).

\subsection{Acknowledgments}
Our deepest gratitude goes to Greg Arone for his encouragement, support, and patience.
We also thank Victor Turchin for his encouragement and for explaining the proof of \refT{formalAD} to us.
We also thank Nathalie Wahl for pointing out some errors and some 
 weaknesses in exposition in an earlier version of this paper.  We are grateful to Paolo Salvatore for pointing out to us the reference 
\cite[Lemma 6.4]{Kon:DQPM}.  
Parts of this paper were written while the first author was visiting the 
University of Virginia, Wellesley College, and the Center for
Deformation and Symmetry at University of Copenhagen and while the
second author was visiting University of{\D} Louvain, Massachusetts
Institute of Technology, \Z and the University of Virginia.  We would like to thank these institutions for their hospitality and support.  Lastly, we wish to thank the referee for a thorough reading of the paper and for helpful suggestions and comments.

\section{Notation, linear orders, weak partitions, and operads}    \label{sec:notation}

In this section we fix some notation, most of which is standard. We also review the notion of linear orders and introduce the notion of a weak ordered partition which is useful in describing the operad structure  maps.

\subsection{Notation}
\label{sec:notationnotation}

$\BK$ will be a commutative ring with unit, often $\BR$.

An integer $N\geq1$  (which gives the  ambient dimension) will be fixed.

For a set ${\ExtVert}$ we denote by $|{\ExtVert}|$ its cardinality.
We denote by $\Perm({\ExtVert})$ the group of permutations of $\ExtVert$.
For a nonnegative integer $n$, we set $\setn{n}=\{1,\dots,n\}$. We will sometimes identify $n$ and the set $\setn{n}$.
\p The set of all functions from a set $X$ to a set $Y$ is denoted
by $Y^X$. \q

When $f\colon X\to Y$ is a map and  $A\subset X$, we denote the restriction
of $f$ to $A$ by $f|A$.

We denote{\D} the one-point space by $*$.

\p We use the notation $x:=\operatorname{def}$ to state that the left hand side is
defined by the right hand side.\q

An extended index of notation is in the Appendix.\previousoldfn{it does not hurt to
  recall it here I guess}

\subsection{Linear orders}
    \label{sec:linord}
%\previousfn{\pl The reason for introducing linearly ordered set is that instead of labelling the internal vertices
%of a diagram by natural numbers, we will just consider them as a linearly ordered set. Many constructions (product of two diagrams, 
%contraction of an edge) become then more natural because we do not have to handle with strange relabelling etc.}

\begin{defin}\label{D:linord}
A \emph{linearly ordered} (or a \emph{totally ordered}) set is a pair $(L,\leq)$ where $L$ is a  set 
and $\leq$ is a reflexive, transitive, and antisymmetric relation on $L$
 such for any   $x,y\in L$ we have $x\leq y$ or $y\leq x$.
 We write $x<y$ when $x\leq y$ and $x\not=y$.

Given two {disjoint} linearly ordered sets $(L_1,\leq_1)$ and $(L_2,\leq_2)$ their
 \emph{ordered sum}
is the linearly ordered set $L_1\ordsum L_2:=(L_1\cup L_2,\leq)$ such that the restriction of $\leq$ to
$L_i$ is the given order $\leq_i$  and such that $x_1\leq x_2$ when $x_1\in L_1$
and  $x_2\in L_2$. 

More generally if $\{{{L}}_p\}_{p\in P}$ is a family of linearly ordered sets indexed by
a linearly ordered set $P$, its \emph{ordered sum}
 \[\underset{p\in P}{\ordsum}{{L}}_p\]
is the disjoint union  $\amalg_{p\in P}{{L}}_p$ equipped with a linear
order $\leq$ \p whose restriction to
each ${{L}}_p$ is the given order on that set and such that $x<y$ when $x\in {{L}}_p$ and $y\in {{L}}_q$ 
with $p<q$ in $P$.

It is clear that the ordered sum $\ordsum$  is associative
but not commutative. 

%Given two linear orders $\leq$ and $\leq'$ on the same set $L$
%we say that they \emph{differ by a single transposition} if
%there exists a transposition $\sigma=(a,b)$ in the group of permutations of $L$,
%for some distinct $a,b\in L$,
%such that $x\leq y$ if and only if  $\sigma(x)\leq'\sigma(y)$, for $x,y\in L$.

We define the \emph{position} function on a linearly ordered finite set $(L,\leq)$ as the unique order-preserving isomorphism
\[\pos\colon L\longrightarrow\{1,\dots,|L|\}.\] 
We write $\pos(x:L)$ for $\pos(x)$ when we want to emphasize the underlying
ordered set $L$.
\end{defin}

\subsection{Weak ordered partitions}

The following terminology will be useful in the description of operad structures  in the next section.
\begin{defin}\label{D:partition}
A \emph{weak partition} of a finite set ${\ExtVert}$ is a  map $\nu\colon {\ExtVert}\to P$, where $P$ is a finite set. 
The preimages $\nu^{-1}(p)$,
for $p\in P$, are the \emph{elements} of the partition. Since we do not ask $\nu$ to be surjective,
some of the elements $\nu^{-1}(p)$ can be empty, and hence the adjective \emph{weak}. The weak partition is \emph{degenerate}
if $\nu$ is not surjective, and \emph{non-degenerate} otherwise. We will simply say \emph{partition} for a non-degenerate weak partition. 
%\previousplfn{ok? }\previousoldfn{\pl commented ``(un)discrete'' and ``proper''}
% \item 
% A partition is \emph{undiscrete} if $P$ is a singleton
% and \emph{discrete} if $\nu$ is a bijection.
% \item
%  A partition
% is \emph{proper} if it is neither discrete, nor undiscrete.
%\item
% A \emph{weak partition} is a map
% $\nu\colon {\ExtVert}\to P$ (not necessarily surjective). 
The (weak) partition $\nu$ is \emph{ordered}
if its codomain $P$ is equipped with a linear order.
The \emph{undiscrete} partition is the partition $\nu\colon A\to\{1\}$
whose only element is $A$.
\end{defin}
\subsection{Operads and cooperads}\label{sec:operad}
Here we review the definition of operads that we will use.
Let $(\calC,\otimes,\unit)$ be a symmetric monoidal category.
Let $\IsoFin$ be the category whose objects are finite sets (including the empty set) and whose morphisms
are bijections between them. This category is equivalent to the category with one object for each integer $n\geq0$
along with the symmetric group $\Sigma_n=\Perm(\setn{n})$ as its set of automorphisms, and no other morphisms.
A \emph{symmetric sequence in $\calC$} is a functor
\[\calO\colon\IsoFin\longrightarrow\calC.\]
Thus a symmetric sequence in $\calC$ is determined
by a sequence $(\calO({n}))_{n\geq0}$ of objects of $\calC$ 
together with an action of  $\Sigma_n$
on $\calO({n})$.

An \emph{operad} $\calO\oprd$ is a symmetric sequence together with a \emph{unit map}
\[u\colon\unit\longrightarrow\calO(1)\]
and, for each ordered weak partition $\nu\colon\ExtVert\to P$, natural \emph{operad structure  maps}
\begin{equation}\label{eq:Thetanu}
\Theta_\nu\colon\calO(P)\otimes\underset{p\in P}{\otimes}\calO(\nu^{-1}(p))\longrightarrow\calO(\ExtVert)
\end{equation}
satisfying the usual associativity, unital, and equivariance conditions.
Here the monoidal product $\underset{p\in  P}{\otimes}$ is taken of
course in the linear{\D} order of $P$.

A \emph{cooperad} is an operad in the opposite category.

Our operads have an object $\calO(0)=\calO(\emptyset)$ in arity $0$. If we were working with operads without
 a nullary term, then we would only need  non-degenerate partitions $\nu$.

When investigating (co)operads, we will often fix the following setting:
\begin{setting}\label{setting:XP}
Fix an ordered weak partition  $\nu\colon {\ExtVert}\to P$,
with ${\ExtVert}$ and $P$ finite, and $P$ linearly ordered. 
We assume that $0\not\in P$ and set
\begin{equation}\label{eq:P*}
P^*:=\{0\}\ordsum P
\end{equation}
where $\ordsum$ is the ordered sum defined in \refS{linord}.
Set ${\ExtVert}_p=\nu^{-1}(p)$ for $p\in P$, and
${\ExtVert}_0=P$. 
\end{setting}
Under this setting the structure maps \refN{eq:Thetanu} become
\[\Theta_\nu\colon\underset{p\in P^*}{\otimes}\calO(\ExtVert_p)\longrightarrow\calO(\ExtVert).\]

\section{CDGA models for operads}\label{sec:CDGAmodel}

In this section we give precise meaning to the notion of a CDGA model for a 
 topological operad or for{\D} a morphism of topological operads. Our definition,
 although not  difficult, is perhaps not so elegant, but it suffices for the
applications we have in mind.
At the end of the section we  sketch an alternative, more concise definition.

Recall that Sullivan \cite{Sul:inf}  (see \cite{BoGu:RHT} or \cite{FHT:RHT}
for a complete development of the theory) constructed a contravariant
functor of piecewise polynomial forms over a field $\BK$ of characteristic $0$,
\[\Apl(-;\BK)\colon\Top\longrightarrow \CDGA\]
which mimics the deRham differential algebra of smooth differential forms on
a manifold. Here $\CDGA$ is the category{\D} of \emph{commutative differential graded $\BK-$algebras} (or CDGA for short) which are non-negatively graded. Sometimes we will also consider \emph{$\BZ$-graded CDGAs} which can be non trivial in negative degree, but those are not the objects of the category $\CDGA$.
 A CDGA $(A,d)$ is a \emph{CDGA model} (over $\BK$) 
for a space $X$ if the CDGAs $(A,d)$ and $\Apl(X;\BK)$ are weakly equivalent, by which we mean
that  there exists a chain of quasi-isomorphisms of CDGAs connecting them:
\[
(A,d)\stackrel{\simeq}\longleftarrow\cdots\quism \Apl(X;\BK).
\]
The main feature of the theory is that when $X$ is
a simply-connected topological space with finite Betti numbers and $\BK=\BQ$,
then any CDGA model for $X$  determines
the rational homotopy type of $X$.
Moreover, many rational homotopy invariants, like the rational cohomology algebra $\Ho^*(X;\BQ)$
or the rational homotopy Lie algebra \p $\pi_*(\Omega X)\otimes\BQ$  can easily
be recovered from the model $(A,d)$.
For fields $\BK$ other than the rationals, we have 
$$\Apl(-;\BK)=\Apl(-;\BQ)\otimes_\BQ\BK,$$ 
and by extension we say that the quasi-isomorphism
type of  $\Apl(X;\BK)$ determines the \emph{$\BK$-homotopy type} of $X$.
We just write $\Apl(X)$ when the field $\BK$ is understood.

Also, if $f\colon X\to Y$ is a map of spaces, we say that 
a CDGA morphism $$\phi\colon(B,d_B)\longrightarrow(A,d_A)$$ is a \emph{CDGA model} for $f$
if there exists a zig-zag of weak equivalences connecting $\phi$ and $\Apl(f;\BK)$,
that is, if there exists a commutative diagram of CDGAs
\[
\xymatrix{
(B,d_B)
\ar[d]_\phi
&
\ar[l]_-\simeq
\bullet
\ar[r]^-\simeq
\ar[d]
&
\cdots
&
\ar[l]_-\simeq
\bullet
\ar[r]^-\simeq
\ar[d]
&
\Apl(Y;\BK)
\ar[d]^{\Apl(f;\BK)}
\\
(A,d_A)
&
\ar[l]_-\simeq
\bullet
\ar[r]^-\simeq
&
\cdots
&
\ar[l]_-\simeq
\bullet
\ar[r]^-\simeq
&
\Apl(X;\BK)
}
\]
in which the horizontal arrows are quasi-isomorphisms.

We would like to define a similar notion of a CDGA model for a topological operad
$\calO\oprd$. A naive definition would be that such a model is a cooperad $\calA\oprd$
of CDGAs that is connected by weak equivalences
of CDGA cooperads to $\Apl(\calO\oprd)$.
However, there is a problem with this definition because the contravariant functor $\Apl$ is not  comonoidal as there is no 
suitable natural map
\begin{equation}\label{eq:Aplcomon}
\Apl(X\times Y)\longrightarrow\Apl(X)\otimes\Apl(Y).
\end{equation}
Therefore it seems that there is no  cooperad structure on $\Apl(\calO\oprd)$
naturally induced from the operad structure  on $\calO$.
On the other hand, $\Apl$ is monoidal through the Kunneth quasi-isomorphism
\begin{equation}\label{eq:Aplkappa}
\kappa\colon
\Apl(X)\otimes\Apl(Y)\quism
\Apl(X\times Y).
\end{equation}
This morphism becomes an isomorphism in the homotopy category, and its inverse  should correspond
to the
homotopy class of the missing
map \refN{eq:Aplcomon}. We would thus like to say that
$\Apl(\calO\oprd)$ is a  cooperad ``up to homotopy''.
However, this sort of  ``up to homotopy'' structure needs to be handled 
with more care than is necessary for our purpose, and so we will not pursue this in detail here and will just give an indication of such a notion at the end of the section. 
Instead we will propose in \refD{CDGAmodeloperad} an ad hoc  definition of a CDGA model for an operad.

There is a second difficulty which we will have do deal with and which comes
from the proof of the formality itself. Namely, in Kontsevich's proof of the weak equivalence between the
(up to homotopy) cooperad $\Apl(\calB\oprd)$ and its cohomology,
 a functor $\ompa$ (to be reviewed in \refS{RHTSAS}) is used.  This functor is weakly equivalent to $\Apl(-;\BR)$ but is defined  
only after restriction to a subcategory of $\Top$, namely the category of
compact semi-algebraic sets. This is analogous to the fact that the deRham CDGA of smooth differential forms $\omDR$
is weakly equivalent to $\Apl(-;\BR)$ after restriction to the subcategory of smooth manifolds.
Consequently, our modeling functors will sometimes be defined on some subcategory $u\colon\calT\hookrightarrow\Top$.

To finally define our notion of a CDGA model for an operad, we will need a few definitions.

Two  cooperads of CDGAs, $\calA\oprd$ and $\calA'\oprd$,
are \emph{weakly equivalent} if they are connected by a chain
of quasi-isomorphism of CDGA cooperads, 
\[ \calA\oprd\stackrel{\simeq}\longleftarrow\cdots\quism \calA'\oprd.
\]

Let $(\calT,\times,\unit)$ be a  symmetric monoidal category and let
\[u\colon\calT\longrightarrow\Top
\]
be a symmetric \emph{strongly} monoidal covariant functor, where by strongly we mean that the natural map
\begin{equation}\label{eq:ucross}
u(X)\times u(Y)\iso  u(X\times Y)
\end{equation}
is an isomorphism and $u(\unit)=*$ is the one-point space.

For us, a contravariant functor
\[\Omega\colon\calT\longrightarrow\CDGA
\]
is \emph{symmetric monoidal} if it is equipped with a natural map
\begin{equation}\label{eq:kappaF}
\kappa\colon \Omega (X)\otimes \Omega (Y)\longrightarrow \Omega (X\times Y)
\end{equation}
satisfying the usual axioms and  such that
$\Omega (\unit)=\BK$. In particular $\Apl\circ u$ is symmetric monoidal.

A \emph{natural monoidal quasi-isomorphism} between two such contravariant
symmetric monoidal functors $\Omega $ and $\Omega '$
is a natural transformation
\[\theta\colon \Omega \longrightarrow \Omega '\]
that induces an isomorphism in homology and that commutes with the monoidal
structure maps.
Two symmetric monoidal contravariant functors are \emph{weakly equivalent}
if they are connected by a chain of natural monoidal quasi-isomorphisms.
If $\Omega $ is weakly equivalent to $\Apl\circ u$ then
the morphism $\kappa$ of \refN{eq:kappaF} is a quasi-isomorphism
because the corresponding one for $\Apl$ in \refN{eq:Aplkappa} is as well and because of
the isomorphism
\refN{eq:ucross}.

Our definition of CDGA models for cooperads is then
\begin{defin}\label{D:CDGAmodeloperad}
A CDGA cooperad $\calA$ is a \emph{CDGA model} for a topological operad $\calO\oprd$
if there exist
\begin{itemize}
\item a CDGA cooperad $\calA'$  weakly equivalent to  $\calA$;
\item a symmetric monoidal category $(\calT,\times,\unit)$;
\item a symmetric strongly monoidal  covariant functor $u\colon\calT\to\Top$;
\item an operad $\calO'\oprd$ in $\calT$ such that $u(\calO'\oprd)$ is weakly equivalent to  $\calO\oprd$;
\item a symmetric monoidal contravariant functor $\Omega $ weakly equivalent to $\Apl\circ u$;
\item for each $n\geq0$ a $\Sigma_n$-equivariant quasi-isomorphism 
\[J_n\colon \calA'(n)\stackrel{\simeq}{\longrightarrow} \Omega (\calO'(n))\]
such that, for each $k\geq0$ and $n_1,\dots,n_k\geq0$ with $n=n_1+\dots+n_k$, the following
diagram commutes:
\[
\xymatrix{
\calA'(n)
\ar[rr]^-{J_n}_-{\simeq}
\ar[dd]_{\Psi}
&  &  
\Omega(\calO'(n))
\ar[d]^{\Omega(\Phi)}
\\
&  & 
\Omega(\calO'(k)\times\calO'(n_1)\times\dots\times\calO'(n_k))
\\
\calA'(k)\otimes\calA'(n_1)\otimes\dots\calA'(n_k)
\ar[rr]^-\simeq_-{J_k\otimes J_{n_1}\otimes\cdots\otimes J_{n_k}}
&  &
\ar[u]_{\kappa}^{\simeq}
\Omega(\calO'(k))\otimes \Omega(\calO'(n_1))\otimes\dots\otimes \Omega(\calO'(n_k))
.}
\]
Here $\Psi$ and $\Phi$ are the (co)operad structure  maps on $\calA'$ and $\calO'$ respectively,
and  the composition
\[\calA'(1)\stackrel{J_1}\longrightarrow \Omega (\calO'(1))\stackrel{\Omega (\eta)}\longrightarrow \Omega (\unit)\cong\BK\]
is required to be the counit of $\calA'$, where $\eta$ is the unit of $\calO'$.
\end{itemize}
\end{defin}
\vspace{5mm}

If $\kappa$ was an isomorphism, then
$\kappa^{-1}\circ\Omega(\Phi)$  would define a cooperad structure on
$\Omega(\calO')$ and  the above diagram{\D} would
simply mean
that the cooperads $\calA'$ and $\Omega(\calO')$ are weakly equivalent.

The main examples of the above that we will consider are:
\begin{itemize}
\item the category $\calT=\CompSemiAlg$ of compact semi-algebraic sets
(\refS{RHTSAS});
\item the forgetful functor $u\colon\CompSemiAlg\to\Top$; \q
\item the functor $\Omega=\ompa$ of semi-algebraic forms (\refS{RHTSAS});
\item the topological operad of little balls $\calO=\calB_N$;
\item the Fulton-MacPherson semi-algebraic operad $\calO'=\Conf[\bullet]$ (\refS{FMoperad}), and
\item its cohomology $\calA=\Ho^*(\Conf[\bullet])$;
\item the cooperad of admissible diagrams $\calA'=\AD$  (Sections \ref{sec:diagrams}-\ref{sec:diagcoop}); and 
\item the Kontsevich configuration space integral $J_n=\IK\colon \AD(n)\to\ompa (\Conf[n])$ (\refS{KCSI}).
\end{itemize}

We will let the reader generalize \refD{CDGAmodeloperad}  in an obvious way to say when a morphism of CDGA cooperads
\[\phi\colon\calA_1\oprd\longrightarrow\calA_2\oprd\]
is a \emph{CDGA model for a morphism of topological operads}
\[f\colon\calO_2\oprd\longrightarrow\calO_1\oprd.\]

\begin{defin}\label{D:formaloperad}
A topological operad is \emph{formal} over $\BK$ if the induced 
cohomology algebra cooperad is a\previousoldfn{\pl I have some trouble of the
  ``its''here and below. Indeed an operad has many CDGA models; so can we write
  ``its model''? -- in french it would be wrong}
  \previousoldfn{\iv Here and elsewhere, there was a strange space that's inserted because of the old footnotes.  I've tried fixing this throughout the paper.} 
  CDGA model for this operad over $\BK$.

A morphism of topological operads  is \emph{formal}
if the induced morphism in cohomology is a CDGA model for this operad morphism.
\end{defin}
%}%dontcompile

This definition, albeit perhaps a bit ad hoc, is good enough for the applications 
we have in mind. A more elegant  definition would have to use a precise notion of 
a (co)operad up to homotopy as follows.

Recall, for example from \cite[\S 1.2]{GiKa:Kos}, that an operad can be seen as a functor on the category of trees.
More precisely let $\Tree$ be the  category whose objects are rooted planar
trees and morphisms compositions of contractions of  non-terminal edges.
 Given trees $S,T_1,\dots,T_k$ where $S$ has $k$ leaves and 
each $T_i$ has $n_i$ leaves, one can build a new tree $S(T_1,\dots,T_k)$ with
$n_1+\dots+n_k$ leaves by grafting the root of each tree $T_i$ to the corresponding leaf of $S$. 
For $n\geq0$ we denote by $\langle n\rangle$  the tree with $n$ leaves and no internal vertex, that is
a tree which is indecomposable with respect to the grafting operation.
Then an operad $\calO$ in a symmetric monoidal category $\calC$ can be seen as a functor
\[O\colon\Tree\longrightarrow\calC
\]
\p where $O(\langle n\rangle)=\calO(n)$, for $n\geq0$.
In order for a functor $O$ to define an operad one asks
for isomorphisms
\[ \alpha_{(S,T_1,\dots,T_k)}\colon O(S(T_1,\dots,T_k))\iso O(S)\otimes\otimes_{i=1}^k O(T_i)
\]
satisfying obvious associativity, unital, and{\D} equivariance relations.%\iv

There is a morphism in $\Tree$ given by 
\[\langle k\rangle(\langle n_1\rangle,\dots,\langle n_k\rangle)\longrightarrow\langle n_1+\dots+n_k\rangle\]
and its image under the functor $O$ composed with the inverse of the isomorphism $\alpha$ 
gives the structure maps of the operad.

An \emph{operad up to homotopy} is an analogous functor $O$ 
except that one only asks 
$\alpha_{(S,T_1,\dots,T_k)}$
to be  a weak equivalence instead of an isomorphism.
Similarly we can define cooperads up to homotopy.

If $\calO$ is a topological operad, 
then $\Apl(\calO)$ naturally becomes a cooperad up to homotopy in this sense, with 
the weak equivalences  $\alpha$ constructed from the Kunneth quasi-isomorphism \refN{eq:Aplkappa}.
There is also an obvious notion of morphisms of (co)operads up to homotopy and of 
weak equivalences. One could\previousoldfn{\pl ``could'' is more honnest than
  ``can'' I guess (see next footnote)} check that if a CDGA cooperad $\calA$ is a CDGA model 
for a topological operad $\calO$ in the sense of \refD{CDGAmodeloperad},
then 
$\calA$ and $\Apl(\calO)$ are also weakly equivalent as cooperads up to homotopy.
This therefore might  give a better definition of an operad model, but it is possible that some further ``$\infty$-version" would be necessary for obtaining something useful.\p\q\pfn{see note
  (16) of the referee)}
\previousoldfn{\pl check \todo \iv  Seems ok, especially since you say ``might", so the reader knows that we're just sketching something.\pl hum... one can check suggest that it was done... But maybe ok like that anyway}
%\previousfn{\pl there was a last paragraph
%about another approach in terms of cofibrant replacements; I guess it was useless and I canceled it}

%One could also try to work with replacements by cofibrant (co)operads in order to avoid some of the problems.
%However there are also difficulties here.
%It is unclear whether there exist good cofibrant replacements in the category of CDGA cooperads.
%Even if one worked in the category of chain complexes instead of CDGA,
%the fact that we consider operads with a non-trivial term in arity $0$ could be an issue for
%the existence of enough cofibrant models.

\section{Real homotopy theory of semi-algebraic sets}
\label{sec:RHTSAS}
In this section we give a brief review of  Kontsevich and Soibelman's  theory 
 of
semi-algebraic differential forms which is outlined{\D} in \cite[\S 8]{KoSo:def}.  In particular we discuss the functor $\ompa$ which is analogous to
the deRham functor $\omDR$ for smooth manifolds.
That functor and{\D}  the way it encodes real homotopy theory of
semi-algebraic sets was developed in full detail
by the authors jointly with Robert Hardt and Victor Turchin in \cite{HLTV:RHTSAS}.

\begin{defin}[\cite{BCR:GAR}]  \previousoldfn{\pl 17/01/11: good idea to make it a definition; I
    added a reference to Coste-Bochniak-Roy} A \emph{semi-algebraic set}  is a subset of $\BR^p$ that is obtained by
 finite unions, finite intersections, and complements of subsets defined by polynomial equations and inequalities.
A \emph{semi-algebraic map} is a continuous map between semi-algebraic sets 
whose graph is a semi-algebraic set.
\end{defin}

We will consider the categories $\SemiAlg$ (and $\CompSemiAlg$)
of (compact) semi-algebraic sets.
Endowed with the cartesian product, this category becomes symmetric monoidal
and 
the obvious forgetful functor
\[u\colon\SemiAlg\longrightarrow\Top\]
is strongly symmetric monoidal because of the natural homeomorphism 
\[u(X)\times u(Y)\iso u(X\times Y).\]

We have for a semi-algebraic set $X$ a functorial chain complex of
\emph{semi-algebraic {\current}s} $\CSA_*(X)$ {\citePADC}, which is
weakly equivalent to{\D} singular chains. A typical element of $\CSA_k(X)$
is represented by a semi-algebraic map $g\colon M\to X$  from a semi-algebraic compact oriented manifold $M$ of dimension
$k$. This element is denoted by $g_*(\bbr{M})\in\CSA_k(X)$. In
particular, taking $g=\id_M$\p, 
\begin{equation}\label{eq:bbrM}
\bbr{M}\in\CSA_k(M)
\end{equation}
represents
a canonically defined fundamental class of the manifold $M$ 
at the level of semi-algebraic chains.  % (i.e.~semi-algebraic currents).
Also, a semi-algebraic map $f\colon X\to Y$ induces a chain map
\begin{equation}\label{eq:CSAf}
f_*\colon \CSA_*(X)\longrightarrow \CSA_*(Y).
\end{equation}

We in addition have a contravariant functor of \emph{minimal forms} {\citePASMinimalForms}
\begin{equation}\label{eq:minimalforms}
\omin\colon\SemiAlg\longrightarrow \CDGA.
\end{equation}
By definition, a minimal form of degree $k$ on $X$ is represented by a linear combination of
\[\mu=f_0\cdot df_1\wedge\dots\wedge df_k
\]
where 
\[f_0,f_1,\dots,f_k\colon X\longrightarrow\BR\]
are semi-algebraic maps.
Even though the $f_i$'s may not be everywhere smooth, 
we can define a differential $d\mu$ which is again a minimal form.
Also for a compact
 semi-algebraic oriented manifold $M$ of dimension $k$ and a semi-algebraic map $g\colon M\to X$, we can 
evaluate the form $\mu$ on  $g_*(\bbr{M})\in\CSA_k(X)$ by the formula
\begin{equation}\label{eq:pairingCOmPA}
\langle\mu\,,\,g_*\bbr{M}\rangle:=\int_M g^*(f_0\cdot df_1\wedge\dots\wedge df_k).
\end{equation}

\pfn{I rephrased the following paragraph}
The convergence of the integral on the right is a consequence of the
semi-algebraicity of $M$. Indeed that integral is the same as
\begin{equation}
\label{eq:intfgM}
\int_{f_*(g_*(M))}x_0\cdot dx_1\wedge\cdots\wedge dx_k
\end{equation}
where $f_*(g_*(M))$ is the image of $M$ in $\BR^{k+1}$ (counted with multiplities)
under the composition \q of $g$ and $f:=(f_0,f_1,\dots,f_k)$.
Thus $f_*(g_*(M))$ is a compact semi algebraic-set of dimension $\leq
k$, which implies that its $k$-volume is finite (this would not be
true for non semi-algebraic compact sets.) Hence the integral in
\refE{intfgM} converges.
%The fact that the integral on the right is  convergent comes from the fact that compact semi-algebraic sets
%are of finite volume 
See \cite[Theorem 2.4 and beginning of Section
3]{HLTV:RHTSAS} for more details.

In this paper, the only minimal forms that we will use are the standard volume form on the sphere and its pullbacks 
along semi-algebraic maps.

The CDGA of minimal forms embeds in that of \emph{PA forms}\previousoldfn{\iv Do
  we say anywhere what PA forms are?  It's sort of strange that we
  don't explain this since we work with the functor $\ompa$ so
  much.\pl I added a phrase saying that PA forms are obtained by
  pushforward of minimal forms (actually I say it twice)} {\citePASSPAForms}
  (``PA" stands for ``piecewise algebraic")
\begin{equation}
\ompa\colon\SemiAlg\longrightarrow\CDGA.\label{eq:ompa}
\end{equation}
Roughly speaking,{\D} PA forms are obtained by integration along the fiber of minimal forms along  oriented semi-algebraic bundles, which are recalled below.
The important feature is the following\previousoldfn{\iv This should be
  changed into a Theorem in the PA paper.\pl\todo I will do that with
  the changes in PA forms paper} 
\begin{thm}[\citePAPmonompaapl]\label{T:ompaapl}
When restricted to the category of compact semi-algebraic sets,
the contravariant symmetric monoidal functors
$\ompa$ and $\Apl(u(-);\BR)$ are weakly equivalent.
\end{thm}
 
Another important feature of minimal and PA forms is that 
classical integration along the fiber for smooth forms can be extended to the semi-algebraic framework.
 To explain, we have from {\citePASSAbundle} the notion of a \q
\emph{semi-algebraic bundle}, or \emph{SA~bundle} for short, which is
the obvious generalization of the usual definition of a locally trivial
{\D} bundle.
% (we leave it to the reader to write this this definition out). 

An SA bundle
\[\pi\colon E\longrightarrow B\]
is \emph{oriented} if its fibers are compact oriented semi-algebraic manifolds,
with orientation which is locally constant in an obvious sense.
For each $b\in B$ we then have{\D} the fundamental class of the fiber over $b$,
\begin{equation}
  \label{eq:fundfiber}
  \bbr{\pi^{-1}(b)}\in\CSA_k(\pi^{-1}(b)),
\end{equation}
where $k$ is the dimension of{\D} the fiber.

\pfn{I added the followin paragraph because from notes (42)( of the
  referee it seems that the concept of fiberwise boundary is obscure}
Given an oriented SA bundle $\pi\colon E\to B$ whose fibers are
compact SA manifolds, there exists a subbundle
\begin{equation}
\label{eq:pipartial}
\pi^\partial \colon E^\partial\to B
\end{equation}
whose fibers are the boundaries of the fibers of $\pi$.
This subbundle is called the\emph{ fiberwise boundary} of $\pi$
(see \cite[Definition 8.1]{HLTV:RHTSAS}).\pfn{check carefully this reference in
  the bibliography.\iv This is the right definition that paper, but that paper's information needs to be updated in the bibliography.}
An example is the map
\[\proj_1\colon E:=[0,1]\times[0,1]\to[0,1]\]
which projects onto the first factor. \q
In this case the fiberwise boundary is $E^\partial=[0,1]\times\{0,1\}$,
 but this is not the boundary of $E$.

For an oriented SA bundle with $k$-dimensional  fiber, there is a linear map
of degree $-k$ {\citePADintegrfiber},
\begin{equation}
\pi_*\colon\omin^{*+k}(E)\longrightarrow\ompa^*(B),\label{eq:pushforward}
\end{equation}
which correponds to integration along the fiber, also{\D} called
\emph{pushforward}. In some sense PA forms are{\D} obtained 
as (generalized) pushforwards of minimal forms \citePADPAForm.
Properties of the pushforward that we will need here are collected in {\citePASpropintfiber}.
They are analogous to the standard properties of integration of smooth differential forms
 along the compact fiber of a smooth bundle. In particular one has a fiberwise Stokes{\D} formula which we will need later.

\section{The Fulton-MacPherson operad}\label{sec:FMoperad}

Fix $N\geq 1$.  In this section we review the Fulton-MacPherson operad
$$\Conf[\bullet]=\{\Conf[n]\}_{n\geq{0}}$$
which is weakly equivalent to the little $N$-disks operad.
As a space, each $\Conf[n]$ is a  compactification of the 
  space $\Conf(n)$ of normalized ordered configurations
of $n$ points in $\BR^N$. It is a  compact  semi-algebraic manifold with boundary,
and so its real homotopy type is encoded by the semi-algebraic analog of deRham theory, $\ompa (\Conf[n])$.
 The operad structure  maps
correspond essentially to inclusions of  various faces of the boundary $\partial\Conf[n]$. 

In this section we will also study \emph{canonical projections} 
$$\pi\colon\Conf[n+l]\longrightarrow\Conf[n]$$
given by forgetting some points of the configuration and will prove that they are SA bundles with compact
manifolds as fibers. This fact will be used in \refS{consGIK} to construct the 
Kontsevich configuration space integral $\IK\colon\AD(n)\to\ompa (\Conf[n])$ of \refN{eq:Iintro}
via integration along the fiber of $\pi$.
We will also study  the interaction of these canonical projections
with the operad structure  in order to later prove that $\IK$ is a map of cooperads.

The plan of this section is the following:
\begin{itemize}
\item[\ref{sec:compconf}:] We define the compactification $\Conf[n]$, compute its dimension, and characterize its 
boundary.
\item[\ref{sec:opC}:] We describe the operad structure  on $\{\Conf[n]\}_{n\geq0}$ and recall 
that this operad\previousoldfn{\pl collection? why not operad?} is equivalent to the operad of little balls.
\item[\ref{sec:canproj}:] We study the canonical projection $\pi\colon\Conf[n+l]\to\Conf[n]$ and 
state that it defines a bundle whose fibers are oriented compact manifolds.
% This will be
%important since the Kontsevich configuration space integral is defined as an integral along the fiber of this bundle.
\item[\ref{sec:bdryCV}:] We  decompose the boundary $\partial\Conf[n]$ into faces which are images of the $\circ_i$ (``circle-$i$'') operad maps.
\item[\ref{sec:singconf}:] We construct \emph{singular} configuration spaces which are variations of spaces $\Conf[n]$ and will
be needed for some technical points. 
\item[\ref{sec:proj-oper}:] In this (long) section, we investigate the pullback of a canonical projection along
 an operad structure map. This will be needed for proving that the Kontsevich configuration space integral respects the (co)operadic
structures. We introduce at the beginning{\D} of this section the notion of a \emph{condensation} which will also be needed for the definition of the cooperad structure on the space of diagrams in \refS{diagcoop}.
\item[\ref{sec:bdryCVpi}:] We describe a decomposition of the fiberwise boundary of the total space $\Conf[n+l]$ of a
canonical projection. This will be used in proving that Kontsevich's
configuration space integral is a chain map.
\item[\ref{sec:orCA}:] We fix an orientation of $\Conf[n]$; this will be important when we integrate  forms over this manifold.
\item[\ref{sec:proofpibdl}:] We prove \refT{projSAbdl}, stated in
  \refS{canproj}, which asserts that the canonical projections
  are{\D} oriented SA bundles. This section also  contains an example
  showing that the canonical projections are not smooth bundles.
\end{itemize}

On a first pass of this section, the reader may just concentrate on Sections \ref{sec:compconf}-\ref{sec:bdryCV} to acquire a good sense of the{\D} Fulton-McPherson operad.
The last  five sections are more technical and are needed only for the details of the proof of certain properties of the Kontsevich configuration space integral in \refS{KCSI}. However, the reader should still look at 
\refD{loc} of a condensation in \refS{proj-oper}, as this will be
needed in \refS{diagcoop} \p to define the cooperadic structure on the
spaces of diagrams $\GD(n)$.

\subsection{Compactification of configuration spaces in $\BR^N$}\label{sec:compconf}

We first recall the Fulton-MacPherson compactification $\Conf[n]$  of the configuration space 
$\Conf(n)$ of $n$ points in $\BR^N$. This compactification (or at least some variation of it) was defined in
\cite{FuMc:com}, with the operad structure given in \cite{GeJo:ope},  and alternatively by Kontsevich
in \cite[Definition 12]{Kon:OMDQ} and \cite[Section 5.1]{Kon:DQPM}. We follow Kontsevich's approach, which was
corrected by Gaiffi in \cite[Section 6.2]{Gai:mod} and developed in detail
by Sinha in \cite{Sin:man} (the equivalence of the Kontsevich and the Fulton-MacPherson definitions follows from Sinha's work as well).

Let ${\ExtVert}$ be a finite set  of cardinality $n$ which will serve as a set of labels for the points of the configurations.
Consider the space 
\begin{equation}\label{eq:Inj}
\Inj({\ExtVert},\BR^N):=\{x\colon {\ExtVert}\hookrightarrow\BR^N\}
\end{equation}
of all injective maps from $A$ to $\BR^N$. 
An element  $x\in\Inj({\ExtVert},\BR^N)$ is an (ordered) configuration $(x(a))_{a\in {\ExtVert}}$ of
$n$ distinct points in $\BR^N$. This space is  topologized as a subspace
of the product $(\BR^N)^{\ExtVert}=\prod_{a\in {\ExtVert}}\BR^N$.

The space $\Inj({\ExtVert},\BR^N)$ is a smooth open manifold of dimension $N\cdot |{\ExtVert}|$.
The group of orientation-preserving similarities $\BR^N\rtimes \BR^+_0$ acts
%\previousfn{\pl is this the correct notation for the affine group} 
by
translation and positive dilation on $\BR^N$, and hence diagonally
on  $\Inj({\ExtVert},\BR^N)$. We denote its orbit space by
\begin{equation}
  \label{eq:defC(X)}
  \Conf({\ExtVert}):= \Inj({\ExtVert},\BR^N)/(\BR^N\rtimes \BR^+_0).
\end{equation}
(This space is denoted by $\widetilde{\Conf_n}(\BR^N)$ in \cite[Definition 3.9]{Sin:man}.)

When $|{\ExtVert}|\geq2$ the action is free and smooth and hence{\D} $\Conf(\ExtVert)$ is a manifold of dimension
\[
\dim\Conf({\ExtVert})=N\cdot|{\ExtVert}|-N-1,\]
\p and when $|{\ExtVert}|\leq1$ then $C({\ExtVert})$ is a one-point space because the action is transitive.

Define the \emph{barycenter} of a map $x\colon A\to\BR^N$ as the point
\begin{equation}\label{eq:bary}
\bary(x)=\bary(x(a):a\in A):=\frac1{|A|}\sum_{a\in A}x(a)
\end{equation}  
and{\D}  its \emph{radius} as the real number
\begin{equation}\label{eq:radius}
\radius(x)=\radius(x(a):a\in A):=\max(\|x(a)-\bary(x)\|:a\in A).
\end{equation}

When $|A|\geq2$, $\Conf({\ExtVert})$  is homeomorphic to
the space of \emph{normalized configurations}
\begin{equation}\label{eq:defC(X)xi}
\Inj_0^1(A,\BR^N):=
\left\{x\in\Inj(A,\BR^N):\bary(x)=0\textrm{ and } \radius(x)=1\right\}.
\end{equation}
We will use{\D}   $\Conf({\ExtVert})$ and $\Inj_0^1(A,\BR^N)$ interchangeably. 
\p \q Most of the time in this paper, a configuration will be denoted by
$x$ or $y$ (maybe with some decoration) and, when seen as an element
of $\Inj_0^1(A,\BR^N)$, its components will be points $x(a)$ for $a$
an element of the set of labels of the components, $\ExtVert$.\pfn{This is to
  answer (58) of the referee; I do not like so much to have bold face for the
  configuration: it would be a huge change and moreover makes no sense
  to have $x(i)$ to be a component of $\mathbf{x}$.\iv I agree.}

Denote by $S^{N-1}$ the unit sphere in $\BR^N$\p.
Given two distinct elements $a,b\in {\ExtVert}$, consider the map
\begin{eqnarray}
  \label{eq:defthetaab}
 \theta_{a,b}\colon\Conf({\ExtVert})&\longrightarrow&S^{N-1}\\
x&\longmapsto&\frac{x(b)-x(a)}{\|x(b)-x(a)\|}\notag
\end{eqnarray}
which gives the direction  between  two points of the configuration.

For three distinct elements $a,b,c\in {\ExtVert}$, also define 
\begin{eqnarray}
  \label{eq:defdeltaabc}
\delta_{a,b,c}\colon\Conf({\ExtVert})&\longrightarrow&[0,+\infty]\\
x&\longmapsto&\frac{\|x(a)-x(b)\|}{\|x(a)-x(c)\|}\notag
\end{eqnarray}
which gives the relative distance of $3$ points of a configuration.

Set
\begin{eqnarray*}
{\ExtVert}^{\{2\}}&=&\{(a,b)\in {\ExtVert}\times {\ExtVert}:a\not= b\}\\
{\ExtVert}^{\{3\}}&=&\{(a,b,c)\in {\ExtVert}\times {\ExtVert}\times {\ExtVert}:a\not= b\not=c\not=a\}
\end{eqnarray*}
and consider the map
\begin{align*}
  \label{eq:iota}
\iota\colon\Conf({\ExtVert})&\longrightarrow(S^{N-1})^{{\ExtVert}^{\{2\}}}\times [0,+\infty]^{{\ExtVert}^{\{3\}}}\\
x & \longmapsto\Big((\theta_{a,b}(x))_{(a,b)\in{{\ExtVert}^{\{2\}}}}\,,\,(\delta_{a,b,c}(x))_{(a,b,c)\in{{\ExtVert}^{\{3\}}}}\Big)\notag.
\end{align*}
Up to translation and dilation, any configuration $x\colon\ExtVert\hookrightarrow \BR^N$
can be recovered from the directions $\theta_{a,b}(x)$ and relative distances $\delta_{a,b,c}(x)$.
Hence $\iota$ is a homeomorphism onto its image \cite[Lemma 3.18]{Sin:man}
and we will identify $\Conf({\ExtVert})$ with $\iota(\Conf({\ExtVert}))$.

\begin{defin}\label{D:FM}
The \emph{Fulton-MacPherson compactification $\Conf[{\ExtVert}]$ of $\Conf({\ExtVert})$} is the topological closure of the image of $\iota$, that is, 
\[\Conf[{\ExtVert}]:=\overline{\iota(\Conf({\ExtVert}))}.\]
\end{defin}

Intuitively, one should think of $x\in \Conf[{\ExtVert}]$ as a ``virtual'' configuration
in which some points are possibly infinitesimally close to each other in such a way
that the direction between any two points and the relative distance between three points is
always well-defined. These directions and relative distances are given by the maps $\theta_{a,b}$ and $\delta_{a,b,c}$,
which obviously extend to $\Conf[{\ExtVert}]$. Moreover an element $x\in\Conf[\ExtVert]$ is completely
 characterized by the values $\theta_{a,b}(x)\in S^{N-1}$
and $\delta_{a,b,c}(x)\in[0,+\infty]$, for distinct $a,b,c\in \ExtVert$. By abuse 
of terminology an element $x\in\Conf[\ExtVert]$ will be called
a configuration and we will talk informally of its components $x(a)\in\BR^N$, for $a\in A$.

The following notation will be useful: For $a,b,c$ distinct  in ${\ExtVert}$ and $x\in\Conf[{\ExtVert}]$,
when $\delta_{a,b,c}(x)=0$ we write
\begin{equation}
  \label{eq:defrel}
  x(a)\simeq x(b)\rel x(c).
\end{equation}
 This happens exactly when the points $x(a)$ and $x(b)$
are infinitesimaly close to each other in comparison to their distance to $x(c)$.   Pictorial interpretations of this situation are given below in Example \ref{Ex:operadicCA}.
In particular  \refF{OperadMapOutput} represents{\D} a configuration
$x\in\Conf[6]$ with $N=2$.

%For example in \refF{operadicCA}
%$x_0\in\Conf[\{1,2,3,4\}]=\Conf[4]$ represents a configuration in the plane $\BR^2$ where the components $x_0(2)$ and $x_0(4)$
%are infinitesimally close to each other in comparaison to their distance to $x_0(1)$ or $x_0(3)$.
%In the same figure $\Phi_\nu(x_0,x_1,x_2,x_3,x_4)$ is another configuration in $\Conf[6]$.
%An example is given at \refF{virtconf9} which represents a configuration $x_0\in\Conf[9]$ of nine points in the plane $\BR^2$.
%The three points labeled $7,8,9$ are pictured in a screen to express that there are infinitesimally close to each other.
%Screens can be nested, so in this picture $x_0(1)\simeq x_0(4)\rel x_0(5)$, but
%also $x_0(1)\simeq x_0(2)\rel x_0(4)$.
%\begin{figure}
%  \centering
%  \fbox{\includegraphics[width=95mm,height=60mm]{virtualconfig9.pdf}}
%  \caption{A  configuration $x_0\in\Conf[9]$}
%  \label{fig:virtconf9}
%\end{figure}

The space $C({\ExtVert})\subset(\BR^N)^{\ExtVert}$ and the map $\iota$ are clearly semi-algebraic, 
therefore so is the closure $C[{\ExtVert}]$. Moreover, by \cite{BoTa:SLK} or \cite{Sin:man}, $\Conf[\ExtVert]$ is a compact manifold with corners.
 It is easy to see that the atlases given in those papers are semi-algebraic, and hence $\Conf[\ExtVert]$ is a compact 
semi-algebraic manifold with boundary (charts are{\D} given in \refL{Phihomeo}).
% \previousfn{\pl this is not  good reference. Markl? \iv Do we mean a compact manifold with corners?
%   Either Bott-Taubes or Sinha ``Manifold-theoretic...'' is good reference. },
% $C[{\ExtVert}]$ is a compact manifold whose interior is
% \[\C[{\ExtVert}]\setminus\partial\Conf[{\ExtVert}]=\Conf({\ExtVert}).\]
% Since $\Conf(\ExtVert)$ is a semi-algebraic manifold we get that
%  $\Conf[{\ExtVert}]$ is a compact semi-algebraic manifold.\previousfn{\pl we need the fact that
% the closure of a semi-algebraic manifold is also a semi-algebraic manifold (it is clear that
% it is a semi-algebraic set). I guess that this should come from the fact that a semi-alg manifold has probably a finite
% atlas.  \iv The details for this, in the particular case of $C[A]$ can either be found in Bott-Taubes or in my "Survey of Bott-Taubes integration".  Both papers have charts for $C[A]$ which are semi-algebraic.  But I'm not sure how much detail you want to get into. }

In conclusion, we have 
\begin{prop}\label{P:dimConf}
For a finite set ${\ExtVert}$, $\Conf[{\ExtVert}]$ is a compact semi-algebraic manifold with interior
 $\Conf({\ExtVert})$ and its dimension is given by
\[\dim(\Conf[{\ExtVert}])=
\begin{cases}
0&\textrm{if }|{\ExtVert}|\leq1;\\
N\cdot|{\ExtVert}|-N-1&\textrm{if }|{\ExtVert}|\geq2.
\end{cases}
\]
\end{prop}

We also have the following important characterization of the boundary 
\begin{prop}
  \label{P:charbdryCX}
For $x\in\Conf[\ExtVert]$, the following are equivalent conditions:
\[x\in\partial\Conf[{\ExtVert}]\quad\Longleftrightarrow\quad(\exists\, a,b,c\in {\ExtVert}\textrm{ distinct }:\,x(a)\simeq x(b)\rel x(c)).
\]
\end{prop}
%\previousfn{\pl add a reference or small arg for that proposition.  \iv I guess I'm not sure what this proposition is saying.  What if all points collide at the same time and at the same rate? \pl This cannot happen since we mod out by scaling. I guess no need for a ref}
For $|A|\leq1$, $\Conf[A]$ is a one-point space; for $|A|=2$, it is homeomorphic to
the sphere $S^{N-1}$. For $n\geq0$ we set $\Conf[n]:=\Conf[\{1,\dots,n\}]$.

\subsection{The operad structure}\label{sec:opC}
We will now define the structure of an operad on 

$$\Conf[\bullet]=\{\Conf[n]\}_{n\geq0}.$$

Recall from Section \ref{sec:notation}  the notion of weak ordered
partitions and how operad structure maps are  associated to them.

Fix a finite set $A$, a linearly ordered finite set $P$, and a weak ordered partition $\nu\colon A\to P$. Set
$$P^*=\{0\}\ordsum P,\ \  \ExtVert_p=\nu^{-1}(p),\ \  \text{and}\  A_0=P$$  as in the setting \ref{setting:XP} from \refS{notation}.
Hence 
\[\prod_{p\in P^*}\Conf[\ExtVert_p]=\Conf[P]\times\prod_{p\in P}\Conf[\nu^{-1}(p)].\]

We now construct an operad structure  map\previousoldfn{\iv What's a ``fiber map"?  Fibration?\pl "operad structure" had globally be changed to "fiber"...}
\begin{equation}\label{eq:Phinu}
\Phi_\nu\colon\prod_{p\in P^*}\Conf[{\ExtVert}_p]\longrightarrow\Conf[{\ExtVert}]
\end{equation}
as follows.
Intuitively the  configuration $x=\Phi_\nu((x_p)_{p\in P^*})$ is obtained by replacing, for each $p\in P$,
the $p$-th component $x_0(p)$ of the configuration $x_0\in \Conf[P]$ 
by the  configuration $x_p\in\Conf[{\ExtVert}_p]$ made infinitesimal.  To illustrate, we first give an example.

\begin{eg}\label{Ex:operadicCA}
Consider $P=\{\alpha, \beta, \gamma, \delta\}$
(with the linear order $\alpha<\beta<\gamma<\delta$)\p,
 $A=\{1, 2, 3, 4, 5, 6\}$, and let
$\nu\colon A\to P$
be given by 
$$
\nu(a)=\begin{cases}
\alpha, &  \text{for } a=1, 2; \\
\beta, &  \text{for } a=3, 4, 5; \\
\delta, &  \text{for } a=6.
\end{cases}
$$

Consider
\begin{align*}
x_0  & \in C[P]\cong C[4]; \\
x_{\alpha}  & \in C[\{1, 2\}]\cong C[2]; \\
x_{\beta}  & \in C[\{3,4, 5\}]\cong C[3]; \\
x_{\gamma}  & \in C[\emptyset]\cong C[0]= \ast; \\
x_{\delta}  & \in C[\{ 6\}]\cong C[1]=\ast
\end{align*}
and suppose{\D} that these configurations are for example as in Figure \ref{fig:OperadMapInput} (with $N=2$). \q

\begin{figure}[h]
\input{OperadMapInput.pstex_t}
\caption{}
\label{fig:OperadMapInput}
\end{figure}

This kind of pictorial representation of compactified configuration spaces first appeared in \cite{Sin:man}.  The plane represents $\mathbb{R}^N$ and the ``funnels" represent infinitesimal configurations.  Thus for example, in the picture of $x_0$, points labeled by{\D} $\alpha$ and $\delta$ are infinitesimally close to each other from the point of view of $\beta$ and $\gamma$. In notation of relation \eqref{eq:defrel},
$x_0(\alpha)\simeq x_0(\delta)\rel x_0(\beta)$ and $x_0(\alpha)\simeq x_0(\delta)\rel x_0(\gamma)$.
 Similarly in the picture of $x$ in Figure \ref{fig:OperadMapOutput}
 below, points (labeled by){\D} 4, 3, and 5 are infinitesimally close
 to each other from the point of view of 6, 1, and 2, but 3 and 5 are
 infinitesimally close to each other from the point of view of 4, as
 are 1 and 2 from the point of view of 6. 

Then the configuration $x=\Psi_{\nu}(x_0, x_{\alpha},
x_{\beta}, x_{\gamma}, x_{\delta})$ 
can be{\D} represented \pfn{in that picture the point 6 should be on
  the same horizontal line as the funnel (12), because $\alpha$ and
  $\delta$ are on the same horizontal line.  \iv  Fixed.} as in Figure \ref{fig:OperadMapOutput}.
  \previousoldfn{\iv In the public version, there is a space between ``represented" and ``as" in the above sentence.  The problem seems to be the $\{\backslash\text{D}\}$ in the tex file (there's a lot of those all over) and I have no idea what that does.  We should maybe get rid of it so we don't get these spaces.}

\begin{figure}[h]
\input{OperadMapOutput.pstex_t}
\caption{}
\label{fig:OperadMapOutput}
\end{figure}

\previousoldfn{\ppfn{\iv Iris saw this picture and said:  ``Are these ice cream
  cones?  Is this ice cream math?"\textbf{pascal} :-)}}

\end{eg}

%An example is given at \refF{operadicCA}\previousfn{\pl the 4 big black
%dots in top of the ficture were just holes in the sheet... ;-) }
%% for $p\in P$ (see \cite[Figure 4.6]{Sin:OKS}).
%where $A=\{a,b,c,d,e,f\}$, $P=\{1,2,3,4\}$, and the weak partition $\nu\colon A\to P$ is defined by $\nu(a)=\nu(b)=\nu(c)=1$,
%$\nu(d)=\nu(e)=2$, and $\nu(f)=4$.
%\begin{figure}\centering
%\fbox{
%\includegraphics[width=140mm,height=80mm]{operadicCA.pdf}
%%{\large \textbf{COMMMENT 14oct09}:} here we should have the picture 
%%``operadicCA.pdf'' which is not found on my computer... 
%%}
%}
%\caption{An example of the value of an operad structure map $\Phi_\nu$}
%\label{fig:operadicCA}
%\end{figure}
 
More precisely, $x=\Phi_\nu((x_p)_{p\in P^*})\in\Conf[A]$
 is characterized by, for  distinct $a,b,c\in A$,
\[
\theta_{a,b}(x)=
\begin{cases}
\theta_{a,b}(x_p),&
\textrm{if }a,b\in {\ExtVert}_p\textrm{ for some }p\in P,\textrm{  that is }\nu(a)=\nu(b)=p;\\
\theta_{\nu(a),\nu(b)}(x_0),&
\textrm{if }\nu(a)\not=\nu(b),
\end{cases}
\]
and
\[
\delta_{a,b,c}(x)=
\begin{cases}
\delta_{a,b,c}(x_p),&
\textrm{if }a,b,c\in {\ExtVert}_p\textrm{ for some }p\in P;\\
\delta_{\nu(a),\nu(b),\nu(c)}(x_0),&
\textrm{if }\nu(a),\,\nu(b),\textrm{ and }\nu(c)\textrm{ are all distinct};\\
0,&\textrm{if }
\nu(a)=\nu(b)\not=\nu(c);\\
1,&\textrm{if }
\nu(a)\not=\nu(b)=\nu(c);\\
+\infty,&\textrm{if }
\nu(a)=\nu(c)\not=\nu(b).\\
\end{cases}
\]

There is an obvious action of the group $\Perm({\ExtVert})$ of permutations of the set ${\ExtVert}$ on $\Conf[{\ExtVert}]$, 
and in particular of the symmetric group $\Sigma_n$ on $\Conf[n]$.
We define the {unit} in $\Conf[{1}]$ as its unique point (or more precisely the unique map $u\colon*\to\Conf[1]$).

The following is straightforward to check (see for example {\cite[Section 4]{Sin:OKS}}). 
\begin{prop}\label{P:operadconf}
The above data endows $\Conf[\bullet]=\{\Conf[n]\}_{n\geq0}$  with the structure of an operad of compact semi-algebraic sets.
\end{prop}
\previousoldfn{\pl add a reference. Sinha? Markl? \iv I added a reference, but Dev does this for the Kontsevich operad and of course doesn't mention semi-algebraic sets, so I'm not sure if that's ok.}

The relevance of the Fulton-MacPherson operad for us is that it is weakly  equivalent to the little balls operad, as
proved by P.~Salvatore: \p\q
\begin{prop}{\cite[Proposition 4.9]{Sal:sum}}\label{P:OperadEquivalence}
The Fulton-MacPherson operad $\Conf[\bullet]$ of configurations in $\BR^N$  and the little $N$-disks operad $\calB\oprd$
are weakly equivalent as topological operads.
\end{prop}

\q For the sake of keeping this paper as self-contained as possible,{\D}\previousoldfn{\iv Paolo apparently also proved this in his ``Configuration spaces with summable labels", Proposition 3.9., if we're looking for a published proof.  But it also seems to me that if we're trying to have a self-contained paper, we should provide the proof of this.  We already wrote it up in an earlier version of the paper (the latest version I have where this is written up is OperadFormality20November06.tex.).  This also might be good if we decided to submit to Memoirs; their lower bound for pages is 80 so we could get more pages by including this proof.  I'll leave it up to you.\pl I changed Paolo reference: actually Paolo told me that it is prop 4.9 (not as in arxiv). His proof is nicer than the one I wrote which was extremely pedestrian, so there is no reason to add it} we summarize Salvatore's proof here.
\begin{proof}[Summary of proof of \refP{OperadEquivalence}]
Recall the $W$ construction of Boardman-Vogt \cite{BoVo:HIA}\previousoldfn{\pl where? \pl Maybe no need to be more precise} which associates to a topological operad 
$\calO({\bullet})$ 
another operad $W\calO\oprd$   consisting of planar rooted trees $\tau$ whose internal edges have length between 0 and 1 and 
 whose internal vertices of valence $i+1$ are decorated by an element
 of $\calO(i)$. 
\p The operad $W\calO$ is a cofibrant replacement of $\calO$.
The main idea of the proof is then to construct a map $R\colon W\calB\oprd\to\Conf[\bullet]$ that sends a decorated tree $\tau$
 to the configuration 
of the centers of the configuration of balls obtained by multicomposition
of all the configurations of balls associated to the vertices of $\tau$ (after rescaling the configuration of balls 
at each internal vertex in a way that depends on the length of the adjacent edge, length $1$ corresponding to an infinitesimal rescaling).
It turns out that $R$ is a homotopy equivalence of operads and, since
 $W\calB\oprd$ is homotopy equivalent to $\calB\oprd$, this proves the proposition.
\end{proof}
In particular, the formality of the little balls operad will follow from that of the Fulton-MacPherson operad.

\subsection{The canonical projections}\label{sec:canproj}

Let $V$ be a finite set containing ${\ExtVert}$ as a subset. Set $I=V\setminus A$.
There is an obvious semi-algebraic map
\begin{equation}\label{eq:pi-canonical}
\pi\colon\Conf[V]\longrightarrow\Conf[{\ExtVert}]
\end{equation}
given by forgetting from the    configuration $y\in\Conf[V]$
all the points labeled by $I$.
This map $\pi$ can also be defined as an operad structure map. Indeed choose an arbitrary linear order on $V$
and consider the inclusion $\iota\colon {\ExtVert}\hookrightarrow V$ as a weak ordered partition.
For $v\in V$,  $\iota^{-1}(v)$ is either empty or a singleton $\{v\}$.
Since $\Conf[\emptyset]$ and $\Conf[\{v\}]$ are both one-point spaces,
the projection on the first factor 
\[
\proj\colon\Conf[V]\times\prod_{v\in V}\Conf[\iota^{-1}(v)]\iso\Conf[V]
\]
gives a homeomorphism
which we  use to identify these two spaces.
Then the operad structure map
\[\Conf[V]= \Conf[V]\times\prod_{v\in V}\Conf[\iota^{-1}(v)]
\stackrel{\Phi_\iota}\longrightarrow\Conf[{\ExtVert}]\]
is exactly the map  $\pi$.

\begin{defin}\label{D:canproj}
The  map $\pi\colon\Conf[V]\to\Conf[\ExtVert]$ of \refN{eq:pi-canonical} is called the
\emph{canonical projection (associated to the inclusion ${\ExtVert}\subset V$).}
\end{defin}
The Kontsevich configuration space integral will be defined through a pushforward of some minimal semi-algebraic forms along such canonical projections.  For this to be possible,  canonical projections have to be oriented SA  bundles (that is, semi-algebraic bundles whose fibers are compact oriented manifolds; see 
\citePADSAbundle):
\begin{thm}\label{T:projSAbdl}
Let  ${\ExtVert}$ be a finite set and let $I$ be a linearly ordered finite set disjoint from ${\ExtVert}$.
The canonical projection
\[\pi\colon\Conf[{\ExtVert}\cup I]\longrightarrow\Conf[{\ExtVert}]
\]
is an oriented SA bundle with fiber of dimension
\[\dim(\fiber(\pi))
\begin{cases}
=N\cdot |I|,&\textrm{if }|{\ExtVert}|\geq2\textrm{ or }I=\emptyset;\\
<N\cdot |I|,&\textrm{otherwise}.
\end{cases}
\]
\pfn{added because the referee wants to know about the fiber}
Assume moreover that $|\ExtVert|\geq2$. Then the fiber of $\pi$ is the
space of configurations of $|I|$ points in $\BR^N\setminus \ExtVert$
compactified by adding a boundary to this open manifold.
\begin{itemize}
\item
When $N$ is odd the orientation of the fiber of $\pi$ depends on the linear order
of $I$. A  transposition of that linear order reverses the orientation.
\item
When $N$ is even the orientation of the fiber is independent of the linear order
on $I$.
\end{itemize}

\end{thm}
\p For example, when $|I|=1$ and $|\ExtVert|\geq2$, the fiber of $\pi$
is a closed $N$-ball with $|\ExtVert|$ disjoint open balls removed
from its interior.

The proof of this theorem is not very difficult but it is long.
Since techniques used in the proof are not used anywhere else in the paper we decided to
delay it until \refS{proofpibdl}.
\previousoldfn{\pl check the the appendix indeed prove all of the theorem }
Notice however that although $\Conf[n]$ are smooth manifolds with corners, it is not true
that the canonical projections are smooth bundles, because their
restrictions to the boundary are usually not  submersions,
as shown{\D} in \refX{notsubm}.
This is the reason why we have to work with semi-algebraic forms instead of smooth forms.

Canonical projections can also be used to construct retractions to the operad structure maps associated to a \emph{non-degenerate} 
partition (see{\D}  \refD{partition}) as in the following easy-to-prove proposition and corollary.
\begin{prop}\label{P:retrPhi}
Let $\nu\colon A\to P$ be an ordered weak partition and set $A_p=\nu^{-1}(p)$ for $p\in P$ as in the setting \ref{setting:XP}.
For $q\in P$ denote by $\pi_q$ the canonical projection associated to the inclusion $A_q\subset A$.
Then the composition
\[\Conf[P]\times\prod_{p\in P}\Conf[A_p]\stackrel{\Phi_\nu}{\longrightarrow}\Conf[A]\stackrel{\pi_q}{\longrightarrow}\Conf[A_q]\]
is the projection on that factor.

Suppose moreover that $\nu$  is  non-degenerate, that is, it is surjective. Use any section of $\nu$ to identify
$P$ as a subset of $A$ and let $\pi_0$ be the associated canonical projection. Then 
the composition
\[\Conf[P]\times\prod_{p\in P}\Conf[A_p]\stackrel{\Phi_\nu}{\longrightarrow}\Conf[A]\stackrel{\pi_0}{\longrightarrow}\Conf[P]\]
is the projection on the first factor.
\end{prop}
\begin{cor}\label{C:retrPhi}
If $\nu\colon A\to P$ is a  non-degenerate ordered partition, then
the operad structure map
\[\Phi_\nu:\Conf[P]\times\prod_{p\in P}\Conf[A_p]\longrightarrow\Conf[A]\]
\p is injective and
admits a continuous semi-algebraic retraction.
\end{cor}
\begin{proof}
A retraction is given by $(\pi_p)_{p\in P^*}$ where $\pi_p$ is as in the previous proposition.
\end{proof}
This corollary is clearly wrong when the weak partition $\nu$ is degenerate.

\subsection{Decomposition of the boundary of $\Conf[n]$ into codimension $0$ faces}\label{sec:bdryCV}

In this section we  show that the boundary of $\Conf[n]$ decomposes as
the union of the images of certain operad structure maps.
Indeed, \refP{bdryCV} below gives a partition of $\partial C[n]$ (up to codimension $1$
intersections) whose pieces are images of ``$\circ_i$'' operations.
Most{\D} of the  operad structure on $\Conf[\bullet]$ can in fact be understood as an explicit
 decomposition of the boundary of $\Conf[n]$ as a union of faces homeomorphic to products of the form
$\Conf[k]\times\Conf[n_1]\times\dots\times\Conf[n_k]$. This is not
true for the nullary{\D} part though.
 
Let $V$ be a finite set.  We will study the boundary\previousoldfn{\pl the
  purpose of this paragraph is to add some intuition for the
  section before going to the gory details. Useful ? \iv Definitely.} of the manifold $\Conf[V]$. Recall that the
elements of that boundary{\D} are  characterized  in \refP{charbdryCX}.
For  a non-empty subset $W$ of $V$, we will consider the
configurations $y\in\Conf[V]$ such that the points $y(w)$ labeled by
$w\in W$ are infinitesimally{\D} closer to each other with respect to
any other point $y(v)$ labeled by $v\in V\setminus  W$. We will show
that these subsets of configurations give a decomposition of
$\partial\Conf[V]$
into codimension $0$ faces (\refP{bdryCV}) when $W$ runs over proper
subsets of cardinality $\geq2$.

For a non-empty subset \p $W\subset V$, let $V/W$ be the quotient set of $V$ in which all the elements
of $W$ are identified to a single element. In particular $|V/W|=|V|-|W|+1$.
Suppose given a linear order on $V/W$
and consider the projection to the quotient
\[q\colon V\longrightarrow V/W\]
which is an ordered non-degenerate partition of $V$.
One then has a structure map 
\[\Phi_{q}\colon\Conf[V/W]\times\prod_{p \in V/W}\Conf[{q}^{-1}(p )]\longrightarrow\Conf[V].
\]
Since ${q}^{-1}(p )$ is either a singleton $\{v\}$ or the subset $W$ and since
$\Conf[\{v\}]$ is a one-point space, we can identify the domain of $\Phi_{q}$ with
 $\Conf[V/W]\times\Conf[W]$. This defines a map
\begin{equation}\label{eq:defPhiW}
\Phi_W:=\Phi_q\colon\Conf[V/W]\times\Conf[W]\longrightarrow\Conf[V]
\end{equation}
that we will denote by $\Phi_W^V$ when we want to emphasize the set $V$.

In terms of operads, the map $\Phi_W$ corresponds to a  ``circle-$i$''
operadic operation $\circ_i$,
 up to some permutation. Indeed, when $V=\{1,\dots,n+k\}=\setn{n+k}$ and $W=\{i,\dots,i+k\}\cong\setn{k+1}$
then $V/W\cong\setn{n}$ and $\Phi_W$ is exactly
\[\circ_i\colon\Conf[{n}]\times\Conf[{k+1}]\longrightarrow\Conf[{n+k}].\]

The image of $\Phi_W$ consists of configurations in $\Conf[V]$
such that the points labeled by $W$ are infinitesimaly
close to each other compared to any  point labeled by $V\setminus W$.
This condition is empty when $V=W$ or when $W$ is a singleton; in other words for such a $W$ the image of $\Phi_W$ is all of
$\Conf[V]$. For proper subsets $W\subset V$ of cardinality $\geq2$, the image of $\Phi_W$ is in the boundary of $\Conf[V]$.
Actually, the next proposition shows that the images of all these $\Phi_W$ supply a decomposition of $\partial\Conf[V]$. The pieces of this
decomposition are indexed by the ``boundary faces" set
\begin{equation}\label{eq:defWV}
\BF(V):=\{W\subset V:W\not= V\textrm{ and }|W|\geq 2\}.
\end{equation}

\begin{prop}\label{P:bdryCV}\ \ 
\begin{enumerate}
\item [(i)] The boundary of $\Conf[V]$ decomposes as
$$\partial\Conf[V]=\underset{W\in\BF(V)}{\bigcup}\im(\Phi_W);$$
\item [(ii)] For $W\in \BF(V)$, \p
\[\dim(\im(\Phi_W))=N\cdot|V|-N-2=\dim(\partial\Conf[V]);\]
\item [(iii)] For $W_1\not= W_2$ in $\BF(V)$, 
$$\dim(\im(\Phi_{W_1})\cap \im(\Phi_{W_2}))<N\cdot|V|-N-2.$$
\end{enumerate}
\end{prop}
\begin{proof}
(i) By \refP{charbdryCX}, $\im(\Phi_W)\,\subset\,\partial\Conf[V]$  for $W\in \BF(V)$.
We will prove that the boundary is contained in the union of the images of the $\Phi_W$.
Let $y\in\partial\Conf[V]$. By \refP{charbdryCX} there exist
\previousoldfn{\pl Add this argument? ``(use for example \cite[Lemma 3.18]{Sin:man})'' \iv I don't quite see how that Lemma is relevant. \pl Resolved: the relevant is \refP{charbdryCX} }
 distinct elements $u_0,v_0,w_0\in V$
 such that
\[ y(v_0)\simeq y(w_0) \rel y(u_0).
\]
Set
\[W=\{w\in V: y(v_0)\simeq y(w) \rel y(u_0)\}.
\]
Then $v_0,w_0\in W$ and $u_0\in V\setminus W$, and hence $W\in \BF(V)$.
Consider the canonical projections
\[\pi_1\colon\Conf[V]\longrightarrow\Conf[(V\setminus W)\cup\{w_0\}]\cong\Conf[V/W]
\quad\textrm{ and }\quad
\pi_2\colon\Conf[V]\longrightarrow\Conf[W].\]
Then $y=\Phi_W(\pi_1(y),\pi_2(y))$. This proves (i).

(ii) For $W\in\BF(V)$, the map $\Phi_W$ is injective (by \refC{retrPhi}) and hence, by compactness,  it 
is a homeomorphism onto its image. Since $|W|\geq2$ and $|V/W|\geq2$,
 \refP{dimConf} implies that
\begin{align*}
\dim(\im\Phi_W) & =
\dim\Conf[V/W]+\dim\Conf[W] \\
& =(N\cdot|V/W|-N-1) +(N\cdot|W|-N-1 \\
&  =N\cdot|V|-N-2.
\end{align*}

(iii)
Let $W_1,W_2\in\BF(V)$ with $W_1\not=W_2$.
We consider three cases.
\begin{itemize}
\item[Case 1:] Suppose that $W_1\cap W_2=\emptyset$.
Then $\im(\Phi_{W_1})\cap \im(\Phi_{W_2})$ is the image
of the composition\previousoldfn{\pl check\todo\pl ok}
\[
\Conf[(V/W_2)/W_1]\times\Conf[W_1]\times\Conf[W_2]\stackrel{\left(\Phi_{W_1}^{V/W_2}\right)\times\id}{\longrightarrow}
\Conf[V/W_2]\times\Conf[W_2]\stackrel{\Phi_{W_2}^V}{\longrightarrow}
\Conf[V]
\]
and an analogous computation as in (ii) implies that this image is of dimension $N\cdot|V|-N-3$.
\item[Case 2:] Suppose that $W_1\subset W_2$ (or the other way around).
Then  $\im(\Phi_{W_1})\cap \im(\Phi_{W_2})$ is the image
of the composition
\[\Conf[V/W_2]\times\Conf[W_2/W_1]\times\Conf[W_1]
\stackrel{\id\times\left(\Phi_{W_1}^{W_2}\right)}{\longrightarrow}
\Conf[V/W_2]\times\Conf[W_2]
\stackrel{\Phi_{W_2}^V}{\longrightarrow}
\Conf[V]
\]
and again this image is of dimension $N\cdot|V|-N-3$.
\item[Case 3:] Suppose that  $W_1\cap W_2\not=\emptyset$, $W_1\not\subset W_2$,  and
 $W_2\not\subset W_1$.
Choose $a\in W_1\cap W_2$, $b\in W_1\setminus W_2$, and
 $c\in  W_2\setminus W_1$.
For $y\in\im(\Phi_{W_1})\cap \im(\Phi_{W_2})$
we simultaneously have 
\[
y(a)\simeq y(b)\rel y(c)\quad\textrm{ and }\quad
y(a)\simeq y(c)\rel y(b),\]
which is impossible. Thus $\im(\Phi_{W_1})\cap \im(\Phi_{W_2})$ is empty.
\end{itemize}
\end{proof}

More generally the operad structure maps 
\[\Phi_\nu\colon\Conf[k]\times\Conf[n_1]\times\dots\times\Conf[n_k]\longrightarrow\Conf[n]\]
 map homeomorphically to faces  of codimension $(k-2)$ in the boundary $\partial\Conf[n]$
when $2\leq k<n$, $n=n_1+\dots+n_k$,  and $n_1,\dots,n_k\geq1$. 
This in fact gives a complete stratification of that boundary, but we will not use this fact.
However, when $n_i=0$ for some $1\leq i\leq k$, then $\Phi_\nu$ is not an inclusion,
and in this case the study of $\Phi_\nu$ can require a more careful treatment 
as will be the case for example in \refS{proj-oper}.

\subsection{Spaces of singular configurations}\label{sec:singconf}

\previousoldfn{\iv  I guess maybe I prefer ``partial configuration spaces", but
  this might be too hard to change in the entire paper.  This is
  because they're compactified hyperplane arrangements where only some
  of the diagonals in a  product of $\BR^N$'s are taken out; so the
  diagonals are ``partially taken out".  Also can think of this in
  terms of the original Fulton-MacPherson compactification as some,
  but not all, diagonals being blown up (maybe this original
  construction should be mentioned somewhere, like after \refD{FM}?).
  But we can leave this as it is, it's fine.\pl So I suggest to keep
  singular configurations then but see footnote below}

\begin{remark}
\previousoldfn{\pl Is this paragraph well located? the idea is to discourage a bit
  the reader to read the following in a first reading. Maybe the
  paragraph should be at the beginning of the next section instead of
  the end of this one? \iv I moved it.} 
This and the next four\previousoldfn{check} sections discuss some of the more technical
properties of the Fulton-MacPherson operad which will be needed for the
corresponding technical parts of the proof of the properties of the Kontsevich
configuration space integral
 in \refS{KCSI}. The reader can thus{\D} safely skip Sections
 \ref{sec:singconf}-\ref{sec:proofpibdl} for the time being and jump
 to \refS{diagrams}\p,
except for the notion of \emph{condensation} in \refD{loc} which is necessary for defining the cooperad structure on the space of diagrams in \refS{diagcoop}.
\end{remark}

At times we will need to consider variations of the configuration spaces $\Conf[V]$
in which  some components of a configuration are allowed to coincide exactly, that is, without extra infinitesimal information to
distinguish the points. The goal of this section is to make this situation precise.

Let $A,I_1, I_2$ be disjoint finite sets.
Set $V_i=A\cup I_i$ for $i=1,2$ and $V=A\cup I_1\cup I_2$. Hence we have a pushout of sets
$V=V_1\cup_AV_2$.  
Consider the following pullback
where $\pi_1$ and $\pi_2$ are  canonical projections:
\begin{equation}\label{eq:pbconfsing}\xymatrix{
\Confsing[V_1,V_2]\ar[r]^-{q_1}\ar[d]_{q_2}\ar@{}[rd]|-{\textrm{pullback}}&\Conf[V_1]\ar[d]^{\pi_1}\\
\Conf[V_2]\ar[r]_{\pi_2}&\Conf[A].}
\end{equation}
Intuitively,
$\Confsing[V_1,V_2]$ can be seen as a compactified singular space of configurations\previousoldfn{\pl it seems to me that it should be ``space of singular configurations'' instead: it is the configuration which is singular, not the space; no?} of points in $\BR^N$ labeled by $v\in V$.
By ``singular" we mean that, for a configuration $y$, the component $y(i_1)$ labeled by $i_1\in I_1$
 may coincide with another component $y(i_2)$ labeled by  $i_2\in I_2$.

Since $V_i\subset V$, we have for $i=1,2$ the canonical projections
\[
\rho_i\colon\Conf[V]\longrightarrow\Conf[V_i].
\]
As $\pi_1\rho_1=\pi_2\rho_2$, we have a surjective map 
\begin{equation}\label{eq:rhoCsing}
\rho\colon\Conf[V]\longrightarrow \Confsing[V_1,V_2]
\end{equation}
to the pullback induced by $(\rho_1,\rho_2)$.
Intuitively, when $y(i_1)$ and $y(i_2)$ are infinitesimally close in $y\in\Conf[V]$, $\rho(y)$ is the singular configuration in which we forget the infinitesimal data associated to those components.

Consider the canonical projections $\pi\colon\Conf[V]\to\Conf[A]$ and
$\pi_{V_1}\colon\Conf[V]\to\Conf[V_1]$,
 and the composition
\[\pi':=q_1\circ\pi_1=q_2\circ\pi_2\colon \Confsing[V_1,V_2]\longrightarrow\Conf[A].\]
%To fix orientations, suppose that $I_1$ and $I_2$ are equipped with linear orders
%and  equip $I_1\cup I_2$ with the linear order of $I_1\ordsum I_2$.
Recall the notation \p$\bbr{M}$ for semi-algebraic {\current}s from \refN{eq:bbrM},
\refN{eq:CSAf}, and \refN{eq:fundfiber}
in \refS{RHTSAS}.
\begin{lemma}\label{L:deg1sing}
There is a  commutative diagram
\[\xymatrix{
\Conf[V]\ar[rr]^-{\rho}\ar[rd]_{\pi_{V_1}}
&&\Confsing[V_1,V_2]\ar[ld]^{q_{1}}\\
&\Conf[{V_1}]
}
\]
where $\pi_{V_1}$ and $q_1$ are orientable SA bundles. If moreover $|V_1|\geq2$, 
then for each $x\in\Conf[V_1]$
\[\rho_*\left(\bbr{{\pi_{V_1}}^{-1}(x)}\right)=\pm\bbr{{q_1}^{-1}(x)}.
\]
In other words, $\rho$ induces a map of degree $\pm1$ between the fibers of $\pi_{V_1}$ and $q_1$.

Similarly there is a commutative diagram
\[\xymatrix{
\Conf[V]\ar[rr]^-{\rho}\ar[rd]_{\pi}
&&\Confsing[V_1,V_2]\ar[ld]^{\pi'}\\
&\Conf[{A}],
}
\]
and if $|A|\geq2$, then 
 $\rho$ induces a map of degree $\pm1$ between the fibers of $\pi$ and $\pi'$.
\end{lemma}
\begin{proof}
\refT{projSAbdl} states that canonical projections are oriented SA bundles, and hence so
are $\pi_{V_1}$ and $\pi_2$. Therefore $q_1$ is also an oriented SA bundle as the pullback of $\pi_2$ along $\pi_1$
\citePAPpbSAbdl.
%Hence $\pi'=q_2\pi_2$ is also an oriented SA bundle as a composition of bundles (\citePAPcompositebdl).
When $|V_1|\geq2$, the fiber ${\pi_{V_1}}^{-1}(x)$ of $\pi_{V_1}$ over any  $x\in \Conf[V_1]$ is a compact manifold whose interior can 
 be identified with the space of injections
\[\Inj(I_2,\BR^N\setminus V_1)=\{y\colon I_2\hookrightarrow \BR^N\setminus V_1\}
\]
where $V_1$ is seen as a fixed subset in $\BR^N$.
From the pullback \refN{eq:pbconfsing} the fiber of $q_1$ is the same as the fiber of $\pi_2$ whose interior can similarly
be  identified with
\[\Inj(I_2,\BR^N\setminus A).\]
Thus $\rho$ maps the interior of the fiber ${\pi_{V_1}}^{-1}(x)$ homeomorphically to a dense subset
of the fiber ${q_1}^{-1}(x)$, and hence induces a degree $\pm1$ map between the fibers
of $\pi_{V_1}$ and $q_1$.

The proof of the second part of the lemma is similar. 
%It is of degree $+1$ because  $\rho$ respects the
%orientations of the fibers.\previousplfn{\todo check}
\end{proof}

\subsection{Pullback of a canonical projection along an operad structure map }
\label{sec:proj-oper}

In \refS{KCSI} we will define the Kontsevich configuration space
integral $\IK$ along the lines of \refN{eq:IGamma} in the{\D} Introduction,  and will want to prove  that
it is a morphism of (almost) cooperads. Since this integral is defined using  
pushforward along a canonical projection,  we  need to investigate
the  pullback of a canonical projection along an operad structure map, as in Diagram \refN{eq:Gpb} below.
This is the aim of this section. The main results are \refP{GGl} (complemented by \refP{GlGm}) and
\refP{pirholambda}.
\previousoldfn{\pl this section needs to be reorganized to emphasize some proposition  \iv What proposition? \pl done}
This section is technical and is only needed{\D} in \refS{GIKalmostcoop}, except for the notion of \emph{condensation} in \refD{loc}, which, as mentioned before, is needed to define the cooperad stucture on the space of diagrams.

Throughout this section we fix a weak ordered partition $\nu\colon A\to P$ and set
 $$P^*=\{0\}\ordsum P,\ \  A_p=\nu^{-1}(p),\ \  \text{and}\  A_0=P$$  as in the setting \ref{setting:XP}. We also have
an associated operad structure map 
\[\Phi_\nu\colon\Conf[P]\times \underset{p\in P}{\prod}\Conf[{\ExtVert}_p]=
\underset{p\in P^*}{\prod}\Conf[{\ExtVert}_p]
\longrightarrow\Conf[{\ExtVert}]
\] from \refN{eq:Phinu}.
We also fix a linearly ordered finite set $I$ disjoint from $A$ and $P$ and set $V=A\cup I$.  Thus we can consider the canonical
projection \[\pi\colon\Conf[V]\longrightarrow\Conf[\ExtVert]
\]
 associated to $A\subset V$ as in \refN{eq:pi-canonical}.
The elements \p of $I:=V\setminus A$ will be called \emph{internal vertices}, the elements of $A$ \emph{external vertices},
and the elements of $V$ \emph{vertices}. As the case $|A|\leq1$ is somewhat degenerate and has to be treated separately,
we will always in this section assume that $|A|\geq2$.

\p Define $\Conf[V,\nu]$ as the pullback
\previousoldfn{\pl in all section change $\Phi'$ to $\Phi'_\nu$ and $\pi'$ into $\pi'_\nu$ \pl done}
\begin{equation}\label{eq:Gpb}
\xymatrix{\Conf[V,\nu]\ar[r]^-{\Phi'_{\nu}}\ar[d]_{\pi'_{\nu}}
\ar@{}[rd]|-{\textrm{{pullback}}}&
\Conf[V]\ar[d]^-{\pi}\\
\underset{p\in P^*}{\prod}\Conf[{\ExtVert}_p]\ar[r]_-{\Phi_\nu}&
\Conf[{\ExtVert}],
}
\end{equation}
where $\pi$ is the canonical projection \refN{eq:pi-canonical} and $\Phi_\nu$ is the operad structure
map \refN{eq:Phinu}.

The main goal of this section is to show that this pullback decomposes as a
union 
\begin{equation}\label{eq:CnuCl}
\Conf[V,\nu]=\underset{\lambda}{\bigcup} \Conf[V,\lambda]
\end{equation}
(\refP{GGl})
such that \D
the restrictions $\Phi'_\lambda:=\Phi'_{\nu}|\Conf[V,\lambda]$ are closely related to some operad structure maps $\Phi'_\lambda$
(\refP{pirholambda}).
Moreover \refN{eq:CnuCl} is ``almost'' a partition,
in the sense that the intersections $\Conf[V,\lambda]\cap\Conf[V,\mu]$ are of lower dimension for $\lambda\not=\mu$
 (\refP{GlGm}).

Let us first give a rough idea of how we will show this. To make it easier, let us temporarily make an additional assumption that
$\nu$ is non-degenerate  (that is, each $A_p$ is non-empty) and that 
$P$ contains at least two elements. 
In that case, the map $\Phi_\nu$ is the inclusion of some part of the boundary of
$\Conf[A]$. More precisely,  $\im(\Phi_\nu)$ consists of all
configurations $x\in\Conf[A]$ such that, for $a,b,c\in A$, if $\nu(a)=\nu(b)\not=\nu(c)$
then $x(a)\simeq x(b)\rel x(c)$. We will say that such a configuration $x\in\Conf[A]$ is \emph{$\nu$-condensed}.
 In other words, a configuration $x\in \im(\Phi_\nu)$
can be thought of as a family indexed by $p\in P$ of clusters of points, where the $p$-th cluster
consist of points $x(a)$ indexed by $a\in A_p=\nu^{-1}(p)$.
\pfn{I tried to add some example to this notion; this paragraph is
  tentative: should be carefully checked}
For example, the
configuration $x\in\Conf[6]$ from \refF{OperadMapOutput} in Section \ref{sec:opC} \q 
is $\nu$-condensed for the partition $\nu$ given at beginning of
Example \ref{Ex:operadicCA}.\p
% It is also $\mu$-condensed for the
% partition defined by $\mu(1)=\mu(2)=\mu(6)=\alpha$,
% $\mu(3)=\mu(5)=\beta$, and $\mu(4)=\gamma$. \pfn{it would be nice to have a picture
%   giving a $\lambda$-localized; or not? \iv I think this last sentence of the paragraph should be deleted.  I think it's too hard for the reader to follow this at this point (but it's ok I think to leave the reference to \refF{OperadMapOutput} since this might click for some readers).}

As $\Phi_\nu$ is an inclusion  (because of our extra assumption), the pullback $\Conf[V,\nu]$ is the subset of $\Conf[V]$
consisting of all configurations
$y\in\Conf[V]$ such that $x:=\pi(y)$ is $\nu$-condensed. Consider such a $y\in\Conf[V,\nu]$.
One can then look at the position of the 
 points $y(i)$, for $i\in I$, with respect to the various clusters of
 points $\{x(a):a\in A_p\}$, for \p  $p\in P$.
Such a point $y(i)$ could be inside or infinitesimally close to some  cluster indexed by $p\in P$,
 in which case we say that, for this configuration, $i$  is \emph{$p$-local};
 or $y(i)$ could be close to none of the clusters in which case we say that $i$ is \emph{global}.
These cases can be encoded by a function $$\lambda\colon I \longrightarrow P^*$$ with $\lambda(i)=p$ if 
$i$ is $p$-local, and $\lambda(i)=0$ if $i$ is global. It is natural to extend $\lambda$ to $V$
by letting $\lambda|A=\nu$. Such a map $\lambda\colon V\to P^*$ will be called a \emph{condensation} (\refD{loc} below),
 and there is a natural 
partition of $\Conf[V,\nu]$ as a union of the subspaces $\Conf[V,\lambda]$ consisting of $\lambda$-condensed
configurations $y$. Moreover, under our extra assumption, each $\Conf[V,\lambda]$ is homeomorphic to the product 
$\prod_{p\in P^*}\Conf[V_p]$ where $V_0=\lambda^{-1}(0)\cup P$ and $V_p=\lambda^{-1}(p)$ for $p\in P$,
and through this homeomorphism the restriction $\Phi'_{\nu}|\Conf[V,\lambda]$ is an operad structure map.

The precise description of the decomposition of $\Conf[V,\nu]$ is a bit more
delicate when the weak partition $\nu$ is degenerate, that is, when our extra assumption does not hold. We now proceed with the details and first define the notion of a condensation.
\begin{defin}\label{D:loc}%\ \ 
Let $A$ be a finite set, $\nu\colon A\to P$ be a weak ordered partition, $I$ be
a finite linearly ordered set disjoint{\D} from $A$, $P^*:=\{0\}\ordsum
P$, and $V:=A\amalg I$.  
\p Set $A_p=\nu^{-1}(p)$ for $p\in P$ and \p $A_0=P$.
 Elements of $V$ are called \emph{vertices} as 
above.
\begin{itemize}
\item
A \emph{condensation of \ $V$  relative to $\nu$}
is a map
\[\lambda\colon V\longrightarrow  P^*\]
such that $\lambda|{\ExtVert}=\nu$.
\item
The set of all such condensations $\lambda$ is denoted by $\Cond(V,\nu)$, or simply $\Cond(V)$
when $\nu$ is understood. 
\item Given a condensation $\lambda\in\Cond(V)$, 
a vertex
$v\in V$ is \emph{$p$-local} if $\lambda(v)=p$ for some $p\in P$, \p
and it is \emph{global} if $\lambda(v)=0$.
\item
A  configuration $y\in\Conf[V]$ is
\emph{$\lambda$-condensed} if for each $u,v,w\in V$ and $p\in P$ such
that
 $u$ and $v$ are $p$-local and $w$ is not $p$-local we have 
$y(u)\simeq y(v)\rel y(w)$.
\item
A condensation $\lambda\in\Cond(V,\nu)$  is \emph{essential}
\previousoldfn{was named \emph{normal} before }
if for each $p\in \lambda(I)$ we have that  $|{\ExtVert}_p|\geq2$.
We denote the set of essential condensations by $\EssCond(V,\nu)$, or simply $\EssCond(V)$
when $\nu$ is understood. 
\end{itemize}
\end{defin}
\p The terminology \emph{condensation} comes  from the idea
that a $\lambda$-condensed configuration $x\in\Conf[V]$
consists of clusters of points condensed together according of the
values
of $\lambda$ on their vertices. 

It is easy to convince oneself that a configuration $y\in\Conf[V]$ is $\lambda$-condensed if and only if
 it is in the image of an operad structure map $\Phi_{\widehat{\lambda}}$,
 where $\widehat{\lambda}$ is some weak partition of $V$ constructed from $\lambda$
(see \refN{eq:lambdahat} and \refN{eq:Phihatl} below for definitions of $\widehat{\lambda}$ and $\Phi_{\widehat{\lambda}}$).

A condensation is essential if there are no internal $p$-vertices when $|A_p|\leq1$, $p\in P$, and
no global (internal) vertices when  $|P|\leq1$.
\p We will see latter that non-essential condensations are in some sense 
\p negligible. For example they are not needed in the decomposition of $ \Conf[V,\nu]$ in
\refP{GGl} below\p and their contribution to the Kontsevich configuration space integral is zero as
we will see in \refL{Inorm}.

%When there are no internal vertices, i.e.~ $I=\emptyset$ hence $V={\ExtVert}$, then $\nu$ itself 
%is an essential condensation of ${\ExtVert}$ relative to $\nu$.  The proof of the following is immediate.
%\begin{lemma}\label{L:imPsi}
%$\im\Phi_\nu=\{x\in\Conf[{\ExtVert}]:x\textrm{ is $\nu$-condensed}\}.$
%\end{lemma}
%\previousplfn{\todo check. Do we need this lemma? Maybe not. In any case it is a corollary of the next one (check) so
%it should be stated after.}

For $\lambda$  a condensation of $V$ relative to $\nu$,
set
\begin{equation}\label{eq:Gl}
\Conf[V,\lambda]:=\{g\in \Conf[V,\nu]:\,\Phi'_{\nu}(g)\textrm{ is $\lambda$-condensed}\}
\end{equation}
where $\Conf[V,\nu]$ and $\Phi'_{\nu}$ are from \refN{eq:Gpb}.

Recall that in this section we assume $|A|\geq2$.
Our first important result is the following decomposition of the pullback $\Conf[V,\nu]$.
\begin{prop}\label{P:GGl}  There is a decomposition
$$\Conf[V,\nu]=\underset{\lambda\in\EssCond(V,\nu)}{\bigcup}\Conf[V,\lambda]$$
where $\lambda$ runs over all essential condensations relative to $\nu$.
\end{prop}
\begin{proof}
\p Recall the pullback $\Conf[V,\nu]$ of diagram \refN{eq:Gpb} and
let $g=(y,(x_p)_{p\in P^*})\in \Conf[V,\nu]$ with $y\in\Conf[V]$, $x_p\in\Conf[A_p]$,
and $\pi(y)=\Phi_\nu\left((x_p)_{p\in P^*}\right)$.
We need to construct an essential condensation $\lambda$ such that $g\in \Conf[V,\lambda]$.
For $i\in I$ and $p\in P$ we say that $i$ is \emph{$p$-local for $g$} if 
\begin{itemize} 
\item[(i)]$|A_p|\geq2$, and
\item[(ii)] $\forall a,b\in A:(\nu(a)=p\textrm{ and }\nu(b)\not=p)\Longrightarrow\left(y(a)\simeq y(i)\rel y(b)\right)$.
\end{itemize}
If $i$ is $p$-local, then it cannot be $q$-local for $q\not=p$ because otherwise there would exist $a\in A_p$ and $b\in A_q$ 
(since $|A_p|,|A_q|\geq2$),  with both $y(a)\simeq y(i)\rel y(b)$ and
$y(b)\simeq y(i)\rel y(a)$, which is impossible.

Define a condensation $\lambda\colon V\to P^*$ by
\[\lambda(v)=
\begin{cases}
  \nu(v),&\textrm{ if }v\in A;\\
  p,     &\textrm{ if $v\in I$ and $v$ is $p$-local for $g$, for some $p\in P$;}\\
  0,     &\textrm{ if $v\in I$ and there is no $p\in P$ for which $v$ is $p$-local for $g$.}
\end{cases}
\]
Let us show that $\lambda$ is essential. If $p\in\lambda(I)\cap P$ then $|A_p|\geq2$
by condition (i). If $0\in\lambda(I)$ then $|P|\geq2$ because otherwise $P$ is a singleton $\{p_1\}$
(as $P$ cannot be empty since $|A|\geq2$), in which case every $i\in I$ is $p_1$-local (except if $|A_{p_1}|<2$ which is 
again impossible
since we assume in this section that $|A|\geq2$).
Therefore $\lambda$ is an essential condensation relative to $\nu$.

We now show that $g\in\Conf[V,\lambda]$. Let $u,v,w \in V$ such that $\lambda(u)=\lambda(v)=p\not=0$ and
$\lambda(w)\not=p$. We need to show that
\begin{equation}\label{yuvwloc}
y(u)\simeq y(v)\rel y(w).
\end{equation}
If $|A_p|\leq1$ then $u,v\in A_p$ since no internal vertex can be $p$-local due to (i).
Hence $u=v$ and \refN{yuvwloc} is obvious. Suppose now that $|A_p|\geq2$.
As $\lambda(w)\not=p$, there are two cases:
\begin{itemize}
\item[(A)] $w\in A_q$ for some $q\not=p$, or
\item[(B)] $w\in I$ and there exists $a\in A_p$ and $b\in A_q$ for some $q\not=p$ such that\previousoldfn{\pl too big
vertical space... \iv We can let the journal people deal with that if it's still a problem after they reformat the paper.}
\[y(a)\not\simeq y(w)\rel y(b).\]
\end{itemize} 
In case (A) we can pick $a\in A_p$ and, whether $u\in I$ or $u\in A_p$, we have $y(a)\simeq y(u)\rel y(w)$. Similarly
 $y(a)\simeq y(v)\rel y(w)$. By transitivity we
get \refN{yuvwloc}.

In case (B) we have $y(a)\not\simeq y(w)\rel y(b)$. Since $y(a)\simeq y(u)\rel y(b)$ we deduce that
$y(a)\simeq y(u)\rel y(w)$. Similarly $y(a)\simeq y(v)\rel y(w)$. Again by transitivity we get \refN{yuvwloc}.\previousoldfn{where do we use in the 
prood the fact that $|A|\geq2$? \pl Done: the hypothesis $|A|\geq2$ used}
\end{proof}

\p The various $\Conf[V,\lambda]$ that appear in the union of the previous proposition are not necessarily
 pairwise disjoint. However,
 their intersection is of  positive codimension as we will see in \refP{GlGm}.\previousoldfn{\pl will this lemma still be there? \pl Yes}
We have assumed in this section that $|A|\geq2$; when $|A|<2$ it is possible that there are
 no essential condensations at all, in which case the decomposition
from the last proposition cannot hold.
\vspace{5mm}

%before

Let \p$\lambda\colon V\to P^*=P\cup\{0\}$ be an essential condensation.
Our next goal is to show that the restriction of $\Phi'_{\nu}\colon \Conf[V,\nu]\to \Conf[V]$ to $\Conf[V,\lambda]$
is closely related to a map
\[\Phi'_{{\lambda}}\colon\prod_{p\in
  P^*}\Conf[V_p]\longrightarrow\Conf[V],\]
\p where $V_p=\lambda^{-1}(p)$ for $p\in P$, 
$V_0=\lambda^{-1}(0)\cup P$,
and the map \p$\Phi'_\lambda$ can be identified with an explicit operad structure map $\Phi_{\widehat\lambda}$.  Here
 $\widehat\lambda$ is some refined partition of $\lambda$.
 In short, this amounts to saying that $\lambda$-condensed configurations are exactly the image of a certain operad structure map. This is the content of \refP{pirholambda}.
\previousoldfn{\pl I was thinking to summarize the section into a diagram; but it is to large to fit the page. A try is made in
the latex file just after this footnote. But actually I guess that this summarizing diagram is unneeded.  \iv I rearranged and finished the diagram, so we can keep it if you want.  Hopefully I got it all.  Constructing it really helped me get through this section.\pl should we keep the diagram? \iv Sure, it helps me as a reference for where to find definitions of maps and results.}

To prove this, we will need to construct various maps that are collected in the following diagram for reader's convenience, along with numbers of equations where they can be found.  The two identifications  are due to the fact that $\prod_{i\in I_0}\Conf[\{i\}]$ is a{\D} one-point space.

%
%%%------ large diagram summarizing this section ------- (still mistakes in the diagram + too large!!!)
%\newpage
%\begin{landscape}
%  
% \[\xymatrix{
% \Conf[V_0]\times\underset{i\in I_0}{\prod}\Conf[\{i\}]\times\underset{p\in P}{\prod}\Conf[V_p]
% \ar@{=}[d]^{\sim}
% \ar[r]
% \ar[rrrd]^-{\Phi_{\widehat{\lambda}}}
% &&
% \Conf[V_0^+]\times\underset{i\in I_0}{\prod}\Conf[\{i\}]\times\underset{p\in P^+}{\prod}\Conf[V_p]
% \ar@{=}[rr]^-{\sim}
% \ar[rrd]^{\Phi_{\lambda^+}}
% &&
% \Conf[V_0^+]\times\underset{p\in P^+}{\prod}\Conf[V_p]
% \ar[d]^{\Phi'_{\lambda^+}}
% \\
% \Conf[V_0]\times\underset{p\in P}{\prod}\Conf[V_p]
% \ar[r]_{\rho_0\times\id}
% \ar[rd]_-{\pi_\lambda=\times_{p\in P}\pi_p}
% \ar[rr]^+{\rho_\lambda}
% \ar[rrrr]^+{\Phi'_\lambda}
% &
% \Confsing[V_0^+,P]\times\underset{p\in P}{\prod}\Conf[V_p]
% \ar[r]^{\cong}_{\Phi''_\lambda}
% &
% \Conf[V,\lambda]
% \ar[ld]^{\pi'_\lambda}
% \ar^{^(->}[r]
% &
% \Conf[C,\nu]
% \ar[r]_{\Phi'_\nu}
% \ar[d]_{\pi'_\nu}
% \ar@{}[rd]|-{\textrm{pullback}}
% &
% \Conf[V]
% \ar[d]^{\pi}
% \\
% &
% \underset{p\in P^*}{\prod}\Conf[A_p]
% \ar@{=}[rr]
% &&
% \underset{p\in P^*}{\prod}\Conf[A_p]
% \ar[r]_{\Phi_\nu}
% &
% \Conf[A]
% }
% \]
%\end{landscape}
%\newpage

\[\xymatrix{
 \Conf[V_0]\times\underset{i\in I_0}{\prod}\Conf[\{i\}]\times\underset{p\in P}{\prod}\Conf[V_p]
 \ar@{=}[d]^{\id}
 \ar[r]^-{\text{proj}}
 \ar
 %@/_18pc/
 %@<-4ex>
 `l[d] `[ddddd] `[r]
 ^-{\Phi_{\widehat{\lambda}}}_-{\eqref{eq:Phihatl}}
  [ddddr]
 &
 \Conf[V_0^+]\times\underset{i\in I_0}{\prod}\Conf[\{i\}]\times\underset{p\in P^+}{\prod}\Conf[V_p]
 \ar@{=}
 @/^3pc/
 @<8ex>
 [dd]^(0.3){\id}
 \ar `r[d] `[ddd] [ddd]
 _-{\Phi_{\lambda^+}}^-{\text{\tiny operad map}}
 \\
 \Conf[V_0]\times\underset{p\in P}{\prod}\Conf[V_p]
 \ar[r]^(0.4){\rho_0\times\id}_(0.4){\eqref{eq:rhoV0+}}
 \ar
 @/_4pc/
 @<-1ex>
 [ddd]^(0.3){\eqref{diag:pirholambda}}_(0.35){\pi_\lambda}
 \ar[d]^-{\eqref{diag:pirholambda}}_-{\rho_\lambda}
 \ar[ddr]^(0.55){\Phi'_\lambda}_(0.55){\eqref{eq:Phi'l}}
 &
 \Confsing[V_0^+,P]\times\underset{p\in P}{\prod}\Conf[V_p]
 \ar[dl]^(0.4){\ \ \cong \ (Lm.\, \ref{L:GlimPsi})}_(0.45){\Phi''_\lambda}
 &
 \\
 \Conf[V,\lambda]
 \ar
 @/_2pc/
 %@<-1ex>
 [dd]^(0.2){\eqref{diag:pirholambda}}_(0.3){\pi'_\lambda}
 \ar
 @/^2.5pc/
 @<-1ex>
 [d]^-{\eqref{eq:Gl}}
 &
 \Conf[V_0^+]\times\underset{p\in P^+}{\prod}\Conf[V_p]
 \ar[d]_-{\Phi'_{\lambda^+}}^-{\eqref{eq:philpplus}}
 \\
 \Conf[V,\nu]
 \ar[r]_-{\Phi'_\nu}
 \ar[d]^-{\pi'_\nu}
 \ar@{}[rd]|-{\textrm{pullback \eqref{eq:Gpb}}}
 &
 \Conf[V]
 \ar[d]^-{\pi}
 \\
 \underset{p\in P^*}{\prod}\Conf[A_p]
 \ar[r]_-{\Phi_\nu}
 &
 \Conf[A]
 \\
  & & 
 }
 \]

First we  show that $\Conf[V,\lambda]$ is homeomorphic to the product of configuration spaces 
$\Conf[V_p]$, $p\in P$, and another, maybe singular, configuration space $\Confsing[V^+_0,P]$.
Let us construct the $V_p$'s.
\p Recall that $A_p=\nu^{-1}(p)=A\cap\lambda^{-1}(p)$ for $p\in P$,
$A_0=P$, and $V=A\amalg I$.
For  $p\in P^*$,
set
\begin{equation}
\label{eq:IpVp}
I_p=I\cap\lambda^{-1}(p) \quad\textrm{ and }\quad V_p={\ExtVert}_p\cup I_p,
\end{equation}
\p so $V_p=\lambda^{-1}(p)$ for $p\in P$, and
$V_0=\lambda^{-1}(0)\cup P$.
\previousoldfn{\iv Wow, there's so much notation in these last two sentences, it took me a long time to trace it all back and figure out what's going on.  This is the part that will drive the reader crazy, but it's actually good that's it's all collected here. }
The linear order of $I$ restricts to linear orders on $I_p$ for $p\in P^*$.
Moreover we order  $V_0$ as\previousoldfn{\pl This order $V_0= I_0\ordsum P$ is strange; it should be in the other ways since
$P$ are external vertices. \iv I see what you mean.  Is this a
problem?  Seems ok to me.} (remember that $A_0=P$)
\[V_0= I_0\ordsum P.\]

Also define the subsets
\[P^+:=\{p\in P:A_p\not=\emptyset\}\,\subset\, P\]
and 
\[V^+_0:=I_0\cup P^+\,\subset\, V_0.\]
Hence we have a pushout of sets $V_0=V^+_0\cup_{P^+}P$.
Consider the following pullback 
\begin{equation}\label{eq:pbV0+P}
\xymatrix{
\Confsing[V^+_0,P]\ar[r]^-{q_{V_0^+}}\ar[d]_{q_P}
\ar@{}[rd]|-{\textrm{pullback}}
&
\Conf[V^+_0]\ar[d]^{\pi^+}\\
\Conf[P]\ar[r]_{\pi'^+}&\Conf[P^+],}
\end{equation}
where $\pi^+$ and $\pi'^+$ are the canonical projections.  This defines a singular configuration space as in
\refS{singconf}.
When $\nu$ is   non-degenerate 
 then
$P^+=P$, $V^+_0=V_0$, and $\Confsing[V^+_0,P]$ is just the
configuration space $\Conf[V_0]$. In any case we have an induced map,
as in \refN{eq:rhoCsing},
\begin{equation}\label{eq:rhoV0+}
\rho_0\colon\Conf[V_0]\longrightarrow\Confsing[V^+_0,P],
\end{equation}
%over $\Conf[P]$, 
which, by \refL{deg1sing}, induces  a degree $\pm1$ map between the fibers.\previousoldfn{do we use at all this remark?}

Define the  weak ordered partition
\begin{align}
\widehat{\lambda}\colon  V & \longrightarrow V_0 \label{eq:lambdahat}\\
v & \longmapsto\widehat{\lambda}(v)=
\begin{cases}
v,&\textrm{if }\lambda(v)=0;\\
\lambda(v),&\textrm{otherwise.}
\end{cases}\notag
\end{align}
There is an associated operad structure map
\begin{equation}\label{eq:Phihatl}
\Phi_{\widehat{\lambda}}\colon\Conf[V_0]\times
\left(\prod_{i\in I_0}\Conf[\{i\}]\times\prod_{p\in P}\Conf[V_p]\right)
\longrightarrow\Conf[V].
\end{equation}
Since $\Conf[\{i\}]$ are one-point spaces,
the domain of $\Phi_{\widehat{\lambda}}$
is homeomorphic to  $\prod_{p\in P^*}\Conf[V_p]$
(through the obvious projection),
and the composition of this homeomorphism with $\Phi_{\widehat{\lambda}}$ gives a map
\begin{equation}\label{eq:Phi'l}
\Phi'_\lambda\colon\prod_{p\in P^*}\Conf[V_p]\longrightarrow\Conf[V].
\end{equation}
We   next show that $\Phi'_\lambda$ factors through the composition of $\Phi'_{\nu}|\Conf[V,\lambda]$ with
 a homeomorphism $\Phi''_\lambda$ between $\Confsing[V_0^+,P]\times\prod_{p\in P}\Conf[V_p]$
and $\Conf[V,\lambda]$.

By definition $P^+=\im(\nu)$ and, since $\lambda$ is essential, $\im(\lambda)\subset P^+\cup\{0\}$.
Therefore the weak partition $\widehat\lambda$ factors as the composition of an ordered  non-degenerate partition 
\[\lambda^+\colon V\longrightarrow V_0^+\]
and the inclusion $V_0^+\hookrightarrow V_0$.

For $p\in P\setminus P^+$ we have $A_p=\emptyset$, and hence $V_p=\emptyset$ because $\lambda$ is essential,
and so $\Conf[V_p]=*$. Also $\Conf[\{i\}]=*$ for $i\in I_0$.
Thus{\D}  the projections induce a homeomorphism 
\[\Conf[V_0^+]\times\prod_{p\in P}\Conf[V_p]\cong \Conf[V_0^+]\times \prod_{i\in I_0}\Conf[\{i\}]\times\prod_{p\in P^+}\Conf[V_p].
\]
The composition of this homeomorphism with the operad structure map $\Phi_{\lambda^+}$ is a map
\begin{equation}\label{eq:philpplus}
\Phi'_{\lambda^+}\colon
\Conf[V_0^+]\times\prod_{p\in P}\Conf[V_p]\longrightarrow\Conf[V].
\end{equation}

Recall $q_{V_0^+}$ and $q_P$ from \refN{eq:pbV0+P}, let  
\[\pi_p\colon\Conf[V_p]\to\Conf[A_p]\]
 be the canonical projections,
 and consider the two maps
\[\Phi'_{\lambda^+}\circ (q_{V_0^+}\times\id)\colon\Confsing[V_0^+,P]\times\prod_{p\in P}\Conf[V_p]\longrightarrow\Conf[V]\]
and
\[q_P\times(\times_{p\in P}\pi_p)\colon\Confsing[V_0^+,P]\times\prod_{p\in P}\Conf[V_p]\longrightarrow\Conf[P]\times\prod_{p\in P}\Conf[A_p]=
\prod_{p\in P^*}\Conf[A_p].
\]
They induce a map 
\begin{equation}
  \label{eq:Phi''l}
\Phi''_\lambda\colon\Confsing[V_0^+,P]\times\prod_{p\in P}\Conf[V_p]\longrightarrow\Conf[V,\nu]  
\end{equation}
into the pullback \refN{eq:Gpb}.

\begin{lemma}\label{L:GlimPsi}
$\Phi''_\lambda$ of \refN{eq:Phi''l} is a homeomorphism onto $\Conf[V,\lambda]\subset\Conf[V,\nu]$.
\end{lemma}
\begin{proof}
Since the domain of $\Phi''_\lambda$ is compact, it is enough to prove that  $\Phi''_\lambda$ is injective and that 
its image is $\Conf[V,\lambda]$.
This is in fact not hard to see using the pictorial interpretations of virtual configurations. 
Here is a more formal proof.

\previousoldfn{\todo prove that $\im(\Phi''_\lambda)\subset\Conf[V,\lambda]$ \pl sort of done/obvious}
For  injectivity, 
let 
$$z=(z_0,(z_p)_{p\in P})\ \  \text{and}\ \  z'=(z'_0,(z'_p)_{p\in P})\in\Confsing[V_0^+,P]\times\underset{p\in P}{\prod}\Conf[V_p]
$$ be such that
$\Phi''_\lambda(z)=\Phi''_\lambda(z')$. 
Since $\pi'_{\nu}\circ\Phi''_\lambda=q_P\times\left(\underset{p\in P}{\times}\pi_p\right)$, we get that
\begin{equation}\label{eq:injGlimPhi1}
q_P(z_0)=q_P(z'_0).
\end{equation}
As $\lambda^+$ is a   non-degenerate  partition, by \refC{retrPhi} $\Phi'_{\lambda^+}$ is injective.
Since $\Phi'_{\nu}\circ\Phi''_{\lambda}=\Phi'_{\lambda^+}\circ(q_{V_0^+}\times\id)$, we deduce that
\begin{equation}\label{eq:injGlimPhi2}
(q_{V_0^+}(z_0),(z_p)_{p\in P})=(q_{V_0^+}(z'_0),(z'_p)_{p\in P}).
\end{equation}
From \refN{eq:injGlimPhi1} and \refN{eq:injGlimPhi2} we deduce that $z=z'$ using the pullback diagram \refN{eq:pbV0+P}.

The image of $\Phi_{\widehat\lambda}$ consists of $\lambda$-condensed
configurations, and{\D} hence $\im(\Phi''_\lambda)\subset\Conf[V,\lambda]$. Let us prove
 surjectivity. Let $$g=(y,(x_p)_{p\in P^*})\in\Conf[V,\lambda].$$ As $y\in\Conf[V]$ is $\lambda$-local it belongs in the image of
$\Phi'_\lambda$. Since $$\Phi'_\lambda=\Phi'_{\lambda^+}\circ(q_{V_0^+}\times\id)\circ(\rho_0\times\id),$$
 we can set 
 $$y=\Phi'_{\lambda^+}(z_0^+,(z_p)_{p\in P})\ \  \text{for some}\ \ 
(z_0^+,(z_p)_{p\in P})\in\Conf[V_0^+]\times\prod_{p\in P}\Conf[V_p].$$ Since $\pi(y)=\Phi_\nu((x_p)_{p\in P^*})$,
 using $\pi^+$ and $\pi'^+$ from \refN{eq:pbV0+P}, we deduce that $\pi^+(z_0^+)=\pi'^+(x_0)$ (this can be seen for example 
by factoring $\nu$ through a  non-degenerate  partition $\nu^+\colon A\to P^+$ and using
\refC{retrPhi}).
Set $z_0=(z_0^+,x_0)\in\Confsing[V_0^+,P]$. Then $g=\Phi''_\lambda(z_0,(z_p)_{p\in P})$.
\end{proof}

Define the product of canonical{\D} projections
\begin{equation}
\label{eq:pil}
\pi_\lambda:=\times_{p\in P^*}\pi_{p}\colon
\prod_{p\in P^*}\Conf[V_p]\longrightarrow\prod_{p\in P^*}\Conf[{\ExtVert}_p]
\end{equation}
and the restriction
\begin{equation}\label{eq:pi'l}
\pi'_\lambda:=(\pi'_\nu|\Conf[V,\lambda])\colon \Conf[V,\lambda]\longrightarrow\prod_{p\in P^*}\Conf[{\ExtVert}_p]
\end{equation}
where $\pi'_\nu$ is from \refN{eq:Gpb}.
Define also
\[\rho_\lambda:=\Phi''_\lambda\circ(\rho_0\times\id)\colon\prod_{p\in P^*}\Conf[V_p]\longrightarrow\Conf[V,\lambda].\]
So we get the following commutative diagram we{\D} are aiming for
\begin{equation}\label{diag:pirholambda}
\xymatrix{
\prod_{p\in P^*}\Conf[V_p]\ar[rr]^-{\rho_\lambda}\ar[rd]_{\pi_\lambda}
&&\Conf[V,\lambda]\ar[ld]^{\pi'_\lambda}\\
&\prod_{p\in P^*}\Conf[{\ExtVert}_p].
}
\end{equation}
\begin{lemma}\label{L:bundlepilambda}
$\pi_\lambda$ and $\pi'_\lambda$ are orientable SA bundles with fibers of dimension $N\cdot |I|$.
\end{lemma}
\begin{proof}
For $\pi_\lambda$, this is a direct consequence of \refT{projSAbdl} since $\pi_\lambda$ is a product of canonical projections.
Since $\lambda$ is essential,  $|A_p|\geq2$ or $I_p=\emptyset$ for each $p\in P^*$, which yields the formula for
the dimension of the fiber.

For $\pi'_\lambda$, the result comes from the homeomorphism
$\Phi''_\lambda$ (\refL{GlimPsi}) 
through which{\D} $\pi'_\lambda$ can be identified with
$q_P\times\underset{p\in P}\times\pi_p$, and from the fact that 
$q_P$ is also an oriented SA bundle as it is the pullback of the canonical projection  $\pi^+$ along 
$\pi'^+$ in Diagram \refN{eq:pbV0+P}. 
\end{proof}

We fix  the orientations of the fibers of $\pi_\lambda$ and $\pi'_\lambda$
as follows.
For the fibers of $\pi_\lambda$, we orient them as the product, in the linear order of $P^*$,
of the fibers
of $\pi_p$ oriented as in \refT{projSAbdl} with respect to the linear order of $I_p$ restricted from that of $I$.
For the fibers of $\pi'_\lambda$, they are connected codimension $0$ submanifolds of the fibers of $\pi'_\nu$, which are
canonically identified with the fibers of $\pi$ because of the pullback \refN{eq:Gpb}. We
orient then the fibers of $\pi'_\lambda$ by the orientation of the fibers of $\pi$ defined in \refT{projSAbdl}
from the given linear order on $I$.

We will see that $\rho_\lambda$ in \refN{diag:pirholambda} induces a change of orientation of the fibers according to the sign
\begin{equation}\label{eq:sigmaIl}
\sigma(I,\lambda):=(-1)^{N\cdot|S(I,\lambda)|}
\end{equation}
where %\previousfn{\pl rassembler ces signes}
\begin{equation}\label{eq:SIl}
S(I,\lambda):=\{(v,w)\in I\times I:v<w\textrm{ and }\lambda(v)>\lambda(w)\}.
\end{equation}

%The composition $\Phi'_\nu\circ\rho_\lambda$ is the map
%\[\Phi'_\nu\circ\rho_\lambda\,=\,\Phi'_\lambda\prod_{p\in P^*}\Conf[V_p]\to\Conf[V]\]
%which can be identified with the operad structure map $\Phi_{\widehat{\lambda}}$.

Recall the fundamental class of the fiber of an oriented SA bundle as
in \refN{eq:fundfiber}.
%Remember from \refS{RHTSAS} that to a compact oriented SA manifold $M$ of dimension $k$ one can associate a
%fundamental semi-algebraic chain $\bbr{M}\in \CSA_k(M)$ of degree $k$. This applies in particular to the fibers of an oriented SA bundle. 
The second main result of this section can then be summarized in the
following
\pfn{IMPORTANT. In the next proposition the item (i) had disappeared
  in your revision. But latter in
  the text there where probably references to that propositoion (i);
  also reference to (ii) should become to (i), and to (iii) should
  become to (ii). I decided to set the (i) back despite
  referees'comment (37). Please chack that there were no change to the
references of the various items of that prop. \iv There are only two instances in the paper where we refer to parts of this proposition and the reference is correct both times.}
\begin{prop}\label{P:pirholambda}
Let $\lambda\in\EssCond(V,\nu)$ be an essential condensation relative
to $\nu$ and consider Diagram \refN{diag:pirholambda} above.
\begin{itemize}
\item[(i)]
$\pi_\lambda$ and $\pi'_\lambda$ are oriented SA bundles
with fibers of dimension $N\cdot|I|$. 
\item[(ii)] $\rho_\lambda$ induces a map of degree $\sigma(I,\lambda)=\pm1$ between the fibers. More precisely, for $x\in\prod_{p\in P^*}\Conf[A_p]$,
\[\rho_{\lambda*}(\bbr{\pi_\lambda^{-1}(x)})=\sigma(I,\lambda)\cdot\bbr{{\pi'_{\lambda}}^{-1}(x)}\]
in $\CSA_*(\pi^{-1}_\lambda(x))$, where $\sigma(I,\lambda)=\pm1$ is defined in \refN{eq:sigmaIl}-\refN{eq:SIl}.
\item[(iii)]
The composition $\Phi'_\nu\circ\rho_\lambda$ is the map
\begin{equation}
\Phi'_\lambda\colon\prod_{p\in P^*}\Conf[V_p]\longrightarrow\Conf[V]
\end{equation}
of \refN{eq:Phi'l}
which can be identified with the operad structure map $\Phi_{\widehat{\lambda}}$ of \refN{eq:Phihatl}.
\end{itemize}
\end{prop}
\begin{proof}
(i) is \refL{bundlepilambda} with the orientations given right{\D}
after it.\p

For (ii)\p, remember that $\rho_\lambda=\Phi''_\lambda\circ(\rho_0\times\id)$ where $\Phi''_\lambda$ is a homeomorphism 
by \refL{GlimPsi}. When $|P|\geq2$,  \refL{deg1sing} implies that $\rho_0$ induces a map of degree $\pm1$
between the fibers over $\Conf[P]$, and when $|P|\leq1$, using that $\lambda$ is essential, $\rho_0$ is the identity map.
Hence in both cases $\rho_\lambda$ 
induces a map of degree $\pm 1$ between the fibers over 
$\underset{p\in  P^*}{\prod}\Conf[A_p]$.
We have fixed the orientations
so that the fibers of $\pi_\lambda$ are oriented according to the linear order of $\ordsum_{p\in P^*}I_p$ 
and the fibers of $\pi'_\lambda$ are oriented according to the linear order
of $I$. The number of transpositions needed to reorder $\ordsum_{p\in P^*}I_p$ as $I$ is exactly the cardinality
of $S(I,\lambda)$. So the sign of the degree of $\rho_\lambda$ on the fibers is a consequence of
the change of orientation rule in \refT{projSAbdl}.

For (iii)\p, the equation  $\Phi'_\nu\circ\rho_\lambda=\Phi'_\lambda$
follows from the construction of $\rho_\lambda$ and $\Phi'_\lambda$.
 The identification of that map with 
the operad structure map $\Phi_{\widehat{\lambda}}$ is through the canonical homeomorphism
\[\Conf[V_0]\times\left(\prod_{i\in I_0}\Conf[\{i\}]\times\prod_{p\in P}\Conf[V_p]\right)
\iso
\prod_{p\in P^*}\Conf[V_p]
\]
induced by the{\D} projection (see \refN{eq:Phihatl} and \refN{eq:Phi'l}).
\end{proof}

Recall  from \refP{GGl}  the decomposition 
\[\Conf[V,\nu]=\underset{\lambda\in\EssCond(V,\nu)}{\bigcup}\Conf[V,\lambda].\]
We just proved that the fibers $\pi'^{-1}_\lambda(x)$ of each  $$\pi'_\lambda\colon\Conf[V,\lambda]\to\prod_{p\in P^*}\Conf[A_p]$$ are of dimension $N\cdot|I|$.
Our next proposition can then be interpreted as saying that the pairwise intersections of the terms
in this union are of codimension $\geq1$.  In other words, the above
union is a partition of $\Conf[V,\nu]$ ``up to codimension $1$''.
\begin{prop}\label{P:GlGm}
If $\lambda\not=\mu$ in $\EssCond(V,\nu)$, then
for each $x\in\prod_{p\in P^*}\Conf[{\ExtVert}_p]$,
\[\dim\left(\pi'^{-1}_{\lambda}(x)\cap \pi'^{-1}_{\mu}(x)\right)<N\cdot|I|.
\]
\end{prop}
\begin{proof}
Let  $x\in\prod_{p\in P^*}\Conf[{\ExtVert}_p]$ and pick $v\in I$  such that $\lambda(v)\not=\mu(v)$.
For concreteness   suppose that $\lambda (v)=k$ where $k=\max(P)$.
Set $V_p=V\cap\lambda^{-1}(p)$, for $p\in P$, and $V_0=\lambda ^{-1}(0)\cup P$.
 Hence $v\in V_k$.

If $(y,x)\in\Conf[V,\lambda]\cap\Conf[V,\mu]$ then, as $y\in\Conf[V]$ is $\lambda$-local, 
\begin{equation}
\label{eq:proxl1l2:1}
\forall a\in A_k,\forall b\in A_p\textrm{ with }p\not=k: y(a)\simeq y(v)\rel y(b)
\end{equation}
and, as $y$ is also $\mu$-local, 
\begin{equation}
\label{eq:proxl1l2:2}
\forall a,a'\in A_k: y(a)\simeq y(a')\rel y(v).
\end{equation}
Consider the following diagram in which $\pi''$ and $\pi_v$ are products of canonical projections,
$\Phi_2\colon\Conf[2]\times\Conf[A_k]\times\Conf[\{v\}]\to\Conf[A_k\cup\{v\}]$ is an operad structure map,
and the upper left square is a pullback:
{\small
\[
\xymatrix{
E
\ar[rr]^-{\widehat{j}}
\ar[d]_-{\widehat{{\pi''}}}
\ar@{}[rd]|-{\textrm{{pullback}}} & &
\prod_{p\in P^*}\Conf[V_p]%\times\Conf[V_k]
\ar[r]^-{\rho_{\lambda }}
\ar[d]^-{\pi''}&
\Conf[V,\lambda]\ar@/^2pc/[ddl]^-{\pi'_\lambda}
\\
\prod_{p<k}\Conf[{\ExtVert}_p]\times\left(\Conf[2]\times\Conf[{\ExtVert}_k]\times\Conf[\{v\}]\right)
\ar@{^(->}[rr]^-{j:=\id\times\Phi_2}& &
\prod_{p<k}\Conf[{\ExtVert}_p]\times\Conf[{\ExtVert}_k\cup\{v\}]
\ar[d]^-{\pi_v}
\\
& &
\prod_{p\in P^*}\Conf[{\ExtVert}_p].
}
\]
}
 
The proximity relations \refN{eq:proxl1l2:1} and \refN{eq:proxl1l2:2}
imply that
\[\Conf[V,\lambda]\cap\Conf[V,\mu]\subset\im(\rho_\lambda\circ \widehat{j}).\]
Therefore
\begin{eqnarray*}
\dim\left(\pi'^{-1}_{\lambda}(x)\cap \pi'^{-1}_{\mu}(x)\right)
&=&\dim(\fiber(\pi'_{\lambda })\cap \Conf[V,\mu])\\
&\leq&\dim(\fiber(\pi_v\circ \pi'')\cap\widehat{j}(E))\\
&=&\dim(\fiber(\pi_v\circ j\circ\widehat{\pi''}))\\
&=&\dim(\fiber(\pi_v\circ j))+\dim(\fiber(\widehat{\pi''}))\\
&=&\dim(\Conf[2])+\dim(\fiber(\pi''))\\
&=&(N-1)+N\cdot(|I|-1)<N\cdot|I|.
\end{eqnarray*}
\end{proof}

\subsection{Decomposition of the fiberwise boundary along a canonical projection}
\label{sec:bdryCVpi}

We turn now to a fiberwise version of the decomposition of the boundary,
extending the results of \refS{bdryCV}. These results will be needed in
\refS{GIKchainmap} to prove that
the Kontsevich configuration spaces integral $\IK$ is a chain map.

Let ${\ExtVert}\subset V$ and consider the canonical projection
\[\pi\colon\Conf[V]\longrightarrow\Conf[{\ExtVert}]
\]
which is a bundle whose fibers are oriented compact manifolds by
\refT{projSAbdl} (\p which we will prove in  \refS{proofpibdl}). 
\c Recall from \refN{eq:pipartial} \p and \cite[Definition
8.1]{HLTV:RHTSAS}
the fiberwise boundary of an oriented SA bundle.
The {fiberwise boundary} of $\pi$ is\pfn{that change was wrong! \iv It's great that you added that example cited below.}
\begin{equation}\label{eq:defCVbdrypi}
\Conf^\partial[V]:=\bigcup_{x\in\Conf[\ExtVert]}\partial(\pi^{-1}(x))
\end{equation}
which is a closed subspace of $\Conf[V]$.
This space is \emph{not} the same as 
\[\partial\Conf[V]\quad\textrm{ or }\quad \underset{x\in
  \Conf[{\ExtVert}]}{\bigcup}\pi^{-1}(x)\cap\partial\Conf[V]\]
(see the example of $[0,1]\times[0,1]\to [0,1]$ right after \refN{eq:pipartial}\p).

We also consider the restriction map
\[\pi^\partial:=(\pi|\Conf^\partial[V])\colon\Conf^\partial[V]\longrightarrow\Conf[{\ExtVert}].\]

Recall from \refN{eq:defWV}  in \refS{bdryCV} the  set $\BF(V)$ indexing  the faces of $\partial\Conf[V]$ and
define 
\begin{equation}\label{eq:WVX}\BF(V,{\ExtVert})=\{W\in\BF(V):{\ExtVert}\subset W\textrm{ or }|W\cap {\ExtVert}|\leq1\}.
\end{equation}
The following is a fiberwise version of \refP{bdryCV}. 
%In it we use the maps $\Phi_W$ defined in \refN{eq:defPhiW}.
\begin{prop}
\label{P:bdryCVpi}  There is a decomposition
\[\Conf^\partial[V]=\underset{W\in\BF(V,{\ExtVert})}{\bigcup}\im(\Phi_W)\]
where $\Phi_W$ are the maps{\D} defined in \refN{eq:defPhiW} of \refS{bdryCV}.
\end{prop}
\begin{proof}
Recall that $\Conf({\ExtVert})$ is the interior of the compact manifold $\Conf[{\ExtVert}]$, that is
\[\Conf({\ExtVert})=\Conf[{\ExtVert}]\setminus\partial\Conf[{\ExtVert}].\]
Then
\[\Conf^\partial[V]\cap\pi^{-1}(\Conf({\ExtVert}))\,=\,
(\partial\Conf[V])\cap\pi^{-1}(\Conf({\ExtVert})).\]
\p Since $\Conf^{\partial}[V]$ is a bundle over $\Conf[\ExtVert]$ and
$\Conf[\ExtVert]=\overline{\Conf(\ExtVert)}$, we get that
\[\Conf^\partial[V]=\overline{\Conf^\partial[V]\cap\pi^{-1}(\Conf({\ExtVert}))}
=
\overline{(\partial\Conf[V])\cap\pi^{-1}(\Conf(\ExtVert))}\]
where by $\overline{E}$ we mean the topological closure of the subspace $E$.\previousoldfn{\pl add this notation in prerequ. \iv I don't think we need to.\pl ok}

For $W\in\BF(V)$, if ${\ExtVert}\not\subset W$ and $|W\cap {\ExtVert}|\geq2$, then $\pi(\im\Phi_W)\subset\partial\Conf[{\ExtVert}]$
because $W\cap {\ExtVert}\in\BF({\ExtVert})$ and $\pi(\im\Phi_W)$ is in the image of
\[\Phi_{W\cap {\ExtVert}}^{\ExtVert}\colon\Conf[{\ExtVert}/(W\cap {\ExtVert})]\times\Conf[W\cap {\ExtVert}]\longrightarrow\Conf[{\ExtVert}].\]
Therefore\previousoldfn{\pl check. \todo. It lacks details \pl ok},
using  \refP{bdryCV}(i),\pfn{is the proof clear enough?\iv Yes, I think so.}
\begin{eqnarray*}
\Conf^\partial[V]&=&
\overline{\partial\Conf[V]\cap\pi^{-1}(\Conf({\ExtVert}))}\\
&=&
\overline{\cup_{W\in\BF(V)}\im(\Phi_W)\cap\pi^{-1}(\Conf({\ExtVert}))}\\
&=&
{\cup_{W\in\BF(V,{\ExtVert})}   \overline{\im(\Phi_W)\cap\pi^{-1}(\Conf({\ExtVert}))}}\\
&=&
\cup_{W\in\BF(V,{\ExtVert})}\im(\Phi_W).
\end{eqnarray*}
\end{proof}

\subsection{Orientation of $\Conf[{\ExtVert}]$}\label{sec:orCA}

In this section,\previousoldfn{\pl this section is sort of technical and not too much interesting ealthough not difficult;
 maybe should be moved latter in the section with the other technivcal
 sections\pl done} we fix an orientation on   $\Conf[{\ExtVert}]$.  This will be important since we
will integrate over this manifold. The orientation will be canonical when $N$ is even and will depend
on a linear order on $A$ when $N$ is odd. We will also fix an orientation on the sphere $S^{N-1}$.

We  first review a few classical facts and fix our conventions about orientation:
\begin{itemize}
\item 
A codimension $0$ submanifold of an oriented manifold inherits that orientation;
\item 
Conversely, the orientation of a connected  manifold is determined by the orientation of any non-empty  codimension $0$ 
connected submanifold;
\item The product $M_1\times M_2$ of two oriented manifolds  has a canonical orientation.
Exchanging the factors preserves or
reverses that orientation according to the sign $$(-1)^{\dim(M_1)\cdot\dim(M_2)};$$
\item $\BR$, and hence
$\BR^N=\BR\times\dots\times\BR$, is equipped with the standard orientation;
\item When $M$ is an oriented smooth manifold and $\omega$ is a smooth differential form with compact support of maximal degree on $M$,
one can consider the integral $$\int_M\omega\in\BR;$$
\item 
The orientation of a non-empty connected smooth manifold $M$
corresponds to an equivalence class of a smooth differential form
$\alpha$  of \D
maximal degree  with connected non-empty bounded  non-vanishing set, so that $\int_M\alpha>0$;
\item 
We orient the boundary of a manifold so that the Stokes' formula holds without a sign, that is, 
\[\int_{\partial M}\omega=\int_Md\,\omega\]
for a smooth differential form $\omega$ with compact support and of maximal degree on the smooth oriented manifold $M$.
\end{itemize}

When $|\ExtVert|\leq1$ then $\Conf[\ExtVert]$ is a one-point space and we choose the positive orientation on it.
Suppose now that $|\ExtVert|\geq2$ and
suppose given a linear order on ${\ExtVert}$. We then  have  a natural orientation on the codimension $0$
submanifold
\[\Inj({\ExtVert},\BR^N)\subset\prod_{a\in {\ExtVert}}{\BR^N}\]
defined in \refN{eq:Inj}, where the product is taken in the linear order of $A$.
A transposition in the linear order of $A$ changes this orientation by a sign $(-1)^N$.

Set
\[\Inj_0({\ExtVert},\BR^N)=\left\{x\in\Inj({\ExtVert},\BR^N):\bary(x)=0\right\}.
\]
This is a manifold without boundary.
We have a diffeomorphism
\begin{eqnarray*}
\Inj_0({\ExtVert},\BR^N)\times\BR^N&\stackrel{\cong}{\longrightarrow}&\Inj({\ExtVert},\BR^N)\\
(x,b)&\longmapsto&x+b
\end{eqnarray*}
{\Z }defined by $(x+b)(a):=x(a)+b$, for $a\in {\ExtVert}$.
We fix the unique orientation on $\Inj_0({\ExtVert},\BR^N)$ for which the above diffeomorphism preserves the orientation.
Consider the codimension $0$ submanifold
\[\Inj_0^{\leq1}({\ExtVert},\BR^N):=\left\{x\in\Inj_0({\ExtVert},\BR^N):\radius(x)\leq1\right\}\subset\Inj_0({\ExtVert},\BR^N) 
\]
with the induced  orientation. This is a manifold with boundary and its boundary inherits the orientation.

Identifying $\Conf(A)$ with $\Inj_0^1(A,\BR^N)$ \Z from \refN{eq:defC(X)xi}, we have
\[
\Conf({\ExtVert})=\partial\Inj_0^{\leq1}({\ExtVert},\BR^N)
\]
and this defines our prefered orientation on $\Conf({\ExtVert})$, and hence on $\Conf[{\ExtVert}]$.

We orient the sphere $S^{N-1}$ so that the map
\[\theta_{a,b}\colon \Conf([\{a,b\}])\iso S^{N-1}
\]
from \refN{eq:defthetaab} is orientation-preserving when the set  $\{a,b\}$ is ordered by $a<b$.

Consider a permutation $\sigma\in \Perm({\ExtVert})$ of the set ${\ExtVert}$.
It induces an obvious automorphism $\Conf[\sigma]$ of the manifold $\Conf[{\ExtVert}]$.
We then have
\begin{prop}
\label{P:permorientCX}
For a permutation $\sigma$ of $\ExtVert$, the induced homeomorphism
$\Conf[\sigma]\colon\Conf[\ExtVert]\to\Conf[\ExtVert]$ is orientation-preserving or orientation-reversing according to the sign
\[(\sign(\sigma))^N
\]
where $\sign(\sigma)=\pm1$ is the signature of the permutation $\sigma$.
\end{prop}

\subsection{Proof of the local triviality of the canonical projections}\label{sec:proofpibdl}\q
\previousoldfn{\iv This really could have been another little separate paper,
  but I guess it's too late for that...  \pl Indeed}
The only aim of this long{\D} section\previousoldfn{\pl to be moved from the appendix to the end of the  the
  section on FultonMcPherson operad.} \p\q 
is to prove \refT{projSAbdl}, which asserts that the canonical projection
$$\pi\colon\Conf[V=A\amalg  I]\longrightarrow\Conf[A]$$ is a semi-algebraic oriented fiber bundle
with fibers of prescribed dimension. 
\p These fibers should be thought of as a compactification of the
configuration space of $|I|$ points in $\BR^N$ with $|V|$ points
removed.
In particular, when $I$ is a singleton, the fiber of $\pi$ is
homeomorphic to a closed ball $D^N$ with $|A|$ disjoint open balls
removed.

That the projection 
 $C(V)\to C(A)$ is a bundle is a classical result due to Fadell and Neuwirth \cite{FaNe:con}.
The proof for the compactified version is more technical because of
the existence of a boundary.
Note \Z 
%also 
that, although the spaces $\Conf[V]$ and $\Conf[A]$ are{\D} smooth
manifolds with corners, \p it is \emph{not true}  that $\pi\colon
\Conf[V]\to\Conf[A]$
is a always a smooth bundle since it is not necessarily a submersion as
the \p
following example shows.
\begin{eg}\label{X:notsubm}
%Here we will show that $\pi$ is not necessary a{\D} smooth bundle.
\p\pfn{the proof of the example was very much changed, and I guess is
  more clear now.  \iv I agree.}
We now show that $\pi\colon\Conf[4]\to\Conf[3]$ is not a smooth bundle. 
 Fix $N=1$, that is, consider configurations of $n$ points on the real line. In that
dimension, $\Conf[n]$ consists of $n!$ copies of the connected component
$\Conf^{\operatorname{incr}}[n]$ corresponding to configurations $(x_1,\dots,x_n)$ with
$x_1<\dots<x_n$. For the remainder of this example we only consider this
connected component but drop the superscript $\operatorname{incr}$ to simplify notation.

The space $\Conf[n]$ is exactly the Stasheff{\D} associahedron $K_{n-2}$ \cite{Sta:wha}. In
particular $\Conf[3]$ is homeomorphic to the interval $[0,1]$ and
$\Conf[4]$ is homeomorphic to a \pfn{planar is meaningless here.
  What we mean is ``plain'' or whathever in the sense of the pentagon
  with its interior; do we need to add an adjective?. \iv No, I don't think so.} pentagon. 
\p Label the $4$ points of a configuration  in $\Conf[4]$ by 
$V=\{a,b,c,d\}$ and set $A=\{a,b,c\}$.
The five vertices of
the pentagon are indexed by all
%possible ways of fully
possible way of
parenthesizing the product $abcd$ in the most refined way as $(ab)(cd)$,  $(a(bc))d$,  etc.
\ppfn{by fully I mean $(a(bc))d$ versus $(abc)d$ for
  example; maybe not correct term. }
\q  Each
parenthetisation  encodes the proximity relations of the points of the
configurations $(a,b,c,d)\in\Conf[4]$, as shown in Figure \ref{fig:associahedron}\ppfn{It would be better for self-contain-ness that this
  figure was drawn here instead of refering to Stasheff paper. Could
  you do that? \iv Done.  Let me know if you think it should be
  modified in any way. \pl either add some adverb like
  ``fully''  as explained in previous note. If this sounds too
  awkward, then maybe  keep in the Figure of
  the pentagon only
the labels of the five vertices but I guess the first option makes
more sense.\iv I see what you mean.  How about saying ``...possible way of parenthesizing the product $abcd$ in the most refined way as $(ab)(cd)$, ..." or ``...possible way of parenthesizing the product $abcd$ so that no more than two factors are in a single set of parentheses, as in $(ab)(cd)$, ..."  I think ``fully" works as well, so if these other options sound weirder, we can keep it.}. 
  
  \begin{figure}[h]
\input{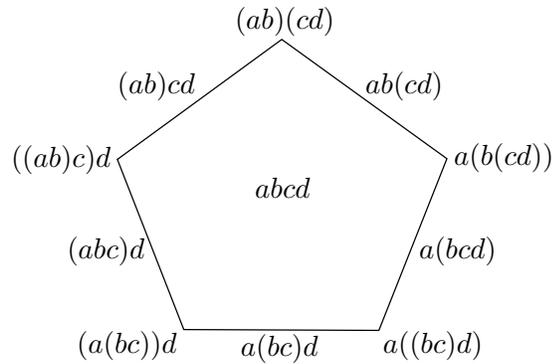}
\caption{Stasheff associahedron $K_2$ depicting the structure of $\Conf[4]$.}
\label{fig:associahedron}
\end{figure}

  Each of these five vertices corresponds to a point on the
boundary of $\Conf[4]$. For example, the vertex $a((bc)d)$ corresponds
to the limit, as $r\to 0^+$, of the configurations
\[(0,1-r-r^2,1-r,1)\,\in\,\Conf(4).\]
%Equivalently one could index the vertices by rooted trees with
%four leaves and two internal vertices as we will explain%
%below.
 Similarly $\Conf[3]$ is an interval whose endpoints are
labeled as $(ab)c$ and $a(bc)$.

\pfn{the following paragraphs  is I guess a clearer proof that in the previous
  version; the figure of the projection is probably not needed anymore (the end of the example
  was fully changed)}
\ppfn{your change in the description of $f$ were not mathematically
  correct, since $f$ does not take value in $\Conf(4)$ for
  example, and $(0,1-r-r^2,1-r,1)$ is not a correct description of an
  element of $\Conf[4]$ when
$r=0$ or $s=0$. What I meant is that $f$ is the unique continutation of
  the map described on $(0,1)\times(0,1)$) by the given formula. I
  tried to rephrase. Maybe still awkward.\iv No, I get it now.  This is good.}

A smooth chart of the manifold with corners $\Conf[4]$ about  the point
$a((bc)d)$
is given by the unique continuous map
\[
f\colon[0,1)\times[0,1)  \longrightarrow\Conf[4] \\
\]
whose restriction to $(0,1)\times(0,1)$ is defined by the map
\begin{align*}
(0,1)\times(0,1) & \longrightarrow\Conf(4)\subset\Conf[4] \\
(r,s)\quad&\mapsto(0,1-r,1-r+rs,1).
\end{align*}

Also, there is a chart
\begin{align*}
g\colon (0,1] & \longrightarrow\Conf[3] \\
\end{align*}
defined, for $0<t<1$, by
\[g(t)=(0,t,1)\in\Conf(3)\subset\Conf[3],\]
and extented continuously to $(0,1]$.

\q
\pfn{the sentence sounds awkward here, no?\iv I changed this and the next one a little but I'm not sure it's better.}
\ppfn{I corrected the domain of $f$ which should not contain $0$. Is this defintiion of $g$
  and $f$ clear?\iv Yes, this is good.}

We then have the following commutative diagram
of smooth maps between manifolds with corners
\[
\xymatrix{
[0,1)\times[0,1)\ar[r]^-{f}\ar[d]_p&\Conf[4]\ar[d]^\pi\\
(0,1]\ar[r]^-g&\Conf[3]
}
\]
where
\[p(r,s)=\frac{1-r}{1-r+rs}.\]
The partial derivatives of $p$ are
\begin{eqnarray*}
\frac{\partial p}{\partial r}(r,s)&=&\frac{s}{(1-r+rs)^2}\\
\frac{\partial p}{\partial s}(r,s)&=&\frac{-r(1-r)}{(1-r+rs)^2}.
\end{eqnarray*}
When $r=0$ and $s=0$, corresponding to the point $f(0,0)=a((bc)d)$,
both these partial derivatives are $0$, showing that $p$ is not a
submersion at
$(0,0)$. Hence $\pi$ is not a submersion at $a((bc)d)$. Therefore
$\pi$ is not a smooth bundle.

\end{eg}

\vspace{4mm}

We now come to the proof  of \refT{projSAbdl}.
The composition of two oriented SA bundles is
again an oriented SA bundle {\citePAPcompositebdl}, and therefore it is enough to prove that
\[\pi\colon\Conf[n+1]\longrightarrow\Conf[n]
\]
is an oriented SA bundle. For $n\leq1$, this is trivial, so we assume
that $n\geq2$.
In that case the fiber $F$ of $\pi$ will be homeomorphic to a disk
$D^N$ with $n$ disjoint open disks removed\p.

We first give a rough idea of the proof in an example. Take $n=9$ and consider the virtual configuration
$x_0\in\Conf[9]$ as in \refF{virtconf9} (see Example \ref{Ex:operadicCA} for an explanation of what such a figure represents).
\p We need to build some neighborhood $V$ of $x_0$ such that the
restriction of $\pi$ over $V$ is equivalent to the projection $V\times
F\to V$.

%\begin{figure}
%  \centering
%  \fbox{\includegraphics[width=95mm,height=60mm]{virtualconfig9.pdf}}
%  \caption{A virtual configuration $x_0\in\Conf[9]$}
%  \label{fig:virtconf9}
%\end{figure}

\begin{figure}[h]
\input{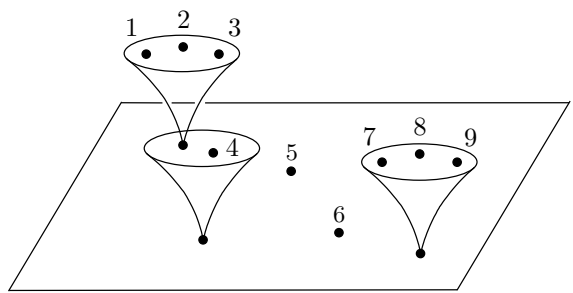}
\caption{A virtual configuration $x_0\in\Conf[9]$.}
\label{fig:virtconf9}
\end{figure}

For this configuration we have proximity relations
such as
\begin{align*}
&  x_0(1)\simeq x_0(2)\rel x_0(4),\\
&x_0(1)\simeq x_0(4)\rel x_0(5),\\
&\textrm{etc.}
\end{align*}
All these  relations are  encoded in  the rooted tree   
$T$ of  \refF{tree9}.
%(for details about trees and stratification of compactified configurations, see \cite[Sections 2 and 3]{Sin:man}).

%\begin{figure}
%  \centering
%  \fbox{\includegraphics[width=95mm,height=60mm]{tree9.pdf}}
%  \caption{The tree $T$ associated to $x_0\in\Conf[9]$}
%  \label{fig:tree9}
%\end{figure}

\begin{figure}[h]
\input{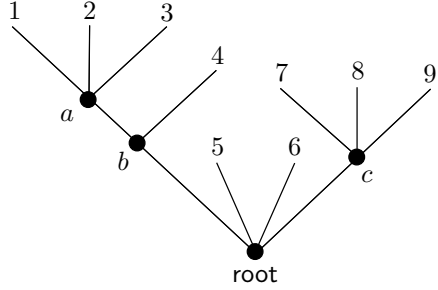}
\caption{The tree $T$ associated to $x_0\in\Conf[9]$.}
\label{fig:tree9}
\end{figure}

To the virtual configuration $x_0$ we associate a configuration of nested balls in $\BR^N$
as in \refF{ball9}, with one ball 
$B_v$ for each vertex $v\in\{1,\dots,9,a,b,c,\theroot\}$ of the tree $T$, so that $B_v\subset B_w$
iff $w$ is below $v$ in the tree and such that
  any two balls  are either disjoint or one is contained in the other.
 The centers of the balls labeled by the leaves define a  configuration $x_1=(x_1(1),\dots,x_1(9))\in\Conf(9)$.
We also assume that each ball is centered at the barycenter of the centers of the balls
immediately contained in that one.  For example, $B_b$ is centered at the barycenter of the centers
of $B_a$ and $B_4$. Also, the largest{\D} ball $B_{\theroot}$ is
centered at the origin\pfn{in \refF{ball9} the bigger balls should be
  labelled by $a,b,c$ and $root$. Could you do this?\iv Done.}.

%\begin{figure}
%  \centering
%  \fbox{\includegraphics[width=95mm,height=60mm]{ball9.pdf}}
%  \caption{A  configuration of nested balls associated to $x_0\in\Conf[9]$}
%  \label{fig:ball9}
%\end{figure}

\begin{figure}[h]
\input{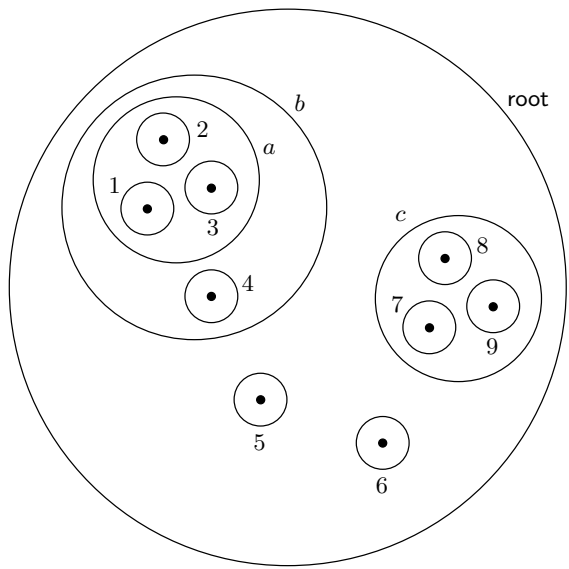}
\caption{A configuration of nested balls associated to $x_0\in\Conf[9]$.}
\label{fig:ball9}
\end{figure}

Consider a self-map
\[\phi_r\colon\BR^N\longrightarrow\BR^N\]
parametrized by $0< r\leq1$ whose effect is to iteratively  shrink each ball
$B_v$ by a homothety of factor $r$ and extend gradually up to the identity map outside of a small neighborhood of
the ball. \p
For $r=1$, $\phi_r$ is just the identity, but as $r\to 0$, the image of the configuration $x_1$ under $\phi_r$  tends to the virtual configuration $x_0$.

Now take a point $z$ anywhere inside the outermost closed ball $B_{\theroot}$ but outside of the innermost
open balls $B_i$ for $1\leq i\leq9$. Let $y_1=(x_1,z)\in \Conf(10)$ be the configuration obtained by adjoining
the point $z$ to the configuration $x_1$. Then the image of $y_1$
under 
$\underset{r\to 0}{\lim}\phi_r$
gives an element in the fiber $\pi^{-1}(x_0)\subset\Conf[10]$.
By choosing  the maps $\phi_r$ with care, we can ensure that  the fiber $\pi^{-1}(x_0)$ is covered by such $z$'s,
giving a homeomorphism $F\cong\pi^{-1}(x_0)$ where $F$ is a closed ball with $9$ small disjoint open
balls removed.

We want to prove the local triviality of $\pi$, so
allow now the centers of the nested  balls to move a bit around their initial value while preserving   the barycentric 
relations. Moreover, bound the shrinking of 
each ball $B_v$ below  by some parameter $\tau(v)\in [0,1]$, for $v$ a vertex other than the root or a leaf $i=1,\dots,9$.
Then applying $\phi_r$ to the configuration of the centers of the balls labeled by the leaves and letting $r\to 0$ describes a neighborhood $V$ of $x_0$ in $\Conf[{9}]$.
A parametrized (by $V$) version of the above construction will then give a trivialization $V\times F\cong\pi^{-1}(V)$
of $\pi$ over $V$. This trivialization can be made semi-algebraic and this will prove that $\pi$ is a semi-algebraic
bundle with an oriented compact generic fiber $F$.

We now proceed with the details of the proof of \refT{projSAbdl}.
\p Our goal is to build a neighborhood $V$ of $x_0$ in $\Conf[n]$ and a
diagram  (Diagram \refN{diag:PhihatPhi})
\[
\xymatrix{
W\times[0,r_1]^{\calV_0^*}\times F
\ar[r]^-{\widehat\Phi}
\ar[d]_{\proj}
&
\Conf[n+1]
\ar[d]^{\pi}
\\
W\times[0,r_1]^{\calV_0^*}
\ar[r]^-{\Phi}_-{\cong}
&
V\subset  \Conf[n].
}
\]
such that $\Phi$ is a semi-algebraic homeomorphism on $V$ and $\widehat\Phi$
 is a semi-algebraic homeomorphism on $\pi^{-1}(V)$.
Here $F$ is the fiber which will be homeomorphic to a unit disk $D^N$
with $n$ open disjoint balls removed. For the domain
$W\times[0,r_1]^{\calV_0^*}$
of the chart $\Phi$, $W$ is  a neighborhood in some products of
configuration spaces, and $\calV_0^*$ is the set of internal vertices
of the tree $T$
associated to $x_0$.

\subsubsection{\textbf{A stratification of $\Conf[n]$}}\

\previousoldfn{\ppfn{I \textbf{bolded} the titles of the subsubsections for clarity
  of the structure}}
\p We first review a classical stratification of the Fulton-MacPherson
configuration spaces indexed by trees (see also
\cite[appendix]{BoTa:SLK} and \cite{Sin:man}).
\begin{defin}
A \emph{rooted tree $T$ with labels in $\setn{n}=\{1,\dots,n\}$} is a tree (that is, an isomorphism class of a simply connected
$1$-dimensional finite simplicial complex) with one distinguished vertex called the \emph{root} of valence $\geq2$ and such 
that none of the other vertices is bivalent. The univalent vertices are called the \emph{leaves} and  are
 in bijection with the set $\setn{n}$.
\end{defin}

An example is given in \refF{tree9} for $n=9$.

Denote by $\calV$ the set of vertices of the tree $T$, including the root and the leaves.
The leaves are identified with the subset $\setn{n}\subset\calV$.
Set 
$$\calV_0:=\calV\setminus\{\theroot\},\ \ \  \calV^*:=\calV\setminus\setn{n},\ \ \  \calV^*_0:=\calV^*\cap\calV_0.$$
Define a partial order on $\calV$ by letting $w\leq v$ when the
shortest path \p in the tree joining $v$ to the root contains $w$.
We write $w<v$ when $w\leq v$ and $w\not=v$. The root is then the minimum of $\calV$.
Two vertices $v_1,v_2$ are \emph{not comparable} if neither $v_1\leq v_2$ nor $v_2\leq v_1$.
For a non-root vertex $v$ we define its \emph{predecessor} 
\[\predec(v):=\max\{w\in\calV:w<v\}.\]
For a non-leaf vertex $w$ we define its\emph{ output set}
\[\outputset(w):=\{v\in\calV:w=\predec(v)\}.\]
The \emph{height function} 
\[{\htree}\colon\calV\longrightarrow\BN\]
is defined inductively by ${\htree}(\theroot)=0$ and ${\htree}(v)={\htree}(\predec(v))+1$ when $v$ is not the root.

For example, in the tree of \refF{tree9} we have:
$ b\leq1$;
$\predec(4)=b$;
$\outputset(\theroot)=\{b,5,6,c\}$;
$b$ and $7$ are not comparable;
 and
${\htree}(a)=2$. For any $w\in\calV^*$, $|\outputset(w)|\geq2$.

For a rooted tree $T$ with leaves labeled by  $\setn{n}$
and set of vertices $\calV$,
consider the product of configuration spaces
\[\Conf_T:=\prod_{w\in\calV^*}\Conf(\outputset(w)).\]

\p We now recall how $\Conf_T$ can be identified, via a  
homeomorphism $h_T$ (see \eqref{eq:treestratumhomeo} below), to a stratum in $\Conf[n]$.\q
Let $\xi=(\xi^w)_{w\in\calV^*}\in\Conf_T$. Thus,{\D} identifying $\Conf(\outputset(w))$ with
$\Inj_0^1(\outputset(w),\BR^N)$ from \refN{eq:defC(X)xi}, for $w\in\calV^*$ we have
\[\xi^w\colon\outputset(w)\hookrightarrow\BR^N,\]
with 
\[
\bary(\xi^w)=0\ \ \  \text{and}\ \ \  
\radius(\xi^w)=1.
\]
For $v\in\calV_0$ we set
\[\xi(v):=\xi^{\predec(v)}(v).\]
For $r>0$ and $v\in\calV$, define
\begin{equation}\label{eq:def-xxirv}
x(\xi,r,v)
:=\sum_{\substack{w\in\calV_0\\w\leq v}}
\xi(w)\cdot 
r^{{\htree}(w)}.
\end{equation}
\p
The latter formula is equivalent to the inductive definition
\[\left\{
  \begin{array}{ll}
x(\xi,r,\theroot)&=\,0\\
x(\xi,r,v)&=\,x(\xi,r,\predec(v))\,+\,r^{\htree(v)}\xi(v).\\
\end{array}
\right.
\]

For $r>0$ small enough,
\[\left(x(\xi,r,i)_{1\leq i\leq n}\right)\]
 determines a 
configuration in $\Conf(n)$.
\p When $r\to0$ this configuration  converges to a virtual
configuration in $\Conf[n]$ whose proximity relations are described
by the tree $T$.
Define
\begin{equation}\label{eq:treestratumhomeo}
h_T\colon\Conf_T\longrightarrow\Conf[n]
\end{equation}
by
\[h_T(\xi)=\lim_{r\to0+}\left(x(\xi,r,i))_{1\leq i\leq n}\right).\]
Then $h_T$ is a homeomorphism onto its image
\previousoldfn{\pl needs arg? \iv I actually think it's clear.}
and the family of $\{\im(h_T)\}$, indexed by all rooted trees $T$ with
labels in $\setn{n}$, gives a stratification of $\Conf[n]$
\cite[Appendix]{BoTa:SLK}
(see also \cite[Sections 2 and 3]{Sin:man})\previousoldfn{\pl citation ok?}\previousfn{\pl
  17Jan2011:a actually because of formula \refN{eq:def-xxirv}  I think
  that Bott-Taubes is a better citation than Sinha here; I do not have
  Bott-Taubes here in Valparaiso. Is it possible to make the change
  with a precise reference inside that Bott-Taubes paper
  (or do you disagree) ? \iv  I cited Bott and Taubes but left Sinha in parentheses since Bott and Taubes don't give too many details.}.
The maximal stratum is $\Conf(n)$, which is the image of $h_{T_0}$ 
where $T_0$ is the tree for which all leaves are of height $1$.

A comment about the notation in this section migth be in order.
Along the rest of the proof, and as it has already appeared above, we will need to consider many configurations in
$\Conf(A)=\Inj_0^1(A,\BR^N)$ or $\Conf[A]$, for some
$A\subset\calV$.   They will sometimes come with various
decorations and arguments, such as
\begin{equation}
\label{eq:fullconfig}
\xi^w,\,\xi_0^w,\,x(\xi,r),\,x_0,\,x_1,\,x(\xi,r),\,\,x(\xi,\tau,r),\, x_1\,, y\,,
y_1\,, \,y_2,\,\textrm{etc.}
\end{equation}
We will also consider the components of these configurations in $\BR^N$, such as 
\begin{equation}
\label{eq:elementconfig}
\xi^w(v),\,\xi(w),\,x(\xi,r,v),\,x(\xi,\tau,r_1,\theroot),\,\xi_0(i),\,x_1(u),\,etc.
\end{equation}
It might therefore sometimes be confusing whether the notation corresponds to
a configuration in $\Conf[A]$ or to one of its components in
$\BR^N$.
As a rule of thumb, the notation will correspond to a point in $\BR^N$ when the last
argument is a vertex, that is an element of $\calV$, like
\[v,\,w,\,i,\,\theroot,\,1,\dots,\,n,\,\,u,\,p,\,q,\,v_1,\,v_2,\,\dots\quad\in\calV,\]
as in \refN{eq:elementconfig}. Otherwise it will be a configuration, 
as in \refN{eq:fullconfig}. In particular, in $\xi^w$ the vertex $w$
is a superscript and not an argument, and indeed $\xi^w$ is a
configuration but $\xi(v)$ is a point in $\BR^N$.

\subsubsection{\textbf{The chart $\Phi$ about $x_0$}}\

\p
Let $x_0\in\Conf[n]$. Our first goal is to
 build a neighborhood $V$ of $x_0$  over which  $\pi$ will be trivial
 \p
and a chart $\Phi$ of that neighborhood.
We have 
\begin{equation}\label{eq:xi0}
x_0=h_T(\xi_0)
\end{equation}
for some tree $T$ and some $\xi_0\in\Conf_T$.

For a finite set $A$ of at least two elements and for \p 
$\zeta\in\Inj^1_0(A,\BR^N)=\Conf(A)$
(see \refN{eq:defC(X)xi}\previousoldfn{\pl check ref}) define
\[\delta(\zeta):=\min\{\|\zeta(a)-\zeta(b)\|:a,b\in A,\,a\not=b\}\in (0,2].\] 
Set
\begin{equation}\label{eq:r1} 
r_1:=\frac{1}{4}\min\{\delta(\xi_0^w):w\in\calV^*\},
\end{equation}
and
set
\[W:=\{\xi\in\Conf_T:\forall\, v\in\calV_0, \|\xi(v)-\xi_0(v)\|\leq r_1^{n+1}\}\]
which is a compact
\previousoldfn{\pl need arg? \iv Clear again, I think.} neighborhood of $\xi_0$ in $\Conf_T$.

Consider now any function
\[\tau\colon\calV_0^*\longrightarrow[0,r_1]\]
that we extend to $\calV$ by  $\tau(\theroot)=0$ and $\tau(i)=0$ for $1\leq i\leq n$.
Define for $\xi\in W$ and $0\leq r\leq r_1$, by induction on the height of
$v\in\calV$,
\begin{equation}\label{eq:xxitaurv}
\left\{
  \begin{array}{ll}
x(\xi,\tau,r,\theroot)&=\,0\\
x(\xi,\tau,r,v)&=\,x(\xi,\tau,r,\predec(v))\,+\,\xi(v)\cdot\underset{\theroot\leq
  u<v}{\prod}
\max(r,\tau(u))\,\,\in\,\BR^n\\
\end{array}
\right.
\end{equation}
%In other words,
%\[x(\xi,\tau,r,v):=\sum_{\substack{w\in\calV_0\\w\leq v}}\left(\xi(w)\,\cdot\,
%\prod_{\substack{u\in\calV\\u\leq w}}\max(r,\tau(u))\right)
%\]
\pfn{in the above product there is a missing factor on the indexing
  set of $u$.  \iv I don't understand what you're saying here.}
Note that $x(\xi,\tau,r,w)$ is the barycenter of the points $x(\xi,\tau,r,v)$
for $v\in\outputset(w)$.
Note also that  when $\tau$ is bounded above by $r$ then
$x(\xi,\tau,r,v)=x(\xi,r,v)$
from \refN{eq:def-xxirv}.

Finally define
\begin{eqnarray}
\label{eq:Phi}\Phi\colon W\times[0,r_1]^{\calV_0^*}&\longrightarrow& \Conf[n]\\
\notag(\xi,\tau)\quad\quad&\longmapsto&\lim_{r\to 0+}\left(x(\xi,\tau,r,i))_{1\leq i\leq n}\right).
\end{eqnarray}

\begin{lemma}\label{L:Phihomeo}
$\Phi$ is a semi-algebraic homeomorphism onto  a compact neighborhood of $x_0$ in  $\Conf[n]$.
\end{lemma}

(Statements similar to this one appear in \cite{Sin:man}, \p  but without the
semi-algebraic condition.)
\previousoldfn{\ivfn{\iv I'm not sure if this comment is useful. \pl  this comment answer referee note (46)}}
\previousoldfn{\ppfn{in the bottom proof I replaced $\phi$ by $\varphi$ to not be
  confused for the further $\phi$ to be defined in \refL{phicre}}}
\begin{proof}
\p Let us first show that $\Phi$ is
semi-algebraic\pfn{referee comment (52)}.
The map
\begin{eqnarray*}
\varphi\colon W\times[0,r_1]^{\calV_0^*}\times(0,r_1]&\longrightarrow&\Conf[n]\times(0,r_1]\\
\quad(\xi,\tau,r)&\longmapsto&((x(\xi,\tau,r,i)_{1\leq i\leq n},r)
\end{eqnarray*}
is semi-algebraic, and hence the graph of $\varphi$ is a semi-algebraic set. 
This map can be continuously extented to 
a function
\[\overline\varphi \colon W\times[0,r_1]^{\calV_0^*}\times[0,r_1]
\longrightarrow\Conf[n]\times[0,r_1]
\]
whose restriction to $r=0$ is the limit function $\Phi$ (after
projection on the first factor).
The graph of $\overline\varphi$ is also semi-algebraic as it is the closure of
a semi-algebraic set. Therefore $\overline\varphi$, and hence $\Phi$ as well, is
semi-algebraic. (This argument is analogous to \cite[Proposition 2.9.1]{BCR:GAR}.)

Next we prove  the injectivity of $\Phi$.
Let  $y$ be  in the image of $\Phi$, that is
\[y=\lim_{r\to 0+}\left(x(\xi,\tau,r,i))_{1\leq i\leq n}\right).\]
We want to show that we can uniquely determine $\xi$ and $\tau$ from $y$.
Define inductively, for $w\in\calV$, $y(w)$ as the (virtual)
barycenter of the points $y(v)$ for $v\in\outputset(w)$. Then 
\[\xi^w=(y(v):v\in\outputset(w))\]
and
the function $\tau$ can be recovered by comparing the radii 
of the various sets $\{y(v):v\in\outputset(w)\}$, $w\in\calV^*$.
This proves the injectivity of $\Phi$.

Since the domain of $\Phi$ is compact, $\Phi$ is a homemorphism onto its image
and it is clear\previousoldfn{\pl need arg? connection with charts of Sinha? \iv I think connecting with Dev's charts would take us too far afield.  But maybe I don't really know what you mean by that.
\pl Need at least a reference to the chart then. \iv In
``Manifold-theoretic compactification...", Section 4.2 gives charts.
But I still don't understand where this reference would go.\pl OK forget} that this image is a neighborhood of $x_0=\Phi(\xi_0,0)$.
\end{proof}
We denote this compact neighborhood of $x_0$ in $\Conf[n]$  by
\begin{equation}
  \label{eq:VimPhi}
V:=\Phi(W\times[0,r_1]^{\calV_0^*}).
  \end{equation}

\subsubsection{\textbf{Shrinking balls to the limit configurations
  $\Phi(\xi,\tau)$}}\

\p
%We will prove the triviality of $\pi:\Conf[n+1]\to\Conf[n]$ over $V$.
We now  build a configuration of nested balls of centers
$x_1(\xi,v)\in\BR^N$ (for $\xi\in W\subset\Conf_T$ and $v$ a vertex in $\calV$) and of suitable radii $\epsilon(v)$, as well as
semi-algebraic  self-maps 
$\phi_r$ of $\BR^N$ which will shrink these balls ($\phi_r$ will depend on $r>0$, but also on $\xi\in W$ and
$\tau\in[0,r_1]^{\calV_0^*}$ not appearing in the notation).

The important features are
\begin{enumerate}
\item Applying the shrinking map $\phi_r$ to
the configuration of centers of innermost balls 
$(x_1(\xi,i))_{1\leq i\leq n}\in\Conf(n)$ gives the configurations 
$x(\xi,\tau,r,i)_{1\leq i\leq n}$ which serves, as $r\to0$, to define
the chart $\Phi$ in \refN{eq:Phi} (\refL{phi=x}).
\item The complement of the innermost balls inside the outermost ball
  will serve as the fiber of the projection $\pi$ (this will appear in
  the next section  and will be based on the properties of \refL{phicre}
  (2).)
\end{enumerate}

% and this will lead to the 
% desired trivialization $\widehat\Phi$ \Z of  $\pi^{-1}(V)$. We denote by $\Ball(c,r)$ and $\Ball[c,r]$ the open and
% the closed balls  in $\BR^N$ of center $c$ and radius $r$.

%%%%% Some section here was moved below 17 jan 2011%%%%%%%

We define first the 
\previousoldfn{\pl 17jan2011: \refL{phicre} (i.e. L:phicre, formerly Lemma 5.25) was moved further to fit  better the logic of the proof}\Z 
centers $x_1(\xi,v)$ and the radii $\epsilon(v)$ of the balls that 
we will consider. Suppose given $\xi\in W$ and recall the map $x$
defined in \refN{eq:def-xxirv} and the radius $r_1>0$ from \refN{eq:r1}.
 For $v\in\calV$, we set
\begin{equation}\label{eq:x1x'1}
x_1(\xi,v):=x(\xi,r_1,v)
\end{equation}
and 
\[\epsilon(v):=4\cdot r_1^{{\htree}(v)+1}.\]
The \Z\previousoldfn{pl 17jan2011:phrase added} \p
balls $\Ball[x_1(\xi,v)\,,\,\epsilon(v)]$ satisfy the following
nesting properties:
\begin{lemma}
\ 

\label{L:epsABC}
\begin{enumerate}
\item
If  $w<v$ in $\calV$  then 
\[\Ball[x_1(\xi,v),\epsilon(v)]\subset  \Ball[x_1(\xi,w),\epsilon(w)/3].
\] 
\item
If  $v_1$ and $v_2$ are not comparable in $\calV$  then
\[\Ball[x_1(\xi,v_1),\epsilon(v_1)]\cap \Ball[x_1(\xi,v_2),\epsilon(v_2)]=\emptyset.
\]
\end{enumerate}
\end{lemma}
\begin{proof}
To simplify notation, we set $r=r_1$ in this proof. Note that $r\leq 1/2$
because $\delta(\xi^w_0)\leq2$.
\previousoldfn{\pl In a previous version there was the following proof of
\[
\|x_1(v)-x'_1(v)\|<\epsilon(v)/3
\]
where $x'_1(v):==x(\xi',r_1,v)$ for some $\xi'\in W$.\\
Proof:
 Using the definition of $W$, 
\begin{multline*}
\left\|x_1(v)-x'_1(v)\right\|=
\left\|
\sum_{\substack{w\in\calV_0\\w<v}}\left(\xi(w)-\xi'(w)\right)\cdot r^{{\htree}(w)}
\right\|
\leq
\left(\sup_{\substack{w<v}} \left\|\xi(w)-\xi'(w)\right\| \right) \cdot
\frac{r}{1-r}\leq
\\
\leq
2 \cdot r^{n+1}\left(\frac{r}{1-r}\right)\quad<\quad(4/3)r^{n+1}\quad\leq\quad \epsilon(v)/3.
\end{multline*}
}%\oldfn

For $w<v$,
\[
x_1(\xi,v)=x_1(\xi,w)+\sum_{w<u\leq v}\xi(u)\cdot r^{{\htree}(u)}.
\]
Therefore
\begin{eqnarray*}
\|x_1(\xi,v)-x_1(\xi,w)\|+\epsilon(v)
&\leq&
\sum_{w<u\leq v}\|\xi(u)\|\cdot r^{{\htree}(u)}+4\cdot r^{{\htree}(v)+1}\\
&\leq&
\left(\sup_{w<u\leq v}\|\xi(u)\|\right)\cdot\frac{r^{{\htree}(w)+1}}{1-r}+4\cdot r^{{\htree}(v)+1}\\
&\leq&r^{{\htree}(w)+1}\left(\frac{1}{1-r}+4\cdot r \right)\\
&\leq&(4/3)\cdot r^{{\htree}(w)+1}\\
&=&\epsilon(w)/3.
\end{eqnarray*}
This proves the first part of the lemma.

For the second part,
 suppose first that $v_1$ and $v_2$ have a common predecessor $w$.
Then
\begin{eqnarray*}
\|\xi(v_2)-\xi(v_1)\|&\geq&
\|\xi_0(v_2)-\xi_0(v_1)\|-\|\xi(v_1)-\xi_0(v_1)\|-\|\xi(v_2)-\xi_0(v_2)\|\\
&\geq&
\delta(\xi_0^w)-2\cdot r^{n+1}\\
&\geq&4\cdot r -2\cdot r^{n+1}\\
&>&2\cdot r.
\end{eqnarray*}
Since  ${\htree}(v_1)={\htree}(v_2)$ we get
\begin{eqnarray*}
\|x_1(\xi,v_1)-x_1(\xi,v_2)\|
&=&
\|\xi(v_1)-\xi(v_2)\|\cdot r^{{\htree}(v_1)}\\
&>& 2\cdot r\cdot r^{{\htree}(v_1)}\\
%&=&2\cdot\epsilon(v_1)\\
&=&\epsilon(v_1)+\epsilon(v_2).
\end{eqnarray*}
This implies the desired formula when $v_1$ and $v_2$ have a common predecessor.

For the general case,
since $v_1$ and $v_2$ are not comparable, there exists
$w_1\leq v_1$ and $w_2\leq v_2$ such that
 $w_1$ and $w_2$ have a common predecessor.
Therefore $\Ball[x_1(\xi,w_1),\epsilon(w_1)]\cap \Ball[x_1(\xi,w_2),\epsilon(w_2)]=\emptyset$.
Combining  this
with the fact that, by the first part of the proposition, $\Ball[x_1(\xi,v_i),\epsilon(v_i)]\subset \Ball[x_1(\xi,w_i),\epsilon(w_i)]$, for $i=1,2$,
we deduce the desired formula.
\end{proof}

%%%%%%%%%%%%%%%%%%%%%%%%%%%%%%%%%%%%%%%%%%%%%%%%%%%
%%%%%%%%This section was moved by Pascal on 17 jan 2011 %%%%%%%%%%%%%
%%%%%%%%%%%%%%%%%%%%%%%%%%%%%%%%%%%%%%%%%%%%%%%%%%%

%We start by defining 
We next define \Z\previousoldfn{\pl17jan2011:  this section was moved further} a suitable morphism shrinking a given ball.

\begin{lemma}\label{L:phicre}
There exists a continuous semi-algebraic map
\begin{eqnarray*}
\phi\colon&
\BR^N\times [0,1]\times[0,2]\times\BR^N&
\longrightarrow\BR^N\\
&(c,r,\epsilon,x)&
\mapsto \phi^{c,\epsilon}_r(x)
\end{eqnarray*}
with the following properties:
\begin{enumerate}
\item the map $x\mapsto\phi^{c,\epsilon}_r(x)$
\begin{enumerate}
\item is radial, centered at $c$;
\item is the identity outside of the ball $\Ball(c,\epsilon)$;
\item restricts on $\Ball[c,\epsilon/3]$ to a homothety of rate $r$;
\item when $r>0$, it is a self-homeomorphism of $\BR^N$;
\item when $r=0$, its restriction  to 
 $\BR^N\setminus\Ball[c,\epsilon/2]$ is a homeomorphism onto $\BR^N\setminus\{c\}$, and
$\phi^{c,\epsilon}_0(\Ball[c,\epsilon/2])=\{c\}$;
\end{enumerate}
\item let $r>0$ and let $x(1),\dots,x(n)$ be $n\geq2$ distinct points in $\Ball[c,\epsilon/3]$; then
\begin{enumerate}
\item
$(\phi_r^{c,\epsilon}(x(1)),\dots,\phi_r^{c,\epsilon}(x(n)))$
determines a configuration in $\Conf(n)$ which  does not depend on
$r$; hence its limit as $r\to 0+$ determines the same configuration;
\item if $z_1,z_2$ are two distinct points in $\Ball(c,\epsilon/2)$ and are different
from the $x(p)$'s for $1\leq p \leq n$, then 
\[y_i:=\lim_{r\to 0+}(\phi_r^{c,\epsilon}(x(1)),\dots,\phi_r^{c,\epsilon}(x(n)),\phi_r^{c,\epsilon}(z_i))\]
determines two different configurations $y_1$ and $y_2$ in
$\Conf(n+1)$ such that  $y_i(p)\not\simeq y_i(q)\rel y_i(n+1)$  for
$1\leq p\not= q\leq n$ and $i=1,2$;
\item if $z\in\BR^N\setminus\Ball(c,\epsilon/2)$ then
\[y:=\lim_{r\to 0+}(\phi_r^{c,\epsilon}(x(1)),\dots,\phi_r^{c,\epsilon}(x(n)),\phi_r^{c,\epsilon}(z))\]
determines a  configuration in $\Conf[n+1]$ such that $y(p)\simeq y(q)\rel y(n+1)$ for $1\leq p,q\leq n$.
\end{enumerate}
\end{enumerate}
\end{lemma}
\begin{proof}
\p\q The proof consists of explicitly constructing the semi-algebraic
function $\phi$.
Define first a  semi-algebraic function
\begin{eqnarray*}
g\colon[0,1]\times\BR_+&\longrightarrow&[0,1]\\
(r,u)\qquad&\longmapsto&g(r,u)
\end{eqnarray*}
by
\[
g(r,u)=
\begin{cases}
r,&\textrm{if }0\leq u\leq1/3;\\
\frac{r}{3-6u},&\textrm{if }1/3\leq u\leq 1/2\textrm{ and }\sqrt{r}\leq3-6u;\\
\sqrt{r},&\textrm{if }1/3\leq u\leq 1/2\textrm{ and }\sqrt{r}\geq3-6u;\\
2\sqrt{r}(1-u)+2u-1,&\textrm{if }1/2\leq u\leq 1;\\
1,&\textrm{if }u\geq1.
\end{cases}
\]
In other words, the function $g$ is determined by the  picture
in \refF{grapheG} where the curve inside the second rectangle
is the parabola $\sqrt{r}=3-6u$.

%\begin{figure}
%  \centering
%  \fbox{\includegraphics[width=85mm,height=50mm]{graphG.pdf}}
%  \caption{Definiton of the function $(r,u)\mapsto g(r,u)$}
%  \label{fig:grapheG}
%\end{figure}

\begin{figure}[h]
\input{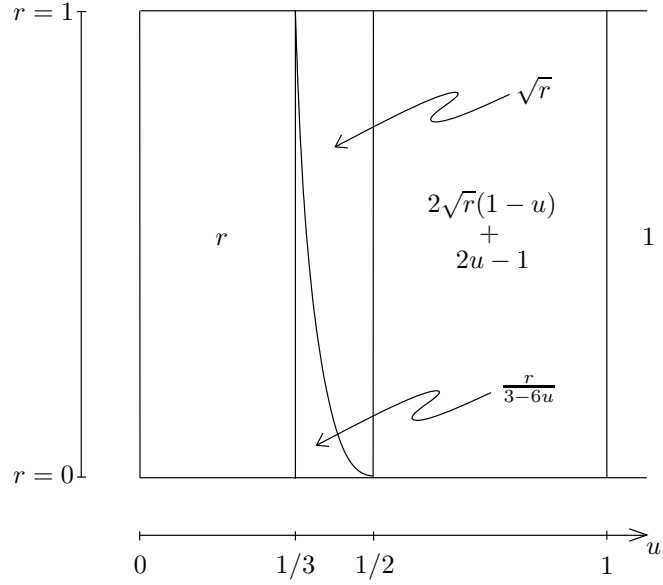}
\caption{Definition of the function $(r,u)\longmapsto g(r,u)$.}
\label{fig:grapheG}
\end{figure}

For $c\in\BR^n$, $\epsilon>0$, and $r\geq0$  define
\begin{eqnarray}
\label{eq:def-phicepsr}
\phi^{c,\epsilon}_r\colon\BR^N&\longrightarrow&\BR^N\\
x&\longmapsto&\phi^{c,\epsilon}_r(x)\,=\,c+(x-c)\cdot
g\left(r\,,\,\frac{\|x-c\|}{\epsilon}\right).\notag
\end{eqnarray}

Properties (1a-e) of $x\mapsto \phi^{c,\epsilon}_r(x)$ are then
immediate.

(2a) follows from (1c).\Z

(2b-c)  are{\D} consequences of the following properties of $g$ as $r\to0$:
\begin{itemize}
\item For $u<1/2$, the map 
\[
u\mapsto \lim_{r\to 0}\frac{g(r,u)}{r}
\]
is the constant $1$ over $[0,1/3]$ and gives a semi-algebraic
homeomorphism
between $[1/3,1/2)$ \Z and $[1,+\infty)$.
\item For $u\geq 1/2$, 
\[\lim_{r\to 0+}\frac{g(r,u)}{r}=+\infty
\]
and the map
\[
u\mapsto \lim_{r\to 0}{g(r,u)}
\]
is a homeomorphism between $[1/2,1]$ and $[0,1]$.
\end{itemize}
\end{proof}

%%%%%%%%%%%%%%%%%%%%%%%%%%%%%%%%%%%%%%%%%%%%%%%%%%%%%%%%%%%%
%%%%%%%%%%% The section above was moved further on 17 jan 2011 %%%%%%%%%%%%%%%%%%%%%%
%%%%%%%%%%%%%%%%%%%%%%%%%%%%%%%%%%%%%%%%%%%%%%%%%%%%%%%%%%%%

Fix 
\previousoldfn{\pl 17Jan2011: the section above starting at \refL{phicre} was
  moved from before \refL{epsABC}}\Z
$\xi\in W$ and  $\tau\in[0,r_1]^{\calV_0^*}$. 
Recall that we extend $\tau$ to $\calV$ by $0$ on the root and the leaves.
For $v\in\calV$ and  $0< r\leq r_1$, set
\[\phi_r^v:=\phi^{x_1(\xi,v),\epsilon(v)}_{\max(r,\tau(v))/r_1}.\]
Note that $\phi_r^v$ depends on $\xi$ and $\tau$ even if this does not appear in the notation.
Then $\phi_r^v$ is a self-map of $\BR^N$ which is the identity
outside of $\Ball[x_1(\xi,v),\epsilon(v)]$, and shrinks the ball 
$\Ball[x_1(\xi,v),\epsilon(v)/3]$ by a homothety of rate $\max(r,\tau(v))/r_1$.
We will compose all these maps $\phi^v_r$ for $v\in \calV$.

If $v_1,v_2\in\calV$ are two distinct vertices  of the same height, then they are non comparable and
Lemmas \ref{L:phicre} (1b) and \ref{L:epsABC} (2) imply that 
\begin{equation}\label{eq:phicomm}
\phi^{v_1}_r\circ\phi^{v_2}_r=\phi^{v_2}_r\circ\phi^{v_1}_r.
\end{equation}

Let $h_{\max}:=\max\{\htree(v):v\in\calV\}$. For $h=1,\dots,h_{\max}$ we define 
\[
\phi^{[h]}_r
:=
\underset{v\in\calV_0, \\\htree(v)=h}{\circ}\phi^v_r\]
which is the composition of the maps $\phi^v_r$ for all vertices $v$ of height $h$,
the order of composition being irrelevant because of \refN{eq:phicomm}.
Finally we set
\[\phi_r:=\phi_r^{[1]}\circ\phi_r^{[2]}\circ\dots\circ\phi_r^{[h_{\max}]}\]
which is a self-map of $\BR^N$ which has the effect of
iteratively shrinking all the  balls
of center $x_1(\xi,v)$ starting with the innermost ones{\D} first.

%%%%%%%%%%%%%%%%%%%%%%%%%%%%%%%%%%%%
%%%% some section below was moved further  ON 17Jan2011%%%%%%%%%
%%%%%%%%%%%%%%%%%%%%%%%%%%%%%%%%%%%%

Using 
\Z\previousoldfn{\pl Some section before the following lemma has been moved
  on 17jan2011 after the proof of the lemma to fit better the logic of
  the proof. Also this phrase has been added to paraphrase the lemma;
  please fix the english as needed.}
$\phi_r$, we can recover the $x(\xi,\tau,r,i)$, which appears
in the chart $\Phi$ from \refN{eq:Phi}, as follows:
\begin{lemma}\label{L:phi=x}For $r>0$ and \Z $1\leq i\leq n$, we have, in $\BR^N$,
\[\phi_r(x_1(\xi,i))
\,=
\,x(\xi,\tau,r,i)
.\]
\end{lemma}
\previousoldfn{\pl check \iv done.}
\begin{proof}
Set, for $w\in\calV$,
\begin{equation}
  \label{eq:xrw}
x_r(w)
:=
  \begin{cases}
\phi_{r}(x_1(\xi,w)), & \textrm{if } w \textrm{ is a leaf};\\
\bary(x_r(v):v\in\outputset(w)), & \textrm{otherwise.}
\end{cases}
\end{equation}
Since $x(\xi,\tau,r,v)$ defined in \refN{eq:xxitaurv} satisfies a barycentric relation
analoguous to \refN{eq:xrw}, the idea of the proof is to compare, by
induction on the height of
 $v\in\calV$,
\[
x_r(v)\textrm{ and } x(\xi,\tau,r,v).
\]
For the sake of the proof, define for $h\geq 1$ and $w\in\calV$
\begin{eqnarray*}
\phi^{\geq[h]}_r
&:=
&\phi^{[h]}_r\circ \phi^{[h+1]}_r\circ\dots\circ \phi^{[h_{\max}]}_r
\\
\phi^{\leq[h]}_{r}
&:=
&\phi^{[1]}_r\circ \phi^{[2]}_r\circ\dots\circ \phi^{[h]}_r
\\
x^{\geq[h]}_r(w)
&:=
&
  \begin{cases}
\phi^{\geq[h]}_{r}(x_1(\xi,w)), &\textrm{if }w\textrm{  is a leaf};\\
\bary(x^{\geq[h]}_r(v):v\in\outputset(w)), &\textrm{otherwise.}
\end{cases}
\end{eqnarray*}
In particular, $\phi_r=\phi_r^{\geq[1]}$ and $x_r(v)=x^{\geq[1]}_r(v)$.

We begin by proving three claims on the relations between these self-maps 
and configurations.

\textbf{Claim 1:} For $w\in\calV^*$, $h=\htree(w)$, and $v>w$, 
\begin{equation}
  \label{eq:claim1}
  x^{\geq[h]}_r(v)=\phi^w_r\left(x^{\geq[h+1]}_r(v)\right).
\end{equation}

The claim is proved by induction on $v>w$. If $v$ is a leaf then
\[
x^{\geq[h]}_r(v)=\phi^{[h]}_r(x^{\geq[h+1]}_r(v))=\phi^{w}_r(x^{\geq[h+1]}_r(v)).
\]
Suppose that the claim has been proved when $v$ is replaced by any
$u\in\outputset(v)$ in \refN{eq:claim1}.
 Then
 \begin{eqnarray}
\notag   x^{\geq[h]}_r(v)
&=&
\bary(x^{\geq[h]}_r(u):u\in\outputset(v))\\
\label{eq:claim1b}
&\stackrel{\textrm{induction}}{=}&
\bary(\phi^{w}_r(x^{\geq[h+1]}_r(u)) :u\in\outputset(v)).
 \end{eqnarray}
Set
\begin{equation}
  \label{eq:Bw}
  \Ball_w:=\Ball[x_1(\xi,w),\epsilon(w)/3].
\end{equation}
By \refL{phicre} (1c), the restriction $\phi^w_r|\Ball_w$ is a
homothety and hence commutes with taking the barycenter.
Since $x^{\geq[h+1]}_r(u)\in\Ball_w$ for $u>w$, we deduce from
\refN{eq:claim1b}  that
\begin{eqnarray*}
    x^{\geq[h]}_r(v)
&=&
\bary(\phi^w_r(x^{\geq[h+1]}_r(u)):u\in\outputset(v))
\\
&=&
\phi^{w}_r(\bary(x^{\geq[h+1]}_r(u) :u\in\outputset(v)))
\\
&
%\stackrel{\refN{eq:claim1b}}
{=}&
\phi^{w}_r(x^{\geq[h+1]}_r(v)).
\end{eqnarray*}
This proves Claim 1.

\textbf{Claim 2:} For $v\in\calV$,
\[x^{\geq[\htree(v)]}_r(v)=x_1(\xi,v).\] 

The proof of Claim 2 is by induction. % on the vertices of the tree $v$.
The claim is clear when $v$ is a leaf.  Suppose that the claim is true
for all $v\in\outputset(w)$ and set $h=\htree(w)$.
Then
\begin{eqnarray*}
  x^{\geq[h]}_r(w)
&=&
\bary( x^{\geq[h]}_r(v):v\in \outputset(w))\\
&\stackrel{\textrm{Claim 1}}{=}&
\bary(\phi^w_r(x^{\geq[h+1]}_r(v)):v\in \outputset(w))\\
&\stackrel{\textrm{induction hyp.}}{=}&
\bary(\phi^w_r(x_1(\xi,v)):v\in \outputset(w))\\
&\stackrel{\phi^w_r|\Ball_w\textrm{ homothety}}{=}&
\phi^w_r(\bary( x_1(\xi,v)):v\in \outputset(w))\\
&=&
\phi^w_r(x_1(\xi,w))\\
&=&x_1(\xi,w),
\end{eqnarray*}
which proves Claim 2.

As a special case we have\Z
\begin{equation}
  \label{eq:xroot}
  x_r(\theroot)=x_1(\xi,\theroot)=0.
\end{equation}

\textbf{Claim 3:} 
For $v\in \calV_0$ and $h=\htree(v)$,
\[x_r(v)=\phi^{\leq[h-1]}_r(x_1(\xi,v))=\phi^{\leq[h]}_r(x_1(\xi,v)).
\]

Let $w=\predec(v)$. Then the restriction of $\phi^{\leq[h-1]}_r$
to $\Ball_w$ is a composition of homotheties and hence commutes with
taking the 
barycenter. Since $x_r(v)$ is defined by iterated barycenters from a
collection of points
\[x_r(i)=\phi^{\leq[h-1]}_r(x^{\geq[h]}_r(i))\]
which belong to the convex $\Ball_w$ (because $i$ are leaves above
$w$),
we deduce that
\begin{eqnarray*}
  x_r(v)
&=&
\phi^{\leq[h-1]}_r(x^{\geq[h]}_r(v))
\\
&\stackrel{\textrm{Claim 2}}{=}&
\phi^{\leq[h-1]}_r(x_1(\xi,v)).
\end{eqnarray*}
Finally, since
\[\phi^{[h]}_r(x_1(\xi,v))=\phi^v_r(x_1(\xi,v))=x_1(\xi,v)\]
we have
\[\phi^{\leq[h-1]}_r(x_1(\xi,v))
=
\phi^{\leq[h]}_r(x_1(\xi,v)).
\]
This proves Claim 3.

We are ready for the proof of the lemma.
Let $w\in\calV^*$. Recall that the restriction of
$\phi^{\leq[\htree(w)]}_r$ to $\Ball_w$
is a composition of homotheties of total rate
\[
R_w:=\underset{\theroot<u\leq w}\prod\frac{\max(r,\tau(u))}{r_1}.
\]
Consider the \emph{normalization} map
\[
N:\Inj(A,\BR^N)\longrightarrow\Inj^1_0(A,\BR^N)=\Conf[A]
\]
that translates the barycenter to the origin and rescales to
$\radius=1$. This map is invariant under homotheties of the arguments.
Therefore
\begin{eqnarray}\label{eq:xrxiw}
\notag
  N(x_r(v):v\in\outputset(w))
&\stackrel{\textrm{Claim 3}}{=}&
N(\phi_r^{\leq[\htree(w)]}(x_1(\xi,v) :v\in\outputset(w))
\\
\notag
&\stackrel{\textrm{homotheties}}{=}&
N(x_1(\xi,v) :v\in\outputset(w))
\\
&=&%%&{\stackrel{\refN{eq:xxitaurv}}{=}}&
\xi^w.
\end{eqnarray}
Also
\begin{eqnarray}\label{eq:radiusxr}
\notag
\radius(x_r(v):v\in\outputset(w))
&=&
\radius(\phi_r^{\leq[\htree(w)]} (x_r(v)):v\in\outputset(w))
\\
\notag
&=&
R_w\cdot \radius((x_1(\xi,v):v\in\outputset(w))
\\
\notag
&=&
R_w\cdot r_1^{\htree(w)}
\\
&=&
\underset{\theroot<u\leq w}{\prod}\max(r,\tau(u)).
\end{eqnarray}
Comparing Equations \refN{eq:xroot}, \refN{eq:xrxiw} and
\refN{eq:radiusxr}
with \refN{eq:xxitaurv}, we deduce that for all $v\in \calV$,
\[x_r(v)=x(\xi,\tau,r,v).\]
The statement of the lemma is the special case when $v$ is a leaf.
\end{proof}

\subsubsection{\textbf{The chart $\widehat\Phi$ of $\pi^{-1}(V)$}}\  

\p
%%%%%%%%%%%%%%%%%%%%%%%%%%%%%%%%%%%%
%%%% MOVED FURTHER ON 17Jan2011%%%%%%%%%
%%>>>>>>>>>>>>>>>>>>>>>>>>>>>>>>>>
 
We 
\Z\previousoldfn{\pl this phrase has been added and the following paragraph was moverd on 17jan2011 from before \refL{phi=x}}
are ready to define the trivialization $\widehat\Phi$ of the canonical
projection $\pi$.
Set 
\[F:=\Ball[0,n+1]\setminus\cup_{i=1}^n\Ball(i,1/4).\]
This is a closed ball with $n$ disjoint{\D} open balls removed  and
will serve as the generic \p\q fiber of $\pi$.
For $\xi\in W$, also set 
\[F_\xi:=\Ball[x_1(\xi,\theroot),\epsilon(\theroot)/2]\setminus\cup_{i=1}^n\Ball(x_1(\xi,i),\epsilon(i)/2).\]
It is easy to build semi-algebraic homeomorphisms
\[\Theta_\xi\colon F\iso F_\xi\]
that depend continuously and semi-algebraically on $\xi\in W$ since $W$ is 		``small".

Recall the homeomorphism 
\[\Phi\colon W\times[0,r_1]^{\calV^*_0}\stackrel\cong\longrightarrow V\subset\Conf[n]
\]
 from \refN{eq:Phi} and
\refL{Phihomeo}. Define
\[\widehat\Phi\colon W\times[0,r_1]^{\calV_0^*}\times F\longrightarrow\Conf[n+1]
\]
by
\begin{equation}
\label{eq:hatPhi}  
\widehat\Phi(\xi,\tau,z_0):=
\lim_{r\to  0+}\left(\phi_r(x_1(\xi,1)),\dots,\phi_r(x_1(\xi,n)),\phi_r(\Theta_{\xi}(z_0))\right).
\end{equation}

\previousoldfn{\pl continuity of $\widehat\Phi$ ?.}
%%<<<<<<<<<<<<<<<<<<<<<<<<<<<<<<
%%%% MOVED FURTHER ON 17Jan2011%%%%%%%%%
%%%%%%%%%%%%%%%%%%%%%%%%%%%%%%%%%%%%%%

By{\D} \refN{eq:Phi}, \refN{eq:hatPhi}, and \refL{phi=x},
 the following diagram commutes:
 \begin{equation}
\label{diag:PhihatPhi}
\xymatrix{
W\times[0,r_1]^{\calV_0^*}\times F
\ar[r]^-{\widehat\Phi}
\ar[d]_{\proj}
&
\Conf[n+1]
\ar[d]^{\pi}
\\
W\times[0,r_1]^{\calV_0^*}
\ar[r]^-{\Phi}_-{\cong}
&
V\subset  \Conf[n].
}
\end{equation}

We want to show that $\widehat\Phi$ is a homeomorphism onto $\pi^{-1}(V)$,
where
$V=\im\Phi$ is from \refN{eq:VimPhi}.
Fix $(\xi,\tau)\in W\times[0,r_1]^{\calV_0^*}$.
It is enough to show that $\widehat\Phi$ restricts to a homeomorphism on
the fibers:
\begin{eqnarray}\label{eq:hatphi}
\hat\phi\colon F_\xi&\iso&\pi^{-1}(\Phi(\xi,\tau))\\
\notag
z&\longmapsto&\widehat\Phi(\xi,\tau,\Theta_{\xi}^{-1}(z)).
\end{eqnarray}

We  first show\previousoldfn{\pl the end of this section has to be read again} that $\hat\phi$ is injective. 
Let $z_1,z_2$ be two distinct elements in $F_\xi$. Set $y_i=\hat\phi(z_i)\in\Conf[n+1]$
for $i=1,2$. We treat different cases.
\begin{itemize}
\item Suppose that there exists a vertex $v\in\calV$ such that
$z_1\in\Ball(x_1(\xi,v),\epsilon(v)/2)$ but 
$z_2\not\in\Ball(x_1(\xi,v),\epsilon(v)/2)$ (or the other way around).
By definition of $F_\xi$,  $v$ is not a leaf. Thus $v$ has at least
 two distinct outputs and we choose two leaves $p$ and $q$
above each of these outputs. Using \refL{phicre} (2b-c), we get
\[y_1(p)\not\simeq y_1(q)\rel y_1(n+1)
\]
but
\[y_2(p)\simeq y_2(q)\rel y_2(n+1).
\]
Thus $y_1\not=y_2$.
\item Suppose that the highest vertex $v\in\calV$ such that $z_1\in\Ball(x_1(\xi,v),\epsilon(v)/2)$
is the same as the highest vertex $w\in\calV$  such that $z_2\in\Ball(x_1(\xi,w),\epsilon(w)/2)$,
that is $v=w$.
\previousoldfn{\iv I don't understand this sentence.\pl clearer?}
Choose again two leaves $p,q$ above two distinct outputs of $v$.
Set 
\[\phi^{\geq v}_r:=\phi_r^{[\htree(v)]}\circ\dots\circ\phi_r^{[\htree_{\max}]}.\]
By  \refL{phicre} (2b), we have that 
\[
\lim_{r\to 0}(\phi^{\geq v}_r(z_i),\phi^{\geq v}_r(x_1(\xi,p)),\phi^{\geq v}_r(x_1(\xi,q)))
\]
defines two distinct configurations in $\Conf(3)$.
Then applying 
\[
\lim_{r\to 0}\phi_r^{[0]}\circ\dots\circ\phi_r^{[\htree(v)-1]},
\]
which is a composition of homotheties of a ball containing the{\D}  configurations,
 still gives two distinct configurations in $\Conf(3)$.
Therefore the images of $y_1$ and $y_2$ under some canonical projection
$\pi\colon\Conf[n+1]\to\Conf[3]$ are distinct. Thus $y_1\not=y_2$.
\item It remains to treat the case when there is no $v\in\calV$
such that $z_i\in\Ball(x_1(\xi,v),\epsilon(v)/2)$ for $i=1$ or $i=2$.
Then $z_1,z_2\in\partial\Ball[x_1(\xi,\theroot),\epsilon(\theroot)/2]$
are in the boundary of the largest ball which is centered at the origin.
In that case
\[\theta_{1,n+1}(y_i)=z_i/\|z_i\|\]
where $\theta_{1,n+1}$ from
\refN{eq:defthetaab} gives the direction between the first and the
last point of the configuration,
and these two directions are distinct. Thus $y_1\not= y_2$.
\end{itemize}
This proves that $\hat\phi$ is injective. For surjectivity,
 since $F_\xi$ and $\pi^{-1}(\Phi(\xi,\tau))$ are compact connected
 manifolds, it is enough to show that $\hat\phi$ is surjective
on the intersection of the boundary with the fiber. This boundary consists of virtual configurations
$y\in\Conf[n+1]$ such that:
\begin{itemize}
\item[(a)]
either for some $1\leq i\leq n$ and for all $j\in\setn{n}\setminus\{i\}$, we have:
$
y(i)\simeq y(n+1)\rel y(j)  
$;
\item[(b)]
or for all $1\leq i,j\leq n$, we have: 
$
y(i)\simeq y(j)\rel y(n+1).
$
\end{itemize}
It is clear that $\hat\phi$ maps
$\partial\Ball[x_1(\xi,i),\epsilon(i)/2]$
surjectively onto the boundaries of type (a)
and
%\previousfn{\pl forced ``line feed''.  \iv I don't know what this means.  Is this some problem with the tex file?}\\
$\partial\Ball[x_1(\xi,\theroot),\epsilon(\theroot)/2]$
onto that of type (b).

This proves that $\hat\phi$ from \refN{eq:hatphi} is a homeomorphism
and hence 
that $\widehat\Phi$ in Diagram \refN{diag:PhihatPhi} is a homeomorphism
onto $\pi^{-1}(V)$.
Thus\Z%, for $n\geq2$, 
\[\pi\colon\Conf[n+1]\longrightarrow\Conf[n]\]
 is an SA bundle.
For $n\geq2$, its 
%generic
 \p  fiber
% (fiber over the interior) 
$F$ is a compact manifold of
dimension $N$ whose interior is homeomorphic to $\BR^N$ with
$n$ points{\D} removed.
 For $n=0$ the fiber is
a point and for $n=1$ the fiber is an $(N-1)$-dimensional sphere $S^{N-1}$.

Thus $\pi:\Conf[V:=A\amalg I]\to\Conf[A]$ is an oriented SA bundle, as
it is the composition of oriented{\D} SA bundles $\pi\colon\Conf[n+1]\to\Conf[n]$   \citePAPcompositebdl.
When $|A|\geq2$,  the interior of the fiber of $\pi$
can be identified with $\Inj(I,\BR^N\setminus A)$ which is of codimension $0$ inside $(\BR^N)^I$.
The latter manifold has a canonical orientation  when $N$ is even or when $N$ is odd and $I$ is linearly ordered, and
in the second case a transposition in the linear order reverses the orientation.

This finishes the proof of \refT{projSAbdl}.

\section{The CDGAs of admissible diagrams}
\label{sec:diagrams}
In this section we introduce the
 CDGA
of \emph{admissible diagrams} $\AD({\ExtVert})$, where $\ExtVert$ is a finite set, for example $\ExtVert=\setn{n}=\{1,\dots,n\}$.
As we will prove later, this differential algebra  is a model for both $\ompa(\Conf[{\ExtVert}])$ and  its cohomology,
 and it will  serve as an intermediate model in the formality proof.
In \refS{diagcoop} we will endow $\AD:=\{\AD({n})\}_{n\geq0}$ with the structure of a cooperad.

The CDGA $\AD({\ExtVert})$ could be defined directly but we 
 will  describe it as a quotient of a larger CDGA of \emph{ diagrams} $\GD({\ExtVert})$
that we will introduce first.
One reason for doing so is that it will be easier to define a cooperad
structure on $\GD:=\{\GD({n})\}_{n\geq0}$ and establish some of its properties, and then induce from
this the cooperad structure for $\AD$.

In this entire section we fix an integer $N\geq2$ which is the 
  ambient dimension {\D} and a unital commutative ring $\BK$.\q\  
% We assume in all this section that $1/2\in\BK$\p\pfn{this is needed to
%   guarantee
% that the space of diagrams is a free module}.
  The case $N=1$ is somewhat special, although trivial, and will be
treated separately in \refS{proofform}.
\subsection{Diagrams}
\label{sec:diagrams1}
Roughly speaking, a  \emph{diagram} is a finite oriented graph
where the vertices come in two flavors,   \emph{external} and \emph{internal}, and where the 
sets of edges and internal vertices  are linearly ordered. An example 
is given in Figure \ref{fig:diagram} and explained in \refX{diag} below.
The  precise definition is as follows.
\begin{defin}\label{D:diag}
A \emph{diagram} is a quintuple $\Gamma=({\ExtVert}_\Gamma,I_\Gamma,E_\Gamma,s_\Gamma,t_\Gamma)$ where \q
\begin{itemize}
\item ${\ExtVert}_\Gamma$ is a finite set;
\item $I_\Gamma$ is a linearly ordered finite set disjoint from ${\ExtVert}_\Gamma$;
\item  $E_\Gamma$ is a   linearly ordered finite set; and
\item $s_\Gamma,t_\Gamma\colon E_\Gamma\to {\ExtVert}_\Gamma\amalg I_\Gamma$ are functions.
\end{itemize}
\previousoldfn{\pl 080611 I removed this condition: ``We require moreover that  $I_\Gamma$ is empty when $|{\ExtVert}_\Gamma|\leq1$.''}
%We say that $\Gamma$ is \emph{diagram on  $A_\Gamma$}.\\

We fix the following terminology and notation:
\begin{itemize}
\item the elements of ${\ExtVert}_\Gamma$ are the \emph{external vertices}, 
 the elements of $I_\Gamma$ are the \emph{internal vertices}, and we set
 $V_\Gamma:={\ExtVert}_\Gamma\amalg I_\Gamma$; this is the set of all \emph{vertices}. We extend the order of $I_\Gamma$ to a partial
order on $V_\Gamma$ by letting $a<i$ when $a\in {\ExtVert}_\Gamma$  and $i$ $\in I_\Gamma$;
\item the elements of $E_\Gamma$ are the \emph{edges};
\item $s_\Gamma(e)$ is the \emph{source} and $t_\Gamma(e)$ is the \emph{target} of the edge $e$;
both are the \emph{endpoints} of the edge;
\item two distinct vertices{\D} are called \emph{adjacent} if they are the endpoints of some edge;
%when endpoints of an edge are distinct, we say that they are \emph{adjacent};
\item 
we  say that the edge $e$ is \emph{oriented} from $s_\Gamma(e)$ to $t_\Gamma(e)$;
\item we partition the set of edges into the following four families:
\begin{itemize}
\item a \emph{loop} is an edge whose endpoints are identical;
\item a \emph{chord} is an edge between two distinct external vertices;
\item a \emph{dead end} is an edge that is  not a loop
and such that at least one if its
endpoints is internal and has only one adjacent vertex\p\pfn{I
  went back to the older  defintion. The referee comment (67) is wrong};
\item a \emph{contractible edge} is an edge that is 
neither a chord, nor a loop, nor a dead end;
\end{itemize}
\item we denote by $\Econtr_\Gamma$ the set of contractible edges of $\Gamma$;
\item the \emph{valence} of a vertex is the number of edges for which the vertex is
an endpoint, with loops adding two to the valence;
\item an edge $e$ is \emph{simple} if there exists no other edge with the same set of endpoints;
\item \emph{double edges} are distinct edges having the same set of endpoints,
that is, a pair $\{ e_1,e_2\}$ 
such that $\{s_\Gamma(e_1),t_\Gamma(e_1)\}=\{s_\Gamma(e_2),t_\Gamma(e_2)\}$;
\item two vertices $v$ and $w$ are \emph{connected} if there exists a path of edges joining them 
(ignoring orientations), that is, if there exists a sequence of edges
$e_1,\dots,e_k$ such that $v\in\{s_\Gamma(e_1),t_\Gamma(e_1)\}$,  $w\in\{s_\Gamma(e_k),t_\Gamma(e_k)\}$, and
$\{s_\Gamma(e_i),t_\Gamma(e_i)\}\cap\{s_\Gamma(e_{i+1}),t_\Gamma(e_{i+1})\}\not=\emptyset$
for $1\leq i<k$;
\item given a finite \Z\previousoldfn{\pl 17jan2011:interchaneg items} set $A$, a \emph{diagram on $A$} is a diagram $\Gamma$ such that $A_\Gamma=A$;
\item a diagram on $A$ \Z is a \emph{unit} if it has no internal vertices or edges. We denote a unit by $\unit$.
In other words $\unit=(A,\emptyset,\emptyset,\emptyset,\emptyset)$;
\previousoldfn{\pl also removed: `` Notice that by definition 
the only diagrams with $0$ or $1$ external vertex are the units.''}\previousoldfn{\pl check latter
the compatibility of this restriction on diagrams on $0$ or $1$ external vertex 
with the cooperad structure.\pl 080611 No need anymore}
\item two  diagrams  $\Gamma$ and  $\Gamma'$ are 
\emph{isomorphic} if ${\ExtVert}_\Gamma={\ExtVert}_{\Gamma'}$ and there exist two order-preserving bijections
$\phi_E\colon E_\Gamma\iso E_{\Gamma'}$
and $\phi_I\colon I_\Gamma\iso I_{\Gamma'}$ (that we extend into a bijection $\phi_V:=\id_{A_\Gamma}\amalg\phi_I\colon V_\Gamma\iso V_{\Gamma'}$)
such that $\phi_V\circ s_\Gamma=s_{\Gamma'}\circ\phi_E$ and $\phi_V\circ t_\Gamma=t_{\Gamma'}\circ\phi_E$.
\end{itemize}
\end{defin}

%We set $V_\Gamma:=A_\Gamma\cup I_\Gamma$ (notice that $A_\Gamma\cap I_\Gamma=\emptyset$.) \\

We will abuse notation by denoting a diagram and its isomorphism
class by the same letter $\Gamma$.

\p\q One should be careful about the definition of a \emph{dead end}. Our
definition is not equivalent to saying that a dead end is an edge with
a univalent internal vertex. Indeed in Example \ref{X:diag} 
and \refF{diagram} below,
the edge $(12,14)_1$ is a dead end (because the vertex $14$ is
internal and has only $12$ as an adjacent vertex), although neither of its
endpoints $12$ and $14$ is of valence $1$\pfn{this comment seems
  needed because of mistake (67) of the referee}. However, when a
diagram has no loops or double edges, dead ends appear only
with univalent internal vertices.\pfn{check the previous paragraph phrasing.\iv Checked and changed slightly.}
The reason for distinguishing dead ends from contractible edges will be given in
Remark \ref{rmk:whydeadend}
\previousoldfn{\pl Nathalie Wahl asked 
why dead ends are not contractible? This is an answer. Is the reference of the remark correct? Yes!}
% \begin{remark}\label{R:pictdiag}
% \textbf{Pictured diagrams.} 
% We will represent certain
% diagrams by a pictured graph with vertices labeled by some obviously ordered
% set (usually the
%  set $\{1,\dots,n+q\}$), with all external vertices drawn on a horizontal line
%  that is not part of the graph and internal vertices outside of that line.
% An example of this is given in Figure \ref{fig:diagram} that we will explain in \refX{diag}.
% By convention we assume that in such a pictured diagram
% each edge 
% is oriented from the lowest to the highest vertex,
% i.e.~$s(e)\leq t(e)$, and that the edges are ordered by the 
% right\previousfn{\pl the reason to choose the right lexicographic 
% order is that it fits well with the definition of the product of diagrams; see below} 
% lexicographic order
% on $\{(s(e),t(e)\}$, i.e.~$e_1\leq e_2$ if $t(e_1)<t(e_2)$ or if $t(e_1)=t(e_2)$
% and $s(e_1)\leq s(e_2)$. A priori this leaves some ambiguity when there are double edges,
% because if $e_1$ and $e_2$ are two edges from $x$ to $y$, then this lexicographic order
% does not determine the strict order between these two edges, but this will be irrelevant because
% transposing the order of these two edges yields to an isomorphic diagram.
% Also we will usually abuse notation by denoting an edge $e$ 
% just by the couple $(s(e),t(e))$ even if there is again some ambiguity in the presence of double
% edges. In case of  need we will     distinguish these double edges by some additional index, like
% $(s(e),t(e))_i$.
% \end{remark}
\previousoldfn{\pl The remark explaining pictured diagram is suppressed 20 june.}

\begin{eg}\label{X:diag}
Consider the diagram  in \refF{diagram}. By
convention, all the external vertices are drawn on a horizontal  line which is not a part of the graph.
This picture
represents a diagram $\Gamma$ with\previousoldfn{\pl the picture of the diagram has to be changed. \pl: done}
\begin{itemize}
\item the set of external vertices ${\ExtVert}_\Gamma=\{1,\dots,5\}$;
\item the set of internal vertices $I_\Gamma=\{6, \dots,15\}$ with its natural order;
\item the set $E_\Gamma$ consists of eighteen edges, each one oriented from the lower to the higher vertex
and  ordered as follows (right lexicographic order):\q
 \vskip 4pt
\begin{tabular}{l}
$(3,4)<(1,6)<(2,6)<(3,7)<(6,7)_1<(6,7)_2<(7,8)<$\\
$<(8,8)<(8,9)<(4,10)<(5,10)<(11,12)<(11,13)<$\\
$<(12,13)<(12,14)_1<(12,14)_2<(14,14)_1<(14,14)_2.$
  \end{tabular}
  
  \vskip 4pt
  \noindent
There are three loops, at vertices $14$ and  $8$;
three dead ends,  $(8,9)$, $(12,14)_1$ and $(12,14)_2$;  a chord $(3,4)$;
double contractible edges $(6,7)_1$ and $(6,7)_2$\previousoldfn{\iv Should $(12,14)_1$ and  $(12,14)_2$ also be listed here or are we just giving examples.  It reads like we're listing everything.\pl fixed};
and nine other simple contractible edges. 
The valence of the vertex $3$ is $2$, that of  $8$ is $4$, that of $14$ is $6$,  that of $15$ is $0$, etc.
\end{itemize}
\end{eg}

%\begin{figure}\centering
%\fbox{\includegraphics[width=80mm,height=45mm]{ExampleDiagram.pdf}}
%\caption{An example of a diagram (see \refX{diag}.)}
%\label{fig:diagram}
%\end{figure}

\begin{figure}[h]
\input{ExampleDiagram.pstex_t}
\caption{An example of a diagram (see \refX{diag}.)}
\label{fig:diagram}
\end{figure}

\begin{remark}\label{R:isodiag}
% The set of external vertices will  usually be taken as  the set $X=\{1,\dots,n\}$.
% Any diagram  with this set of external vertices is isomorphic to a diagram $\Gamma$
% for which $I_\Gamma=\{n+1,\dots,n+q\}$ where $q=|I_\Gamma|$ is the number of internal vertices.\\
% On the other hand,
Given\previousoldfn{\pl The first part of the remark is commented 20june}
two diagrams $\Gamma_1$ and $\Gamma_2$ 
with the same set of external vertices, we can always find a diagram $\Gamma'_2$
isomorphic to $\Gamma_2$  such that 
the sets $I_{\Gamma_1}$and  $I_{\Gamma'_2}$, and $E_{\Gamma_1}$ and $E_{\Gamma'_2}$ respectively,
are disjoint. This will be used in the definition of the 
product of two (isomorphism classes of) diagrams in \refS{proddiag}.
\p Also, if $A$ and $P$ are disjoint sets and $\Gamma$ is a diagram on
$A$,
we can assume (after maybe replacing diagram $\Gamma$ by an isomorphic
one) 
that $P$ and $I_\Gamma$ are disjoint\p.
\end{remark}

\subsection{The module $\GD({\ExtVert})$ of  diagrams}
\label{sec:modulediagrams}
We will define the $\BK$-module generated by isomorphism classes of diagrams modulo some
signed relations when the linear order of internal vertices or edges is permuted, or the 
orientation of some edge is reversed. To make this precise, we need the following:

\begin{defin}\label{D:differs}
  Let $\Gamma$ and $\Gamma'$ be two diagrams with the same  set of external
  vertices.
  \begin{itemize}
  \item  $\Gamma$ and $\Gamma'$ \emph{differ  by an
      inversion of an edge} if, up to isomorphism, these two diagrams
have the same ordered sets of  internal vertices and edges, there
exists an edge $e$ such that $s_{\Gamma'}(e)=t_{\Gamma}(e)$ and $t_{\Gamma'}(e)=s_{\Gamma}(e)$,
and $s_{\Gamma}$ and $s_{\Gamma'}$ (respectively, $t_{\Gamma}$ and $t_{\Gamma'}$)
agree on all the other edges.
  \item $\Gamma$ and $\Gamma'$ \emph{differ  by a
      transposition in the linear order of internal vertices}, if, up to isomorphism, 
 they    have the same ordered set of edges, the same {underlying} set
    of internal vertices $I$, the same source and target functions,
    and there exists a transposition $\sigma=(a,b)$ in the group of
    permutations of the set $I$, for some pair of distinct internal
    vertices $a$ and $b$, such that for all internal vertices
    $i_1,i_2\in I$ we have that $i_1\leq_{I_\Gamma} i_2$ if and only
    if $\sigma(i_1)\leq_{I_{\Gamma'}}\sigma(i_2)$.
  \item $\Gamma$ and $\Gamma'$ \emph{differ  by a
      transposition in the linear order of the edges}, if, up to isomorphism, 
 they have the
    same ordered set of internal vertices, the same {underlying}
    set of edges $E$, the same source and target functions, and there
    exists a transposition $\sigma=(a,b)$ in the group of permutations
    of the set $E$, for some pair of distinct edges $a$ and $b$, such
    that, for all edges $e_1,e_2\in E$, $e_1\leq_{E_\Gamma}
    e_2$ if and only if $\sigma(e_1)\leq_{E_{\Gamma'}}\sigma(e_2)$.
  \end{itemize}
\end{defin}

\begin{defin}\label{D:spacediag}  Fix an integer $N\geq1$.
The \emph{space of  diagrams on a set  $\ExtVert$} is the free $\BK$-module 
$\GD({\ExtVert})$  generated by the isomorphism classes of diagrams with the set of external vertices
 ${\ExtVert}$, modulo the equivalence relation $\simeq$ generated by the following:
 \begin{itemize}
 \item $\Gamma\simeq(-1)^{N}\Gamma'$  if $\Gamma$ and  $\Gamma'$ differ
by an inversion of an edge;
\item $\Gamma\simeq(-1)^N\Gamma'$ if $\Gamma$ and  $\Gamma'$ differ
by  a transposition in the linear order of internal vertices;
\item $\Gamma\simeq(-1)^{N+1}\Gamma'$  if $\Gamma$ and  $\Gamma'$ differ
by   a transposition in the linear order of  edges.
 \end{itemize}
When we want to emphasize the ambient dimension $N$, we will denote the space of diagrams
by $\GD_N({\ExtVert})$.

By abuse of notation we will denote by the same symbol a diagram and its
equivalence class in $\GD({\ExtVert})$.
\end{defin}

Because of the relations, when $N$ is odd (respectively even) and $1/2\in\BK$\p, a diagram with a loop (respectively a double edges) vanishes in the space of diagrams. Other symmetries \D
 of a diagram can also make it  vanish.
Also because of the{\D} relations, when $N$ is even the orientation of
the edges and the linear order on internal vertices are irrelevant;
when $N$ is odd it is the  linear order on the  edges which is irrelevant.
When \p $1/2\in\BK$,  $\GD(A)$ is a free $\BK-$module generated by a
suitable collection of diagrams (if $1/2\not\in\BK$
the relation $\simeq$  produces $2$-torsion.)

\begin{defin}\label{D:degdiag}
The \emph{degree} of a diagram $\Gamma$ is defined to be
\[\deg(\Gamma)=|E_\Gamma|\cdot(N-1)-|I_\Gamma|\cdot N\]
where $|E_\Gamma|$ is the number of edges and $|I_\Gamma|$ is the number of internal vertices.
\end{defin}
\previousoldfn{\pl maybe there is a problem with de definition of degree when $|{\ExtVert}|\leq1$. What about the diagram $\Gamma$
with one ext vert, one int vert, and one edge? Is its degree $0$ or $-1$. Maybe OK because when $|{\ExtVert}|\leq1$ we assume that there is no edge. \iv In the previous version, I think we only defined degree by this formula on admissible diagrams because of this problem.  I any case, I believe we treated the degree of non-admissible diagrams separately because the degree of the form we get from Kontsevich integration doesn't satisfy this formula for non-admissible diagrams.\pl I guess all is ok now. The diagram $\Gamma$
 is of degree -1
but is not admissible. $\GIK(\Gamma)=0$ by definition since $|A|\leq1$.}

The motivation for defining $\deg(\Gamma)$ as such
comes from the fact that in \refS{KCSI} we will construct
a differential form $\IK(\Gamma)\in\ompa(\Conf[A])$ whose degree is exactly that.
The integration producing this form is what motivates the 
signs in \refD{spacediag}.

The degree is compatible with the equivalence relation $\simeq$, and so $\GD({\ExtVert})$ becomes a graded $\BK$-module.

\subsection{Product of diagrams}\label{sec:proddiag}
Let $\Gamma_1$ and $\Gamma_2$  be two isomorphism classes 
of diagrams on the same  set ${\ExtVert}$.
By \refR{isodiag}, we can assume  that 
the sets $I_{\Gamma_1}$and  $I_{\Gamma_2}$, and $E_{\Gamma_1}$ and $E_{\Gamma_2}$ respectively,
are disjoint. Remember the sum of  linearly ordered sets $\ordsum$ defined in \refS{linord}.
Define the product diagram $\Gamma=\Gamma_1\cdot\Gamma_2$
by
\begin{itemize}
\item ${\ExtVert}_\Gamma:={\ExtVert}$;
\item $ I_{\Gamma}:=I_{\Gamma_1}\ordsum I_{\Gamma_2}$;
\item $ E_{\Gamma}:=E_{\Gamma_1}\ordsum E_{\Gamma_2}$;
\item $s_\Gamma|E_{\Gamma_i}=s_{\Gamma_i}$ and  $t_\Gamma|E_{\Gamma_i}=t_{\Gamma_i}$.
\end{itemize}

\begin{eg}\label{X:proddiag}
An example of a product of two isomorphism classes of  diagrams is represented in 
\refF{prodgendiag}.
% where we use the pictural convention of \refR{pictdiag}. Notice that 
% our convention of the \emph{right} lexicographic order on the edges
% makes that the multiplication can always 
% be represented as follows: Given two diagrams on the same set $\ExtVert$ of external vertices and
% represented by pictured graphs, their product is represented by the pictured 
% graph obtained by gluing the pictures of the two diagrams identifying their external vertices
% and by relabelling the internal vertices of the second graph so that they appear  in the linear order after those of the first graph.
In each picture the edges are oriented from the lower-labeled to the higher-labeled vertex and are ordered by the right lexicographic order
as in \refX{diag}.

%\begin{figure}
%\fbox{\includegraphics[width=90mm,height=32mm]{ProductDiagrams.pdf}}
%\caption{Example of a product of two diagrams}\label{fig:prodgendiag}
%\end{figure}

\begin{figure}[h]
\input{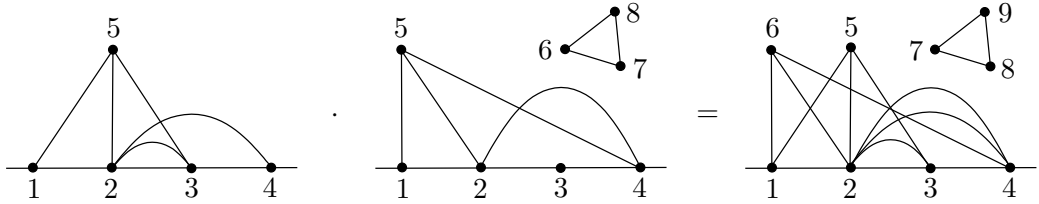}
\caption{Example of a product of two diagrams.}
\label{fig:prodgendiag}
\end{figure}

\end{eg}

\previousoldfn{\pl To add here: that there is a unit and that it is graded commutative; see NB20march2008 page 48. \pl done}
  
\begin{prop}\label{P:GD-CDGA}
The above product extends to a degree $0$ linear map
\[\GD({\ExtVert})\otimes\GD({\ExtVert})\longrightarrow\GD({\ExtVert})\] 
 which endows 
$\GD({\ExtVert})$ with the structure of  a commutative $\BZ$-graded algebra.
\end{prop}
\begin{proof}
 The multiplication has been defined on generators and we extend it bilinearly. This multiplication is compatible with the equivalence relation $\simeq$ on diagrams. It is also clearly associative and
$$\deg(\Gamma_1\cdot\Gamma_2)=\deg(\Gamma_1)+\deg(\Gamma_2).$$
The unit diagram $\unit=({\ExtVert},\emptyset, \emptyset, \emptyset, \emptyset)$
is of degree $0$ and is indeed a unit for the product.\\
It remains to check that the multiplication is graded-commutative.
Let $\Gamma_i=({\ExtVert},I_i,E_i,s_i,t_i)$, for $i=1,2$, be two diagrams. We distinguish two cases.
\begin{itemize}
\item Suppose that $N$ is odd. The diagrams $\Gamma_1\cdot\Gamma_2$ and
$\Gamma_2\cdot\Gamma_1$ differ by the order of the edges, which is irrelevant in
this case, and the order of internal vertices. The number of pairs of transposed vertices is
$|I_1|\cdot|I_2|$. Since $N$ is odd, $|I_i|\equiv\deg(\Gamma_i)\mod2$.
Therefore 
$\Gamma_2\cdot\Gamma_1=
(-1)^{\deg(\Gamma_1)\cdot\deg(\Gamma_2)}\Gamma_1\cdot\Gamma_2$.
\item Suppose that $N$ is even. The argument is the same as for $N$ odd after exchanging the
roles  of the linear orders of edges and the internal vertices.
\end{itemize}
\end{proof}

\subsection{A differential on the space of diagrams}
\label{sec:ddiag}

We define now a differential on the $\BK$-module $\GD(A)$ by ``contracting edges'' on diagrams.
Recall from \refD{diag} the notion of a contractible edge in a diagram.
Also remember that in \refD{diag} we extended the linear order on internal vertices{\D} into a partial order on the set of all vertices by making external vertices precede internal ones. In particular, if $e$ is a contractible edge
of a diagram $\Gamma$ then the pair $\{s_\Gamma(e),t_\Gamma(e)\}$ is
linearly ordered.
\begin{defin}
\label{D:Gamma/e}
Let $\Gamma$ be a diagram and let $e$ be a contractible edge of $\Gamma$.
The diagram obtained from $\Gamma$ by \emph{contraction of the edge $e$} is the diagram
${\Gamma/e}$
% denoted by $\Gamma/e$
\previousoldfn{\iv Why are we giving it two notations?\pl I removed one notation $\overline\Gamma$}
 defined as follows:
\begin{itemize}
\item ${\ExtVert}_{{\Gamma/e}}:={\ExtVert}_\Gamma$
\item $I_{{\Gamma/e}}:=I_\Gamma\setminus\{\max(s_\Gamma(e),t_\Gamma(e))\}$
\item $E_{{\Gamma/e}}:=E_\Gamma\setminus{\{e\}}$
\item $s_{{\Gamma/e}}:=q\circ s_\Gamma$ and $t_{{\Gamma/e}}:=q\circ t_\Gamma$ where $q$ is defined by:
\begin{align*}
q\colon V_{\Gamma} & \longrightarrow V_{{\Gamma/e}} \\
 v & \longmapsto\left\{ \begin{array}{ll}
      \min(s_\Gamma(e),t_\Gamma(e)),&\textrm{if }v=\max(s_\Gamma(e),t_\Gamma(e));\\
      v,&\textrm{otherwise},
    \end{array}\right.
\end{align*}
\end{itemize}
where the linear orders on  $I_{{\Gamma/e}}$ and $E_{{\Gamma/e}}$ are the restrictions of those
on  $I_{\Gamma}$ and $E_{\Gamma}$.
\end{defin}
Notice that ${\Gamma/e}$ is well-defined because $\max(s_\Gamma(e),t_\Gamma(e))$ is internal since
$e$ is not a chord, and $\min(s_\Gamma(e),t_\Gamma(e))\not=\max(s_\Gamma(e),t_\Gamma(e))$
since $e$ is not a loop.

When $e'$ is an edge distinct from a contractible edge $e$, we will denote by $\overline{e'}$
the edge of $\Gamma/e$ corresponding to $e'$ in $\Gamma$ through the inclusion 
$E_{{\Gamma/e}}\hookrightarrow E_\Gamma$.

\begin{eg}
 \label{X:contrdiag}
An example of contraction of an edge is given in Figure \ref{fig:contrdiag}
(where we omit precise ordering and orientation of the edges).
% Notice that a sign $\pm$ appears because our convention on representing
% diagrams by pictured graphs is not necessarily
%  compatible with the contraction of an edge.\previousfn{\pl problem?}
The edge $\overline{(8,7)}$ is the edge $(8,4)$ in the diagram after
contraction of the edge $(4,7)$.\previousoldfn{\pl check with precise \refD{Gamma/e}}
\previousoldfn{\pl \todo the picture need to be changed as in corrected draft of 16June page 8 left. \pl Done on 24 june.}
%\begin{figure}
%\fbox{\includegraphics[width=110mm,height=48mm]{ContrDiagram.pdf}}
%\caption{Contraction of an edge}\label{fig:contrdiag}
%\end{figure}

\begin{figure}[h]
\input{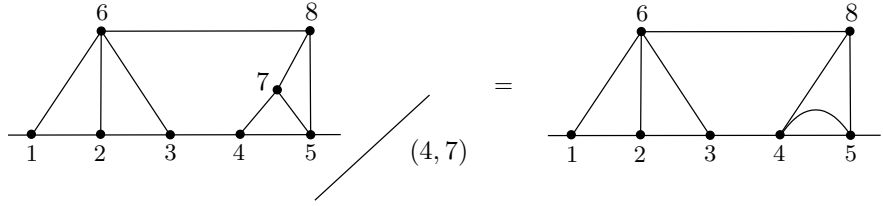}
\caption{Contraction of an edge.}
\label{fig:contrdiag}
\end{figure}

\end{eg}

%before

Before defining 
 the differential $d$, we need to introduce some signs.
Recall from \refS{linord} the position function   $\pos$ associated to  linearly ordered sets such as $I_\Gamma$ and $E_\Gamma$.
Define
$\epsilon(\Gamma,e)=\pm1$ according to the following table (where
$s:=s_\Gamma$ and $t:=t_\Gamma$ are the source and the target of edges)\p
\previousoldfn{\iv I don't understand the notation in the exponent in the table.\pl I added a recall of $\pos$. Is it better? \iv yes.}
\previousoldfn{\pl see NB7E: page 26 and NB:8C page 26}
\begin{center}
\begin{tabular}{|c|c|}
\hline
\multicolumn{2}{|c|}{Value of $\epsilon(\Gamma,e)$ }\\
\hline
$N$ odd&$N$ even\\
  \hline
$\quad (-1)^{\pos(\max(s(e),t(e)):I_\Gamma)}$ if $s(e)<t(e)$&$(-1)^{\pos(e:E_\Gamma)}$\\
$-(-1)^{\pos(\max(s(e),t(e)):I_\Gamma)}$ if $s(e)>t(e)$&\\
\hline
\end{tabular}
\end{center}
Let $\Gamma$ be a diagram on a set of external vertices ${\ExtVert}$.
Define its differential $d(\Gamma)\in\GD({\ExtVert})$ by the formula
\begin{equation}
\label{eq:d}
d(\Gamma):=\sum_{e\in\Econtr_\Gamma}\epsilon(\Gamma,e)\cdot\Gamma/e
\end{equation}
where the sum runs over all contractible edges $e$ in $\Gamma$ and the sign
$\epsilon(\Gamma,e)$ is from the{\D} above table.
An example of this  is the diagram $\Gamma$ in \refF{killer3terms} of the Introduction
for which $d(\Gamma)$ is the diagram of \refF{3terms} with $k=3$, $(i,j,l)=(1,2,3)$, and with a suitable orientation
and ordering of the edges.

\begin{lemma}\label{L:d}
Formula \refN{eq:d} induces a linear map
$d\colon\GD({\ExtVert})\to\GD({\ExtVert})$.
\end{lemma}

\begin{proof}
\q The proof that $d$ is
compatible with the equivalence relation $\simeq$ of
\refD{spacediag} is straightforward (but tedious). We give the sketch of the argument in one case (the one for which the
proof is more complicated) and leave the others to the
reader. Proofs of analogous statements can be find in 
\cite[Theorem 4.2]{CCRL:Vas}.

\p\pfn{this longer proof answer comment (68) of referee. \iv This is essentially Theorem 4.2 in Cattaneo et. al., but I think it's good to have this proof here as well, especially if the referee asked for it.}
  
Assume that $N$ is odd  and let $\Gamma$ and $\Gamma'$ be two diagrams
that
 differ by a transposition in the linear order of internal vertices.
Then $\Gamma\simeq-\Gamma'$.

 Let 
 $a$ and $b$ be the transposed vertices.
Since any transposition of the linear order of internal vertices is  obtained as a composition of transpositions
of adjacent vertices, we can assume without loss of generality  that $a$ and $b$ are consecutive in $V_\Gamma$,
so 
$$\pos(b:I_\Gamma)=\pos(a:I_\Gamma)+1.$$ 
Moreover, since it is easy to check that the differential is
compatible with inversion of orientations of edges,
we can assume that each edge $e$ of $\Gamma$ is oriented so that
 $s_\Gamma(e)\leq t_\Gamma(e)$ in $V_\Gamma$.

 Let $e$ be a contractible edge. 
We need to show that 
$$\epsilon(\Gamma,e)\cdot\Gamma/e\simeq-\epsilon(\Gamma',e)\cdot\Gamma'/e.$$
 We distinguish three cases.
\begin{enumerate}
\item Suppose that $\{s_\Gamma(e),t_\Gamma(e)\}=\{a,b\}$.
Since $s(e)$ and $t(e)$ are permuted, 
we have that  $\epsilon(\Gamma,e)=-\epsilon(\Gamma',e)$. On the other hand, $\Gamma/e\simeq\Gamma'/e$ since one of the two 
consecutive vertices $a$ or $b$ disappears.
\item Suppose that $\{s_\Gamma(e),t_\Gamma(e)\}\cap\{a,b\}=\emptyset$.
In that case $\epsilon(\Gamma,e)=\epsilon(\Gamma',e)$ but  $\Gamma/e\simeq-\Gamma'/e$. 
\item Suppose that $\{s_\Gamma(e),t_\Gamma(e)\}\cap\{a,b\}$ is a singleton.
We have assumed that $a$ and $b$ are consecutive in $V_\Gamma$ and that $s_\Gamma(e)\leq t_\Gamma(e)$.
Moreover $s_\Gamma(e)\not= t_\Gamma(e)$ since $e$ is contractible.
Therefore we have four possibilities for the order of $a,b,s_\Gamma(e),t_\Gamma(e)$ in $V_\Gamma$:
\begin{enumerate}
\item $a=s_\Gamma(e)<b<t_\Gamma(e)$.
Then $\epsilon(\Gamma,e)=\epsilon(\Gamma',e)$ and $\Gamma/e\simeq-\Gamma'/e$.
\item $a<s_\Gamma(e)=b<t_\Gamma(e)$.
Then  $\epsilon(\Gamma,e)=\epsilon(\Gamma',e)$ and $\Gamma/e\simeq-\Gamma'/e$.
\item $s_\Gamma(e)<a<b=t_\Gamma(e)$.
Then  $\epsilon(\Gamma,e)=-\epsilon(\Gamma',e)$ and $\Gamma/e\simeq\Gamma'/e$.
\item $s_\Gamma(e)<a=t_\Gamma(e)<b$. 
Then $\epsilon(\Gamma,e)=-\epsilon(\Gamma',e)$ and $\Gamma/e\simeq\Gamma'/e$.
\end{enumerate}
\end{enumerate}
In all cases we have $\epsilon(\Gamma,e)\cdot\Gamma/e=-\epsilon(\Gamma',e)\cdot\Gamma'/e$.
This proves that $d(\Gamma)=-d(\Gamma')$ as desired.
\optproof[\refL{d}]{
({\pl See NB8C page 38-43 et p.27.})\\
It is clear that $d(\Gamma)$ depends only on the isomorphism class of $\Gamma$. We
need to check that it is compatible with the equivalence relation $\simeq$ of \refD{spacediag}.
Let $\Gamma$ and $\Gamma'$ be two diagrams.
\\({\pl Maybe to shorter the proof
we should split first among the cases $N$ odd/even as in the proof of \refL{d20} etc.})\\
\begin{enumerate}
\item {\bf Suppose that $\Gamma$ and $\Gamma'$ differ by an inversion of an edge.}
\begin{itemize}
\item Suppose $N$ even.
Then $\Gamma\simeq\Gamma'$. Also for each contractible edge $e$, we have $\Gamma/e\simeq\Gamma'/e$ and 
$\epsilon(\Gamma,e)=\epsilon(\Gamma',e)$. Therefore $d(\Gamma)=d(\Gamma')$.
\item Suppose $N$ odd.
Then   $\Gamma\simeq-\Gamma'$. Let $e_0$ be the edge that is reversed.
If $e_0$ is contractible, we get that $\Gamma/e_0\simeq\Gamma'/e_0$ but $\epsilon(\Gamma,e_0)=-\epsilon(\Gamma',e_0)$.
On the other hand for a contractible edge $e\not=e_0$ we have $\Gamma/e\simeq-\Gamma'/e$ and
$\epsilon(\Gamma,e)=\epsilon(\Gamma',e)$.
Therefore, for any contractible edge $e$ we have $\epsilon(\Gamma,e)\cdot\Gamma/e=-\epsilon(\Gamma',e)\cdot\Gamma'/e$,
and hence $d(\Gamma)=-d(\Gamma')$ as expected.
\end{itemize}
\item {\bf Suppose that $\Gamma$ and $\Gamma'$ differ by a transposition in the linear order of internal vertices.} Let 
 $a$ and $b$ be the transposed vertices.\\({\pl should change $a$ $b$ to $i$ $j$})
Since any transposition of the linear order of internal vertices is  obtained as a composition of transpositions
of adjacent vertices, we can assume without loss of generality  that $a$ and $b$ are consecutive in $V_\Gamma$,
$\pos(b:I_\Gamma)=\pos(a:I_\Gamma)+1$. Also, by Case (1) above in the proof,
we can assume that each edge $e$ of $\Gamma$ is oriented so that
 $s_\Gamma(e)\leq t_\Gamma(e)$ in $V_\Gamma$.
\begin{itemize}
\item Suppose $N$ even.
Then $\Gamma\simeq\Gamma'$.
Let $e$ be a contractible edge. As $N$ is even,  $\epsilon(\Gamma,e)=\epsilon(\Gamma',e)$.
Also $\Gamma/e$ is the same as $\Gamma'/e$ except maybe for the order of vertices, and hence
$\Gamma/e\simeq\Gamma'/e$. Thus $d(\Gamma)=d(\Gamma')$.
\item Suppose $N$ odd.
Then  $\Gamma\simeq-\Gamma'$. Let $e$ be a contractible edge. 
We need to show that $\epsilon(\Gamma,e)\cdot\Gamma/e\simeq-\epsilon(\Gamma',e)\cdot\Gamma'/e$.
 We distinguish three cases.
\begin{enumerate}
\item Suppose that $\{s_\Gamma(e),t_\Gamma(e)\}=\{a,b\}$.
Since $s(e)$ and $t(e)$ are permuted, 
we have that  $\epsilon(\Gamma,e)=-\epsilon(\Gamma',e)$. On the other hand $\Gamma/e\simeq\Gamma'/e$ since one of the two 
consecutive vertices $a$ or $b$ disappears.
\item Suppose that $\{s_\Gamma(e),t_\Gamma(e)\}\cap\{a,b\}=\emptyset$.
In that case $\epsilon(\Gamma,e)=\epsilon(\Gamma',e)$ but  $\Gamma/e\simeq-\Gamma'/e$. 
\item Suppose that $\{s_\Gamma(e),t_\Gamma(e)\}\cap\{a,b\}$ is a singleton.
We have assumed that $a$ and $b$ are consecutive in $V_\Gamma$ and that $s_\Gamma(e)\leq t_\Gamma(e)$.
Moreover $s_\Gamma(e)\not= t_\Gamma(e)$ since $e$ is contractible.
Therefore we have four possibilities for the order of $a,b,s_\Gamma(e),t_\Gamma(e)$ in $V_\Gamma$:
\begin{enumerate}
\item $a=s_\Gamma(e)<b<t_\Gamma(e)$.
Then $\epsilon(\Gamma,e)=\epsilon(\Gamma',e)$ and $\Gamma/e\simeq-\Gamma'/e$.
\item $a<s_\Gamma(e)=b<t_\Gamma(e)$.
Then  $\epsilon(\Gamma,e)=\epsilon(\Gamma',e)$ and $\Gamma/e\simeq-\Gamma'/e$.
\item $s_\Gamma(e)<a<b=t_\Gamma(e)$.
Then  $\epsilon(\Gamma,e)=-\epsilon(\Gamma',e)$ and $\Gamma/e\simeq\Gamma'/e$.
\item $s_\Gamma(e)<a=t_\Gamma(e)<b$. 
Then $\epsilon(\Gamma,e)=-\epsilon(\Gamma',e)$ and $\Gamma/e\simeq\Gamma'/e$.
\end{enumerate}
\end{enumerate}
In all cases we have $\epsilon(\Gamma,e)\cdot\Gamma/e=-\epsilon(\Gamma',e)\cdot\Gamma'/e$.
This proves that $d(\Gamma)=-d(\Gamma')$ as expected.
\end{itemize}
\item {\bf Suppose that $\Gamma$ and $\Gamma'$ differ by a transposition in the linear order of edges.}
Without loss of generality we can assume that the two edges that are permuted, $e_1$ and $e_2$,  are consecutive in $E_\Gamma$,
$\pos(e_2:E_\Gamma)=\pos(e_1:E_\Gamma)+1$.
\begin{itemize}
\item Suppose $N$ even.
Then $\Gamma\simeq-\Gamma'$. Let e be a contractible edge of $\Gamma$. 
\begin{enumerate}
\item If $e<e_1$ or $e>e_2$ then $\epsilon(\Gamma,e)=\epsilon(\Gamma',e)$ and $\Gamma/e\simeq-\Gamma'/e$.
\item If $e=e_1$ or $e=e_2$ then the linear order of the edges in $\Gamma/e$ and $\Gamma'/e$ are the same, so $\Gamma/e\simeq\Gamma'/e$.
But in that case $\epsilon(\Gamma,e)=-\epsilon(\Gamma',e)$ because $\pos(e:E_\Gamma)=\pos(e:E_{\Gamma'})\pm1$.
\end{enumerate}
In both cases $\epsilon(\Gamma,e)\cdot\Gamma/e=-\epsilon(\Gamma',e)\cdot\Gamma'/e$, and hence $d(\Gamma)=-d(\Gamma')$ as expected.
\item Suppose $N$ odd. 
In that case the order of edges is irrelevant. Therefore $\Gamma\simeq\Gamma'$. Also,
for each contractible edge $e$,  $\Gamma/e\simeq\Gamma'/e$,
and $\epsilon(\Gamma,e)=\epsilon(\Gamma',e)$ because $\epsilon$ does not depends on the order of the edges.
Thus  $d(\Gamma)=d(\Gamma')$.
\end{itemize}
\end{enumerate}
}
\end{proof}
\begin{lemma}\label{L:degred}
$d$ is homogeneous of degree $+1$.
\end{lemma}
\begin{proof}\previousoldfn{\pl NB 8Cp.43}
This is clear from \refD{degdiag} of degree since, for a contractible edge $e$ of a diagram $\Gamma$, the diagram $\Gamma/e$
has one fewer edge and one fewer internal vertex than $\Gamma$.
\end{proof}
\begin{lemma}\label{L:dLeibniz}
$d$ satisfies the Leibniz rule, that is, 
$$d(\Gamma\cdot\Gamma')=d(\Gamma)\cdot\Gamma'+(-1)^{\deg(\Gamma)}\Gamma\cdot d(\Gamma').$$
\end{lemma}
\begin{proof}
Recall that $E_{\Gamma\cdot\Gamma'}=E_\Gamma\ordsum E_{\Gamma'}$. It is clear
than an edge  is contractible in $\Gamma$ or  $\Gamma'$ if and only if it is
contractible in ${\Gamma\cdot\Gamma'}$. Moreover
if $e$ is a contractible edge of $\Gamma$ then $(\Gamma\cdot\Gamma')/e=(\Gamma/e)\cdot\Gamma'$,
and
if $e'$ is a contractible edge of $\Gamma'$ then $(\Gamma\cdot\Gamma')/e'=\Gamma\cdot(\Gamma'/e')$.
It remains to study the signs $\epsilon$ which appear in the differentials, which is straightforward.
\optproof[\refL{dLeibniz}]{
We distinguish two cases.
\begin{itemize}
\item {\bf $N$ even.} 
For a contractible edge $e$ of $\Gamma$ we have 
$\pos(e:E_{\Gamma} \ordsum E_{\Gamma'})=\pos(e:E_{\Gamma})$,
hence 
\begin{equation}\label{eq:dLeibnizEven1}
\epsilon(\Gamma\cdot\Gamma',e)=\epsilon(\Gamma,e).
\end{equation}
On the other hand, for
a contractible edge $e'$ of $\Gamma'$,
\[\pos(e':E_\Gamma\ordsum E_{\Gamma'})=\pos(e':E_{\Gamma'})+|E_\Gamma|.\]
Therefore
\[\epsilon(\Gamma\cdot\Gamma',e')=(-1)^{\pos(e':E_{\Gamma'})}\cdot(-1)^{|E_\Gamma|}\epsilon(E_{\Gamma'},e').\]
Since $N$ is even,  $\deg(\Gamma)\equiv|E_\Gamma|\mod 2$, and the latter equation
becomes 
\begin{equation}\label{eq:dLeibnizEven2}
\epsilon(\Gamma\cdot\Gamma',e')=(-1)^{\deg(\Gamma)}\cdot\epsilon(\Gamma',e').
\end{equation}
From Equations $(\ref{eq:dLeibnizEven1}-\ref{eq:dLeibnizEven2})$ we deduce immediately
the Leibniz formula.
\item {\bf $N$ odd.}
Without loss of generality we can assume that each edge of $\Gamma$ and of $\Gamma'$
 is oriented so that  its source is smaller than
its target.
In that case the target of a contractible edge is always an internal vertex. Also, for a contractible
edge $e$ of $\Gamma$,
\[\epsilon(\Gamma,e)=(-1)^{\pos(t_\Gamma(e):V_\Gamma)},\]
and similarly for $\Gamma'$ and $\Gamma\cdot\Gamma'$.
If $e$ is a contractible edge of $\Gamma$ we deduce that 
\begin{equation}\label{eq:dLeibnizEven3}
\epsilon(\Gamma\cdot\Gamma',e)=\epsilon(\Gamma,e).
\end{equation}
If $e'$ is a contractible edge of $\Gamma'$ then
$\epsilon(\Gamma\cdot\Gamma',e')=(-1)^{|I_\Gamma|}\cdot\epsilon(\Gamma',e')$,
which implies that
\begin{equation}\label{eq:dLeibnizEven4}
\epsilon(\Gamma\cdot\Gamma',e')
=
(-1)^{\deg(\Gamma)}\cdot\epsilon(\Gamma',e')
\end{equation}
because $\deg(\Gamma)\equiv|I_\Gamma|\mod 2$.
From Equations $(\ref{eq:dLeibnizEven3}-\ref{eq:dLeibnizEven4})$ we get
the Leibniz formula.
\end{itemize}
}%\optproof
\end{proof}
\begin{lemma}\label{L:d20}
$d^2=0$.
\end{lemma}
\begin{proof}
\previousoldfn{\pl Add the proof. cf NB7E 010108 pages 33-36.\pl Done in 080607}
\previousoldfn{\pl The proof should be rewritten more concisely by first separating  the two cases $N$ odd and even,
and inside each case treat the various possiblilities for $\simeq$.  \iv I think it's ok this way as well.} 
Let $\Gamma$ be a diagram and let $e_1$ and $e_2$ be distinct edges.
If $e_1$ is contractible, denote by $\overline{e_2}$ the edge in $\Gamma/e_1$
corresponding to $e_2$. It is easy to check that $\overline{e_2}$ is contractible 
in $\Gamma/e_1$ if and only the following two conditions hold:
\begin{itemize}
\item $e_1$ and  $e_2$ are contractible in $\Gamma$, and 
\item  $e_1$ and $e_2$ do not have the same endpoints, and  if $e_1$ and $e_2$ have one endpoint in common, then another
endpoint of $e_1$ or $e_2$ is an internal vertex. 
\end{itemize}
Since these conditions are symmetric, we deduce
that $\overline{e_2}$ is contractible in $\Gamma/e_1$
if and only if $\overline{e_1}$ is contractible  in $\Gamma/e_2$,
where $\overline{e_1}$ is the edge in $\Gamma/e_2$ corresponding to $e_1$ in $\Gamma$.
Moreover, in that case $(\Gamma/e_1)/\overline{e_2}$ is isomorphic to
$(\Gamma/e_2)/\overline{e_1}$. Therefore 
\begin{equation}\label{eq:d2G}
d^2(\Gamma)=\sum_{e_1<e_2}
\{\epsilon(\Gamma,e_1)\cdot\epsilon(\Gamma/e_1,\overline{e_2})
+
\epsilon(\Gamma,e_2)\cdot\epsilon(\Gamma/e_2,\overline{e_1})
\}\,\cdot\,
(\Gamma/e_1)/\overline{e_2},
\end{equation}
where the sum runs over each pair $e_1,e_2$ of distinct contractible edges of $\Gamma$
such that $e_1<e_2$ and the other condition above making $\overline{e_2}$
contractible in $\Gamma/e_1$ holds.
It is straightforward to check that the brackets in this sum vanish.
\optproof[\refL{d20}]{
\begin{itemize}
\item{\bf $N$ is even.}
Since $e_1<e_2$ we have
\begin{eqnarray*}
\pos(\overline{e_2}:E_{\Gamma/e_1})&=&\pos(e_2,E_\Gamma)-1\quad\textrm{and}\\
\pos(\overline{e_1}:E_{\Gamma/e_2})&=&\pos(e_1,E_\Gamma),
\end{eqnarray*}
which imply that the brackets are zero in \refE{d2G}.
\item{\bf $N$ is odd.}
Without loss of generality we can assume that for each  edge $e$ of $\Gamma$ we have
$s_\Gamma(e)\leq t_\Gamma(e)$. We can also assume that the edges $e$ of $\Gamma$ are ordered
by the left lexicographic order on $(s_\Gamma(e),t_\Gamma(e))$. Since $e_1<e_2$ this implies that
$s_\Gamma(e_1)\leq s_\Gamma(e_2)$ and if  $s_\Gamma(e_1)= s_\Gamma(e_2)$ then
$t_\Gamma(e_1)\leq t_\Gamma(e_2)$. For $\{i,j\}=\{1,2\}$, we set $s_i=s_\Gamma(e_i)$,
$t_i=t_\Gamma(e_i)$, $\overline{s_i}=s_{\Gamma/e_j}(\overline{e_i})$,
and $\overline{t_i}=t_{\Gamma/e_j}(\overline{e_i})$. 
% Also by $\pos(v)$, for some vertex $v$, we mean its position inside the
% set  of vertices to which $v$ naturally belongs, more precisely: In $V_\Gamma$ for $s_i$ and $t_i$,
% in $V_{\Gamma/e_1}$ for $\overline{s_2}$ and  $\overline{t_2}$, and
% in $V_{\Gamma/e_2}$ for $\overline{s_1}$ and  $\overline{t_1}$. 
Notice also that $\max(\overline{s_i},\overline{t_i})$ is internal in ${\Gamma/e_j}$.
%\oldfn
\private{\pl Check the end of the proof, from here. \pl Done}

By inspection we fill the following table which depends on the order of $t_1$ and $t_2$ in $V_\Gamma$.\\
\vspace{2mm}

\begin{tabular}{|c||c|c||c|c||}
\hline
&\multicolumn{2}{c||}{\quad\quad in $\Gamma/e_1$}&\multicolumn{2}{c||}{\quad\quad in $\Gamma/e_2$}\\
&$\pos(\max(\overline{s_2},\overline{t_2}):I_{{\Gamma/e_1}})$& Is $\overline{s_2}<\overline{t_2}$ ? &
$\pos(\max(\overline{s_1},\overline{t_1}):I_{{\Gamma/e_2}})$& Is $\overline{s_1}<\overline{t_1}$ ? \\
\hline
\hline
$t_1<t_2$&$\pos(t_2:I_{\Gamma})-1$&yes&$\pos(t_1:I_{\Gamma})$&yes\\
$t_1=t_2$&$\pos(s_2:I_{\Gamma})$&no&$\pos(s_2:I_{\Gamma})$&yes\\
$t_1>t_2$&$\pos(t_2:I_{\Gamma})$&yes&$\pos(t_1:I_{\Gamma})-1$&yes\\
\hline
\end{tabular}\\
\vspace{2mm}

From this table we deduce that the brackets in \refE{d2G} vanish.
\end{itemize}
}%\optproof
\end{proof} 

\begin{thm}
$(\GD({\ExtVert}),d)$ is a  commutative differential $\BZ$-graded algebra.
\end{thm}
\begin{proof}
This is a consequence of Proposition \ref{P:GD-CDGA}  and
Lemmas \ref{L:d}--\ref{L:d20}.
%\previousfn{\pl Make sure in all versions that the numbering of these
%three lemmas is indeed consecutive}
\end{proof}
\subsection{The CDGA $\AD(A)$ of admissible diagrams}
\begin{defin}\label{D:admdiag}
A diagram is \emph{admissible} if it contains no loops, no double edges,
no internal vertices of valence $\leq 2$, and if each of its internal vertices is connected to some external vertex.
Otherwise a diagram is \emph{non-admissible}.
We denote by $\NAI({\ExtVert})$ the graded submodule of $\GD({\ExtVert})$ generated by the non-admissible diagrams.
\end{defin}
An admissible diagram does not have dead ends
either, because \p a dead end implies the existence of an internal
vertex of valence $1$, or of a loop, or of a double edge\p. Hence \p an
admissible diagram consists  only of simple chords and simple
contractible edges.\previousoldfn{\pl At some point I was thinking that we need that
  the only admissible diagram with $0$ or $1$ vertex would be the unit
  diagram. But I do not see why anymore (advantage it makes the
  cooperad $\AD$ reduced.) If this is needed then we obtain this by
  saying  that a diagram is not
  admissible also if some of the internal vertices become disconnected
  from the external vertices 
  when one remove a single external vertex. With such a definition,
  then the only admissible diagram with $1$ external vertex is the
  unit diagram. This seems to be still compatible with \refL{ideal-NAI} and
  \refL{PhiNAI}, as well as \refT{formalAD} and \refP{GIKNAI}.}

\begin{lemma}\label{L:ideal-NAI}
The module of non-admissible diagrams $\NAI({\ExtVert})$ is a differential ideal\previousoldfn{\pl this would be wrong if we allow to contract dead ends. Should make a remark about that.\pl Done: rmk \ref{rmk:whydeadend}}
 of $\GD({\ExtVert})$.
\end{lemma}
\begin{proof}
It is easy to check that $\NAI({\ExtVert})$ is an ideal of the algebra $\GD({\ExtVert})$.

We show that $\NAI({\ExtVert})$ is preserved by the differential $d$.
 Let $\Gamma$ be a non-admissible diagram.
\previousoldfn{\pl cancel??? maybe not. \iv What do you mean by ``cancel''?  Delete this proof?  I guess I'd leave it.}
\begin{itemize}
\item If $\Gamma$ contains a loop or a dead end then the same is true
  for each term of $d(\Gamma)$ \p since loops and dead ends are not contracted.
\item If $\Gamma$ contains a double edge then  each term of $d(\Gamma)$ contains a double edge
or a loop (when one of the double edges is contracted).
\item  If $\Gamma$ contains a path component with all vertices internal, then the same is true for each 
term of $d(\Gamma)$.
\item If $\Gamma$ contains an internal vertex $i$ of valence $2$ but no double edges or dead ends,
then for  most of the terms
of $d(\Gamma)$, $i$ is still a bivalent internal vertex, except
 for the two terms obtained by contracting  each of the two edges with endpoint $i$. 
 These two terms cancel each other
 \previousoldfn{\pl\todo check \iv Done.}.
 \item If $\Gamma$ has an internal vertex of valence $1$ then it has a
   dead end, and \p this case is treated in the first bullet above.
 \item If $\Gamma$ has an internal vertex of valence $0$ then it has a connected component
 with all vertices internal, and \p this case is treated in the third bullet above.
 \end{itemize}
 This proves that $d(\NAI({\ExtVert}))\subset\NAI({\ExtVert}))$.
 \end{proof}
\begin{remark}\label{rmk:whydeadend}
The previous lemma would be wrong if in the definition of the differential $d$ we allowed contractions of dead ends. 
This is why dead ends are not defined as  contractible edges in \refD{diag}.
\end{remark}
\begin{defin}\label{D:AD}
 The \emph{$\BZ$-graded CDGA of admissible diagrams} is the quotient $$\AD({\ExtVert}):=\GD({\ExtVert})/\NAI({\ExtVert}).$$
 We write $\AD_N({\ExtVert})=\AD({\ExtVert})$ when we want to emphasize the ambient dimension $N$.
 \end{defin}
 By abuse of notation we will denote by the same symbol a diagram on ${\ExtVert}$, its equivalence class in
 $\GD({\ExtVert})$, and its larger equivalence class in $\AD({\ExtVert})$. The context should always clear up any 
 ambiguity.
 \previousoldfn{In  the previous version there was some  $\AD'(A)$ which was pretented to be a
subCDGA of $\AD(A)$. Nathalie Wahl pointed that it was wrong...}
As a $\BK$-module, $\AD(A)$ is generated by admissible diagrams.
% Denote by $\AD'({\ExtVert})$ the submodule of $\GD({\ExtVert})$ generated by admissible diagrams.
%  Consider 
%  \[
%  \xymatrix{
%  \AD'({\ExtVert})\ar@{^(->}[r]^-\iota&\GD({\ExtVert})\ar[r]^-\pi&\AD({\ExtVert})
%  }
%  \]
%  where $\iota$ is the inclusion and $\pi$ is the projection on the quotient.
% \begin{prop}\label{P:ADAD'}
% $\AD'({\ExtVert})$ is a subCDGA of $\GD({\ExtVert})$ and the composition $\pi\circ\iota$ is an isomorphism
% of CDGAs, so $\AD'({\ExtVert})\cong\AD({\ExtVert})$.
% \end{prop}
% \begin{proof}
% The unit diagram is admissible. It is clear that the product of two admissible diagrams is again admissible.
% Also, if $\Gamma$ is an admissible diagram and if $e$ is a contractible edge then $\Gamma/e$ is
% also admissible.  Therefore $d(\Gamma)$ is in $\AD'({\ExtVert})$. This proves that $\AD'({\ExtVert})$
%  is a subCDGA of $\AD({\ExtVert})$.
 
%  It is immediate that $\pi\circ\iota$ is an isomorphism since $\GD({\ExtVert})=\AD'({\ExtVert})\oplus\NAI({\ExtVert})$
%  \end{proof}
% Because of this proposition, we will from now on identify $\AD'({\ExtVert})$ and $\AD({\ExtVert})$.
A (co)chain complex is said to be \emph{connected} if it is concentrated in non-negative degrees and is isomorphic to $\BK$ in degree $0$.
 \begin{prop}\label{P:connAD} 
 If $N\geq 3$, then $\AD_N({\ExtVert})$ is  a connected CDGA.

 \end{prop}
 \begin{proof}
 Let $\Gamma=({\ExtVert},I,E,s,t)$ be an admissible diagram different from the unit.
 We think of an edge of $\Gamma$ as the union of two half-edges, each with one endpoint
  which is a vertex of $\Gamma$. Since $\Gamma$ is not the unit and since internal vertices are connected
   to some
  external one, there is at least one half-edge whose endpoint is an external vertex.
  Since each internal vertex is of valence $\geq 3$, there are at least $3\cdot|I|$ other half-edges.
  Therefore $ |E|\geq \frac{1}{2} (1+3 |I|) $. We deduce that
 \begin{eqnarray*}
 \deg(\Gamma)
  &=& 
 |E|\cdot(N-1)-|I|\cdot N
\\
 &\geq&
 \frac{1}{2}
 (1+3 |I|) \cdot(N-1)
 -
 |I| \cdot N
\\
&=&
\frac{N-1}{2}+|I|\cdot\frac{N-3}{2}>0.
 \end{eqnarray*}
 \end{proof}

 \begin{remark}
\label{Rmk:connD} It is in fact true that  $\AD_N({\ExtVert})$ is $(N-3)$-connected
  for $N\geq 3$. 
 Indeed if $|I|=0$ then  $\deg(\Gamma)=|E|\cdot(N-1)\geq N-1$ and if $|I|\geq1$ then
 the inequalities in the above proof show that $\deg(\Gamma)\geq\frac{N-1}{2}+1\cdot\frac{N-3}{2}=N-2$.
 When $N\geq 4$ we can refine this argument by treating separately the cases $|I|=1$ and $|I|=2$  and
 deducing that  $\AD_N({\ExtVert})$  is in fact $(N-2)$-connected and
 of finite type\previousoldfn{\iv Where does the finite type from?  Just because
   it has the homology of configuration space? \pl Refinning the
   argument shows that for $N\geq4$ in a given degree the number of
   internal vertices  is bounded, and hence there are finitely
   many admissuible diagrams in that degree}. On the other hand, for $N=3$ it is not true that it is $1$-connected as can be seen from an easy example with $2|E|=1+3|I|$.  Also $\AD_N({\ExtVert})$ cannot be $(N-1)$-connected
 for $|{\ExtVert}|\geq2$ since its homology is the homology of configuration spaces in $\BR^N$, as we will see in
\refT{formalAD}.

However, for $N=2$, $\AD_2(A)$ is not concentrated in non-negative degrees, and is thus not a CDGA 
which is suitable for{\D}  modeling a rational homotopy type, even if its cohomology is non-negatively graded (since it is the cohomology of a configuration space).  
 \end{remark}

\section{Cooperad structure on the spaces of  (admissible)  diagrams}
\label{sec:diagcoop}

In this section we will endow the sequence of CDGAs  $\{\AD({n})\}_{n\geq0}$ with the structure of a cooperad.
We will do this by first  endowing  $\{\GD({n})\}_{n\geq0}$
 with the structure of a cooperad of graded $\BK$-algebras (not differential!).
We fix an ambient dimension $N\geq2$.
 
The plan is as follows. First we construct in section
\ref{sec:PsiConstr}  the cooperad structure maps $\widehat\Psi$ and $\Psi$ on $\GD$ and $\AD$ using
the notion of \emph{condensation} from \refD{loc}. Then we prove in section 
\ref{sec:Phi-alg} that these are morphisms of algebras, and in section 
\ref{sec:Phichain}  we show that  $\Psi$ is a chain map (this is not the case for $\widehat\Psi$\previousoldfn{\pl really?}.)
Finally we prove in section \ref{sec:PsiOper} that this defines  the structure of a cooperad of CDGAs on $\AD$;
this is our main result, \refT{coopGDAD}.

For the  several following subsections, fix a weak ordered partition $\nu\colon A\to P$ and set
 $P^*=\{0\}\ordsum P$, $A_p=\nu^{-1}(p)$, and $A_0=P$,  as in the setting \ref{setting:XP}.

\subsection{Construction of the cooperad structure maps $\widehat\Psi_\nu$ and $\Psi_\nu$}\label{sec:PsiConstr}

In this section we build maps
\begin{eqnarray*}
\widehat\Psi_\nu\colon\GD({\ExtVert})\longrightarrow\GD(P)\otimes\underset{p\in P}{\otimes}\GD({\ExtVert}_p)&\quad\textrm{and}\\
\Psi_\nu\colon\AD({\ExtVert})\longrightarrow\AD(P)\otimes\underset{p\in P}{\otimes}\AD({\ExtVert}_p)
\end{eqnarray*}
which will serve as cooperad  structure maps.
Of course, the tensor product over $p\in P$ is taken in the order fixed on $P$.
Since  ${\ExtVert}_0=P$ we have
\[\GD(P)\otimes\underset{p\in P}{\otimes}\GD({\ExtVert}_p)=\underset{p\in P^*}{\otimes}\GD({\ExtVert}_p).\]

Let us first describe roughly the idea of $\widehat\Psi_\nu$. Let $\Gamma$ be a diagram on $A$ with the set of vertices $V_\Gamma$. Recall from \refD{loc} that a condensation
of $V_\Gamma$ relative{\D} to $\nu$ is a map
\[\lambda\colon V_\Gamma\to P^*\]
such that $\lambda|A=\nu$.
For each condensation $\lambda\in\Cond(V_\Gamma,\nu)$, we first construct 
diagrams $\Gamma(\lambda,0)\in\GD(P)$
 and $\Gamma(\lambda,p)\in\GD(A_p)$ as follows. For $p\in P$, the diagram $\Gamma(\lambda,p)$ on 
 $A_p$  is the full subgraph of $\Gamma$ whose vertices are the $p$-locals
 (that is, those in $\lambda^{-1}(p)$). 
The diagram  $\Gamma(\lambda,0)$  on $P$  is obtained from $\Gamma$ by shrinking each subgraph $\Gamma(\lambda,p)$ into
a single external vertex $p$, for $p\in P$.
Then $\widehat\Psi_\nu(\Gamma)\in\GD(P)\otimes\otimes_{p\in P}\GD(A_p)$ is defined to be
\[\widehat\Psi_\nu(\Gamma)=\sum_{\lambda\in\Cond(V_\Gamma,\nu)}\Gamma(\lambda)
\quad\textrm{ where }\quad\Gamma(\lambda)=\pm\Gamma(\lambda,0)\otimes\underset{p\in P}{\otimes}\Gamma(\lambda,p).
\]
The precise formulas are in equations \refN{eq:defGammal} and
\refN{eq:defPsihat}.

Here is an example illustrating this.

%\eg\label{eg:Psinu}
%Take $A=\{a_1,\dots, a_7\}$, $P=\{1,2,3\}$ and the ordered partition $\nu\colon A\to P$ defined by
%$\nu(a_1)=\nu(a_2)=1$, $\nu(a_3)=\nu(a_4)=\nu(a_5)=2$, and $\nu(a_6)=\nu(a_7)=3$.
%Consider $\Gamma\in\GD(P)$ as in \refF{Gammalambda} with set of internal vertices $I=\{i_1,i_2\}$.
%Consider the two condensation $\lambda$  and $\mu$ relative to $\nu$ and characterized by $\lambda(i_1)=2$ and $\lambda(i_2)=3$,
%and $\mu(i_1)=1$ and $\mu(i_2)=0$. Then $\Gamma(\lambda)$ and $\Gamma(\mu)$ are as in \refF{ }.

\begin{eg}\label{Ex:GammaLambda}
Suppose $A=\{ 1, 2, 3, 4, 5\}$ and $I=\{ 6, 7\}$, so $V=\{ 1, 2, 3, 4, 5, 6, 7\}$.  Let $\Gamma$ be as in Figure \ref{fig:GammaForLambda}.

\begin{figure}[h]
\input{GammaForLambda.pstex_t}
\caption{}
\label{fig:GammaForLambda}
\end{figure}

Let $P^*=\{0, \alpha, \beta\}$ and let $\lambda\colon V\to P^*$ be defined by
$$
\lambda(v)=\begin{cases}
\alpha, & \text{for } v=1, 2, 3;\\
\beta, & \text{for } v=4, 5, 7;\\
0, & \text{for } v=6.
\end{cases}
$$

Then\previousoldfn{\pl the picture of $\Gamma(\lambda,0)$ is wrong: it should have a
  double edges $(\alpha,6)$} $\Gamma(\lambda, 0)$, $\Gamma(\lambda, \alpha)$, and $\Gamma(\lambda, \beta)$ are given in Figure \ref{fig:GammaLambdas}.

\begin{figure}[h]
\input{GammaLambdas.pstex_t}
\caption{}
\label{fig:GammaLambdas}
\end{figure}

\end{eg}

Here is another heuristic description of $\Gamma(\lambda)$%
\ppfn{Here is a last request: it would be
  nice to add a figure here explaining the idea. The figure would be
  almost the same as in the one giving $\Gamma$ in Example 7.1 except
  that the vertices $1,2,3$ on one hand, $7,4,5$ on other hand would
  be pictured very close to each other (maybe inside  dashed circles
  labeled $\alpha$ and $\beta$ respectively) and $6$ would be far from
  these two circles. If we hasd such a picture we could  replace the
  beginning of next sentence in the text by: ``Picture $\Gamma$ as in
  Figure xxx, so that all...''\iv Done.  Let me know if you think I should change something.} .
Picture $\Gamma$ as in Figure \ref{fig:GammaForLambdaClusters}, so that all the $p$-local vertices and their connecting edges are drawn infinitesimally close
to each other, and all the global vertices and the various clusters of $p$-local vertices are drawn far from each other.
Then $\Gamma(\lambda,0)$ is the diagram $\Gamma$ seen from far away, 
and each $\Gamma(\lambda,p)$ for $p\in P$  is that diagram seen through a microscope
centered at the $p$th cluster (forgeting the edges getting out of range.)
This interpretation will correspond, through the Kontsevich configuration space integral,
to what happens to the configurations of points in the Fulton-MacPherson operad.  See the discussion afterÊ \refN{eq:coopI1} as to why the Kontsevich configuration space integral
commutes with the cooperadic structures.\previousoldfn{\pl add a reference about this to some point later in
the paper?}

\begin{figure}[h]
\input{GammaForLambdaClusters.pstex_t}
\caption{}
\label{fig:GammaForLambdaClusters}
\end{figure}

%
%\begin{figure}
%  \centering
%  \fbox{\includegraphics[width=120mm,height=90mm]{GammaLambda.pdf}}
%  \caption{Two terms of $\widehat\Psi_\nu(\Gamma)$ from Example \ref{eg:Psinu}}
%\label{fig:Gammalambda}
%\end{figure}

\begin{defin}
\label{D:locG}
Let $\nu\colon A\to P$ be a weak ordered partition, let
$P^*=\{0\}\ordsum P$,
let $\Gamma$ be a diagram on ${\ExtVert}$ and assume that $I_\Gamma\cap P=\emptyset$.
\begin{itemize}
\item  A \emph{condensation} $\lambda$
on $\Gamma$ is a condensation of $V_\Gamma$ relative to $\nu$
as in \refD{loc}, \p that is, it is a map
$\lambda\colon V_\Gamma\to P^*$ such that $\lambda|A=\nu$.
We consider the set $\Cond(\Gamma,\nu):=\Cond(V_\Gamma,\nu)$ of \p all
condensations on $\Gamma$ relative to $\nu$, and write 
 $\Cond(\Gamma)$ when $\nu$ is understood.
\item
The \emph{extension to the edges} of the condensation $\lambda$  on $\Gamma$ is the map
\[\lambda_E\colon E_\Gamma\longrightarrow  P^*\]
defined by
\[
\lambda_E(e)=\left\{
\begin{array}{ll}
\lambda(s_\Gamma(e)),&\textrm{if }\lambda(s_\Gamma(e))=\lambda(t_\Gamma(e)),\\
0,&\textrm{otherwise.}
\end{array}
\right.
\]
\item Given a condensation $\lambda$ of $\Gamma$, a vertex $v$ (respectively an edge $e$)
is \emph{$p$-local}, for $p\in P$, if $\lambda(v)=p$ (respectively $\lambda_E(e)=p$).
It is \emph{global}  if  $\lambda(v)=0$ (respectively $\lambda_E(e)=0$).
\end{itemize}
\end{defin}

The \p terminology \emph{condensation} is motivated in the case of
diagrams (as it was in \refD{loc} for configurations) by the idea
explained right before \refD{locG} that the diagram should be pictured
with its vertices condensed into clusters depending on the values of $\lambda$.

Clearly the set of condensations on $\Gamma$ is in bijective correspondence with the set of maps from $I_\Gamma$ to $ P^*$,
since the value of a condensation $\lambda$ on an external vertex $a$
is determined by $\lambda(a)=p$ for $a\in {\ExtVert}_p$, that is, $\lambda(a)=\nu(a)$.
An edge is $p$-local if and only if both of its endpoints are.  Otherwise it is global. Also, a global vertex 
is always internal but a global edge can be a chord.

Let $\Gamma$ be a diagram on ${\ExtVert}$ and let
$\lambda\in\Cond(\Gamma)$. \p Without loss of generality we can
assume that $I_\Gamma\cap P=\emptyset$ (see \refR{isodiag}\p).
For $p\in P^*$ we define a diagram
\begin{equation}\label{eq:Glp}
\Gamma(\lambda,p):=({\ExtVert}_p,I_p,E_p,s_p,t_p)
\end{equation}
with
\begin{itemize}
\item $I_p=I_\Gamma\cap\lambda^{-1}(p)$;
\item $E_p=\lambda_E^{-1}(p)$;
\item 
\begin{itemize}
\item For $p\in P$, $s_p$ and $t_p$ are the restrictions of $s_\Gamma$ and $t_\Gamma$
to $E_p$;
\item For $p=0$, $s_0=\widehat{\lambda}\circ s_\Gamma$ and $t_0=\widehat{\lambda}\circ t_\Gamma$
where 
\[\widehat{\lambda}\colon V_\Gamma\longrightarrow P\cup I_0\]
is defined by $\widehat{\lambda}(v)=v$ if $\lambda(v)=0$, and $\widehat{\lambda}(v)=\lambda(v)$ otherwise.
\end{itemize}
\end{itemize}

% In other words, for $p\in P$, the diagram $\Gamma(\lambda,p)$ is the diagram obtained from $\Gamma$ 
% by only keeping its $p$-local vertices and the edges between them. The diagram $\Gamma(\lambda,0)$ 
% should be thought of as the diagram obtained from $\Gamma$ by drawing all $p$-local vertices
% infinitesimaly close to each other so that they become a single external vertex $p$ for each $p\in P$,
% and keeping all the global internal vertices and all  the edges which did not became infinitesimal.
% \previousfn{\pl Maybe add a pictured example \todo \iv  I agree that this would be a good idea.  I think we had a picture in the older version of this paper.}

The set of edges (respectively of internal vertices) of $\Gamma$ is the disjoint union
for $p\in P^*$
of the set of edges (respectively of internal vertices)  
of the diagrams $\Gamma(\lambda,p)$.
Even if $\Gamma$ is admissible,{\D} $\Gamma(\lambda,p)$
may not be. Note also that
$\widehat\lambda$ above is the same as in \refN{eq:lambdahat}.

The equivalence class of $\Gamma(\lambda,p)$ in $\GD({\ExtVert}_p)$
with respect to the relation $\simeq$ of \refD{spacediag}, or even $\otimes_{p\in P^*}\Gamma(\lambda,p)$, is \emph{not}
an invariant of the equivalence class of $\Gamma$ in
$\GD({\ExtVert})$\previousoldfn{\pl Ass example of NB8D p.66?\pl No}. To
correct\previousoldfn{all the signs $\sigma$ should be collected here recalling
  sign $\sigma(I,\lambda)$ already introduced} this
we introduce some signs.\previousoldfn{\pl Do we need a number for this  equation? \pl 080613 No.\pl080620 YES!}
%\[
%\begin{equation}\label{eq:defsigmaGl}
%\sigma(\Gamma,\lambda):=(-1)^{(N-1)\cdot|S(E_\Gamma,\lambda)|-N\cdot|S(I_\Gamma,\lambda)|}
%\end{equation}
%\]
 Define for $I=I_\Gamma$ and $E=E_\Gamma$,
\[S(I,\lambda):=
\{(v,w)\in I\times I:v<w\textrm{ and
}\lambda(v)>\lambda(w)\}
\]
(which was already{\D} introduced in \refN{eq:SIl}), and
\[
%\label{eq:defSEl}
S(E,\lambda):=
\{(e,f)\in E\times E:e<f\textrm{ and }\lambda_E(e)>\lambda_E(f)\}.
\]Define the signs
\begin{eqnarray}
\label{eq:defsigmaIl}
\sigma(I,\lambda)&:=&(-1)^{N\cdot|S(I,\lambda)|}\,,\\
\label{eq:defsigmaEl}
\sigma(E,\lambda)&:=&(-1)^{(N-1)\cdot|S(E,\lambda)|}\,, \textrm{and}\\
\label{eq:defsigmaGl}
\sigma(\Gamma,\lambda)&:=&\sigma(I,\lambda)\cdot
\sigma(E,\lambda)
\end{eqnarray}
(The sign $\sigma(I,\lambda)$ was already defined in \refN{eq:sigmaIl}.)
The proof of the following is straightforward.
\begin{lemma}\label{L:locsign}\previousoldfn{\pl NB8E p.23}
For a diagram $\Gamma$ and a condensation $\lambda$ on $\Gamma$, the element
\[\sigma(\Gamma,\lambda)\cdot
\underset{p\in P*}{\otimes}\Gamma(\lambda,p)
\quad\in\quad\underset{p\in P*}{\otimes}\GD({\ExtVert}_p)
\] depends only
on the equivalence class of $\Gamma$ in $\GD({\ExtVert})$.
\end{lemma}

\optproof[\refL{locsign}]{
\begin{proof}%\previousfn{\pl cancel?}
\begin{itemize}
\item Suppose $N$ even. In that case the order of internal vertices and the orientation of edges are irrelevant
both for a representative of the equivalence class of $\Gamma$ and for the function $\rho(\Gamma,\lambda)$.

Let $\Gamma'$ be a diagram that differs from $\Gamma$ by transposing the order of two consecutive edges
$e_1$ and $e_2$. Then $\Gamma'=-\Gamma$ in $\GD({\ExtVert})$ and $\lambda$ is also a condensation of $\Gamma'$.
\begin{itemize}
\item Suppose that $\lambda_E(e_1)=\lambda_E(e_2)$. Then
 $\rho(\Gamma,\lambda)=\rho(\Gamma',\lambda)$. Also
 $\Gamma(\lambda,\lambda_E(e_1))=-\Gamma'(\lambda,\lambda_E(e_1))$,
 and $\Gamma(\lambda,p)=\Gamma'(\lambda,p)$ for $p\not=\lambda_E(e_1)$.
 \item Suppose that $\lambda_E(e_1)\not=\lambda_E(e_2)$. 
 Then  $\rho(\Gamma,\lambda)=-\rho(\Gamma',\lambda)$ but
  $\Gamma(\lambda,p)=\Gamma'(\lambda,p)$ for each $p\in P^*$.
  \end{itemize}
  In both cases we get 
$\rho(\Gamma,\lambda)\cdot\otimes_{p\in  P^*}\Gamma(\lambda,p)=-
 \rho(\Gamma',\lambda)\cdot\otimes_{p\in  P^*}\Gamma'(\lambda,p)$ as expected.
 \item Suppose $N$ odd. In that case the order of edges is irrelevant.\\

Let $\Gamma'$ be a diagram that differs from $\Gamma$ by an inversion of an edge $e$.
Then $\Gamma'=-\Gamma$. Also $\Gamma'(\lambda,\lambda_E(e))=-\Gamma(\lambda,\lambda_E(e))$ and
$\Gamma'(\lambda,p)=-\Gamma(\lambda,p)$ for $p\not=\lambda_E(e)$. On the other hand 
$\rho(\Gamma,\lambda)=\rho(\Gamma',\lambda)$. Therefore
$\rho(\Gamma,\lambda)\cdot\otimes_{p\in  P^*}\Gamma(\lambda,p)=-
 \rho(\Gamma',\lambda)\cdot\otimes_{p\in  P^*}\Gamma'(\lambda,p)$ as expected.\\

Let $\Gamma'$ be a diagram that differs form $\Gamma$ by the transposition of two consecutive
  internal vertices.
An argument analogous to the one for the case $N$ even yields to the desired equation.
\end{itemize}
\end{proof}
}%\optproof

For a diagram  $\Gamma$ on ${\ExtVert}$  and a condensation $\lambda$ of $\Gamma$
we set 
\begin{equation}
  \label{eq:defGammal}
\Gamma(\lambda):=\sigma(\Gamma,\lambda)\cdot
\underset{p\in P*}{\otimes}\Gamma(\lambda,p)\quad\quad
\in\quad\underset{p\in P*}{\otimes}\GD({\ExtVert}_p),
\end{equation}
where $\sigma(\Gamma,\lambda)=\pm1$ is from \refN{eq:defsigmaGl}
and $\Gamma(\lambda,p)$ is from \refN{eq:Glp}.
By \refL{locsign}  we get a linear map \Z
\[\widehat\Psi_\nu\colon\GD({\ExtVert})\longrightarrow
\underset{p\in P*}{\otimes}\GD({\ExtVert}_p)
\]
defined on generators by
\begin{equation}\label{eq:defPsihat}
\widehat\Psi_\nu(\Gamma):=\sum_{\lambda\in\Cond(\Gamma,\nu)}\Gamma(\lambda).
\end{equation}
                                                    
Recall $\NAI({\ExtVert}_p)\subset\GD({\ExtVert}_p)$, the differential ideal of non-admissible diagrams (\refD{admdiag} and \refL{ideal-NAI}).
Set
\begin{equation}
\label{eq:Nnu}
\NAI(\nu)
:=
\sum_{p\in P^*}
\underset{\substack{q\in P^*\\q<p}}{\otimes}\GD({\ExtVert}_q)\otimes
\NAI({\ExtVert}_p)\otimes
\underset{\substack{q\in P^*\\q>p}}{\otimes}\GD({\ExtVert}_q),
\end{equation}
which is a differential ideal in 
$\otimes_{p\in  P^*}\GD({\ExtVert}_p)$.
Since $\AD({\ExtVert}_p)=\GD({\ExtVert}_p)/\NAI({\ExtVert}_p)$, 
we have an isomorphism of CDGAs ($\BZ$-graded if $N=2$)
\begin{equation}\label{eq:xADisoxD/xNAI}
\left(\underset{p\in P*}{\otimes}\GD({\ExtVert}_p)\right)/\NAI(\nu)\,\cong\,\underset{p\in P*}{\otimes}\AD({\ExtVert}_p)
\end{equation}

\begin{lemma}\label{L:PhiNAI}
$\widehat\Psi_\nu(\NAI({\ExtVert}))\subset\NAI(\nu)$.\previousoldfn{\pl NB:E p.25}
\end{lemma}
\begin{proof}
Let $\Gamma$ be a non-admissible diagram on ${\ExtVert}$ and let $\lambda\in\Cond(\Gamma)$.
\begin{itemize}
\item
If $\Gamma$ has a loop at a vertex $v$, then $\Gamma(\lambda,\lambda(v))$ also has a loop.
 \item If $\Gamma$ has double edges $e_1$ and $e_2$, then so does
 $\Gamma(\lambda,\lambda_E(e_1))$.
\item If $\Gamma$ has an internal vertex $v$ of valence $\leq 2$, then the same
is true for $\Gamma(\lambda,\lambda(v))$ because the valence of $v$ can only decrease.
%\item If $\Gamma$ has a dead end, then it has  an internal vertex of valence $1$, or a loop, or a double edge.
\item If, for some $p\in P$, $\Gamma$ has an internal $p$-local vertex that is not connected to any external vertex, then
the same is true for $\Gamma(\lambda,p)$.
If  $\Gamma$ has a connected component consisting only of internal global vertices, then
the same is true for $\Gamma(\lambda,0)$.
 \end{itemize}
In all cases we see that if $\Gamma$ is not admissible, then the same is true for $\Gamma(\lambda,p)$
for some $p\in P^*$. Therefore $\Gamma(\lambda)\in\NAI(\nu)$ and
$\widehat\Psi_\nu(\NAI({\ExtVert}))\subset\NAI(\nu)$.
 \end{proof}

\begin{prop}\label{P:Phi}
$\widehat\Psi_\nu$ defined in \refN{eq:defPsihat} induces a linear map
\[\Psi_\nu\colon\AD({\ExtVert})\longrightarrow\AD(P)\otimes\underset{p\in P}{\otimes}\AD({\ExtVert}_p).\]
\end{prop}
\previousoldfn{\iv Here I'm confused.  Should this be $P$ or $P^*$? \pl I do not
  undesrtand why it is confusing. It is $p\in P$ here since the extra
  factor corresponmding to $p=0$ appears already as $\AD(P)$.}
\begin{proof}
This is an immediate consequence of the isomorphism \refN{eq:xADisoxD/xNAI}
and  \refL{PhiNAI}.
\end{proof}

Thus, for an  admissible diagram $\Gamma$, $\Psi_\nu(\Gamma)$ is
obtained as the sum \refN{eq:defPsihat} in which non-admissible terms
  are set to zero.
Actually, there are many condensations $\lambda$ for which
$\Gamma(\lambda)$
is not admissible and therefore{\D} does not contribute to
$\Psi_\nu(\Gamma)$.
In particular, only \emph{admissible condensations} (to be defined in
\refD{locreg}) can contribute to the sum, and hence we can use the sum
\refN{eq:PsiReg} below, which has many fewer terms, to define $\Psi_\nu(\Gamma)$.

%Since the weak ordered partition $\nu\colon A\to P$ has been fixed, for the next few sections%we set $\widehat\Psi:=\widehat\Psi_\nu$ and $\Psi:=\Psi_\nu$.
\subsection{$\widehat\Psi_\nu$ and $\Psi_\nu$ are morphisms of algebras}\label{sec:Phi-alg}
The aim of this section is to prove
\begin{prop}\label{P:Phi-alg}
$\widehat\Psi_\nu$ and $\Psi_\nu$ are morphisms of algebras.
\end{prop}

\begin{proof}
%[Proof of \refP{Phi-alg}]
\previousoldfn{\pl NB:8D p.65 and 8E p.3}
We first prove the statement for $\widehat\Psi_\nu$.
Let $\Gamma_1$ and $\Gamma_2$ be two diagrams on ${\ExtVert}$ and suppose that $I_{\Gamma_1}$ and 
$I_{\Gamma_2}$, $E_{\Gamma_1}$ and 
$E_{\Gamma_2}$ respectively, are disjoint.

Define the function
\[\Cond(\Gamma_1)\times\Cond(\Gamma_2)\longrightarrow\Cond(\Gamma_1\cdot\Gamma_2),\quad\quad(\lambda_1,\lambda_2)
\longmapsto\lambda_1\cdot\lambda_2\]
 by $(\lambda_1\cdot\lambda_2)(v)=\lambda_i(v)$ when $v\in V_{\Gamma_i}$ for $i=1,2$.
This map is well-defined because if $v\in V_{\Gamma_1}\cap V_{\Gamma_2}$
then $v$ is external and $\lambda_1(v)=\lambda_2(v)=\nu(v)$. Moreover, it is a bijection 
whose inverse is given by $\lambda\mapsto(\lambda|V_{\Gamma_1},\lambda|V_{\Gamma_2})$.

Since 
\[\Gamma_1(\lambda_1,p)\cdot\Gamma_2(\lambda_2,p)=
(\Gamma_1\cdot \Gamma_2)(\lambda_1\cdot \lambda_2,p)
\]
it is easy to see that
\[
\underset{p\in P*}{\otimes}(\Gamma_1\cdot\Gamma_2)(\lambda_1\cdot\lambda_2,p)
=
\eta(\Gamma_1,\lambda_1,\Gamma_2,\lambda_2)
\cdot
\left(
\underset{p\in P*}{\otimes}\Gamma_1(\lambda_1,p)
\right)
\cdot
\left(
\underset{q\in P*}{\otimes}\Gamma_2(\lambda_2,q)
\right).
\]
where 
\[
\eta(\Gamma_1,\lambda_1,\Gamma_2,\lambda_2):=
(-1)^s
\quad\textrm{with}\quad
s=
{
%\left\lbrace
\sum_{\substack{p,q\in P^*\\q<p}}
\deg(\Gamma_1(\lambda_1,p)) 
\cdot
\deg(\Gamma_2(\lambda_2,q)).
%\right\rbrace
}
\]

We have
\begin{eqnarray}
\widehat\Psi_\nu(\Gamma_1\cdot\Gamma_2)
&=&
\sum_{\lambda\in\Cond(\Gamma_1\cdot\Gamma_2)}\sigma(\Gamma_1\cdot\Gamma_2,\lambda)\cdot
\underset{p\in P*}{\otimes}(\Gamma_1\cdot\Gamma_2)(\lambda,p)\notag\\
&=&\sum_{\lambda_1\in\Cond(\Gamma_1)}\sum_{\lambda_2\in\Cond(\Gamma_2)}
\sigma(\Gamma_1\cdot\Gamma_2,\lambda_1\cdot\lambda_2) \cdot
\underset{p\in P*}{\otimes}(\Gamma_1\cdot\Gamma_2)(\lambda_1\cdot\lambda_2,p)\notag\\
&=&\sum_{\lambda_1\in\Cond(\Gamma_1)}\sum_{\lambda_2\in\Cond(\Gamma_2)}
\Big\{
\sigma(\Gamma_1\cdot\Gamma_2,\lambda_1\cdot\lambda_2) \cdot
\eta(\Gamma_1,\lambda_1,\Gamma_2,\lambda_2) 
\label{eq:Phialg1}
\\
&& \quad\quad\cdot 
\left(
\underset{p\in P*}{\otimes}\Gamma_1(\lambda_1,p)
\right)
\cdot
\left(
\underset{q\in P*}{\otimes}\Gamma_2(\lambda_2,q)
\right)
\Big\}
\notag
\end{eqnarray}
On the other hand
\begin{eqnarray}
\widehat\Psi_\nu(\Gamma_1)\cdot\widehat\Psi_\nu(\Gamma_2)&=&
\sum_{\lambda_1\in\Cond(\Gamma_1)}\sum_{\lambda_2\in\Cond(\Gamma_2)}
\Gamma_1(\lambda_1)\cdot\Gamma_2(\lambda_2)\notag\\
&=&
\sum_{\lambda_1\in\Cond(\Gamma_1)}\sum_{\lambda_2\in\Cond(\Gamma_2)}
\Big\{
\sigma(\Gamma_1,\lambda_1)\cdot\sigma(\Gamma_2,\lambda_2)
\label{eq:Phialg2}
\\ 
&&\quad\quad\cdot
\left(
\underset{p\in P*}{\otimes}\Gamma_1(\lambda_1,p)
\right)
\cdot
\left(
\underset{q\in P*}{\otimes}\Gamma_2(\lambda_2,q)
\right)
\Big\}
\notag
\end{eqnarray}

It remains to check that  the signs of \refN{eq:Phialg1} and  \refN{eq:Phialg2}
agree, which is straightforward.\previousoldfn{this is the content of \refL{rho-prod} below in the endnotes}.

For $\Psi_\nu$, the statement is a consequence of the definition of $\Psi_\nu$ in \refP{Phi} and of the fact that
\refN{eq:xADisoxD/xNAI} is an isomorphism of algebras.
\end{proof}

\optproof[\refP{Phi-alg}]{
\begin{lemma}\label{L:rho-prod}
For $i=1,2$, let $\Gamma_i$ be diagrams on ${\ExtVert}$ and let $\lambda_i\in\Cond(\Gamma_i)$.
\[\rho(\Gamma_1\cdot\Gamma_2,\lambda_1\cdot\lambda_2)=
\rho(\Gamma_1,\lambda_1)\cdot \rho(\Gamma_2,\lambda_2)\cdot
\eta(\Gamma_1,\lambda_1,\Gamma_2,\lambda_2)
\]
\end{lemma}
\begin{proof}\private{\pl NB8E page 2}
\\
\begin{itemize}
\item Suppose $N$ odd.
Set, for $i=1,2$, 
\[R_i:=\{(v,w)\in I_{\Gamma_i}\times I_{\Gamma_i}:v<w
\textrm{ and }\lambda_i(v)>\lambda_i(w)\}
\]
and set
\[R_{12}:=\{(v_1,v_2)\in I_{\Gamma_1}\times I_{\Gamma_2}:\lambda_1(v_1)>\lambda_2(v_2)\}.
\]
It is easy to see that
\[\rho(\Gamma_1\cdot\Gamma_2,\lambda_1\cdot\lambda_2)=
(-1)^{|R_1|+|R_2|+|R_{12}|}=\rho(\Gamma_1,\lambda_1)\cdot\rho(\Gamma_2,\lambda_2)\cdot
(-1)^{|R_{12}|}.
\]
For $i=1,2$, and $p\in P^*$,
\[\deg(\Gamma_i(\lambda_i,p))\equiv\left|\{v\in I_{\Gamma_i}:\lambda_i(v)=p\}\right|\mod2.\]
Therefore
\begin{eqnarray*}
|R_{12}|&=&\sum_{q<p}
\left|\{(v_1,v_2)\in I_{\Gamma_1}\times I_{\Gamma_2}:
         \lambda_1(v_1)=p,\lambda_2(v_2)=q\}\right|
\\&\equiv&
\sum_{q<p}\deg(\Gamma_1(\lambda_1,p))\cdot\deg(\Gamma_2(\lambda_2,q))\mod2,
\end{eqnarray*}
which implies the formula of the lemma.
\item Suppose $N$ even. The argument is completely analogous to the case $N$ odd,
 with the edges playing the role of  the internal vertices.
\end{itemize}
\end{proof}
}%\oldproof

\subsection{$\Psi_\nu$ is a chain map}\label{sec:Phichain}

This section is devoted to the proof of the following
\begin{prop}\label{P:Phi-d}
 $\Psi_\nu$ commutes with the differentials.
 \end{prop}
\chaug{}
%Note  that it is \emph{not}
%true that $\widehat\Psi_\nu$ commutes with the differential 
The analog \Z for $\widehat\Psi_\nu$ is \emph{not} true, as illustrated in the following example.
\previousoldfn{\pl To check. Find an example! \pl done but still problem :-/ \iv problem with what?  mod 2? \pl Aug2010: the problem with the previous example
was that with the diagram
with 3 internal vertices and the 3 edges connecting tem as a triangle, actually the differential were $0$ because symmetries and the relations
$\simeq$ on diagrams (for $N$ odd or even). Since it is important to convince ourelf and the reader that there could not be an easy proof that
$Psi$ is a chain map, I prefer to make the example quite explicit}
\previousoldfn{
The following counterexample was wrong because the relations $\simeq$ imply that the difference 
$\widehat\Psi_\nu(d(\Gamma))-d(\widehat\Psi_\nu((\Gamma))$ was actually $0$.
A counterexample is the diagram $\Gamma$ on $A=\{a\}$ with three internal vertices
and three edges joining these internal vertice to make a triangle. Considering the unique
partition $\nu\colon\{a\}\to P:=\{1\}$, one computes%\previousoldfn{\pl check \iv seems ok}
 that
$\widehat\Psi_\nu(d(\Gamma))\not=d(\widehat\Psi_\nu((\Gamma))$. This computation can even be done over $\BK=\BZ/2$ 
if one does not want to be bothered by the signs.
}%}\oldfn

\begin{eg} 
Here we will show that $\widehat\Psi_\nu$ is not a chain
 \Z\previousoldfn{\pl At the beginning
  of an example I like to have a phrase summarizing  what this example
  is about}
map.
Consider the diagram $\Gamma$ given in \refF{NonCommuting}, 
with the set of external vertices $A=\{a\}$.
The internal vertices $\{1,2,3,4\}$ have their natural order and each edge is oriented from the lower to the higher vertex.
Suppose also that $N$ is odd, and hence the order of edges is irrelevant.
Set $P=\{\alpha\}$ and consider the unique
ordered partition
\[\nu\colon \{a\}\longrightarrow \{\alpha\}.\]

The boundary of $\Gamma$ is given by the diagram pictured on the top right of \refF{NonCommuting} (obtained as an alternating
sum of three such diagrams).

Then $\Psi_\nu(d\Gamma)$ is a sum of eight terms corresponding to the eight condensations
\[\lambda\colon\{1,2,3\}\longrightarrow P^*=\{0,\alpha\}.\]
Using that dead ends are not contractible and  diagrams with loops vanish when $N$ is odd
(because of the relations $\simeq$ of \refD{spacediag}), one computes
that \Z
\[d(\widehat\Psi_{\nu}(d\Gamma))\in \GD(\{\alpha\})\otimes\GD(\{a\})\]
consists of a single term, corresponding to the condensation
\Z\previousoldfn{\pl 17/01/2011: I believe that the example as phrased in your
  version  was not correct, more precisely the condensation $\lambda$
  should be as I changed it just below. Could
    you check it again please. \\ Indeed it seems to me that the non trivial
    term of
$\widehat\Psi_\nu(d\Gamma)$ corresponds to the condensation with $1$
$\alpha$-local and $2$ and $3$ global, so that the corresponding term
\[(d\Gamma)(\lambda)=(d\Gamma)(\lambda,0)\otimes
(d\Gamma)(\lambda,\alpha)\]
is such that $(d\Gamma)(\lambda,0)$ has two internal vertices $2$ and
$3$, an external vertex $\alpha$, an edge $(\alpha,2)$ and double
edges $(2,3)$; and $(d\Gamma)(\lambda,\alpha)$ has an external vertex
$a$, an internal vertex $1$ and no edge. Taking the differential of
$(d\Gamma)(\lambda)$
consists only in contracting the only non dead-end  $(\alpha,2)$ to
get indeed the third picture of Figure 13. I guess you're right.  I think I just never changed the tex file from an older example we had where the labeling was incorrect.}
 $\lambda(1)=\alpha\,,\,\lambda(2)=\lambda(3)=0$,
and   represented by\pfn{the comment (72) of the referee is correct:
  it would help to have a $3$ on the internal vertex of left factor of
  $d(\widehat\Psi_\nu)$ and $1$ ont the internal isolated vertex of the
  right factor in the bottom of f \refF{NonCommuting}. Could you fix this? \iv Done.} the bottom\previousoldfn{\pl 17/01/11: the third picture of Figure 13 is not correct. The
  second factor of the tensor product of $d(\widehat\Psi_\nu(d\Gamma))$
  should have its external vertex labeled as $a$ and not as $\alpha$.}
picture  of \refF{NonCommuting}.

\begin{figure}[h]
\input{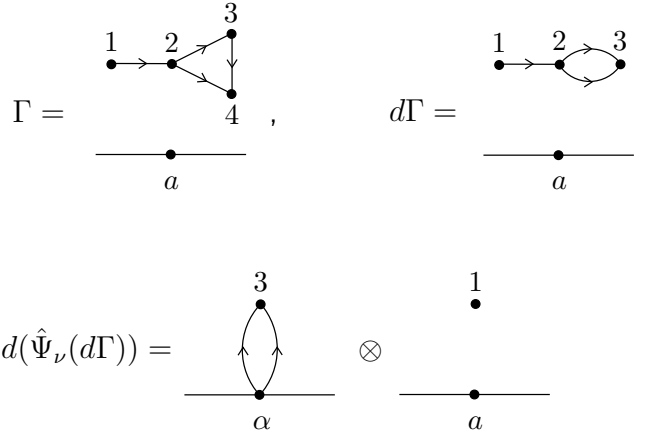}
\caption{A diagram $\Gamma$ for which $d(\widehat\Psi_\nu d(\Gamma))\not=0$.}
\label{fig:NonCommuting}
\end{figure}

\previousoldfn{\todo\\
Here there should be a figure with three pictures, side by side:
\begin{itemize}
\item
the left picture is a diagram with one external vertex $a$, four internal vertices
$1,2,3,4$ and edges $12,13,23,34$ (hence it is a triangle with a hair, and the external vertex is isolated).
This diagram is named $\Gamma$
\item
the middle picture is a diagram with  one external vertex $a$, three internal vertices
$1,2,3$ and three edges:  double edges $12,12$ and a single edge $23$.
This diagram is named $d\Gamma$
\item the right picture is actually a tensor picture of two pictures: \emph{pict1} $\otimes$ \emph{pict2}
where 
\begin{itemize}
\item the left factor \emph{pict1} is a diagram with one external vertex  $I$, one internal vertex with no name
and  double edges joiniung the two vertices;
\item the  right factor \emph{pict2} is a diagram with one external vertex  $I$, one internal vertex with no name
and no edges.
\end{itemize}
this tensor product is named $d\widehat\Psi_\nu(d\Gamma)$
\end{itemize}}

Thus
\[d(\widehat\Psi_\nu d(\Gamma))\not=0.\]
On the other hand
\[d(d\widehat\Psi_\nu(\Gamma))=0\]
since $d^2=0$. Therefore  $\widehat\Psi_\nu$ is not a chain map.
\end{eg}

In fact, the equality $\widehat\Psi_\nu d=d\widehat\Psi_\nu$
is not really expected  to hold. Indeed, if $\Gamma$ is a diagram with $l$ internal vertices then
there are $l^{|P|+1}$ terms in $\Psi_\nu(\Gamma)$, corresponding to \D
the various condensations. 
On the other hand,  for each contractible edge $e$ in $\Gamma$, $\Psi_\nu(\Gamma/e)$ has only
 $(l-1)^{|P|+1}$ terms. Thus there is no clear correspondence between the terms of the sums
$\widehat\Psi_\nu d(\Gamma)$ and of $d\widehat\Psi_\nu(\Gamma)$, and hence no evidence that these two
sums should be equal.

This explains why the proof below that $\Psi_\nu$ is a chain map is quite elaborate.
The idea is to restrict to  condensations of $\Gamma$ for which $\Gamma(\lambda)$
is admissible and to establish \refL{omegabij}, which amounts to
exhibiting 
  a $1$--$1$ correspondence
between condensations of $\Gamma$ and of $\Gamma/e$.

Let $\Gamma$ be a diagram
 \Z\previousoldfn{\pl 17/01/11 \refD{lcontr}  below has
  been moved before statement of \refL{lcontr} and the defintion of
  the function $\lambda/e$ has been moved after the proof of
  \refL{lcontr}.} 
on ${\ExtVert}$.

\previousoldfn{\pl the former version of \refL{lcontr} ``lemma:lcontr'' characterizing $\lambda$-contractible edges was wrong.}

 \begin{defin}\label{D:locreg}\previousoldfn{\pl NB8E:p. 4}
A condensation $\lambda$ of $\Gamma$ is \emph{admissible}
 if for each internal vertex $i$ and each $p\in P$ there is an equivalence
 \[\lambda(i)=p\quad\Longleftrightarrow 
 \quad i \textrm{ admits at least two distinct adjacent $p$-local  vertices}.\]
 Denote by  $\AdmCond(\Gamma)$ the set of admissible condensations on $\Gamma$.
 \end{defin}
 This terminology is motivated by the following:
 \begin{lemma}\label{L:notregloc}
 If $\lambda$ is not admissible then $\Gamma(\lambda)\in\NAI(\nu)$.
 \end{lemma}
 \begin{proof}
 \previousoldfn{\pl NB:8D p.62}
 Suppose that $i$ is a $p$-local internal vertex for some $p\in P$,  and
 suppose that it does not have two adjacent $p$-local vertices.
 Then $i$ is internal of valence $<2$ in $\Gamma(\lambda,p)$, and hence
 $\Gamma(\lambda,p)\in\NAI({\ExtVert}_p)$.
 
Suppose that  $i$ is an internal vertex that is not $p$-local but that 
has   two adjacent $p$-local vertices for some $p\in P$.
 Then in $\Gamma(\lambda,0)$,  the external vertex $p$ is connected
 by a double edge to either $i$
 (if $\lambda(i)=0)$ or to the external vertex $q$ (if $\lambda(i)=q\in P\setminus\{p\}$).  Thus
 $\Gamma(\lambda,0)\in\NAI(P)$.
 \end{proof}
 Lemma \ref{L:notregloc}  implies that, in $\otimes_{p\in P^*}\AD({\ExtVert}_p)$ and for $\Gamma$ admissible,
 \begin{equation}\label{eq:PsiReg}
\Psi_\nu(\Gamma)=\sum_{\lambda\in\AdmCond(\Gamma, \nu)}\Gamma(\lambda).
\end{equation}

 \begin{lemma}\label{L:unicityregular}
 Let $\lambda_1,\lambda_2$ be two admissible condensations on $\Gamma$. If $\lambda_1$ and $\lambda_2$ coincide on all vertices except possibly on one,
 then $\lambda_1=\lambda_2$.
 \end{lemma}
 \begin{proof}
Let $u$ be a vertex of $\Gamma$ such that $\lambda_1(v)=\lambda_2(v)$ for $v\not=u$.
 If $u$ is external then the values of $\lambda_i(u)$ are determined by $\nu$,
and hence $\lambda_1=\lambda_2$.
 Suppose that $u$ is internal.
 If $u$ has two adjacent vertices that are $p$-local (for both $\lambda_1$ and $\lambda_2$) for some $p\in P$,
  then $\lambda_i(u)=p$ by admissibility. Otherwise $\lambda_i(u)=0$, again by admissibility.
 \end{proof}

For a condensation\Z\previousoldfn{\pl 17/01/11 this has been moved and a bit changed}
$\lambda$ of $\Gamma$,
recall the extension to vertices $\lambda_E$ from \refD{locG} and $\Gamma(\lambda, p)$ from Equation \eqref{eq:Glp}.
\begin{defin}\label{D:lcontr}
An edge $e$ of $\Gamma$ is \emph{$\lambda$-contractible}
 if it is contractible in $\Gamma(\lambda,\lambda_E(e))$.
\end{defin}

\chaug
\begin{lemma}\label{L:lcontr}
Assume that $\Gamma$ is admissible and let $\lambda$ be an admissible condensation.\\
 An edge $e$ of $\Gamma$ is \emph{$\lambda$-contractible} if and only if
 the following conditions hold:
 \begin{enumerate}
 \item  $e$ is contractible in $\Gamma$, and
 \item $\lambda(s_\Gamma(e))=\lambda(t_\Gamma(e))$
 or $\min(\lambda(s_\Gamma(e)),\lambda(t_\Gamma(e)))=0$.
 \end{enumerate}
 \end{lemma}
Condition (2) of the lemma is equivalent to having either both endpoints of $e$
$p$-local for some $p\in P$ or some endpoint being global. 
\begin{proof}
\chaug
First we show that (1) and (2) are necessary conditions.
If (1) does not hold then $e$ is a chord in $\Gamma$ (because $\Gamma$ is admissible).
Therefore $e$ is also a chord in $\Gamma(\lambda,\lambda_E(e))$, and hence it is not
$\lambda$-contractible.
If (2) does not hold then $e$ is a chord in $\Gamma(\lambda,0)$
 joining the external vertices $\lambda(s_\Gamma(e))$ and  $\lambda(t_\Gamma(e))$.
Therefore $e$ is not $\lambda$-contractible.

Suppose now that (1) and (2) hold. We will prove that $e$
is a contractible edge in $\Gamma(\lambda,\lambda_E(e))$.
Let $v$ and $w$ be the endpoints of $e$. As $\Gamma$ is admissible, $e$ is not a loop and $v\not=w$.
 We distinguish three cases.

\begin{itemize}
\item{Case 1: $\lambda(v)=\lambda(w)=p\in P$.}\\ 
Since $e$ is contractible, $v$ or $w$ is internal in $\Gamma$, and hence the same is true in
 $\Gamma(\lambda,p)$. Clearly $e$ is neither a chord, nor a loop there.
Since $\lambda$ is admissible, if $v$ is internal, then it has an adjacent $p$-local vertex, distinct from  $w$.
Therefore $v$ also has another adjacent vertex in $\Gamma(\lambda,p)$.
The same is true for $w$ if it is internal.  Therefore $e$ is not a dead end  in $\Gamma(\lambda,p)$.
Thus $e$ is $\lambda$-contractible.
\item{Case 2: $\lambda(v)=\lambda(w)=0$.}\\
Then $e$ is an edge of $\Gamma(\lambda,0)$ joining two distinct internal vertices $v$ and $w$.
In particular $e$ is neither a chord nor a loop.
Since $e$ is not a dead end in $\Gamma$, there exists a vertex $v'\not=w$ that is adjacent to $v$ in
$\Gamma$. If $v'$ is global then it is also a vertex of $\Gamma(\lambda,0)$ and if $v'$ is $p$-local, for some $p\in P$,
then it becomes an external vertex $p$ in  $\Gamma(\lambda,0)$. In both cases in $\Gamma(\lambda,0)$, $v$ has an adjacent vertex
distinct from $w$. Similarly $w$ has an adjacent vertex in $\Gamma(\lambda,0)$ distinct from $v$.
Therefore $e$ is not a dead end in $\Gamma(\lambda,0)$  and $e$ is $\lambda$-contractible.
\item{Case 3: $\lambda(v)=0$ and $\lambda(w)=p\in P$ (or the other way).}\\
Then $e$ is an edge of $\Gamma(\lambda,0)$ joining the internal vertex $v$ to the external vertex $p$.
Since $e$ is not a dead end in $\Gamma$, 
there is a vertex $v'$ adjacent to $v$ in $\Gamma$ and distinct from $w$.
Since $\lambda$ is admissible, we have that $\lambda(v')\not=p$, and hence $v'$ is either global,
or $q$-local for some $q\not=p$. Then, in $\Gamma(\lambda,0)$,
 $v'$  either is an internal vertex or becomes the external vertex $q$. In both cases it is a vertex distinct from $p$ and adjacent to $v$. 
This proves that $e$ is not a dead end. 
Thus $e$ is $\lambda$-contractible.
\end{itemize} 
\end{proof}

Let $e$ be a
 \Z\previousoldfn{\pl 17/01/11 this paragraph has been moved} 
  contractible edge in $\Gamma$. Let $v$ and $w$ be the endpoints of $e$
 with $v<w$. Thus $V_{\Gamma/e}=V_\Gamma\setminus\{w\}$. Define the function
 \begin{equation}
\label{eq:lovere}
\lambda/e\colon V_{\Gamma/e}\longrightarrow P^*
\end{equation}
by
 \[(\lambda/e)(z)=\left\{
 \begin{array}{ll}
 \lambda(z),&\textrm{if $z\not=v$ or $z=v$ is external;}\\
 \max(\lambda(v),\lambda(w)),&\textrm{if $z=v$ is internal.}
 \end{array}
 \right.
 \]
 It is clear that $\lambda/e$ is a condensation of $\Gamma/e$. Notice also that if $e$ is
 $\lambda$-contractible then $(\lambda/e)(v)=\max(\lambda(v),\lambda(w))$.

 Assume that \Z $\Gamma$ is an admissible diagram.
  Consider the sets
 
 \begin{align*}
 \Omega & =
 \Big\{(e,\lambda):e\in E_\Gamma,\,\lambda\in\AdmCond(\Gamma),\,e\textrm{ $\lambda$-contractible},\,
 (\Gamma/e)(\lambda/e)\not=0\textrm{ in }\underset{p\in P^*}{\otimes}\AD({\ExtVert}_p)\Big\},\\
  \overline{\Omega} & =
 \Big\{(e,\bar\lambda):e\in E_\Gamma,\,e\textrm{ contractible},\,
 \bar\lambda\in\Cond(\Gamma/e),\,
 (\Gamma/e)(\bar\lambda)\not=0\textrm{ in }\otimes_{p\in P^*}\AD({\ExtVert}_p)\Big\},
 \end{align*}
 
  and the map
\begin{align*} 
\omega\colon\Omega & \longrightarrow\overline{\Omega}\\
    (e,\lambda)  & \longmapsto(e,\lambda/e).
\end{align*}
  \begin{lemma}\label{L:omegabij}
  $\omega$ is a bijection.
  \end{lemma}
  \begin{proof}\previousoldfn{\pl NB8E p.15}
  We first show that $\omega$ is injective.
  Let $e$ be a contractible edge and, for $i=1,2$,
   let $\lambda_i$  be  admissible condensations
  of $\Gamma$ such that $e$ is $\lambda_i$-contractible
   and $(\Gamma/e)(\lambda_i/e)\not=0$.
  Assume that  $\lambda_1/e=\lambda_2/e$. We will show that $\lambda_1=\lambda_2$.
  
  Set $\bar\lambda=\lambda_1/e=\lambda_2/e$.  This is an admissible condensation because $(\Gamma/e)(\bar\lambda)\not=0$ in $\otimes_{p\in P^*}\AD({\ExtVert}_p)$ and because 
of \refL{notregloc}.
  Let $v$ and $w$ be the endpoints of $e$ with $v<w$.
  Thus $V_{\Gamma/e}=V_\Gamma\setminus\{w\}$. We know that $\lambda_i$ 
  agrees with $\bar\lambda$ on $V_\Gamma\setminus\{v,w\}$, and therefore
  we only need to show that $\lambda_1(v)=\lambda_2(v)$ and $\lambda_1(w)=\lambda_2(w)$.
  Moreover, since each $\lambda_i$ is admissible, by \refL{unicityregular}
  it is enough to prove only one of these two equations.
  
If $v$ is external, then $\lambda_1(v)=\lambda_2(v)$ is determined and hence $\lambda_1=\lambda_2$.
Suppose that $v$ is internal.
 If $\bar\lambda(v)=0$ then, since  $\bar\lambda(v)=\max(\lambda_i(v),\lambda_i(w))$,
 we get $\lambda_i(v)=\lambda_i(w)=0$ for $i=1,2$, and hence $\lambda_1=\lambda_2$. Suppose that $\bar\lambda(v)=p\in P$.
By admissibility of $\bar\lambda$, there exist two vertices $x$ and $y$ other than $v$ and $w$
that are adjacent to $v$ in $\Gamma/e$ and with $\bar\lambda(x)=\bar\lambda(y)=p$.
This implies that $\lambda_i(x)= \lambda_i(y)=p$ for $i=1,2$.
Then in $\Gamma$, either $x$ and $y$ are both adjacent to $v$ (respectively to $w$), or
$x$ is adjacent to $v$ and $y$ is adjacent to $w$ (or the other way around).
In the first case
we get by admissibility of $\lambda_i$  that $\lambda_1(v)=\lambda_2(v)=p$ (respectively $\lambda_1(w)=\lambda_2(w)=p$), 
and hence $\lambda_1=\lambda_2$ by \refL{unicityregular}. In the second case, since 
$p=\bar\lambda(v)=\max(\lambda_1(v),\lambda_1(w))$,
we get that $\lambda_1(v)=p$ or $\lambda_1(w)=p$. Let us say that $\lambda_1(v)=p$, the other case
being analogous. Then
$w$ is adjacent to $v$ and to either $x$ or $y$, and thus $w$ is adjacent to two $p$-local vertices (for
the condensation $\lambda_1$),
and hence we also have $\lambda_1(w)=p$ by admissibility. The same argument shows that $\lambda_2(v)=\lambda_2(w)=p$.
This proves injectivity of $\omega$.
\vspace{3mm}

To show $\omega$ is surjective, let $e$ be a contractible edge of $\Gamma$ and let $\bar\lambda$ be a condensation of $\Gamma/e$
such that $(\Gamma/e)(\bar\lambda)\not=0$.
We will construct an admissible condensation $\lambda$ of $\Gamma$ such that $e$ is $\lambda$-contractible
and $\lambda/e=\bar\lambda$. 
Let $v$ and $w$ again be the endpoints of $e$ with $v<w$.
For $z\in V_\Gamma\setminus\{v,w\}$, set $\lambda(z)=\bar\lambda(z)$.
We need to define $\lambda(v)$ and $\lambda(w)$ and to check that $\lambda$ has the desired properties.
We consider the following cases.
\begin{enumerate}
\item Suppose $v$ is external. Then $\lambda(v)=\bar\lambda(v)=p\in P$
  is prescribed by $\nu$.
\begin{enumerate}
\item Suppose that there exists a vertex $x$ in $\Gamma$ different from $v$ and adjacent to $w$ such that
$\bar\lambda(x)=p$. In that case, set $\lambda(v)=\lambda(w)=p$.
Then $\lambda$ is admissible at the vertex $w$ because it has two $p$-local adjacent vertices $v$ and $x$.
It is easy to check that $\lambda$ is also admissible at the other internal vertices of $\Gamma$, using the fact
that $\bar\lambda$ is. Moreover $e$ is $\lambda$-contractible by \refL{lcontr} since it is contractible and $\lambda(v)=\lambda(w)$.
\item Suppose that $w$ is not adjacent to any $p$-local vertex other than $v$. 
 In that case, set $\lambda(v)=p$ and $\lambda(w)=0$.
The vertex $w$ is not adjacent in $\Gamma$ to two vertices $x$ and $y$ such that
$\bar\lambda(x)=\bar\lambda(y)=q\in P$ with $q\not=p$ because
 otherwise $(\Gamma/e)(\bar\lambda,0)$ would contain a double edge joining $p$ and $q$, and hence
$\Gamma/e\in\NAI(\nu)$, contrary to our hypothesis. This proves that $\lambda$ is admissible at $w$ and the admissibility at other vertices
 is a consequence of the admissibility of $\bar\lambda$.
Also $e$ is $\lambda$-contractible.
 \end{enumerate}
 \item Suppose that $v$ is internal. Then $w$ is also internal since $v<w$.
 \begin{enumerate}
 \item
 Suppose that $\bar\lambda(v)=0$. In that case set $\lambda(v)=\lambda(w)=0$.
 By admissibility of $\bar\lambda$, there do not exist
 two vertices $x,y\in V_\Gamma\setminus\{v,w\}$ adjacent in $\Gamma$ to either $v$ or $w$
with $\bar\lambda(x)=\bar\lambda(y)\in P$. 
 It is easy to see that $\lambda$ is admissible and $e$ is $\lambda$-contractible. 
 \item Suppose that $\bar\lambda(v)=p\in P$.
 By admissibility of $\bar\lambda$ there exist two distinct  vertices 
 $x,y\in V_\Gamma\setminus\{v,w\}$ adjacent in $\Gamma$ to either $v$ or $w$
 such that $\bar\lambda(x)=\bar\lambda(y)=p$.
 \begin{itemize}
 \item  If $v$ is not adjacent to any vertices in $(V_\Gamma\setminus\{v,w\})\cap\bar\lambda^{-1}(p)$
 then set $\lambda(v)=0$ and $\lambda(w)=p$.
  \item  If $w$ is not adjacent to any vertices in $(V_\Gamma\setminus\{v,w\})\cap\bar\lambda^{-1}(p)$
 then set $\lambda(v)=p$ and $\lambda(w)=0$.
 \item If both $v$ and $w$ are adjacent to some vertices in $(V_\Gamma\setminus\{v,w\})\cap\bar\lambda^{-1}(p)$
 then set $\lambda(v)=\lambda(w)=p$.
 \end{itemize}
 In each case it is easy to see that $\lambda$ is admissible and that $e$ is $\lambda$-contractible.
 \end{enumerate}
 \end{enumerate}
 This proves surjectivity of $\omega$.
   \end{proof}
   
    \begin{lemma}\label{L:G/e(l/e)}
 If $\Gamma$ is admissible, if $\lambda$ is admissible, and if $e$ is a $\lambda$-contractible edge 
 of $\Gamma$, then,
 for $p\in P^*$, we have in $\AD({\ExtVert}_p)$:
 \[
 (\Gamma/e)(\lambda/e,p)=
 \left\{
 \begin{array}{ll}
  \Gamma(\lambda,p)/e,&\textrm{if $p=\lambda_E(e)$;}\\
 \Gamma(\lambda,p),&\textrm{otherwise.}\\
 \end{array}
 \right.
 \]
 \end{lemma}
 \begin{proof}
 Let $v$ and $w$ be the endpoints of $e$ with $v<w$. Then $V_{\Gamma/e}=V_\Gamma\setminus\{w\}$
 and $E_{\Gamma/e}=E_\Gamma\setminus\{e\}$.
 It is easy to see that the equations to prove are equivalent to
 \[
 \left\{
 \begin{array}{rcl}
 \lambda/e&=&\lambda|(V_{\Gamma}\setminus\{w\})\\
 (\lambda/e)_E&=&\lambda_E|(E_{\Gamma}\setminus\{e\})
 \end{array}
 \right.
 \]

 Since $e$ is $\lambda$-contractible, by \refL{lcontr}, $\lambda(v)=\lambda(w)$ or
 $\min(\lambda(v),\lambda(w))=0$. If $\lambda(v)<\lambda(w)$  then $\lambda(v)=0$ which implies
 that $v$ is internal, and the same for $w$ because $v<w$, in which case we can transpose the order
 of $v$ and $w$ to get an equivalent diagram (up to sign) in which the roles of $v$ and $w$ are exchanged.
 Therefore, without loss of generality we
 can always assume that $\lambda(v)\geq\lambda(w)$. 
 This implies that $(\lambda/e)(v)=\lambda(v)$. Also for $z\not=v,w$ we have 
 $(\lambda/e)(z)=\lambda(z)$. Thus $\lambda/e=\lambda|V_{\Gamma}\setminus\{w\}$.
 
 It remains to prove that  $(\lambda/e)_E=\lambda_E|(E_{\Gamma}\setminus\{e\})$.
 Let $f\not=e$ be an edge of $\Gamma$. If $w$ is not an endpoint of $f$, then,
  since $\lambda/e=\lambda|V_{\Gamma}\setminus\{w\}$, 
   $(\lambda/e)_E(f)=\lambda_E(f)$. Suppose that $w$ is an endpoint of $f$.
 If $\lambda(w)=\lambda(v)$ then $(\lambda/e)_E(f)=\lambda_E(f)$.
Otherwise $\lambda(w)=0$ and $\lambda(v)=p\in P$, and hence $f$ is global in $\Gamma$.
As $\Gamma$ is admissible and $f\not=e$, the other endpoint of $f$ is not $v$.
Since $\lambda$ is admissible and since $w$ is not $r$-local but is adjacent to the $p$-local vertex $v$,
we get that the other endpoint of $f$ is not $p$-local. This implies that $f$ is global in $\Gamma/e$,
and hence $(\lambda/e)_E(f)= \lambda_E(f)=0$. This proves that 
 $(\lambda/e)_E=\lambda_E|(E_{\Gamma}\setminus\{e\})$.
 \end{proof}

%before

 \begin{proof}[Proof of \refP{Phi-d}]
 \previousoldfn{\pl The signs need to worked out properly. \pl 2july Done but to be checked}
 Let $\Gamma$ be an admissible diagram on $A$.
 For $p\in P^*$ and for a condensation $\lambda$ of $\Gamma$,
 define the sign
 \[\eta(\Gamma,\lambda,p):= (-1)^s
\quad\textrm{with}\quad
s={ \sum_{\substack{q\in P^*\\q<p}}\deg(\Gamma(\lambda,q))}.
 \]
We have
 \begin{align}
 d(\Psi_\nu(\Gamma))&\stackrel{\textrm{\refN{eq:PsiReg}}}{=}
 d\left(\sum_{\lambda\in\AdmCond(\Gamma)}\Gamma(\lambda)\right)\notag\\
 &=
 \sum_{\lambda\in\AdmCond(\Gamma)} 
 \sum_{p\in P^*}
\sigma(\Gamma,\lambda)\cdot
\eta(\Gamma,\lambda,p)\cdot
\underset{q<p}{\otimes} \Gamma(\lambda,q)\otimes 
 d\left(\Gamma(\lambda,p)\right)
 \otimes\underset{q>p}{\otimes}\Gamma(\lambda,q)\notag\\
  &=
  \sum_{\lambda\in\AdmCond(\Gamma)} 
 \sum_{p\in P^*}
 \sum_{e\in\Econtr_{\Gamma(\lambda,p)}}
\sigma(\Gamma,\lambda)\cdot
\eta(\Gamma,\lambda,p)\cdot
 \epsilon(\Gamma(\lambda,p),e)\cdot\notag\\
&\quad\quad\quad\quad\quad\quad\quad\quad\quad\cdot
\underset{q<p}{\otimes} \Gamma(\lambda,q)\otimes 
\Gamma(\lambda,p)/e
 \otimes\underset{q>p}{\otimes}\Gamma(\lambda,q)\notag\\
   &\stackrel{\textrm{\refL{G/e(l/e)}}}{=}
  \sum_{(e,\lambda)\in\Omega}
   \sigma(\Gamma,\lambda)\cdot
\eta(\Gamma,\lambda,\lambda_E(e))\cdot
 \epsilon(\Gamma(\lambda,\lambda_E(e)),e)\cdot
\left(
   \underset{p\in P^*}{\otimes}(\Gamma/e)(\lambda/e,p)
\right).
\label{eq:Phi-d1}
\end{align}
On \Z the other hand,
 \begin{align}
\Psi_\nu(d(\Gamma))
&=\Psi_\nu\left(
\sum_{e\in\Econtr_\Gamma}\epsilon(\Gamma,e)\cdot
\Gamma/e\right)
\notag\\
&=
 \sum_{(e,\bar\lambda)\in\overline{\Omega}}
\epsilon(\Gamma,e)
\cdot
\sigma(\Gamma/e,\bar\lambda)
\cdot
\left(
\underset{p\in P^*}{\otimes}
(\Gamma/e)(\bar\lambda,p)
\right)
\notag\\
&\stackrel{\textrm{\refL{omegabij}}}{=}
 \sum_{(e,\lambda)\in\Omega}
\epsilon(\Gamma,e)
\cdot
\sigma(\Gamma/e,\lambda/e)
\cdot
\left(
\underset{p\in P^*}{\otimes}
(\Gamma/e)(\lambda/e,p)
\right)
\label{eq:Phi-d2}
\end{align}
It remains to check that the signs of \refN{eq:Phi-d1} and \refN{eq:Phi-d2} agree,
which is straighforward.
%\previousfn{\pl\todo NB8I p.38seq}
\end{proof}

\subsection{Proof that the cooperad structure is well-defined}
\label{sec:PsiOper}

We show in this section that $\{\widehat\Psi_\nu\}$ and $\{\Psi_\nu\}$
endow $\GD$ and $\AD$ with the cooperadic structure (the former in the category of vector spaces and the latter in the category of chain complexes), 
when $\nu$ runs over all weak ordered partitions.

First let us show the associativity of the structure maps. Suppose given a
 weak ordered partition $\nu\colon A\to P$ as before.
Suppose moreover that ${\ExtVert}$ is itself linearly ordered, that $\nu$ is increasing,
and that $P^*\cap {\ExtVert}=\emptyset$.
Let $\xi\colon B\to {\ExtVert}$ be an ordered weak partition of a finite set $B$.
Set $B_a=\xi^{-1}(a)$ for $a\in {\ExtVert}$.
Also set $A^*=\{0\}\ordsum A$.

We then have a natural bijection
\[\underset{a\in {\ExtVert}}{\amalg}B_a\cong\underset{p\in P}{\amalg}\underset{a\in {\ExtVert}_p}{\amalg}B_a.
\]
For $p\in P$, the partition $\xi$ restricts to a weak ordered partition
\[\xi_p\colon\underset{a\in {\ExtVert}_p}{\amalg}B_a\longrightarrow {\ExtVert}_p.
\]
The associativity of $\widehat\Psi$ amounts to the following lemma whose proof is straightforward.
\begin{lemma}\label{L:assocphi}
The following diagram is commutative:
\begin{equation}\label{eq:assocphi}
\xymatrix{
\GD\left(\underset{p\in P}{\amalg}\left(\underset{a\in {\ExtVert}_p}{\amalg}B_a\right)\right)
\ar@{=}[r]
\ar[d]^-{\widehat\Psi_{\nu\circ\xi}}
&
\GD\left(\underset{a\in {\ExtVert}}{\amalg}B_a\right)
\ar[d]^-{\widehat\Psi_{\xi}}
\\
\GD(P)\otimes\underset{p\in P}{\otimes}\GD\left(\underset{a\in {\ExtVert}_p}{\amalg}B_a\right)
\ar[d]^-{\id\otimes\underset{p\in P}{\otimes}\widehat\Psi_{\xi_p}}
&
\GD({\ExtVert})\otimes\underset{a\in {\ExtVert}}{\otimes}\GD({B_a})
\ar[d]^-{\widehat\Psi_{\nu}\otimes\id}
\\
\GD(P)\otimes\underset{p\in P}{\otimes}
\left(\GD({\ExtVert}_p)\otimes\underset{a\in {\ExtVert}_p}{\otimes}\GD({B_a})
\right)
\ar[r]^-{\cong}_-{\tau}
&
\left(\GD(P)\otimes
\underset{p\in P}{\otimes}\GD({\ExtVert}_p)\right)
\otimes
\underset{a\in {\ExtVert}}{\otimes}\GD({B_a})
}
\end{equation}
The horizontal bottom isomorphism $\tau$ is the obvious reordering 
of factors (with the usual Koszul sign).
\end{lemma}

\optproof[\refL{assocphi}]{
Let $\Gamma$ be a diagram in $\GD(\amalg_{p\in P}\amalg_{a\in {\ExtVert}_p}B_p)$.
The image of $\Gamma$ by the left vertical arrows in \refN{eq:assocphi} is
\begin{equation}\label{eq:assphiL}
\begin{split}
\sum_{\lambda\in\Cond(\Gamma,P)}&
\sum_{(\mu_p\in\Cond(\Gamma(\lambda,p),{\ExtVert}_p))_{p\in P}}
\rho(\Gamma,\lambda)\cdot
\Gamma(\lambda,0)\otimes\\
&\otimes
\otimes_{p\in P}\left(\rho(\Gamma(\lambda,p),\mu_p)
\cdot
(\Gamma(\lambda,p))(\mu_p,0)\otimes 
\otimes_{a\in {\ExtVert}_p}(\Gamma(\lambda,p))(\mu_p,a)
\right).
\end{split}
\end{equation}
The image of $\Gamma$  by the right vertical arrows  in \refN{eq:assocphi} is
\begin{equation}\label{eq:assphiR}
\begin{split}
\sum_{\mu\in\Cond(\Gamma,{\ExtVert})}&\sum_{\omega\in\Cond(\Gamma(\mu,0),P)}
\rho(\Gamma,\mu)\cdot\rho(\Gamma(\mu,0),\omega)\cdot(\Gamma(\mu,0))(\omega,0)\otimes
\\
&\otimes\otimes_{p\in P}(\Gamma(\mu,0))(\omega,p)
\otimes\otimes_{a\in {\ExtVert}}\Gamma(\mu,a).
\end{split}
\end{equation}
Consider the sets
\begin{eqnarray*}
\Theta'\,&:=&\{(\lambda,(\mu_p)_{p\in P}):\lambda\in\Cond(\Gamma,P),\,\mu_p\in\Cond(\Gamma(\lambda,p),{\ExtVert}_p)\textrm{ for }p\in P\},\\
\Theta''&:=&\{(\mu,\omega):\mu\in\Cond(\Gamma,{\ExtVert}),\,\omega\in\Cond(\Gamma(\mu,0),P)\}
\end{eqnarray*}
that index the doubles sums in  Equations \refN{eq:assphiL} and \refN{eq:assphiR}.

For the sake of the proof define a \emph{bicondensation} 
as a map
\[\theta\colon I_\Gamma\longrightarrow\{0\}\cup P\cup {\ExtVert}
\]
and denote by $\Theta$ the set of all bicondensations.
We can associated to $\theta\in\Theta$ the following condensations:
\begin{itemize}
\item A condensation $\lambda\in\Cond(\Gamma,P)$ defined by, for $v\in I_\Gamma$,
\[\lambda(v)=\left\{
\begin{array}{ll}
\theta(v)&\textrm{if }\theta(v)\in P^*\\
\pi(\theta(v))&\textrm{if }\theta(v)\in {\ExtVert},
\end{array}\right.
\]
where 
\[\pi\colon {\ExtVert}\longrightarrow P\]
is defined by $\pi(a)=p$ for $a\in {\ExtVert}_p$;
\item condensations $\mu_p\in\Cond(\Gamma(\lambda,p),{\ExtVert}_p)$ defined by,
for $p\in P$ and $v\in I_\Gamma\cap\lambda^{-1}(p)$,
\[\mu_p(v)=\left\{
\begin{array}{ll}
\theta(v)&\textrm{if }\theta(v)\in {\ExtVert}\\
0&\textrm{otherwise };
\end{array}\right.
\]
\item a condensation $\mu\in\Cond(\Gamma,{\ExtVert})$ defined by,
for $v\in I_\Gamma$,
\[\mu(v)=\left\{
\begin{array}{ll}
\theta(v)&\textrm{if }\theta(v)\in {\ExtVert}\\
0&\textrm{otherwise };
\end{array}\right.
\]
\item a condensation $\omega\in\Cond(\Gamma(\mu,0),P)$ defined by,
for $v\in I_\Gamma\cap\mu^{-1}(\{0\})$,
\[\omega(v)=\theta(v).\]
\end{itemize}

These definitions are summarized in
 the following table:\\
\begin{center}
\begin{tabular}{|ll||c|c|c|}
\hline
$v\in I_\Gamma$&&
$\theta(v)=0$&
$\theta(v)\in P$&
$\theta(v)\in {\ExtVert}$
\\
\hline
\hline
$\lambda(v)\in P^*$&&
$0$&
$\theta(v)$&
$\pi(\theta(v))$
\\
\hline
$\mu_p(v)\in {\ExtVert}^*$,& for $p=\lambda(v)$&
not defined&
$0$&
$\theta(v)$
\\
\hline
$\mu(v)\in {\ExtVert}^*$&&
$0$&
$0$&
$\theta(v)$
\\
\hline
$\omega(v)\in P^*$&&
$0$&
$\theta(v)$&
not defined
\\
\hline
\end{tabular}
\end{center}

It is easy to check that these associations define maps
\begin{eqnarray*}
\Theta\longrightarrow\Theta'&\,,\,&\theta\longmapsto(\lambda,(\mu_p)_{p\in P})\\
\Theta\longrightarrow\Theta''&\,,\,&\theta\longmapsto(\mu,\omega)
\end{eqnarray*}
that are bijections.\private{\pl\todo check}.
Moreover, if $(\lambda,(\mu_p)_{p\in P})$ and $(\mu,\omega)$ are the images of the same bicondensation $\theta$
then  \private{\pl\todo check}
\begin{equation}\label{eq:Gammaassphi}
\left\lbrace
\begin{aligned}
\Gamma(\lambda,0)&=(\Gamma(\mu,0))(\omega,0)\\
(\Gamma(\lambda,p))(\mu_p,0)&=(\Gamma(\mu,0))(\omega,p)\\
(\Gamma(\lambda,p))(\mu_p,a)&=\Gamma(\mu,a)
\end{aligned}
\right.
\end{equation}
for $p\in P$ and $a\in {\ExtVert}_p$.

Let $\eta(\Gamma)=\pm1$ be Koszul sign associated to the reordering map
 $\tau$ on the image of $\Gamma$ through the vertical maps.
It is easy to check \private{\pl\todo} that
\begin{equation}\label{eq:signassphi}
\rho(\Gamma,\lambda)\cdot\prod_{p\in P}\rho(\Gamma(\lambda,p),\mu_p)\,=\,
\rho(\mu,\Gamma)\cdot\rho(\omega,\Gamma(\mu,0))\cdot\eta(\Gamma).
\end{equation}
By \refN{eq:Gammaassphi} and \refN{eq:signassphi} we deduce that
expressions \refN{eq:assphiL} and \refN{eq:assphiR} are equal up to the sign $\eta(\Gamma)$ and
this achieves to prove the commutativity of Diagram \refN{eq:assocphi}.
}%\optproof

%\subsection{Action of the symmetric group on the space of diagrams}\label{sec:coopsym}
%Let ${\ExtVert}$ be a linearly ordered finite set. 
We next define an action of the group $\Perm(A)$ of permutation of the finite set $A$ on  diagrams on $A$.
Given a permutation  $\sigma\in\Perm({\ExtVert})$ and a diagram $\Gamma=({\ExtVert},E,I,s,t)$,
we define a new diagram
\[\sigma\cdot\Gamma=({\ExtVert},E,I,\sigma\circ s,\sigma\circ t)\]
where the bijection $\sigma\colon {\ExtVert}\iso {\ExtVert}$ is extended to
 all vertices by $\sigma(i)=i$ for $i\in I$.
The following is immediate.
\begin{prop}\label{P:permAD}
There is an induced action of ($\BZ$-graded) CDGA of $\Perm({\ExtVert})$ on $\GD({\ExtVert})$ and $\AD({\ExtVert})$.
\end{prop}
\optproof[\refP{permAD}]{
\begin{proof}
It is easy to check that this action is compatible with the equivalence relation $\simeq$ of
\refD{spacediag}. Hence it induces a linear action on $\GD({\ExtVert})$. Since $\NAI({\ExtVert})$ is stable by the
action, we also have an action on $\AD({\ExtVert})$. The action is compatible with the multiplicative structure. To see that it is compatible with the differential, the action commutes
with contraction of edges, and that the sign $\epsilon(\Gamma,e)$ defined before
\refE{d} is unaffected by the action because the highest endpoint of an edge is always internal and
so its position is not changed by the action. 
\end{proof}
}%\optproof

%\subsection{The cooperad structure on $\{\GD(n)\}_{n\geq0}$ and on  $\{\AD(n)\}_{n\geq0}$}
%\label{sec:coopend}
% For an integer $n\geq0$,  set $\GD(n)=\GD(\setn{n})$ and $\AD(n)=\AD(\setn{n})$.

% Let $k$ and $n_1,\dots,n_k$ be non-negative integers and set $n=n_1+\dots+n_k$.
% We have an obvious bijection 
% \[\setn{n}\cong\amalg_{i\in\setn{k}}\,\setn{n_i}=\setn{n_1}\amalg\dots\amalg\setn{n_k}.
% \]
% which correponds to an obvious ordered partition $\nu\colon\setn{n}\to\setn{k}$
% which is an increasing function.
% Therefore we have maps
% \begin{align*}
% \widehat\Psi\colon\GD(n)\to\GD(k)\otimes\GD(n_1)\otimes\dots\otimes\GD(n_k)&\textrm{\quad and}\\
% \Psi\colon\AD(n)\to\AD(k)\otimes\AD(n_1)\otimes\dots\otimes\AD(n_k)&.
% \end{align*}

% We also have by \refS{coopsym} an action of the symmetric group $\Sigma_n$ on $\GD(n)$ and $\AD(n)$.

To define the counits of the cooperad structure, consider the CDGA maps
\[
\hat\eta\colon\GD(1)\longrightarrow\BK\textrm{\quad and \quad}
\eta\colon\AD(1)\longrightarrow\BK
\]
defined by $\hat\eta(\unit)=1$ and $\hat\eta(\Gamma)=0$ for a diagram other than the unit, and similarly for $\eta$.

\begin{thm}\label{T:coopGDAD}
The structure maps $\widehat\Psi_\nu$ and $\Psi_\nu$, for all weak ordered
partitions $\nu$, the symmetric action, and the counits $\hat\eta$ and $\eta$
described above define:
\begin{itemize}
\item the structure of a cooperad of $\BZ$-graded $\BK$-algebras on $\GD\oprd$,\previousoldfn{\pl\todo still to check whether $\GD$ is not a cooperad
of CDGAs} and
\item the structure of a cooperad of CDGAs on $\AD\oprd$ ($\BZ$-graded
  if $N=2$).
\end{itemize}
\end{thm}
\begin{proof}
The associativity of the structure maps $\widehat\Psi$ required for a cooperad structure is exactly
\refL{assocphi}. We have the corresponding associativity for $\Psi$ since, by \refP{Phi}, 
that structure map  is induced by $\widehat\Psi$. It is easy to check that $\hat\eta$ and $\eta$ are counits.
The equivariance is also easy to check. 
\end{proof}

Note that the cooperad structures developed here are related to cooperad strucures on the category of sets (as developed in \cite{Sin:OKS}).

\section{Equivalence of the cooperads $\AD\oprd$ and $\Ho^*(\Conf[\bullet])$}
\label{sec:HoAD}

We show in this section that the CDGA cooperad $\AD\oprd$ of admissible diagrams is weakly equivalent 
to the cohomology algebra of the Fulton-MacPherson{\D} cooperad $\Ho^*(\Conf[\bullet];\BK)$ for any commutative ring with unit
$\BK$ and ambient dimension $N\geq2$.

Fix a finite set $A$.
We first recall the computation of the algebra $\Ho^*(\Conf[A];\BK)$ due to
F. Cohen  \cite{Coh:CBS}.
%\previousfn{\pl Is there a better reference? Earlier I mean.}.
Denote by $[\dvol]\in \Ho^{N-1}(S^{N-1};\BK)$ the  orientation class
of the sphere. For $a,b$ which are distinct in $A$, recall the map
$\theta_{ab}\colon\Conf[A]\to S^{N-1}$ from \refN{eq:defthetaab} which
gives the direction between two points of the configuration, and set 
\begin{equation}
\label{eq:gab}
g_{ab}:=\theta_{ab}^*([\dvol])\in \Ho^{N-1}(\Conf[A];\BK).
\end{equation}
Then as graded algebras we have
\[
\Ho^*(\Conf[A];\BK)=\frac{\wedge\left(\{g_{ab}:\,a,b\in A,\,a\not=b\}\right)}
{\left(\text{$3$-term relation}\,;\,(g_{ab})^2\,;\,g_{ab}-(-1)^Ng_{ba}\right)}
\]
where $\wedge(\{g_{ab}\})$ is the
free  commutative graded $\BK$-algebra generated by the $g_{ab}$'s, and
 the \emph{$3$-term relation}  is 
\[g_{ab}g_{bc}+g_{bc}g_{ca}+g_{ca}g_{ab}\]
for all distinct $a,b,c\in A$.
Here we follow the standard conventions\previousoldfn{\pl was moved from
  \refS{notation} to here} in{\D} rational homotopy theory and denote by $\wedge Z$ the free commutative graded algebra
generated by a graded vector space $Z$. This is thus the tensor product of the symmetric algebra on $Z^{\textrm{even}}$ and
the exterior algebra on $Z^{\textrm{odd}}$.

For $a,b$ distinct in $A$,
denote by
\begin{equation}
  \label{eq:Gab}
\Gamma\langle a,b\rangle  
\end{equation}
the diagram on $A$ with no internal vertices and whose only edge is a chord
from $a$ to $b$. This is an admissible cocycle of degree $N-1$.

We endow the cohomology algebra with a zero differential to make it a CDGA.
\begin{thm}\label{T:formalAD}
For $N\geq2$, there is a quasi-isomorphism of CDGAs ($\BZ$-graded if $N=2$)
\[\bIK\colon\AD(A)\quism\left(\Ho^*(\Conf[A];\BK),0\right)\]
characterized by
\[
\begin{cases}
\bIK(\Gamma\langle a,b\rangle)=g_{ab},&\textrm{for $a,b$ distinct in $A$};\\
\bIK(\Gamma)=0,&\textrm{for a diagram $\Gamma$ with internal vertices.}
\end{cases}
\]
Moreover $\bIK$ is a weak equivalence of cooperads.
\end{thm}

The rest of the section is devoted to the proof of this theorem.

Consider the submodule $\AD^{(0)}(A)$ of $\AD(A)$ generated by admissible diagrams without internal vertices.
Then
\[\AD^{(0)}(A)\,=\,\frac{\wedge\left(\{\Gamma\langle a,b\rangle:\,a,b\in A,\,a\not=b\}\right)}
{\left((\Gamma\langle a,b\rangle)^2\,;\,\Gamma\langle a,b\rangle-(-1)^N\Gamma\langle b,a\rangle\right)}.
\]
Therefore we have a surjective algebra map
\[
\bIK_0\colon\AD^{(0)}(A)\longrightarrow\Ho^*(\Conf[A];\BK)
\]
defined by $\bIK_0(\Gamma\langle a,b\rangle)=g_{ab}$.

\begin{lemma}
\[\bIK_0(\AD^{(0)}(A)\cap d(\AD(A)))=0.\]
\end{lemma}
\begin{proof}[Proof]
\pfn{this proof has been revised; the previous one was flawed}
It is enough to prove that $\bIK_0(d\Gamma)=0$ when $\Gamma$ is an
admissible diagram consisting of 
one internal vertex $i$ and $n$ edges connecting it to
the external vertices $a_1,\dots,a_n$. In that case,
\[
\bIK_0(d\Gamma)=
\sum_{k=1}^n(-1)^kg_{a_1a_k}g_{a_2a_k}\dots g_{a_{k-1}a_k}
 g_{a_{k}a_{k+1}}\dots g_{a_ka_n}.
\]
The right hand side is a generalization of the 3-term relation 
\[g_{a_1a_k}g_{a_2a_k}=g_{a_1a_2}(g_{a_2a_k}-g_{a_1a_k}).
\]
from \eqref{eq:arnold} (which
corresponds to the case $n=3$) and can be proved to vanish by an easy
induction on $n\geq3$ using only this relation.
\pfn{needed to be developped? \iv No, I don't think so. You set it up nicely by writing the 3-term relation this way.}

When $\BK=\BR$,
an alternative non-computational proof is possible: The Kontsevich configuration space integral 
\[\IK\colon\AD(n)\rightarrow\ompa(\Conf[n])\]
(to be defined in \refS{KCSI}) commutes with the differential 
(\refP{GIKchainmap}), and hence
\[\bIK_0(d\Gamma)=[\IK(d\Gamma)]=[d\IK(\Gamma)]=0\]
in $\Ho(\ompa(\Conf[n]))\cong\Ho^*(\Conf[n];\BR)$.
\previousoldfn{\ivfn{\iv This is a nice observation.}}
% Since $\bIK_0$ is a morphism of algebras, it is enough to check
%  the vanishing of the indecomposable coboundaries in $\AD^{(0)}(A)$.
%  Those correspond to the image under $d$ of admissible diagrams $\Gamma$
% consisting of one internal vertex $i$ and edges connecting it to various external vertices.
% Let us say that these external vertices are $1,\dots,n$ and that the edges are $(a,i)$
% for $1\leq a\leq n$. By admissibility, $n\geq3$. Then $\bIK_0(d(\Gamma))=0$ because of the following easy consequence of the
%  $3$-term relation:
% \[g_{12}g_{23}\dots g_{n-1,n}+g_{23}\dots g_{n-1,n}g_{n,1}+\dots+ g_{n,1}g_{12}\dots g_{n-2,n-1}=0,
% \]
% which can be proved for example by working with a standard
% $\BK$-module basis of $\Ho^*(\Conf[A];\BK)$
% consisting of some{\D} monomials in the $g_{ab}$'s.
\end{proof}

This lemma implies that we can define the CDGA morphism $\bIK$ by
\begin{equation}
\bIK(\Gamma)=
\begin{cases}
\bIK_0(\Gamma),&\text{if $\Gamma$ has no internal vertices;}\\
0,&\text{otherwise}.
\end{cases}\label{eq:defbarIK}
\end{equation}

It is straightforward to check that this induces a morphism of cooperads.

Since $\bIK$ induces a surjection in homology, in order to prove that it is a quasi-isomorphism we only need to
establish the following
\begin{lemma}\label{L:HADX}
The graded  $\BK$-modules
$\Ho_*(\AD({\ExtVert}))$ and $\Ho_*(\Conf[{\ExtVert}];\BK)$ are isomorphic.
\end{lemma}
The proof of this lemma will take up the rest of this section.

A diagram $\Gamma$ on ${\ExtVert}$ induces a partition of ${\ExtVert}$
into its path-connected components, and we denote this partition by $\nu_\Gamma$. 
In other words, two external vertices $a$ and $b$ belong to the same
element $C\in\nu_\Gamma$
(see \refD{partition} for definitions{\D}  regarding partitions)
if and only if they are connected by a path of unoriented edges in $\Gamma$.
For a partition $\nu$ of ${\ExtVert}$, denote
by
\[\AD({\ExtVert})\langle\nu\rangle\]
the submodule of $\AD({\ExtVert})$ generated by admissible diagrams $\Gamma$
 whose partition of connected components is $\nu$.
It is clear that $\AD({\ExtVert})\langle\nu\rangle$ is a subcomplex of $\AD({\ExtVert})$.
In the particular case of the indiscrete partition $\nu=\{{\ExtVert}\}$,
we get the subcomplex of \emph{connected admissible diagrams}
\[\cAD({\ExtVert}):=\AD({\ExtVert})\langle\{{\ExtVert}\}\rangle.\]
We have an isomorphism of complexes
\begin{equation}\label{eq:splitD}
\AD({\ExtVert})\cong\underset{\nu}{\oplus}\underset{C\in\nu}{\otimes}\cAD(C)
\end{equation}
where the sum runs over all partitions $\nu$ of the set ${\ExtVert}$.

The Poincar\'e series of the homology of the configuration space $\Conf[{\ExtVert}]$ is
given by \cite{Coh:CBS}
\begin{equation}\label{eq:PoincareC}
(1+t)(1+2t)\dots(1+(|{\ExtVert}|-1)t)
\end{equation} 
with $t$ of degree $N-1$.  In particular the top degree Betti number is
\begin{equation}\label{eq:topBettiC}
\dim\Ho^{(N-1)(|{\ExtVert}|-1)}(\Conf({\ExtVert});\BK)=(|{\ExtVert}|-1)!
\end{equation} 
In view of the isomorphism \refN{eq:splitD} and formulas \refN{eq:PoincareC} and \refN{eq:topBettiC},
\refL{HADX} will be a direct consequence of the following
\begin{lemma}\label{L:HcADX}
For ${\ExtVert}$ non-empty,
\[
\dim\Ho^i(\cAD({\ExtVert}))=\left\{
\begin{array}{ll}
(|{\ExtVert}|-1)!,&\textrm{if }i=(N-1)\cdot(|{\ExtVert}|-1),\\
0,&\textrm{otherwise.}
\end{array}
\right.
\]
\end{lemma}

Before
 proving this lemma, we introduce further submodules.
Fix an element $a\in {\ExtVert}$ and consider the following submodules of $\cAD({\ExtVert})$:
\begin{itemize}
\item $\calU_0$ is the submodule  generated by connected admissible diagrams
with $a$ of valence $1$ and such that the only edge with endpoint $a$ is contractible;
\item $\calU_1$ is the submodule generated by connected admissible diagrams
with $a$ of valence $\geq2$;
\item $\cAD'({\ExtVert})$ is the submodule generated by all connected admissible diagrams that 
are neither in $\calU_0$ nor in{\D} $\calU_1$.
\end{itemize}
It is clear that
$\cAD'({\ExtVert})$ is a subcomplex of $\cAD({\ExtVert})$.
\begin{lemma}\label{L:cAD'cAD}
The inclusion 
\[ \cAD'({\ExtVert})\hookrightarrow\cAD({\ExtVert})\]
is a quasi-isomorphism.
\end{lemma}
\begin{proof}
Consider the quotient complex $\calU:=\cAD({\ExtVert})/\cAD'({\ExtVert})$.
 We need to show that $\calU$ is acyclic.

Identify $\calU$ in the obvious way with the graded $\BK$-module $\calU_0\oplus\calU_1$,
and define an increasing filtration on $\calU$ where
elements of filtration $\leq p$ are the linear combinations of diagrams in $\calU_0$ with less than $p$ edges and diagrams
in $\calU_1$ with less than $p-1$ edges. The differential preserves the filtration.
Consider the spectral sequence associated to this filtration  and which converges to the homology of $\calU$.
 The differential at the $0$th page
\[d^0\colon\calU_0\longrightarrow\calU_1\]
consists of contracting the only edge with endpoint $a$.
It is an isomorphism because there is an  inverse  given by ``blowing up'' the vertex $a$ of a diagram $\Gamma\in\calU_1$
into a contractible edge $(a,a')$ as in \refF{blowupU0U1}.
Therefore the page $E^1$  of the spectral sequence is trivial and hence $\calU$ is acyclic.
\end{proof}

%\begin{figure}\centering
%\fbox{\includegraphics[width=80mm,height=45mm]{BlowupU1U0.pdf}}
%\caption{Example of blowing up vertex $a=1$ into a contractible edge $(a,a')$.}
%\label{fig:blowupU0U1}
%\end{figure}\previousfn{\pl the edge $(x,x')$ should be renammed $(a,a')$ in the picture.}

\begin{figure}[h]
\input{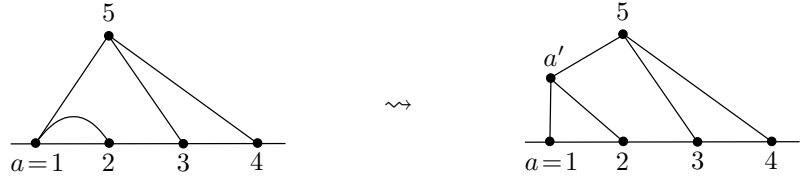}
\caption{Example of blowing up vertex $a=1$ into a contractible edge $(a,a')$.}
\label{fig:blowupU0U1}
\end{figure}

 We are now ready for the
\begin{proof}[Proof of \refL{HcADX}]
The proof is by induction on the cardinality of $\ExtVert$.

If ${\ExtVert}$ is a singleton then $\cAD'({\ExtVert})=\BK\cdot\unit$, where $\unit$ is the unit diagram with a single
external vertex and no internal vertices or edges.  \refL{HcADX} is then a consequence of
\refL{cAD'cAD}.

Let ${\ExtVert}$ be of cardinality $k\geq2$ and suppose that the lemma has been proved for $<k$ external vertices.
Fix $a\in {\ExtVert}$.
Any diagram in $\cAD'({\ExtVert})$ has exactly one edge with endpoint $a$ and it is a chord.
We have an isomorphism of complexes
\[\cAD'({\ExtVert})\cong\underset{b\in {\ExtVert}\setminus\{a\}}{\oplus} \Gamma\langle a,b\rangle\cdot\cAD({\ExtVert}\setminus\{a\}).
\]
Using \refL{cAD'cAD} we conclude that
\[\dim \Ho^i(\cAD({\ExtVert}))=(|{\ExtVert}|-1)\cdot\dim\Ho^{i-(N-1)}(\cAD({\ExtVert}\setminus\{a\}))
\]
and deduce the desired conclusion using the induction hypothesis.
\end{proof}
We now finish the
\begin{proof}[Proof of \refL{HADX}]
An\previousoldfn{\pl this proof has to be improved. \iv I think it's fine like this.}
 elementary computation\previousoldfn{\pl this proof has to be clarified}\previousoldfn{\pl find a reference. \iv What is an elementary computation?  Find reference for what?\pl no reference needed I guess. Is 
the elementary computation clearer now?}
by induction on $|A|$
 using isomorphism \refN{eq:splitD}
%, formulas \refN{eq:PoincareC} and \refN{eq:topBettiC},
 and \refL{HcADX}
shows that the Poincar\'e series of $\Ho^*(\AD(A))$ is exactly 
\refN{eq:PoincareC}, and this is also the Poincar\'e series of $\Ho^*(\Conf[A])$.
\end{proof}

This finishes the proof of \refT{formalAD}.

\section{The Kontsevich configuration space integrals}\label{sec:KCSI}\previousoldfn{\iv We can't call this ``the Kontsevich integral''
 because that's what everyone calls Kontsevich's construction of the universal finite type invariant, so it would be very confusing if
 we starting calling this thing the Kontsevich integral.}

In the previous section we built a quasi-isomorphism
\[\bIK\colon\AD(n)\quism \Ho^*(\Conf[n])
\]
of cooperads.
The goal of this section is to  construct a CDGA morphism
\[\IK\colon\AD(n)\longrightarrow \ompa (\Conf[n])
\]
%called the \emph{Kontsevich integrals},
which will turn out to be a quasi-isomorphism as well as \D
``almost'' a morphism of  cooperads
(see \refP{GIKalmostcoop}{\D} below for the precise meaning of this.)
In this entire section the ground ring is the field of real numbers $\BK=\BR$. We also fix an integer $N\geq2$
which is the dimension of the euclidean space $\BR^N$ on which we consider the configuration spaces $\Conf[n]$,
as well as the underlying dimension of the space of admissible diagrams $\AD=\AD_N$.

We will throughout use many constructions related to semi-algebraic forms that we quickly reviewed in
\refS{RHTSAS} and which are fully developed in \cite{HLTV:RHTSAS}.

The plan of this section \Z is as follows. 

\begin{itemize}
\item[\ref{sec:consGIK}:] We  construct a linear map
\[\GIK\colon\GD(n)\longrightarrow \ompa (\Conf[n]).
\]
\item[\ref{sec:GIKalg}:] We prove that $\GIK$ is a map of algebras.
\item[\ref{sec:vanishing-gik-NAI}:] We show that $\GIK$ induces the
  desired map $\IK$
on $\AD(n)$
by showing that it vanishes on non-admissible diagrams.
\item[\ref{sec:GIKchainmap}:] We prove that $\GIK$, and hence $\IK$,
commutes with the differentials. 
\item[\ref{sec:GIKalmostcoop}:]
We prove that $\GIK$ and $\IK$ are almost morphisms of cooperads.
\end{itemize}

\subsection{Construction of the Kontsevich configuration space integral {$\GIK$}}
\label{sec:consGIK}

Fix a finite set ${\ExtVert}$.
We construct  a linear map
\[\GIK\colon \GD({\ExtVert})\longrightarrow\ompa(\Conf[{\ExtVert}])
\]
as follows.

Let $\Gamma$ be a diagram on ${\ExtVert}$. 
Let $\dvol$ be the standard normalized
 volume form on the sphere $S^{N-1}\subset\BR^N$ defined as
 %\previousfn{\pl check the formula}
\begin{equation}\label{eq:dvol}
\dvol=\kappa_N\cdot
\sum_{i=1}^N(-1)^it_i\,dt_1\wedge\dots\wedge\widehat{dt_i}\wedge\dots\wedge dt_N
\end{equation}
where $t_1,\dots,t_N$ are the standard coordinates in  $\BR^N$, 
$\widehat{dt_i}$ means $dt_i$ is omitted, and $\kappa_N\in\BR$ is a normalizing constant
such that
\[
\int_{S^{N-1}}\dvol=1.\]
Since \p all the functions in
\refN{eq:dvol}
 are polynomials, and are hence semi-algebraic, $\dvol \in\omin^{N-1}(S^{N-1})$ is what was called
in \refS{RHTSAS} a
\emph{minimal form}{\D}\q. 
More generally, for any linearly ordered finite set $E$, 
consider the product of spheres
 \[(S^{N-1})^E=\prod_{e\in E}S^{N-1},\] and 
denote  by $\dvol_E$ the top volume form
in that product, that is, 
\begin{equation}
\dvol_E:=\times_{e\in E}\dvol_e\in\omin((S^{N-1})^E)
\label{eq:volE}
\end{equation}
where the products are taken in the order of $E$ and $\dvol_e$ is the
standard normalized volume
form on the $e$th factor.

For  $v$ and $w$  two distinct vertices in $V_{\Gamma}$,
 recall from \refN{eq:defthetaab} the map
 %\previousfn{\pl remember from where?}
\previousoldfn{\pl this part is to rephrase after
making precise the abuse of notion for  configurations. \iv I'm not sure I understand what you mean by this.\pl nevermind}
\[\theta_{v,w}\colon\Conf[V_\Gamma]\longrightarrow S^{N-1}\]
which associates to a configuration $x$ the direction from $x(v)$ to $x(w)$.
By convention, when $v=w$, we set $\theta_{v,v}$ to be the constant map to a fixed basepoint of the sphere. For an edge $e$ of $\Gamma$ we set
$\theta_e=\theta_{s_\Gamma(e),t_{\Gamma}(e)}$ and we define
\[\theta_\Gamma:=\left(\theta_e\right)_{e\in E_\Gamma}\colon\Conf[V_\Gamma]\longrightarrow(S^{N-1})^{E_\Gamma}.\]
\\
Recall the definition of a minimal form from Equation \eqref{eq:minimalforms}.  We then have such a  form 
\begin{equation}
  \label{eq:thetaGamma}
\theta_\Gamma^*(\dvol_{E_\Gamma})\in\omin(\Conf[V_\Gamma])  
\end{equation}

which is of degree $l=|E_\Gamma|\cdot(N-1)$.

By \refT{projSAbdl},\previousoldfn{\pl missing yet. \pl done 23june}
the canonical projection 
\begin{equation}
\pi_\Gamma\colon\Conf[V_\Gamma]\longrightarrow\Conf[{\ExtVert}]\label{eq:piGamma}
\end{equation}
is an oriented SA bundle\previousoldfn{\pl should fix the coice between the terminologies 'semi-algebraic bundle' vs 'SA bundle'
 in all the paper}. When $|{\ExtVert}|\geq 2$,
the fiber of $\pi_{\Gamma}$ is of dimension $N\cdot|I_\Gamma|$ and
  integration along the fiber\previousoldfn{\pl add a reference to RHTSAS.  \iv
    It's Definition 8.3 there.  I didn't add this since you have some
    special way of adding references to that paper.} 
\citePADintegrfiber
\ gives \Z a
pushforward  map
 \begin{equation}\label{eq:pushpiG}
 (\pi_\Gamma)_*\colon
 \omin^{l}(\Conf[V_\Gamma])
 \longrightarrow
 \ompa^{l-N\cdot|I_\Gamma|}(\Conf[{\ExtVert}]).
\end{equation}

 When $|{\ExtVert}|\geq{2}$,
 define $\GIK(\Gamma)$ as the pushfoward
 \begin{equation}\label{eq:defGIK2}
 \GIK(\Gamma):=(\pi_\Gamma)_*(\theta_\Gamma^*(\dvol_{E_\Gamma}))
 \,\in\,
 \ompa(\Conf[{\ExtVert}]).
 \end{equation}
 For example when $\Gamma$ is the diagram from \refF{killer3terms} in \D
 the Introduction, $\GIK(\Gamma)$ corresponds to formula
\refN{eq:IGamma}.

If ${\ExtVert}$ is empty or a singleton\previousoldfn{\pl The reason to treat the case $|{\ExtVert}|\leq1$ separately
is that in that case the fiber has not the right dimension. \pl I added a remark in the text, 20 june.}
 we just set
\begin{equation}
  \label{eq:defGIK1}
  \GIK(\Gamma):=
  \begin{cases}
    1,&\text{if $\Gamma$ is the unit diagram;}\\
    0,&\text{otherwise.}
  \end{cases}
\end{equation}
The reason we treat the case $|{\ExtVert}|\leq1$ separately is that 
 the dimension of the
fiber of $\pi_\Gamma$ is then smaller than expected when there are internal vertices
(see \refT{projSAbdl}). Therefore we should in those cases
consider the pushforward $\pi_{\Gamma*}$ of  \refN{eq:pushpiG} to be $0$. Formula \refN{eq:defGIK1} is a clean way to do this. 
\previousoldfn{\pl In the definiton of normal condensation there will be a lemma treating that.}

 \begin{lemma}\label{L:defGIK}
For any finite set ${\ExtVert}$, 
 formulas (\ref{eq:defGIK2}) and  (\ref{eq:defGIK1}) induce a degree $0$ linear map
 \[\GIK\colon\GD({\ExtVert})\longrightarrow\ompa(\Conf[{\ExtVert}]).\]
 \end{lemma}
 \begin{proof}
 This is clear for $|{\ExtVert}|\leq1$.
 Suppose that $|{\ExtVert}|\geq2$.
 It is easy\previousoldfn{\pl check \iv looks good} to check that  \refN{eq:defGIK2} is compatible with the 
 equivalence relation $\simeq$ of \refD{spacediag}
 (it is the compatibility with $\GIK$ which is the motivation for the
definition of $\simeq$). We extend it by linearity.
  It  is clear that $\GIK$ is of degree $0$ (recall \refD{degdiag} of
  the degree of a diagram). 
 \end{proof}
 
\subsection{$\GIK$ is a morphism of algebras}\label{sec:GIKalg}
In this section we prove
 \begin{prop}\label{P:GIKalg}
$\GIK$ is a morphism of algebras.
\end{prop}
\begin{proof}
If $|{\ExtVert}|\leq1$ then the proposition is obvious. Suppose now that $|{\ExtVert}|\geq2$.
Let $\Gamma_1$ and $\Gamma_2$ be two diagrams on ${\ExtVert}$ and
suppose, without loss{\D} of generality, that
they have  disjoint sets of internal vertices
and of edges. Notice that $V_{\Gamma_1\cdot\Gamma_2}=V_{\Gamma_1}\cup_A V_{\Gamma_2}$ and
 consider the pullback 
\begin{equation}\label{eq:pbGLIalg}
\xymatrix{
\Confsing[V_{\Gamma_1},V_{\Gamma_2}]\ar[d]_{q_1}\ar[r]^-{q_2}
\ar@{}[rd]|-{\textrm{{pullback}}}
&\Conf[V_{\Gamma_2}]\ar[d]^{\pi_2}\\
\Conf[V_{\Gamma_1}]\ar[r]^{\pi_1}&\Conf[{\ExtVert}]
}
\end{equation}
which defines a singular configuration space as in \refS{singconf}.

Set $\pi'=\pi_i\circ q_i\colon \Confsing[V_{\Gamma_1},V_{\Gamma_2}]\to\Conf[{\ExtVert}]$.
Consider  the  canonical  projections
\[\pi\colon\Conf[V_{\Gamma_1\cdot\Gamma_2}]\longrightarrow\Conf[A]\]
and
\[\rho_i\colon \Conf[V_{\Gamma_1\cdot\Gamma_2}]\longrightarrow \Conf[V_{\Gamma_i}]\]
for $i=1,2$, and the induced map to the pullback
\[\rho\colon\Conf[V_{\Gamma_1\cdot\Gamma_2}]\longrightarrow \Confsing[V_{\Gamma_1},V_{\Gamma_2}].\]
%such that $q_i\circ\rho=\rho_i$.
%  The space $P$ can be interpreted as the space of
% ``singular''  configurations of points labeled by $V_{\Gamma_1\cdot\Gamma_2}$,
% where by singular we mean that some points labeled by $I_{\Gamma_1}$ are allowed to coincide
% with points labeled by $I_{\Gamma_2}$.
% Then $\rho$ maps a  configuration in $\Conf[V_{\Gamma_1\cdot\Gamma_2}]$
% to the corresponding singular configuration in $P$ where the information about the directions
% between two components labeled by $I_{\Gamma_1}$ and $I_{\Gamma_2}$ respectively is lost,
% and similarly for some information about relative distance of three components.
% When $|{\ExtVert}|\geq2$, then $\rho$
% is  a homeomorphism onto its dense image\previousfn{\pl check}
%  when restricted to the complement  of the codimension 1  subspace of $\Conf[V_{\Gamma_1\cdot\Gamma_2}]$
% consiting of  configurations where some component labeled by $I_{\Gamma_1}$ is
% infinitesimally close to a component labeled by $I_{\Gamma_2}$.

% The canonical projection 
% \[\pi'\colon\Conf[V_{\Gamma_1\cdot\Gamma_2}]\to\Conf[{\ExtVert}]\]
% is exactly the map $\pi'=\pi\circ\rho$.
% Since $\pi$ and $\pi'$ are oriented SA bundles over $\Conf[{\ExtVert}]$,
% the fibers $\pi^{-1}(x)$ and  $\pi'^{-1}(x)$ are compact oriented manifolds for $x\in\Conf[{\ExtVert}]$.
By the second part of \refL{deg1sing},
$\pi$ and $\pi'$ are oriented SA bundles and $\rho$ induces a map of degree $\pm1$ between their fibers.
It is easy to check that it is actually of degree $+1$ because it preserves their orientations (which depend, when
$N$ is odd, on the linear order of $I_{\Gamma_1}\ordsum I_{\Gamma_2}$).
Therefore by {\citePAPnatpush}, for any minimal form $\mu\in\omin(P)$,
we have
\begin{equation}
  \label{eq:pp'rho}
\pi'_*(\mu)=\pi_*(\rho^*(\mu)).  
\end{equation}
We  have then
\begin{eqnarray*}
\GIK(\Gamma_1\cdot\Gamma_2)
&=&
\pi_*(\theta^*_{\Gamma_1\cdot\Gamma_2}(\dvol_{E_{\Gamma_1}\ordsum E_{\Gamma_2}}))
\\
&=&
\pi_*(\rho^*(q_1^*\theta^*_{\Gamma_1}(\dvol_{E_{\Gamma_1}})\wedge q_2^*\theta^*_{\Gamma_2}(\dvol_{E_{\Gamma_2}})))
\\
&\stackrel{\text{\refE{pp'rho}}}{=}&
\pi'_*\left(q_1^*\theta^*_{\Gamma_1}(\dvol_{E_{\Gamma_1}})\wedge q_2^*\theta^*_{\Gamma_2}(\dvol_{E_{\Gamma_2}})\right)
\\
&\stackrel{\text{{\citePAPprodpushPB}}}{=}&
\pi_{1*}(\theta^*_{\Gamma_1}(\dvol_{E_{\Gamma_1}}))\wedge\pi_{2*}(\theta^*_{\Gamma_2}(\dvol_{E_{\Gamma_2}}))
\\
&=&
\GIK(\Gamma_1)\cdot\GIK(\Gamma_2).
\end{eqnarray*}
\end{proof}
\subsection{Vanishing of $\GIK$ on non-admissible diagrams}\label{sec:vanishing-gik-NAI}
Recall from \refD{admdiag} the ideal $\NAI(A)$ of non-admissible diagrams.
In this section we prove
\begin{prop}\label{P:GIKNAI}
$\GIK(\NAI({\ExtVert}))=0$.
\end{prop}
\begin{remark}\p\q\pfn{answer comment (86) of referee. \iv Very nice.}The ideal $\NAI({\ExtVert})$ is not the entire kernel of $\GIK$
since
there are admissible diagrams of arbitrarily high degrees but
$\ompa^*(\Conf[\ExtVert])$ is bounded above.
\end{remark}

Since $\AD(A)=\GD(A)/\NAI(A)$, we deduce the following
\begin{cor}\label{C:defIK}
$\GIK$ induces a map of algebras
\[\IK\colon\AD({\ExtVert})\longrightarrow\ompa(\Conf[{\ExtVert}]).\]
\end{cor}
\begin{defin}\label{D:KCSI}
The maps
\[
\GIK\colon\GD({\ExtVert})\longrightarrow\ompa(\Conf[{\ExtVert}])
\,\,\textrm{ and }\,\,
\IK\colon\AD({\ExtVert})\longrightarrow\ompa(\Conf[{\ExtVert}]).
\]
are called the \emph{Kontsevich configuration space integrals}.
\end{defin}
The proof of \refP{GIKNAI} consists of Lemmas \ref{L:GIK0loop}--\ref{L:GIK0bivalent}.
\begin{lemma}
\label{L:GIK0loop}
$\GIK$ vanishes on diagrams with loops.
\end{lemma}
\begin{proof}
If $|{\ExtVert}|\leq1$ the lemma is obvious. Suppose that $|{\ExtVert}|\geq2$ and
let $\Gamma$ be a diagram with a loop. One of the components of the map $\theta_\Gamma$
to the product $(S^{N-1})^{E_\Gamma}$ is a constant map. Therefore $\theta_\Gamma$
factors through a space\previousoldfn{dont we need that $N\geq3$ or alike???
  \pl We assume here $N\geq2$ which is enough} of dimension $<(N-1)\cdot|E_\Gamma|$. By {\citePAPPAdegdimz} we deduce that
the pullback of the maximal degree form $\dvol_{E_\Gamma}$ by $\theta_\Gamma$ is zero,
and hence the same is true for $\GIK(\Gamma)$. \
\end{proof}
\begin{lemma}
\label{L:GIK0double}
$\GIK$ vanishes on diagrams with double edges.
\end{lemma}
\begin{proof}
If $|{\ExtVert}|\leq1$ the lemma is obvious. Suppose that $|{\ExtVert}|\geq2$ and
let $\Gamma$ be a diagram with double edges. The two components of the map $\theta_\Gamma$
corresponding to the double edges factor through the diagonal map
\[\Delta\colon S^{N-1}\longrightarrow S^{N-1}\times S^{N-1}.\]
Therefore $\theta_\Gamma$
factors through a space of dimension $<(N-1)\cdot|E_\Gamma|$. The conclusion is the same as in the proof
of \refL{GIK0loop}.
\end{proof}
\begin{lemma}\label{L:GIK0nonconn}
$\GIK$ vanishes on diagrams containing an internal vertex not connected to any external vertices.
\end{lemma}
\begin{proof}
The lemma is trivial if $|{\ExtVert}|\leq1$.
Assume that $|{\ExtVert}|\geq 2$. Let $\Gamma$ be a diagram as in the statement.
We have a factorization  $\Gamma=\Gamma_1\cdot\Gamma_2$ where $\Gamma_1$ is a diagram with at
 least one internal vertex and such 
that all edges are between internal vertices. Since $\GIK$ is a morphism of algebras,
it is enough to prove that $\GIK(\Gamma_1)=0$. So without loss of generality we assume that $\Gamma=\Gamma_1$.

The canonical projection $\pi_\Gamma$ factors as
\[
\Conf[V_\Gamma]\stackrel{\rho}\longrightarrow\Conf[I_\Gamma]\times\Conf[{\ExtVert}]\stackrel{q}\longrightarrow\Conf[{\ExtVert}]
\]
where $\rho$ is induced by the canonical projections on each factors, and $q$ is the projection on the second factor.
Since we have assumed that the  edges  of $\Gamma$ are only between internal vertices,
there is a factorization $\theta_\Gamma=\theta'\circ\rho$ for some map
\[\theta'\colon \Conf[I_\Gamma]\times\Conf[{\ExtVert}]\longrightarrow(S^{N-1})^{E_\Gamma}.\]

Since $\Gamma$ contains at least one      internal vertex, \refP{dimConf} implies that 
\[\dim(\Conf[I_\Gamma])\leq N\cdot|I_\Gamma|-N.\]
 Therefore for  $x\in\Conf[{\ExtVert}]$ 
we have
\[\dim(q^{-1}(x))<N\cdot |I_\Gamma|=\dim(\pi_{\Gamma}^{-1}(x)))\]
and {\citePAPpushfactcodim} implies that
\[\GIK(\Gamma)=\pi_{\Gamma*}(\theta_\Gamma(\dvol_{E_\Gamma}))=\pi_{\Gamma*}(\rho^*(\theta'^*(\dvol_{E_\Gamma})))=0.\]
\end{proof}
\begin{lemma}
\label{L:GIK0univalent} 
$\GIK$ vanishes on diagrams containing a univalent internal vertex.
\end{lemma}
\begin{proof}
If $|{\ExtVert}|\leq1$, lemma is trivial. Suppose that $|{\ExtVert}|\geq2$.
Let $\Gamma$ be a diagram with an internal vertex 
\previousoldfn{\pl change $v$ to $i$? \iv I don't see the need to do that.\pl done}
 $i$ of valence $1$ and let $v$ be the only vertex adjacent to $i$.
Then $V_\Gamma$ has at least three vertices.
Consider the projection
\[\rho\colon\Conf[V_\Gamma]\longrightarrow\Conf[\{i,v\}]\times\Conf[V_\Gamma\setminus\{i\}]\]
induced by the canonical projections on each factor.
Since $(i,v)$ is the only edge with endpoint $i$ we have a factorization $\theta_\Gamma=\theta'\circ\rho$
for some map 
\[\theta'\colon\Conf[\{i,v\}]\times\Conf[V_\Gamma\setminus\{i\}]\longrightarrow(S^{N-1})^{E_\Gamma}.
\]
Since $i$ is internal we get a map
\[q\colon \Conf[\{i,v\}]\times\Conf[V_\Gamma\setminus\{i\}]\longrightarrow\Conf[{\ExtVert}]\]
obtained as  the projection on the second factor followed by the canonical projection,
and $\pi_\Gamma=q\circ\rho$.
For  $x\in \Conf[{\ExtVert}]$,
\[\dim(q^{-1}(x))<\dim(\pi_\Gamma^{-1}(x)).\]
Then {\citePAPpushfactcodim} implies that
\[\GIK(\Gamma)=\pi_{\Gamma*}(\theta_\Gamma(\dvol_{E_\Gamma}))=\pi_{\Gamma*}(\rho^*\theta'^*(\dvol_{E_\Gamma}))=0.\]
\end{proof}
\begin{lemma}\label{L:GIK0bivalent}
$\GIK$ vanishes on diagrams containing a bivalent internal vertex.
\end{lemma}
\begin{proof}
Lemma is trivial when $|{\ExtVert}|\leq1$.
Assume that $|{\ExtVert}|\geq 2$. We will use Kontsevich's trick from \cite[Lemma 2.1]{Kon:Fey}. 
Let $\Gamma$ be a diagram with an internal vertex $i$ of valence $2$ and
let $v$ and $w$ be its adjacent vertices.  The key idea will be to
consider the automorphism of $\Conf[V_\Gamma]$ which replaces the point
labeled by $i$ by a point symmetric to it with respect to the barycenter of the
points labeled by $v$ and $w$, and to use this symmetry to show that
$\GIK(\Gamma)$ is equal to its{\D} negative.
For concreteness, suppose  that the two
edges at $i$ are oriented as $(v,i)$ and $(w,i)$, and  ordered by $(v,i)<(w,i)$
as the last two edges of the ordered set $E_\Gamma$.

To give the idea of the proof suppose first
that the diagram consists only of 
  these  three vertices and two edges, with $v$ and $w$ external.
Set $\theta=(\theta_{v,i},\theta_{w,i})$, which in this special case is exactly $\theta_\Gamma$,
and set $\pi=\pi_\Gamma$.

Consider the continuous involution
\[\chi\colon \Conf[\{v,w,i\}]\iso\Conf[\{v,w,i\}]
\]
defined on $\Conf(\{v,w,i\})$  by
\[\chi(y)=(y(v)\,,\,y(w)\,,\,y(v)+y(w)-y(i))\]
where $y(v)+y(w)-y(i)$ is  the point symmetric to $y(i)$
with respect to the barycenter $y(v)$ and $y(w)$. This is a
semi-algebraic automorphism of degree $(-1)^N$.

Let \[A\colon S^{N-1}\longrightarrow S^{N-1}\] be the antipodal map and let
\begin{equation}\label{eq:interchange}
\tau\colon S^{N-1}\times S^{N-1}\longrightarrow S^{N-1}\times S^{N-1}
\end{equation}
be the interchange of factors which is of degree $(-1)^{N-1}$. By construction of $\chi$, the following 
diagram commutes
\begin{equation}\label{eq:commtauAA}
\xymatrix{
C[\{v,w,i\}]\ar[d]_\chi\ar[r]^-{\theta}& S^{N-1}\times S^{N-1}\ar[d]^{\tau\circ(A\times A)}\\
C[\{v,w,i\}]\ar[r]^-{\theta}& S^{N-1}\times S^{N-1}
}
\end{equation}
By symmetry of $\dvol$, we have $A^*(\dvol)=\pm\dvol$, so
\begin{equation*}  \label{eq:sign-tauAA}
(\tau\circ(A\times A))^*(\dvol\times\dvol)=  (-1)^{N-1}(\dvol\times\dvol)
\end{equation*}
and hence 
\begin{equation}
  \label{eq:tildechidvol}
\chi^*\theta^*(\dvol_{E_\Gamma})=(-1)^{N-1}\theta^*(\dvol_{E_\Gamma}).
\end{equation}

On the other hand the restriction of $\chi$ to each fiber $\pi^{-1}(x)$,  $x\in\Conf[{\ExtVert}]$,
is an SA homeomorphism of degree $(-1)^N$. By {\citePAPnatpush}, 
\begin{equation}
  \label{eq:pushpitildechi}
  \pi_{*}(\chi^*(\theta^*(\dvol_{E_\Gamma})))=(-1)^N\pi_{*}(\theta^*(\dvol_{E_\Gamma})).
\end{equation}
We deduce that

\begin{eqnarray*}
  \GIK(\Gamma)&=&
\pi_*(\theta^*\dvol_{E_\Gamma})\\
&\stackrel{\text{\refE{pushpitildechi}}}=&
(-1)^N\pi_*(\chi^*\theta^*\dvol_{E_\Gamma})\\
&\stackrel{\text{\refE{tildechidvol}}}=&
(-1)^{N-1}(-1)^N\pi_*(\theta^*\dvol_{E_\Gamma})\\
&=&
-\GIK(\Gamma),
\end{eqnarray*}
and hence $\GIK(\Gamma)=0$.

For the case of a general diagram,
consider the fiber product
\begin{equation}\label{eq:pb3}
\xymatrix{
P\ar[d]\ar[r]
\ar@{}[rd]|-{\textrm{pullback}}
&\Conf[\{v,w,i\}]\ar[d]^{\pi_1}\\
\Conf[V_{\Gamma}\setminus\{i\}]\ar[r]_{\pi_2}&\Conf[\{v,w\}]
}
\end{equation}
where $\pi_1$ and $\pi_2$ are the canonical projections.
Since $\pi_1\circ\chi=\pi_1$,
the automorphism $\chi$ of $\Conf[\{v,w,i\}]$ can be mixed with the identity map on
on  $\Conf[V_{\Gamma}\setminus\{i\}]$ to give an automorphism of $P$ that we also denote by $\chi$.
The canonical projections 
\[\Conf[V_\Gamma]\longrightarrow\Conf[V_\Gamma\setminus\{i\}]\quad\textrm{ and }\quad
\Conf[V_\Gamma]\longrightarrow\Conf[\{v,w,i\}]
\]
induce a map $\rho\colon \Conf[V_\Gamma]\to P$. We have a factorization $\pi_\Gamma=\pi\circ\rho$ for some
map $\pi\colon P\to\Conf[{\ExtVert}]$ which is an oriented SA bundle.

Since the only edges with endpoint $i$ are $(v,i)$ and $(w,i)$,
there is a factorization  $\theta_\Gamma=\theta\circ\rho$ for some map
\[
\theta\colon P\longrightarrow \left(S^{N-1}\right)^{E_\Gamma\setminus\{(v,i),(w,i)\}}\times S^{N-1}\times S^{N-1}.
\]
For each $x\in \Conf[{\ExtVert}]$ the restriction of $\rho$ to the interior
of $\pi_\Gamma^{-1}(x)$ is an  oriented homeomorphism onto a dense image in the fiber $\pi^{-1}(x)$.
By naturality of integration along the fiber {\citePAPnatpush},
\begin{equation}
  \label{eq:pushpirhopi'}
  \pi_{\Gamma*}(\theta_\Gamma^*(\dvol_{E_\Gamma}))=\pi_{*}(\theta^*(\dvol_{E_\Gamma})).
\end{equation}
As for the Diagram \refN{eq:commtauAA}, we have  $\theta\circ\chi=(\id\times\tau\circ(A\times A))\circ\theta$.
The rest of the proof is the same as in the special case treated above, starting with \refE{tildechidvol}.
\end{proof}
\begin{proof}[Proof of \refP{GIKNAI}]
A non admissible diagram satisfies the hypothesis of one of
Lemmas \ref{L:GIK0loop}--\ref{L:GIK0bivalent}.
%\previousfn{\pl keep checking that the numerotation of these lemmas is consecutive}
\end{proof}
%%%% 

\subsection{$\GIK$ and $\IK$ are chain maps}\label{sec:GIKchainmap}

This section is devoted to the proof of the following.
\begin{prop}\label{P:GIKchainmap}  The Kontsevich configuration space integrals commute with the differential, that is, 
$$\GIK \,d=d\,\GIK\ \ \  \text{and}\ \ \   \IK \,d=d\,\IK.$$
\end{prop}

Let ${\ExtVert}$ be a finite set and let $\Gamma$ be a diagram on ${\ExtVert}$.
We will prove that $\GIK(d(\Gamma))=d(\GIK(\Gamma))$,
which by \refC{defIK} implies the analogous result for $\IK$.
If $|{\ExtVert}|\leq1$ then this is obvious.
Also if $\Gamma$ is non-admissible, then by \refP{GIKNAI} and \refL{ideal-NAI}
we have
\[d(\GIK(\Gamma))=0=\GIK(d(\Gamma)).\]
So now we  assume  that $|{\ExtVert}|\geq2$ and that $\Gamma$ is admissible.

From now on we will drop $\Gamma$ from the notation when it appears as an index, so $\Gamma=({\ExtVert},I,E,s,t)$, $V:=V_\Gamma$,
 $\pi:=\pi_{\Gamma}$, etc.
Also, to easily define  orientations of certain configuration spaces
we assume that ${\ExtVert}$ is equipped with an arbitrary linear order and that
$V={\ExtVert}\ordsum I$.

On one side, by definition of $d(\Gamma)$ in \refN{eq:d},
\begin{equation}\label{eq:dIG5}
\GIK(d(\Gamma))=\sum_{e\in\Econtr}\epsilon(\Gamma,e)\cdot\GIK(\Gamma/e).
\end{equation}

To develop the other side $d(\GIK(\Gamma))$, we will need the results from Sections \ref{sec:bdryCV} and \ref{sec:bdryCVpi} on the decomposition of the 
fiberwise boundary of $\Conf[V]$ into faces which are the images of
operadic maps $\Phi_W$ defined in \refN{eq:defPhiW}, mainly Propositions
\ref{P:bdryCV} and \ref{P:bdryCVpi}.
Recall from \refN{eq:defCVbdrypi} the fiberwise boundary of $\pi$,
\[\pi^\partial\colon\Conf^{\partial}[V]\longrightarrow\Conf[{\ExtVert}].
\]
Since $\GIK(\Gamma)=\pi_*(\theta^*(\dvol_E))$ and 
 $\theta^*(\dvol_E)$ is a cocycle, 
 the fiberwise Stokes formula of {\citePAPdpush} 
implies that
\begin{equation}\label{eq:dIG0}
d(\GIK(\Gamma))=(-1)^{\deg(\Gamma)}\cdot
\pi_*^{\partial}\left((\theta^*\dvol_E)|\Conf^{\partial}[V]\right),
\end{equation}
where $(\theta^*\dvol_E)|\Conf^{\partial}[V]$ denotes the restriction of the form $\theta^*\dvol_E$ to that subspace.
Set 
\[\mu:=(\theta^*\dvol_E)|\Conf^{\partial}[V]\in\omin^*(\Conf^{\partial}[V]).\]
Using the decomposition of the fiberwise boundary of $\Conf[V]$
from   \refP{bdryCVpi} and \refP{bdryCV} (ii)-(iii), we get,
by additivity of integration along the fiber\previousoldfn{\pl reference in
  RHTSAS \todo} \citePAPaddpush\previousfn{\pl this prop has still do be added
  in RHTSAS},
\begin{equation}\label{eq:dIG1}
\pi^\partial_*(\mu)
=
\sum_{W\in\BF(V,{\ExtVert})}(\pi^\partial|\im\Phi_W)_*(\mu)
\end{equation}
with the notation from Sections \ref{sec:bdryCV} and
\ref{sec:bdryCVpi}.
Recall in particular that $\BF(V,A)$ is the indexing set of some faces
of the fiberwise boundary and consist of some subsets $W\subset V$.

The core of the proof of \refP{GIKchainmap} consists of computing the terms of the sum in
\refN{eq:dIG1}. We will prove that they all vanish except when $W$ 
is the pair of endpoints of a contractible edge $e$ of $\Gamma$,
and in that case 
\[(\pi^\partial|\im\Phi_W)_*(\mu)=\pm \GIK(\Gamma/e),\]
which are exactly
the terms of $\GIK(d(\Gamma))$ in \refN{eq:dIG5}.\previousoldfn{\pl should we explain more of the geometric idea behind this? How? \iv I don't think we need to.}

Let $W\in\BF(V,{\ExtVert})$, that is:
 $W\subsetneq V$, $|W|\geq2$, and
either ${\ExtVert}\subset W$ or $|W\cap {\ExtVert}|\leq1$ (see
\refN{eq:WVX}). 
Consider the projection to the quotient set
\[q\colon V\longrightarrow V/W.\]
The composite
\begin{equation}
  \label{eq:orderV/W}
(V\setminus W)\cup\{\min(W)\}\hookrightarrow V\stackrel{q}\longrightarrow V/W
\end{equation}
is a bijection
and we use it to transport the linear order of $V$ to $V/W$.

In order to compute $(\pi^\partial|\im\Phi_W)_*(\mu)$ in \refN{eq:dIG2} below, 
we first associate to $\Gamma$ and $W$ two diagrams: $\Gamma'$ which is the full subgraph
of $\Gamma$ with set of vertices $W$, and $\Gbar$ which is the quotient of $\Gamma$ by the
subgraph $\Gamma'$.
More precisely, $\Gamma':=({\ExtVert}',I',E',s',t')$ where
\begin{itemize}
\item ${\ExtVert}':={\ExtVert}\cap W$;
\item $I':=I\cap W$;
\item $E':=E\cap s^{-1}(W)\cap t^{-1}(W)$;
\item $s'=s|E'$ and $t'=t|E'$,
\end{itemize}
and  $\Gbar:=({\ExtVertbar},\Ibar,\Ebar,\sbar,\tbar)$ with
\begin{itemize}
\item ${\ExtVertbar}:=q({\ExtVert})$
\item $\Ibar:=(V/W)\setminus q({\ExtVert})$;
\item $\Ebar:=E\setminus E'$;
\item $\sbar=q\circ(s|\Ebar)$ and $\tbar=q\circ(t|\Ebar)$.
\end{itemize}
Hence $V_{\Gamma'}=W$ and $V_{\overline\Gamma}=V/W$.

Set $\thetabar:=\theta_{\Gbar}$ and $\theta':=\theta_{\Gamma'}$.
Set also the minimal forms
$\mubar=\thetabar^*(\dvol_{\Ebar})$
and $\mu'=\theta'^*(\dvol_{E'})$.
The following diagram is commutative 
\begin{equation}
\label{eq:diagdI=Id}
\xymatrix{
(S^{N-1})^{\Ebar}\times (S^{N-1})^{E'}
\ar[rr]^-{\tau_W\cong}
&
&
(S^{N-1})^E
\\
\Conf[V/W]\times\Conf[W]
\ar[u]^-{\thetabar\times\theta'}
\ar[r]^-{\Phi_W}
\ar[rd]_{\pi^{\partial}\circ\Phi_W}
&
\Conf^{\partial}[V]
\ar[d]^-{\pi^{\partial}}
\ar@{^(->}[r]
&
\Conf[V]
\ar[u]^-{\theta}
\ar[dl]^-{\pi}
\\
&
\Conf[{\ExtVert}]
}
\end{equation}
Here $\tau_W$ is the obvious reordering of factors which is a homeomorphism since $E=\Ebar\amalg E'$.

Since $W\in\BF(V,A)$, there{\D} are two cases:
\begin{enumerate}
\item ${\ExtVert}\subset W$. Then  we have a canonical projection
\[\pi'\colon \Conf[W]\longrightarrow\Conf[{\ExtVert}],\]
and $\pi^\partial\circ\Phi_W=\pi'\circ\proj_2$ where
$\proj_2\colon\Conf[V/W]\times\Conf[W]\to\Conf[W]$ is the projection on
the second factor.
\item $|W\cap {\ExtVert}|\leq 1$. Then the composite
\[{\ExtVert}\hookrightarrow V\stackrel{q}\longrightarrow V/W
\]
is injective, and we have an associated canonical projection
\[\pibar\colon\Conf[V/W]\longrightarrow\Conf[{\ExtVert}].\]
Further, $\pi^\partial\circ\Phi_W=\pibar\circ\proj_1$ where
$\proj_1$ is the projection
on the first factor.
\end{enumerate}
In both cases, $\pi^\partial\circ\Phi_W$ is the composition of two oriented SA bundles, and hence is 
itself an oriented{\D} SA bundle \citePAPcompositebdl.\previousoldfn{\pl reference to RHTSAS.  \iv It's Proposition 8.5.}

The linear orders on $V/W$ and $W$ give $\Conf[V/W]\times\Conf[W]$  a natural orientation,
as well as to the fibers of $\pi^\partial\circ\Phi_W$.
Define the sign
\[\sign(\Phi_W)=\pm1\]
according to whether
\[\Phi_W\colon\Conf[V/W]\times\Conf[W]\longrightarrow\Conf^{\partial}[V],\]
which is a homeomorphism onto its image of codimension $0$,
preserves or reverses orientation. Then $\Phi_W$ induces the same change
of orientation between the fibers over any $x\in\Conf[{\ExtVert}]$.
Define also  $\sign(\tau_W)=\pm1$ by
\[\tau^*_W(\dvol_E)=\sign(\tau_W)\cdot(\dvol_{\Ebar}\times\dvol_{E'}).
\]

The Diagram \refN{eq:diagdI=Id} and {\citePAPnatpush}
imply that
\begin{equation}\label{eq:dIG2}
(\pi^\partial|\im\Phi_W)_*(\mu)
=
\sign(\Phi_W)\cdot\sign(\tau_W)\cdot\left((\pi^\partial\circ\Phi_W)_*(\mubar\times\mu')\right).
\end{equation}
Our computation of $(\pi^\partial|\im\Phi_W)_*(\mu)$ goes through the following lemma,
in which we use the notation $\langle\omega,\bbr{M}\rangle$ to denote the evaluation on 
a compact oriented  semi-algebraic manifold $M$ of a PA form
$\omega\in\ompa(M)$ (see equations \refN{eq:bbrM} and \refN{eq:pairingCOmPA}); in other words
$$\langle\omega,\bbr{M}\rangle=\int_M\omega.$$
\previousoldfn{\iv I guess it's still not clear to me why the switch in notation is necessary.  Maybe a sentence justifying this?}
\begin{lemma}\label{L:pimumu}
\begin{equation}\label{eq:pimumu}
(\pi^\partial\circ\Phi_W)_*(\mubar\times\mu')=
\begin{cases}
\quad\pibar_*(\mubar)\cdot\langle\mu'\,,\,\bbr{\Conf[W]}\rangle,
&
\textrm{if }|W\cap {\ExtVert}|\leq1;
\\
\pm\pi'_*(\mu')\cdot
\langle\mubar\,,\,\bbr{\Conf[V/W]}\rangle,&
\textrm{if }{\ExtVert}\subset W.
\end{cases}
\end{equation}
\end{lemma}
\begin{proof}
If $|W\cap {\ExtVert}|\leq1$ then 
$\pi^\partial\circ\Phi_W=\pibar\circ\proj_1$
and the desired formula is a consequence of
the double pushforward formula of
{\citePAPdoublepush}.

If ${\ExtVert}\subset W$ then 
$\pi^\partial\circ\Phi_W=\pi'\circ\proj_2$
and the desired formula is again a consequence of
the double pushforward formula,
 with an extra sign because of the interchange of factors \previousoldfn{\p\ppfn{``from
   Equation \eqref{eq:interchange}'' has nothing to do with this !?
   \iv I'm not sure what you're talking about here.  This equation
   isn't in the old or the new proof as far as I can tell. \pl I was
   refering to your change in january answering comment (89): you
   refered to what was then equation (101) (now (114)) but this is not
 the interchange we are talking about here}}
in $\Conf[V/W]\times\Conf[V]$ to apply the double pushforward formula.
\end{proof}

Our next task is to show that in the right hand side of
\refN{eq:pimumu}, the expressions
\[\langle\mubar\,,\,\bbr{\Conf[V/W]}\rangle\quad\textrm{and}\quad
\langle\mu'\,,\,\bbr{\Conf[W]}\rangle
\]
vanish, except when $W$ is the pair of endpoints of a contractible edge.
This is the content of Lemmas \ref{L:mubar0}--\ref{L:mu'e}.
To prove them we first establish the following general vanishing lemma.
\begin{lemma}\label{L:W30}
Let $\Gamma_0$ be a  diagram with at least $3$ vertices. Then
\begin{equation}\label{eq:W30}
\langle\theta^*_{\Gamma_0}(\dvol_{E_{\Gamma_0}})\,,\,\bbr{\Conf[V_{\Gamma_0}]}\rangle\,=\,0.
\end{equation}
\end{lemma}
\begin{proof}
In this proof we drop $\Gamma_0$ from the notation when it appears as an index, 
so here $V:=V_{\Gamma_0}$,  $E:=E_{\Gamma_0}$, and  $\theta:=\theta_{\Gamma_0}$.
By hypothesis, $|V|\geq3$.

We can assume that
\begin{equation}\label{eq:W30-degdim}
\deg\theta{}^*(\dvol_{E{}})=\dim\Conf[{V{}}]
\end{equation}
because otherwise the left side of \refN{eq:W30}  vanishes for degree reasons.

If $\Gamma_0$ has an isolated vertex $v$ then $\theta{}$ factors through $\Conf[{V{}}\setminus\{v\}]$.
Since
\[\dim\Conf[{V{}}\setminus\{v\}]<\dim\Conf[{V{}}],
\]
the  left side of \refN{eq:W30} again vanishes for degree reasons.

If $\Gamma_0$ has a univalent vertex and $|{V{}}|\geq 3$ 
then the  left side of \refN{eq:W30} vanishes by the same argument 
as in the main part of the  proof of \refL{GIK0univalent} (where the relevant hypothesis is that
 there are at least three vertices).

If $\Gamma_0$ has a bivalent vertex then the vanishing follows by the same argument as
in \refL{GIK0bivalent}.

Finally, suppose that all the vertices of $\Gamma_0$ are at least trivalent.
%We treat separately the cases $N=1$, $N=2$, and $N\geq3$.

%If $N=1$ then $\theta^*(\dvol_E)$ is of degree $0$ and
%\refN{eq:W30-degdim} cannot hold  when $|V|\geq3$. 

If $N=2$ then $|E|\geq3$ and the statement is exactly that of
\cite[Lemma 6.4]{Kon:DQPM}.\previousoldfn{\pl check\todo}

Suppose that $N\geq3$. Since all the vertices are at least trivalent,
$|E|{}\geq\frac{3}{2}|{V{}}|$.
Therefore
\[\deg(\theta{}^*(\dvol_{E{}}))=(N-1)\cdot|E{}|
\geq\frac{3(N-1)}{2}|{V{}}|
=N\cdot|V{}|+\frac{N-3}{2}\cdot|V{}|
\geq N\cdot|{V{}}|
>
\dim\Conf[{V{}}]
\]
which contradicts \refE{W30-degdim}. 
\end{proof}
% \begin{remark}\label{R:N=2break}
% It is unclear wheter this lemma is true when $N=2$. This is the only place in the proof
%  where we need the hypothesis $N\neq2$.\previousfn{\pl I tried some computation. The first diagram I could think
% of where the lemma can break was the 1-skeletton of the prisma $\Delta^2\times I$} 
% \end{remark}

% From now on we assume that $N\not=2$.
\q\begin{remark}
The above proof is essentially the one given in \cite[Appendix A.3]{CCRL:Vas}.
  However, the context is different in that situation since the configuration space integrals produce differential forms on the spaces of knots rather then on configuration spaces, as is the case here.
\end{remark}

\begin{lemma}\label{L:mubar0}
If ${\ExtVert}\subset W$, then $\langle\mubar\,,\,\bbr{\Conf[V/W]}\rangle=0$.
\end{lemma}
\begin{proof}
If $|V/W|\geq3$ then we apply  \refL{W30}  to $\Gamma_0=\Gbar$.

Otherwise $|V/W|=2$ and $V=W\cup\{v\}$ for some internal vertex $v$ of $\Gamma$.
Since $\Gamma$ is admissible, $v$ is at least trivalent and its adjacent vertices are 
in $W$. Therefore $\Gbar$ has double edges (even triple) and the conclusion is the same  as in the proof of \refL{GIK0double}. 
\end{proof}

\begin{lemma}\label{L:mu'0}
If $|W|\geq3$ or if $W$ is a pair of non-adjacent vertices of $\Gamma$,
then 
\[\langle\mu'\,,\,\bbr{\Conf[W]}\rangle\,=\,0.\]
\end{lemma}
\begin{proof}
If  $|W|\geq3$, apply \refL{W30} to $\Gamma_0=\Gamma'$.

If $W$ is a pair of non adjacent vertices, then $\Gamma'$ has no edges
and hence $\mu'=1\in\omin^0(\Conf[W])$. As $N>1$, 
$\deg(\mu')=0<\dim\Conf[W]$ and the statement follows.
\previousoldfn{The following was removed since we assume now $N>1$:\\ If $N=1$, then our choice of orientation from \refS{orCA} implies that 
the $0$-dimensional manifold 
$\Conf[W]$ consists of two points of opposite orientations. Therefore the
constant form $\mu'$ integrates to $0$ on it.}
\end{proof}

We are finally left\previousoldfn{the end of this section has to be read
  again. post-it ``ici'' in the draft} with the case when
$W$ is a pair of adjacent vertices of $\Gamma$ and $|W\cap {\ExtVert}|\leq1$. Then the edge $e$ 
connecting these two vertices is contractible because  at most one of the endpoints is external and it is not a loop nor a dead end since 
$\Gamma$ is admissible.
Moreover in that case we have  
\[
\Gbar=\Gamma/e\quad\textrm{and}\quad\pibar_*(\mubar)=\GIK(\Gamma/e).\]
(The order of internal vertices in $\Gbar$ is the same as for $\Gamma/e$ 
because the ordering \refN{eq:orderV/W} is
compatible with that of $I_{\Gamma/e}$ from \refD{Gamma/e}.)
Define the sign
\begin{equation}\label{eq:defeta}
\eta(e)=\begin{cases}
+1,&\textrm{if $N$ is even or $s(e)<t(e)$}\\
-1,&\textrm{otherwise.}
\end{cases}
\end{equation}

\begin{lemma}\label{L:mu'e}
If $W$ is a pair of vertices connected by a contractible edge $e$ of $\Gamma$ then
\[\langle\mu'\,,\,\bbr{\Conf[W]}\rangle\,=\,\eta(e).
\]
\end{lemma}
\begin{proof}
$\Gamma'$ consists of a single edge and we have a homeomorphism
\[
\theta'=\theta_{s(e),t(e)}\colon\Conf[\{s(e),t(e)\}]\longrightarrow S^{N-1}
\]
which preserves or reverses orientation according to the sign
$\eta(e)$.
Thus
\[
\langle\mu'\,,\,\bbr{\Conf[W]}\rangle\,=\,\eta(e)\cdot\int_{S^{N-1}}\dvol\,=\,\eta(e).
\]
\end{proof}

Also set  $\Phi_e=\Phi_W$ and $\tau_e=\tau_W$ in that case.

Collecting 
\refN{eq:dIG0}, \refN{eq:dIG1}, \refN{eq:dIG2}, \refL{pimumu}, and
Lemmas \ref{L:mubar0}--\ref{L:mu'e}, 
%\previousfn{\pl check consecutiveness}
we get
\begin{equation}\label{eq:dIG4}
d(\GIK(\Gamma))=\sum_{e\in\Econtr}
(-1)^{\deg(\Gamma)}\cdot
\sign(\Phi_e)\cdot
\sign(\tau_e)\cdot
\eta(e)\cdot
\GIK(\Gamma/e).
\end{equation}

Comparing this to the formula \refN{eq:dIG5} for $\GIK(d(\Gamma))$, 
it remains to compare the signs of the terms in \refN{eq:dIG4} and \refN{eq:dIG5}.
Let $e$ be a contractible edge of $\Gamma$.

\begin{lemma}\label{L:signtaue}
$\sign(\tau_e)=(-1)^{(N-1)\cdot(\pos(e:E)+|E|)}$.
\end{lemma}
\begin{proof}
If $e$ is the last edge in the order of $E$ then $\tau_e$ is the identity map,
and hence $\sign(\tau_e)=+1$ which is the expected value since $\pos(e:E)=|E|$.

When one transposes $e$ with a consecutive  edge in the linear order of $E$ then both
$\sign(\tau_e)$ and $(-1)^{(N-1)\cdot(\pos(e:E)+|E|)}$ change by a factor of $(-1)^{N-1}$.
This proves the  lemma in full generality.
\end{proof}
\begin{lemma}\label{L:signPhie}
$\sign(\Phi_e)=(-1)^{N\cdot(\pos(\max(s(e),t(e)):I)+|I|)}$.
\end{lemma}
\begin{proof}
Suppose first that $t(e)$ is the last and $s(e)$ the second to the 
last vertex in the linear order of $A\ordsum I$.
Then it is easy to see\previousoldfn{\pl check \iv ok} that
\[\Phi_e\colon\Conf[V\setminus\{t(e)\}]\times\Conf[\{s(e),t(e)\}]\longrightarrow\partial\Conf[V]
\]
is orientation-preserving, and hence
\[\sign(\Phi_e)=+1=(-1)^{N\cdot(|I|+|I|)}\]
as expected.

Consider now a permutation $\sigma$
of the set of vertices and its induced action
on the following diagram
\[
\xymatrix{
\Conf[V\setminus\{\max(s(e),t(e))\}]
\times\Conf[\{s(e),t(e)\}]\qquad
\ar[r]^-{\Phi_{\{s(e),t(e)\}}}
\ar[d]_{\sigma\times\sigma}
&
\qquad\partial\Conf[V]
\ar[d]^-{\sigma}
\\
\Conf[V\setminus\{\max(\sigma(s(e)),\sigma(t(e)))\}]
\times
\Conf[\{\sigma(s(e)),\sigma(t(e))\}]\qquad
\ar[r]^-{\Phi_{\sigma(s(e)),\sigma(t(e))\}}}
&
\qquad\partial\Conf[V].
}\]
Inspecting the changes of signs through this diagram, 
it is straighforward to check that the formula is true in  general.
\end{proof}

By Lemmas \ref{L:signtaue}--\ref{L:signPhie}
%\previousfn{\pl check consecutive}
we get that the expressions at \refN{eq:dIG4} and \refN{eq:dIG5} are equal.
This finishes the proof of \refP{GIKchainmap} showing that $\GIK$ and $\IK$
are chain maps.

\subsection{$\GIK$ and $\IK$ are  almost  morphisms of cooperads}\label{sec:GIKalmostcoop}
Ideally,  $\GIK$ and $\IK$ would be morphisms of cooperads.
However, as we explained in Section \ref{sec:CDGAmodel}, this is not true since
$\ompa(\Conf[\bullet])$ is not a cooperad because $\ompa$ is not comonoidal. 
However, these maps are almost morphisms of cooperads
in the following sense.

\begin{prop}\label{P:GIKalmostcoop}
The Kontsevich configuration space integrals $\GIK$ and $\IK$ are compatible
 with the cooperad structures on $\GD$ and $\AD$ as well as with the structure induced on $\ompa(\Conf[\bullet])$ by the
operad structure  on $\Conf[\bullet]$.  Namely, we have 
\begin{enumerate}
\item  Given a weak ordered partition $\nu\colon A\to P$,  set
 $P^*=\{0\}\ordsum P$, $A_p=\nu^{-1}(p)$, and $A_0=P$  as in the
 setting \ref{setting:XP}.  
Recall the (co)operad{\D} structure maps
\begin{align*}
\Phi_\nu\colon & \underset{p\in  P^*}\prod\Conf[A_p] \longrightarrow\Conf[A]\,,\\
\widehat\Psi_\nu\colon & \GD(A) \longrightarrow\underset{p\in  P^*}\otimes\GD(A_p).
\end{align*}
Then  the following diagram is commutative:
\[
\xymatrix{
\GD({\ExtVert})\ar[r]^-{\GIK}\ar[dd]_-{\widehat\Psi_\nu}&
\ompa(\Conf[{\ExtVert}])\ar[d]^-{\Phi^*_\nu=\ompa(\Phi_\nu)}\\
&
\ompa(\prod_{p\in P^*}\Conf[{\ExtVert}_p])\\
\otimes_{p\in P^*}\GD({\ExtVert}_p)\ar[r]_-{\otimes_{p\in P^*}\GIK}&
\otimes_{p\in P^*}\ompa(\Conf[{\ExtVert}_p])\ar[u]^-{\simeq}_-{\times}
}
\]
%Here $\widehat\Psi_\nu$ and $\Phi_\nu$ are the (co)operad structure maps
%associated to the given ordered weak partition of ${\ExtVert}$;
Here the right vertical quasi-isomorphism $\times$ is the standard Kunneth quasi-isomorphism on forms;
\item $\GIK$ is equivariant with respect to the action of the
  permutations of ${\ExtVert}$;
\item $\GIK$ commutes with the counits $\hat\eta\colon\GD(1)\to\BR$ and $\ompa(C[1])\iso\BR$.
\end{enumerate}
 $(1)$--$(3)$ are also true when we replace
$\GIK$ by $\IK$,
 $\GD$ by $\AD$, $\widehat\Psi_\nu$ by $\Psi_\nu$,
and $\hat\eta$ by $\eta$.
\end{prop}
The rest of the section is devoted to the proof of this proposition.
Statements (2) and (3) are easy, as is (1) when $|A|\leq1$,
and that the statements pertaining to $\IK$ follow from those pertaining to $\GIK$.

We now focus on the proof of (1) for a given weak ordered partition $\nu\colon A\to P$
such that $|A|\geq2$. Let $\Gamma$ be a diagram on $A$.
We need to prove that
\begin{equation}\label{eq:coopI1}
(\times_{p\in P^*}\GIK)(\widehat\Psi_\nu(\Gamma))
=
\Phi^*_\nu(\GIK(\Gamma)).
\end{equation}
To understand why this formula holds, remember the discussion of
condensations of configurations starting soon after  \refN{eq:CnuCl}  and ending at \refD{loc}.
Morally, the right hand side of the formula is
the restriction of  the form $\GIK(\Gamma)$ to the part of the
boundary{\D}   of $\Conf[A]$ consisting of $\nu$-condensed configurations
(assuming that $\nu$ is non-degenerate.)
When performing integration along the fiber of $\pi_\Gamma$ over a
$\nu$-condensed configuration $x\in \Conf[A]$, the points of the
configuration $y\in\pi_\Gamma^{-1}(x)\subset\Conf[V]$ labeled by internal vertices can
be differently condensed with respect to the various clusters of points
in $x$, and this corresponds exactly to the different condensations
$\lambda$ relative to $\nu$. Thus the integral
$\Phi^*_\nu(\GIK(\Gamma))$ is obtained by summing over various subdomains 
$\Conf[V,\lambda]\subset\Conf[V]$ indexed by condensations, and the cooperad structure map
$\widehat\Psi_\nu$ that appears on the left side of \refN{eq:coopI1} is
precisely the
 sum over these condensations.
\previousoldfn{\pl add a sligh geometric idea of the proof? cf Kontsevich
  paper\todo. \pl I tried to do it; is it meaningful? \iv definitely.}

We now  proceed with the details. To simplify notation, we will drop
  $\Gamma$ from the notation
when it appears as an index, 
so $I:=I_\Gamma$, $\pi:=\pi_\Gamma$, $E:=E_\Gamma$, $\theta=\theta_\Gamma$, etc.
Also for a given condensation $\lambda$ of $V$ relative to $\nu$ and
for $p\in P^*$ we will replace the index $\Gamma(\lambda,p)$ by $p$, as in 
$V_p:=V_{\Gamma(\lambda,p)}$, $\theta_p:=\theta_{\Gamma(\lambda,p)}$, etc.
%We also set  $\widehat\Psi:=\widehat\Psi_\nu$ and
%$\Phi:=\Phi_\nu$.

The proof of \refE{coopI1} relies on the decomposition of the pullback
of the operad structure map $\Phi_\nu$ along the canonical projection $\pi$
that we have investigated in \refS{proj-oper}.  We will use the notation
and results from that section.
Thus consider the pullback $\Conf[V,\nu]$ of $\Phi_\nu$ along $\pi$ from
Diagram \refN{eq:Gpb}. Recall from \refP{GGl} that we have a
decomposition
\[\Conf[V,\nu]=\cup\Conf[V,\lambda]\]
where $\lambda$ runs over all (essential){\D} condensations $\lambda$ of $\nu$, and that this \D
decomposition is ``almost'' a  partition (\refP{GlGm}).
Fix such a condensation $\lambda$ and
consider the following diagram, where the bottom left triangle
is Diagram \refN{diag:pirholambda}, the right bottom pullback is
Diagram \refN{eq:Gpb},
% $\Phi'_\lambda$ is the map \refN{eq:Phi'l} which
%by \refP{pirholambda} is identified to an operadic map
%and coincide with $\Phi'_\nu\circ\rho_\lambda$,
 and $\tau_\lambda$ is the obvious
interchange of factors:
\begin{equation}\label{eq:diagIcoop}
\xymatrix{
\prod_{p\in P^*}(S^{N-1})^{E_{p }}
\ar[rrr]_-{\cong}^-{\tau_\lambda}&
&
&
(S^{N-1})^{E}
\\
\prod_{p\in P^*}\Conf[V_p]
\ar[u]^-{\times_{p\in P^*}\theta_{p }}
\ar[r]^-{\rho_\lambda}
\ar[rd]_{\pi_\lambda=\times_{p\in P^*}\pi_p}&
\Conf[V,\lambda]
\ar@{}[ld]|{  {\quad\quad \quad\quad  \quad\quad  (\ref{diag:pirholambda})}  }
\ar@{}[rd]|{  {  (\ref{eq:pi'l}) \quad\quad \quad   }  }
\ar@{^(->}[r]
\ar[d]^-{\pi'_\lambda}&
\Conf[V,\nu]
\ar[r]^-{\Phi'_\nu}
\ar[dl]^-{\pi'_\nu}
\ar@{}[rd]|-{\textrm{pullback }(\ref{eq:Gpb})\quad}
&
\Conf[V]
\ar[u]_-{\theta}
\ar[d]^-{\pi}
\\
&
\prod_{p\in P^*}\Conf[{\ExtVert}_p]
\ar[rr]_-{\Phi_\nu}&&
\Conf[{\ExtVert}].
}
\end{equation}
The top rectangle in this diagram is also commutative. Indeed by 
\refP{pirholambda} (iii), $\Phi'_\nu\circ\rho_\lambda=\Phi'_\lambda$ and this can be
identified with an operadic map (see \refN{eq:Phi'l} and \refN{eq:Phihatl}).
From this it follows easily that the rectangle commutes.
It is exactly in the commutativity of that rectangle that the compatibility between 
$\Phi_\nu$ and $\widehat\Psi_\nu$ appears.
Recall also from
%\previousfn{\pl in diagram above add the map $\Phi'_\lambda$ as a curved map} that by 
\refP{pirholambda} that
 $\pi_\lambda$ and $\pi'_\lambda$ are oriented SA bundles  and that
 $\rho_\lambda$ induces 
a map of degree $\sigma(I,\lambda)=\pm1$ between the fibers.

The idea of the proof of Equation \refN{eq:coopI1} is to use this diagram to
relate  the left side
of \refN{eq:coopI1}, which is the sum over all condensations $\lambda$ of
\[
\times_{p\in P^*}\GIK(\Gamma(\lambda,p))
=
\times_{p\in P^*}\pi_{p*}(\theta_p^*\dvol_{E_p}),
\]
to the right side of  \refN{eq:coopI1},
which is the pullback through $\Phi^*_\nu$ of
\[
\GIK(\Gamma)=\pi_*(\theta^*(\dvol_E)).
\]

To make this precise, recall the sign
\[
\sigma(E,\lambda)=\pm1
\]
defined in \refN{eq:defsigmaEl} just before \refL{locsign}.
%and set the sign\previousfn{\pl all these signs should be deined together around \refN{eq:defsigmaGl}}
%\begin{equation}\label{eq:defsigmaEl}
%\sigma(E,\lambda):=(-1)^{(N-1)\cdot|S(E,\lambda)|}.
%\end{equation}
\begin{lemma}\label{L:signtaul}
$\tau^*_\lambda(\dvol_E)=\sigma(E,\lambda)\cdot\left(\times_{p\in P^*}\,\dvol_{E_p}\right).$
\end{lemma}
\begin{proof}
Switching two factors of $S^{N-1}$  is a map of degree $(-1)^{N-1}$.
The factors of $\prod_{p\in P^*}(S^{N-1})^{E_{p }}$
are ordered as
$\ordsum_{p\in P^*}E_{p }$
and the number of transpositions needed to reorder this set as $E$ is
the cardinality of $S(E,\lambda)$.
\end{proof}

Recall from \refD{loc} the notion of an essential condensation,
which is a condensation $\lambda$ such that, for each $p\in P^*$,
\p$I_p=\emptyset$ (that is, $I_\Gamma\cap\lambda^{-1}(p)=\emptyset$)
%$\lambda^{-1}(p)\cap{I_p}=\emptyset$
 when $|A_p|\leq2$, and let
\[\EssCond(\Gamma)=\EssCond(V_\Gamma,\nu)\]
be the set of essential condensations of the diagram $\Gamma$.
\begin{lemma}\label{L:Inorm}
Let $\lambda$ be a condensation of $\Gamma$.
\begin{enumerate}
\item[(i)]
If $\lambda$ is essential, then for each $p\in P^*$
\[\GIK(\Gamma(\lambda,p))=\pi_{p*}(\theta_{p}^*(\dvol_{E_{p }})).\]
\item[(ii)]
If $\lambda$ is not essential, then
\[\left(\times_{p\in P^*}\GIK\right)(\Gamma(\lambda))=0.\]
\end{enumerate}
\end{lemma}
\begin{proof}
Suppose that $\lambda$ is essential. Then for each $p\in P^*$, either $|{\ExtVert}_p|\geq2$, in which case
$\GIK(\Gamma(\lambda,p))$ is given by the pushforward \refN{eq:defGIK2} as expected, or $I_p=\emptyset$ in which
case formulas \refN{eq:defGIK1} and \refN{eq:defGIK2} agree because
$\pi_p$ is the identity map and $\Conf[A_p]=*$.

If $\lambda$ is not essential then for some $p\in P^*$ we have $|{\ExtVert}_p|\leq1$ and $I_p\not=\emptyset$,
in which case $\GIK(\Gamma(\lambda,p))=0$ by \refN{eq:defGIK1}.
\end{proof}
We can now prove the commutativity of the diagram in
\refP{GIKalmostcoop}, part (1), which amounts to showing \refE{coopI1}.
By inspection of Diagram \refN{eq:diagIcoop} we have the following
sequence of equalities

\begin{tabular}{lcl}
  &  \tiny{by definition of $\GIK$}  &  \\
$\Phi_\nu^*(\GIK(\Gamma))$ & = & $\Phi_\nu^*(\pi_{*}(\theta^*(\dvol_{E})))$ \\
  &  \tiny{by pullback formula of the pushforward}  &  \\
  &  \tiny{ \citePAPpushpullback}  &  \\
  &    =   & $\pi'_{\nu*}(\Phi'^*_\nu(\theta^*\dvol_{E}))$ \\
  &  \tiny{by Propositions \ref{P:GGl} and \ref{P:GlGm}}, \\ 
  &  \tiny{and additivity of the pushforward}  &  \\
  &  \tiny{\citePAPaddpush}  &  \\
  &  =   & $\displaystyle{\sum_{\lambda\in\EssCond(\Gamma)}\pi'_{\lambda*}((\Phi'_\nu{}|\Conf[V,\lambda])^*\theta^*(\dvol_{E}))}$  \\
  &  \tiny{by \refP{pirholambda} (ii)}  & \\
  &  \tiny{and naturality of the pushforward}  &  \\
  &  \tiny{\citePAPnatpush}  & \\
  &  = & $\displaystyle{\sum_{{\lambda\in\EssCond(\Gamma)}}\sigma(I,\lambda)\cdot\pi_{\lambda*}(\rho^*_\lambda(\Phi'_\nu|\Conf[V,\lambda])^*\theta^*(\dvol_{E}))}$ \\
 % &  \tiny{by \refP{pirholambda} (iii)} & \\
%  & = & $\displaystyle{\sum_{{\lambda\in\EssCond(\Gamma)}}\sigma(I,\lambda)\cdot\pi_{\lambda*}(\Phi'^*_\lambda\theta^*\dvol_E)}$ \\
  &  \tiny{by commutativity of \refN{eq:diagIcoop}}  & \\
  & = & $\displaystyle{\sum_{{\lambda\in\EssCond(\Gamma)}}\sigma(I,\lambda)\cdot\pi_{\lambda*}\left((\times_{p\in P^*}\theta_{p })^*(\tau^*_\lambda(\dvol_{E})\right)}$ \\
  &  \tiny{by \refL{signtaul}}  & \\
  & = & $\displaystyle{\sum_{{\lambda\in\EssCond(\Gamma)}}\sigma(I,\lambda)\cdot\sigma(E,\lambda)\cdot\pi_{\lambda*}\left((\times_{p\in P^*}\theta_{p })^*(\times_{p\in P^*}\dvol_{E_{p }})\right)}$ \\
  &  \tiny{by definition of $\sigma(\Gamma,\lambda)$ in \refN{eq:defsigmaGl},}  & \\
  &  \tiny{definition of $\pi_\lambda$ in \refN{eq:pil},}  & \\
  &  \tiny{and \refL{Inorm} (i)}  & \\
  & = & $\displaystyle{\sum_{{\lambda\in\EssCond(\Gamma)}}\sigma(\Gamma,\lambda)\cdot\times_{p\in P^*}\left(\GIK(\Gamma(\lambda,p))\right)}$ \\
  &  \tiny{by definition of $\Gamma(\lambda)$ in \refN{eq:defGammal}}  & \\
  &  \tiny{and \refL{Inorm}(ii)}  & \\
  & = & $\displaystyle{\sum_{{\lambda\in\Cond(\Gamma)}}(\times_{p\in P^*}\GIK)(\Gamma(\lambda))}$ \\
  &  \tiny{by definition of $\widehat\Psi_\nu$ in \refN{eq:defPsihat}}  & \\
  & = & $\left(\times_{p\in P^*}\GIK\right)(\widehat\Psi_\nu(\Gamma)).$
\end{tabular}

This\previousoldfn{\pl need a  references for additivity of pushforward here above and RHTSAS?????\todo} finishes the proof 
of \refP{GIKalmostcoop}, showing that $\GIK$ and $\IK$ are almost morphisms of cooperads.

\section{Proofs of the formality theorems}
\label{sec:proofform}

In this section we  prove all the formality theorems given in the
Introduction. Here $\BK$ is the field of real numbers $\BR$.
%Suppose $N\neq2$.

For (non-relative) \p formality, the case of ambient dimension $N=1$ is trivial
because the little intervals operad is weakly{\D} equivalent to the
associative operad which is clearly formal. Assume that $N\geq 2$.
Let us show first that
\[\IK\colon\AD(A)\longrightarrow\ompa(\Conf[A])\]
is a weak equivalence. It is a CDGA map by \refC{defIK} and 
\refP{GIKchainmap}.
 The map induced in cohomology is surjective because,
for $a,b$ distinct in $A$,
the single-chord diagrams $\Gamma\langle a,b\rangle$ defined in \refN{eq:Gab}
are sent to 
$\theta_{ab}^*(\dvol)$ which correpond clearly to the generators $g_{ab}$ 
of the cohomology algebra of the configuration space (see \refS{HoAD}).
Since by \refT{formalAD} $\Ho(\AD(A))\cong\Ho^*(\Conf[A]))$, we deduce that 
$\IK$ is a quasi-isomorphism.

As reviewed in \refS{RHTSAS}, by {\citePAPmonompaapl} $\ompa$ and $\Apl(u(-);\BR)$ are weakly equivalent
symmetric monoidal contravariant functors where
\begin{equation}
u:\CompSemiAlg\longrightarrow\Top\label{eq:u}
\end{equation}
is the forgetful functor which is symmetric strongly monoidal.
In view of \refD{CDGAmodeloperad}, all of this combined with \refT{formalAD} and \refP{GIKalmostcoop}
implies that, for $N\geq3$, $\Ho(\Conf[\bullet];\BR)$ is a CDGA model for the operad $\Conf[\bullet]$,
and hence the same is true for the little $N$-disks operad. This establishes\previousfn{\pl There is a lie here. Actually for $N\leq2$
the CDGA $\AD$ is not connected, even not bounded below as a chain complex. This is not so nice. I try to fix this but could not yet.
For $N=1$ the formality has no content, so its ok (but watch out for relative formality below when $m=1$.
For $N=2$ I still dont know how to fix this. So we are talking of CDGA
which are unbounded below... To be settled!\todo \pl Well I could not
fix it properly. So the case $N=2$ remains exceptionnal.}
\refT{unstableformality}, that is, the formality of the little balls operad over $\BR$ in the
sense of Definitions \ref{D:CDGAmodeloperad} and \ref{D:formaloperad}.

When $N=2$, the above argument does not prove the formality because 
 $\AD_2$ is only a cooperad of $\BZ$-graded CDGAs (see end of Remark
\ref{Rmk:connD}) and is therefore not suitable for modeling rational (or real) homotopy theory.
However, we do have a zig-zag of quasi-isomorphisms of $\BZ$-graded{\D} CDGA
 (almost) cooperads
between $\Ho^*(\Conf_2[\bullet];\BR)$ and $\ompa^*(\Conf_2[\bullet])$.
Moreover, if we replace the zeroth term of the little disks operad (corresponding to
operations in arity $0$) by the empty space and replace $\AD(0)$ by
$0$, 
then it is possible to truncate the cooperad $\AD$ by an acyclic
operadic ideal to make it a connected CDGA and recover
formality. However, we will not pursue this here.\previousoldfn{\pl is this
  paragraph useful? understandable? \iv I think it's pretty useful and understandable.}

We  now deduce the stable formality of the operad, which is the
formality in the category of operads of chain complexes\p. We assume $N\geq2$.
 Recall from {\citePADC} the chain complex of semi-algebraic {\current}s
\[\CSA_*\colon\SemiAlg\longrightarrow\ChZ,\]
which is monoidal.
We define the $\BR$-dual of  a graded real vector space or of a graded $\BZ$-module $V$ as
\begin{equation}
V^\vee:=\Hom(V,\BR),\label{eq:Vdual}
\end{equation}
and denote the dual of a linear map $f\colon V\to W$ by $f^\vee\colon W^\vee\to V^\vee$.
There is a natural  pairing
\begin{align}
\label{eq:pairing}
\langle-,-\rangle\colon\ompa(X)\otimes\CSA(X)& \longrightarrow\BR \\
\omega\otimes\gamma\quad\quad& \longmapsto\langle\omega,\gamma\rangle\notag
\end{align}
and, by {\citePAPevalmostmon},
the evaluation map
\begin{align*}
\ev\colon\CSA_*(X)\otimes\BR & \stackrel{\simeq}{\longrightarrow}\left(\ompa(X)\right)^\vee \\
\gamma & \longmapsto\langle-,\gamma\rangle
\end{align*}
 is a monoidal symmetric weak equivalence when $X$ is a compact semi-algebraic set.

Fix  a weak ordered partition $\nu\colon A\to P$ and set
 $P^*=\{0\}\ordsum P$, $A_p=\nu^{-1}(p)$, and $A_0=P$  as in  the setting \ref{setting:XP}.
Consider the following diagram (in which 
we write $\otimes_{P^*}$ for $\otimes_{p\in P^*}$) 
\[
\xymatrix{
\underset{P^*}{\otimes}\CSA_*(\Conf[A_p])\otimes\BR
\ar[r]^-{\otimes_{P^*}\ev}_-{\simeq}
\ar[dd]_{\times}^-{\simeq}
&
\underset{P^*}{\otimes}(\ompa(\Conf[A_p]))^\vee
\ar[d]^-\simeq
\ar[r]^-{\otimes_{P^*}(\IK^\vee)}_-\simeq
&
\underset{P^*}{\otimes}\left(\AD(A_p)\right)^\vee
\ar[d]^-\simeq
\\
&
\left(\underset{P^*}{\otimes} \ompa(\Conf[A_p])\right)^\vee
\ar[r]^-{(\otimes_{P^*}\IK)^\vee}_-\simeq
&
\left(\underset{P^*}{\otimes}\AD(A_p)\right)^\vee
\ar[dd]^-{(\Psi_\nu)^\vee}
\\
\CSA_*\left(\underset{p\in P^*}{\prod}\Conf[A_p]\right)\otimes\BR
\ar[r]^-{\ev}_-{\simeq}
\ar[d]_-{\CSA_*(\Phi_\nu)}
&
\left(\ompa((\underset{p\in P^*}{\prod}\Conf[A_p])\right)^\vee
\ar[u]^-\simeq_-{(\times)^\vee}
\ar[d]^-{(\ompa(\Phi_\nu))^\vee}
\\
\CSA_*(\Conf[A])\otimes\BR
\ar[r]^-{\ev}_-{\simeq}
&
\left(\ompa((\Conf[A])\right)^\vee
\ar[r]^-{\IK^\vee}_-\simeq
&
(\AD(A))^\vee
}
\]
This diagram is commutative by {\citePAPevalmostmon}
and \refP{GIKalmostcoop}.

Note that $\AD^\vee$, as the dual of the cooperad of $\BZ$-graded differential vector spaces $\AD$,
is an operad.
The above diagram implies that the operad $\CSA_*(\Conf[\bullet])\otimes\BR$ is weakly equivalent to
$\AD^\vee$.
By  \refT{formalAD}, the latter is weakly equivalent to $H_*(\Conf[\bullet])\otimes\BR$.
By {\citePAPmonCHSing}, the symmetric monoidal functors of semi-algebraic chains $\CSA_*$ and 
of singular chains $\Sing_*$ are weakly equivalent.
This proves \refT{stableformality}, the stable formality of the little $N$-disks operad which says that the chains and the homology of the little balls operad are quasi-isomorphic.
% when $N\not=2$.
%For $N=2$, this has been proved by Tamarkin \cite{Tam:for}.

We now arrive to the proof of the relative formality.
Let $1\leq m< N$ be integers. Suppose given a linear isometry
\[\epsilon\colon\BR^m\longrightarrow\BR^N.\]
For an integer $d\geq1$ and a finite set $A$, denote by $\Conf_d[A]$ the Fulton-MacPherson space of
configurations in $\BR^d$. 
Define  the map
\[
\Conf_\epsilon[A]\colon\Conf_m[A]\longrightarrow\Conf_N[A]
\]
which sends a configuration in $\BR^m$ to its image under $\epsilon$ in
$\BR^N$.
Clearly this map induces  a morphism of operads and is equivalent to the morphism induced by $\epsilon$ between
the little balls operads. 

Define the morphism between CDGAs of admissible diagrams in dimensions $N$ and $m$
\[\AD_\epsilon\colon\AD_N(A)\longrightarrow\AD_m(A)\]
by, for a diagram $\Gamma$ in $\AD_N(A)$,
\[
\AD_\epsilon(\Gamma)=
\begin{cases}
\unit,&\text{if $\Gamma$ is the unit diagram;}\\
0,&\text{otherwise}.
\end{cases}
\]
The case $m=1$, however, is special. We set
\[\AD_1(A):=\Ho^*(\Conf_1[A])\cong\BR[\Perm(A)],\]
hence $\AD_1$ is the{\D} associative cooperad. The unit $\unit\in\AD_1(A)$
is the constant cohomology class $1\in\Ho^0(\Conf[A])$, and
$\AD_\epsilon$ is defined in the same way as above.
\begin{lemma}
\label{L:formADeps}
$\AD_\epsilon$ is a morphism of CDGA cooperads and it is weakly equivalent to
\[\Ho^*(\Conf_\epsilon[\bullet])\colon\Ho^*(\Conf_N[\bullet])\longrightarrow\Ho^*(\Conf_m[\bullet]).\]
\end{lemma}
\begin{proof}
It is clear that $\AD_\epsilon$ is a morphism of CDGA cooperads.
Since $m<N$, $\Ho^*(\Conf_\epsilon;\BR)$
is the trivial map, that is, it is zero in positive degrees and maps 
 the unit of the cohomology algebra to the unit. The same is true for $\AD_\epsilon$.
 This, combined with \refT{formalAD} (which is{\D} tautological for $m=1$)
implies the result.
\end{proof}

We now want to prove that $\AD_\epsilon$ is weakly equivalent to
$\ompa(\Conf_\epsilon)$. For this we will use the Kontsevich
configuration space integral that we extend in ambient dimension
$m=1$
as the linear \Z map
\[\IK=\IK_1\colon\AD_1(A)=\Ho^*(\Conf_1[A])\longrightarrow\ompa^*(\Conf_1[A])\]
that sends a (degree $0$) cohomology class to the corresponding locally
constant function on $\Conf_1[A]$. It is clearly a weak equivalence of
``almost'' cooperads as in \refP{GIKalmostcoop}.

We need the following
\begin{lemma}
Assume that $m\geq1$ and $N\geq2m+1$.
Then the following diagram commutes:
\begin{equation}\label{eq:diagrelform}
\xymatrix{
\AD_N(A)
\ar[r]^-{\IK_N}_-{\simeq}
\ar[d]_{\AD_\epsilon}
&
\ompa (\Conf_N[A])
\ar[d]^{\ompa (\Conf_\epsilon)}\\
\AD_m(A)
\ar[r]^-{\IK_m}_-{\simeq}&
\ompa (\Conf_m[A]).
}
\end{equation}
\end{lemma}
\begin{proof}  The case $A=\emptyset$ is clear. Assume that $|A|\geq1$. 
Let $\Gamma$ be an admissible diagram in $\AD_N(A)$.
We have to show that
\begin{equation}\label{eq:commrelfor}
\ompa (\Conf_\epsilon)\left(\IK_N(\Gamma)\right)=
\IK_m\left(\AD_\epsilon(\Gamma)\right),
\end{equation}
where the right hand side is zero except when $\Gamma$ is the unit diagram.
For the unit diagram this is clear, so 
  assume that $\Gamma$ is not a unit and let us show that the left hand side of
\refN{eq:commrelfor} vanishes. 
Denote $\Gamma$'s set of edges  by $E$ and by $I$ its set of internal vertices.
Suppose first that moreover each external vertex of $\Gamma$ is an 
endpoint of some edge. Since internal vertices are at least trivalent,
 this implies
that
\[
|E|\geq\frac{1}{2}(|A|+3\cdot|I|).
\]
Since $N\geq2m+1\geq3$, we deduce that 
\begin{eqnarray*}
\deg(\Gamma)&=&(N-1)\cdot|E|-N\cdot|I|\\
&\geq&\frac{N-1}{2}(|A|+3\cdot|I|)-N\cdot|I|\\
&=&\frac{N-3}{2}\cdot|I|+\frac{N-1}{2}\cdot|A|\\
&\geq&|A|\cdot m\\
&>&\dim(\Conf_m[A]),
\end{eqnarray*}
and hence the left hand side of \refE{commrelfor} vanishes for degree reasons.

Consider now a general admissible non-unit diagram $\Gamma$ on $A$
and let $B\subset A$ be the set of external vertices that are the
endpoints of some edge of $\Gamma$.
We have an obvious associated map
\begin{align*}
\iota \colon\AD_N(B) & \longrightarrow\AD_N(A) \\
\Gamma' & \longmapsto\iota(\Gamma')
\end{align*}
defined by adding to a diagram $\Gamma'$ in $\AD_N(B)$ 
isolated external vertices labeled by $A\setminus B$.
Thus $\Gamma$ is the image under $\iota $ of some 
diagram $\Gamma'\in\AD_N(B)$.
The map $\iota $ can easily be described in terms of cooperadic operations
analogously to the operadic description before \refD{canproj}
of the canonical projection
\[\pi:\Conf[A]\longrightarrow\Conf[B].
\]

The following
diagram commutes
\[\xymatrix{
\AD_N(B)
\ar[r]^-{\IK_N}_-{\simeq}
\ar[d]_{\iota}
&
\ompa (\Conf_N[B])
\ar[d]^{\ompa (\pi)}
\\
\AD_N(A)
\ar[r]^-{\IK_N}_-{\simeq}&
\ompa (\Conf_N[A]).
}
\]
Since each external vertex of $\Gamma'$ is the endpoint of an edge,
we get by the discussion above that
\[
\ompa (\Conf_\epsilon)\IK_N(\Gamma')=\IK_m\AD_\epsilon(\Gamma').
\]
The commutativity of the last diagram and naturality imply then that
\refE{commrelfor} holds
for $\Gamma=\iota(\Gamma')$. This  proves the lemma.
\end{proof}
Finally we prove the last statement of the Introduction:
\begin{proof}[Proof of \refT{relfor}]
The last two lemmas clearly imply 
 formality of the morphism $\Conf_\epsilon$, and hence the same for the corresponding 
map between operads of little balls (when $m\not=2$).

The stable formality of the morphism of operads $\Conf_\epsilon[\bullet]$ is deduced from the unstable  formality above
 exactly as in the absolute case.
\end{proof}

\previousoldfn{\iv References updates:  If a paper has appeared, I don't think we
  should say ``available on arXiv".  That seems strange.  This is in
  references [2] and [3].  References [17] and [22] should have the
  accent on my last name going the other way.  Reference [17] should
  say ``submitted", and should say ``available at arXiv:0806.0476".
  Reference [22] should say ``available {\it at} arXiv..."\pl I
  updated. About arXiv, I suggest to leave it like this (so the
  referee can easierly find the refences). In the accepted version we
  will clean this.}
\newpage

\appendix
\section{Index of notation}
\label{sewc:idxnot}

\p\pfn{I decided that an index of notation would be nice. 
This index should be cross-checked. It could also be better arranged
  from the typographical viewpoint. \iv I made a few changes, but this is a fantastic idea.}
\previousoldfn{\ppfn{I made a few changes and additions in the index}}

For the convenience of the reader, we include a short index of the most
important notation. Each entry is followed by a short description
and a reference to where the notation is defined
\previousoldfn{\ppfn{instead of:
  first occurs}} 
in the paper.

\vskip 6pt
\begin{center}
\textbf{Names and latin letters}
\end{center}
\vskip 6pt

{\small 

\begin{longtable}{p{4.5cm}p{9.5cm}}
$\ExtVert$, $A_\Gamma$  \dotfill &
set of external vertices in  a configuration space or in a diagram; \refS{FMoperad}, \refD{diag} (see also $A_p$, $I$, $V$)\\
$A_p$ \dotfill &
$A_p=\nu^{-1}(p)$; Setting \ref{setting:XP} 
\\
$A_0$ \dotfill &
alternative notation for the codomain $P$ of a weak partition
$\nu\colon A\to P$; Setting \ref{setting:XP}\\
$\AdmCond(\Gamma)$ \dotfill &
set of admissible condensations on a diagram $\Gamma$; \refD{locreg}
\\
$\Apl$ \dotfill &
Sullivan functor of polynomial forms; \refS{CDGAmodel}\\
$\calB_N(n)$, $\calB_N(\bullet)$, $\calB(\bullet)$ \dotfill &
little $N$-disk
operad; \refS{intro}\\
$\bary(x)$ \dotfill &
barycenter of a configuration in $\BR^N$; \refN{eq:bary}
\\
$\BF(V)$ \dotfill&
indexing set of the boundary faces of $\Conf[V]$; \refN{eq:defWV}
\\
$\BF(V,A)$ \dotfill &
indexing set of the boundary faces of the fiberwise boundary
$\Conf^\partial[V]$; \refN{eq:WVX}
\\
$\CSA_*(X)$, $\CSA_k(X)$ \dotfill &
semi-algebraic chains on a semi-algebraic set $X$; \refS{RHTSAS}\\
$\Conf(A)$, $\Conf(n)$ \dotfill 
&
space of normalized configurations, identified with
$\Inj_0^1(A,\BR^N)$; \refN{eq:defC(X)}, \refN{eq:defC(X)xi}
\\
$\Conf[A]$, $\Conf[V]$, $\Conf[\bullet]$ \dotfill &
Fulton-MacPherson compactification of configuration spaces and corresponding operad; \refS{FMoperad}, \refD{FM}
\\
$\Conf[V,\nu]$ \dotfill &
set of $\nu$-condensed configurations; \refN{eq:Gpb}
\\
$\Conf^\partial[V]$ \dotfill &
fiberwise boundary of $\pi\colon\Conf[V]\to\Conf[A]$;
\refN{eq:defCVbdrypi}
\\
$\Confsing(V_1,V_2)$ \dotfill &
singular configuration space; \refN{eq:pbconfsing}
\\
$\Cond(V,\nu)$, $\Cond(V)$ \dotfill & set of
condensations on $V$ relative to a weak partition $\nu\colon A\to P$; 
\refD{loc}
\\
$\Cond(\Gamma,\nu)$, $\Cond(\Gamma)$ \dotfill &
set of condensations on $V_\Gamma$ relative to a weak partition
$\nu\colon A\to P$; \refD{locG}
\\
$d$ \dotfill &
differential of a diagram; \refN{eq:d}
\\
$\AD(n)$, $\AD_N(n)$, $\AD(A)$, $\AD(\bullet)$ \dotfill &
spaces and cooperad of admissible
diagrams; \refD{AD} 
\\
$\GD(n)$, $\GD_N(n)$, $\GD(A)$,  $\GD(\bullet)$ \dotfill &
spaces and cooperad of 
diagrams; \refD{spacediag}
\\
$\deg(\Gamma)$ \dotfill &
degree of a diagram; \refD{degdiag}
\\
$E_\Gamma$ \dotfill &
ordered set of edges of a diagram $\Gamma$; \refD{diag}
\\
$\Econtr_\Gamma$ \dotfill &
set of contractible edges of a diagram; \refD{diag}
\\
$\EssCond(V)$, $\EssCond(V,\nu)$ \dotfill &
set of essential condensations; \refD{loc}
\\
$g_{ab}$ \dotfill &
standard generator of the cohomology algebra of $\Conf[A]$;
\refN{eq:gab}.
\\
$\IK$, $\GIK$ \dotfill &
Kontsevich configuration space integrals; \refS{KCSI}, \refN{eq:defGIK2},
\refN{eq:defGIK1},
\refC{defIK}
\\
$\overline{\IK}$ \dotfill &
quasi-isomorphism between $\AD(A)$ and $\Ho^*(\Conf[A])$;
\refT{formalAD}, \refN{eq:defbarIK}\\
$I$, $I_\Gamma$ \dotfill &
set of vertices on a configuration space or 
ordered sets of internal vertices of a diagram; beginning of
\refS{canproj}, \refD{diag} (see also $A$ and $V$)
\\
$I_p$ \dotfill& set of $p$-local internal vertices (or global if $p=0$),
$I_p=I\cap\lambda^{-1}(p)$; \refN{eq:IpVp}, after \refN{eq:Glp} (see
also $A_p, V_p$) 
\\
$I_0$ \dotfill& set of global internal vertices,
$I_0=\lambda^{-1}(0)$; \refN{eq:IpVp}, after \refN{eq:Glp} (see
also $A_0, V_0$) 
\\
$\Inj(A,\BR^N)$ \dotfill &
space of injections of $A$ into $\BR^N$; \refN{eq:Inj}
\\
$\Inj_0^1(A,\BR^N)$ \dotfill &
space of injections of $A$ into $\BR^N$
with barycenter at the origin and radius $1$, identified with
$\Conf(A)$; 
\refN{eq:defC(X)xi}
\\
$\BK$ \dotfill &
ground unital ring (often $\BK=\BR)$; \refS{notation}
\\
$N$ \dotfill &
fixed positive integer giving the 
ambient dimension of the little disks operad
or the configuration space; \refS{notation}
\\
$\NAI(A)$ \dotfill 
&ideal of non-admissible diagrams on $A$; \refD{admdiag}
\\
$\NAI(\nu)$ \dotfill &
ideal of non-admissible diagrams associated to a weak partition $\nu$;
\refN{eq:Nnu}
\\
$P$ \dotfill &
codomain of a weak (ordered) partition; \refD{partition}, Setting
\ref{setting:XP}
\\
$P^*$ \dotfill &
extended codomain $P^*:=\{0\}\ordsum P$ of a partition; Setting \ref{setting:XP}\\
$\Perm(A)$ \dotfill &
set of permutations of a set $A$; \refS{notation}
\\
$\pos$, $\pos(x:L)$ \dotfill &
position function; \refS{linord}\\
$\radius(x)$ \dotfill &
radius of a configuration; \refN{eq:radius}\\
$s_\Gamma(e)$ \dotfill &
source of an edge; \refD{diag}\\
$S^{N-1}$ \dotfill &
unit sphere in $\BR^N$; \refN{eq:defthetaab}\\
$S(I,\lambda)$, $S(E,\lambda)$ \dotfill &
sets used to define signs $\sigma(I,\lambda)$, $\sigma(E,\lambda)$;
\refN{eq:SIl}, preceeding \refN{eq:defsigmaIl}\\
$t_\Gamma(e)$ \dotfill & target of an edge; \refD{diag}\\
$u$ \dotfill &
forgetful functor $\CompSemiAlg\to \Top$; \refS{RHTSAS}, \refN{eq:u}\\
$V$, $V_\Gamma$ \dotfill &
set of vertices of a configuration or a diagram, $V=A\amalg I$, $V_\Gamma=A_\Gamma\amalg
I_\Gamma$; \refD{diag}
\\
$V_p$ \dotfill& set of $p$-local internal vertices (or global if $p=0$),
$V_p=A_p\cup I_p$; \refN{eq:IpVp}, after \refN{eq:Glp} (see
also $A_p, I_p$) 
\\
$V_0$ \dotfill& set of global internal vertices,
$V_0=\lambda^{-1}(0)\cup P$; \refN{eq:IpVp}, after \refN{eq:Glp} (see
also $A_0, I_0$) 
\\
$\dvol$ \dotfill &
symmetric volume form on the unit sphere; \refN{eq:dvol}
\\
$\dvol_E$ \dotfill &
volume form on a product of a family of spheres indexed by $E$;
\refN{eq:volE}
\end{longtable}

}

%\newpage
\vskip 6pt
\begin{center}
\textbf{Greek letters}
\end{center}
\vskip 6pt

{\small

\begin{longtable}{p{4.5cm}p{9.5cm}}
$\Gamma$\dotfill &
a diagram, or an isomorphism class of diagram, or an equivalence class
of
diagram; \refD{diag}\\
$\Gamma\langle a,b\rangle$ \dotfill &
a diagram consisting of a single chord joining the external vertices $a$
and $b$; \refN{eq:Gab}\\
$\Gamma(\lambda)$, $\Gamma(\lambda,p)$ \dotfill &
used to define the cooperadic structure on diagrams; 
\refN{eq:Glp}, \refN{eq:defGammal}
\\
$\delta_{a,b,c}$ \dotfill &
relative distance between three points of a configuration;
\refN{eq:defdeltaabc}\\
$\epsilon(\Gamma,e)$ \dotfill &
sign associated to the contraction of an edge $e$ in a diagram; preceeding \refN{eq:d}
\\
$\theta_{a,b}$, $\theta_{ab}$ \dotfill &
map $\Conf[A]\to S^{N-1}$ giving the direction between two points of a
configuration; \refN{eq:defthetaab}\\
$\theta_e$ \dotfill&$\theta$-function associated to an edge $e$,
$\theta_e=\theta_{s_{\Gamma}(e),t_{\Gamma}(e)}$;
before \refN{eq:thetaGamma}\\
$\theta_\Gamma$ \dotfill &
product of maps $\theta_{e}$ indexed by the edges $e$ of $\Gamma$;
\refN{eq:thetaGamma}\\
$\kappa$ \dotfill &
Kunneth quasi-isomorphism; \refN{eq:Aplkappa}, \refN{eq:kappaF}\\
$\lambda$ \dotfill &
condensation $\lambda\colon V\to P^*$; \refD{loc}\\
$\lambda_E$ \dotfill &
extension of the condensation $\lambda$ to edges; \refD{locG}, following \refN{eq:Glp}\\
$\widehat\lambda$ \dotfill &
condensation associated to $\lambda$; \refN{eq:lambdahat}, below \refN{eq:Glp}\\
$\nu$ \dotfill &
weak partition $A\to P$; \refD{partition}, Setting \ref{setting:XP} \\
$\pi$ \dotfill & 
canonical projection $\Conf[V]\to\Conf[A]$ for $A\subset V$;
\refN{eq:pi-canonical}, \refD{canproj}
\\
$\pi_*$ \dotfill &
pushforward or integration along the fiber
$\omin^{k+*}(E)\to\ompa^*(B)$;
\refN{eq:pushforward}
\\
$\pi_\Gamma$ \dotfill &
canonical projection $\Conf[V_\Gamma] \to\Conf[A\Gamma]$;
\refN{eq:piGamma}.\\
$(\pi_\Gamma)_*$ \dotfill &
pushforward along the canonical projection $\pi_\Gamma$; \refN{eq:pushpiG}
\\
$\pi^\partial\colon E^\partial\to B$ \dotfill &
fiberwise boundary of an SA bundle $\pi$; \refN{eq:pipartial}
\\
$\sigma(I,\lambda)$, $\sigma(E,\lambda)$, $\sigma(\Gamma,\lambda)$ \dotfill &
signs; \refN{eq:sigmaIl}, \refN{eq:defsigmaIl}, \refN{eq:defsigmaEl},  \refN{eq:defsigmaGl}
\\$\Phi_\nu$ \dotfill &
operadic structure map in $\Conf[\bullet]$
associated to a weak partition $\nu$; \refN{eq:Phinu}\\
$\Phi_W^V=\Phi_W$ \dotfill &
operadic structure map corresponding to a circle operation
$\Conf[V/W]\times\Conf[W]\to\Conf[V]$
for $W\subset V$; \refN{eq:defPhiW}
\\
$\Psi_\nu$ \dotfill &
cooperadic structure map on $\AD(\bullet)$ 
associated to a weak partition $\nu$; \refS{PsiConstr}, \refP{Phi}\\
$\widehat\Psi_\nu$ \dotfill &
cooperadic structure map on $\GD(\bullet)$ 
associated to a weak partition $\nu$; \refS{PsiConstr}, \refN{eq:defPsihat}\\
$\omin$ \dotfill &
functor of minimal forms on semi-algebraic sets;
\refN{eq:minimalforms}, \refN{eq:pushforward}\\
$\ompa$ \dotfill &
functor of PA forms on semi-algebraic sets; \refN{eq:ompa},
\refT{ompaapl}, \refN{eq:pushforward}\\
\end{longtable}

}

%\newpage
\vskip 6pt
\begin{center}
\textbf{Other symbols}
\end{center}
\vskip 6pt

{\small

\begin{longtable}{p{4.5cm}p{9.5cm}}
$x(a)\simeq x(b)\rel x(c)$ \dotfill &
proximity relation in $\Conf[A]$; \refN{eq:defrel}\\
$\Gamma\simeq\pm\Gamma'$ \dotfill &
equivalence relation of diagrams; \refD{spacediag}\\
$\bbr{M}$, $g_*(\bbr{M})$, $\bbr{\pi^{-1}(b)}$ \dotfill &
semi-algebraic chain represented by a compact semi-algebraic manifold
$M$, its image by a semi-algebraic map $g$, or semi-algebraic chain
represented
by a fiber of an oriented SA bundle; \refN{eq:bbrM}, \refN{eq:fundfiber} \\
$\underline{n}$ \dotfill &
set $\{1,\dots,n\}$; \refS{notationnotation}\\
$f|A$ \dotfill&restriction of a function to a subdomain; \refS{notationnotation}\\
$L_1\ordsum L_2$, $\ordsum_{p\in P}L_p$ \dotfill &
ordered sum; \refS{linord}\\
$\langle-,-\rangle$ \dotfill &
evaluation of a form; \refN{eq:pairingCOmPA}, \refN{eq:pairing}\\
$Y^X$ \dotfill &
set of functions from $X$ to $Y$; \refS{notationnotation}
\\
$|\ExtVert|$ \dotfill &
cardinality of a set $\ExtVert$; \refS{notationnotation}
\\
$\overline{E}$ \dotfill&
closure of a subset $E$ in a topological space
\\
$V/W$ \dotfill &
quotient of a set $V$  by a subset $W\subset V$; \refS{bdryCV}\\
$\Gamma/e$ \dotfill &
contraction of an edge in a diagram; \refD{Gamma/e}\\
$\lambda/e$ \dotfill &
condensation induced on a contracted diagram; \refN{eq:lovere}.
\\
$V^\vee$, $f^\vee$ \dotfill&
linear dual of a vector space or of a linear map; \refN{eq:Vdual}\\
$\|x\|$ \dotfill&
Euclidean norm of $x\in\BR^N$\\
$\wedge Z$ \dotfill&
free commutative graded algebra generated by the graded vector space
$Z$;
after \refN{eq:gab}
\end{longtable}

}

%% =====BIBLIOGRAPHY ============

\previousoldfn{\ivfn{\iv The bibliography entry for the paper with Bob Hardt needs to
  be updated. \pl done}}

\bibliographystyle{plain}
\bibliography{/Users/pascal/Library/texmf/tex/latex/bibliographies/PLbiblio}

\def\cprime{$'$}
\begin{thebibliography}{10}

\bibitem{Arn:coh}
V.~I. Arnol{\cprime}d.
\newblock The cohomology ring of the group of dyed braids.
\newblock {\em Mat. Zametki}, 5:227--231, 1969.

\bibitem{ALTV:cof}
Greg Arone, Pascal Lambrechts, Victor Turchin, and Ismar Voli{\'c}.
\newblock Coformality and rational homotopy groups of spaces of long knots.
\newblock {\em Math. Res. Lett.}, 15(1):1--14, 2008.
\newblock Available at arXiv:math.AT/0701350.

\bibitem{ALV:HQE}
Gregory Arone, Pascal Lambrechts, and Ismar Voli{\'c}.
\newblock Calculus of functors, operad formality, and rational homology of
  embedding spaces.
\newblock {\em Acta Math.}, 199(2):153--198, 2007.
\newblock Available at arXiv:math/0607486.

\bibitem{BoVo:hom}
J.~M. Boardman and R.~M. Vogt.
\newblock {\em Homotopy invariant algebraic structures on topological spaces}.
\newblock Springer-Verlag, Berlin, 1973.
\newblock Lecture Notes in Mathematics, Vol. 347.

\bibitem{BoVo:HIA}
J.~M. Boardman and R.~M. Vogt.
\newblock {\em Homotopy invariant algebraic structures on topological spaces}.
\newblock Lecture Notes in Mathematics, Vol. 347. Springer-Verlag, Berlin,
  1973.

\bibitem{BCR:GAR}
J.~Bochnak, M.~Coste, and M.-F. Roy.
\newblock {\em G\'eom\'etrie alg\'ebrique r\'eelle}, volume~12 of {\em
  Ergebnisse der Mathematik und ihrer Grenzgebiete (3) [Results in Mathematics
  and Related Areas (3)]}.
\newblock Springer-Verlag, Berlin, 1987.

\bibitem{BoTa:SLK}
Raoul Bott and Clifford Taubes.
\newblock On the self-linking of knots.
\newblock {\em J. Math. Phys.}, 35(10):5247--5287, 1994.
\newblock Topology and physics.

\bibitem{BoGu:RHT}
A.~K. Bousfield and V.~K. A.~M. Gugenheim.
\newblock On {${\rm PL}$} de {R}ham theory and rational homotopy type.
\newblock {\em Mem. Amer. Math. Soc.}, 8(179):ix+94, 1976.

\bibitem{CCRL:Vas}
Alberto~S. Cattaneo, Paolo Cotta-Ramusino, and Riccardo Longoni.
\newblock Configuration spaces and {V}assiliev classes in any dimension.
\newblock {\em Algebr. Geom. Topol.}, 2:949--1000 (electronic), 2002.

\bibitem{Coh:CBS}
Fred Cohen.
\newblock Cohomology of braid spaces.
\newblock {\em Bull. Amer. Math. Soc.}, 79:763--766, 1973.

\bibitem{FaNe:con}
Edward Fadell and Lee Neuwirth.
\newblock Configuration spaces.
\newblock {\em Math. Scand.}, 10:111--118, 1962.

\bibitem{FHT:RHT}
Yves F{\'e}lix, Stephen Halperin, and Jean-Claude Thomas.
\newblock {\em Rational homotopy theory}, volume 205 of {\em Graduate Texts in
  Mathematics}.
\newblock Springer-Verlag, New York, 2001.

\bibitem{FuMc:com}
William Fulton and Robert MacPherson.
\newblock A compactification of configuration spaces.
\newblock {\em Ann. of Math. (2)}, 139(1):183--225, 1994.

\bibitem{Gai:mod}
Giovanni Gaiffi.
\newblock Models for real subspace arrangements and stratified manifolds.
\newblock {\em Int. Math. Res. Not.}, 12:627--656, 2003.

\bibitem{GeJo:ope}
E.~Getzler and J.~D.~S Jones.
\newblock Operads, homotopy algebra, and iterated integrals for double loop
  spaces.
\newblock Preprint arXiv:hep-th/9403055v1.

\bibitem{GiKa:Kos}
Victor Ginzburg and Mikhail Kapranov.
\newblock Koszul duality for operads.
\newblock {\em Duke Math. J.}, 76(1):203--272, 1994.

\bibitem{GoWe:EI2}
Thomas~G. Goodwillie and Michael Weiss.
\newblock Embeddings from the point of view of immersion theory. {II}.
\newblock {\em Geom. Topol.}, 3:103--118 (electronic), 1999.

\bibitem{GNPR:mod}
F.~Guill{\'e}n~Santos, V.~Navarro, P.~Pascual, and A.~Roig.
\newblock Moduli spaces and formal operads.
\newblock {\em Duke Math. J.}, 129(2):291--335, 2005.

\bibitem{HLTV:RHTSAS}
Robert Hardt, Pascal Lambrechts, Victor Turchin, and Ismar Voli{\'c}.
\newblock Real homotopy theory of semi-algebraic sets.
\newblock {\em Algebr. Geom. Topol.}, 11(5):2477--2545, 2011.
\newblock Available as arXiv:0806.0476v3.

\bibitem{Kon:Fey}
Maxim Kontsevich.
\newblock Feynman diagrams and low-dimensional topology.
\newblock In {\em First European Congress of Mathematics, Vol.\ II (Paris,
  1992)}, volume 120 of {\em Progr. Math.}, pages 97--121. Birkh\"auser, Basel,
  1994.

\bibitem{Kon:OMDQ}
Maxim Kontsevich.
\newblock Operads and motives in deformation quantization.
\newblock {\em Lett. Math. Phys.}, 48(1):35--72, 1999.
\newblock Mosh\'e Flato (1937--1998).

\bibitem{Kon:DQPM}
Maxim Kontsevich.
\newblock Deformation quantization of {P}oisson manifolds.
\newblock {\em Lett. Math. Phys.}, 66(3):157--216, 2003.

\bibitem{KoSo:def}
Maxim Kontsevich and Yan Soibelman.
\newblock Deformations of algebras over operads and the {D}eligne conjecture.
\newblock In {\em Conf\'erence Mosh\'e Flato 1999, Vol. I (Dijon)}, volume~21
  of {\em Math. Phys. Stud.}, pages 255--307. Kluwer Acad. Publ., Dordrecht,
  2000.

\bibitem{LTV:HQLK}
P.~Lambrechts, V.~Turchin, and I.~Voli\'c.
\newblock The rational homology of spaces of long knots in codimension $>2$.
\newblock {\em Geom. Topol.}, 14:2151--2187 (electronic), 2010.
\newblock Available at arXiv:math.AT/0703649.

\bibitem{NeMi:for}
Joseph Neisendorfer and Timothy Miller.
\newblock Formal and coformal spaces.
\newblock {\em Illinois J. Math.}, 22(4):565--580, 1978.

\bibitem{Sal:sum}
Paolo Salvatore.
\newblock Configuration spaces with summable labels.
\newblock In {\em Cohomological methods in homotopy theory ({B}ellaterra,
  1998)}, volume 196 of {\em Progr. Math.}, pages 375--395. Birkh\"auser,
  Basel, 2001.

\bibitem{SeWi:EFL}
Pavol {\v{S}}evera and Thomas Willwacher.
\newblock Equivalence of formalities of the little discs operad.
\newblock {\em Duke Math. J.}, 160(1):175--206, 2011.

\bibitem{Sin:man}
Dev~P. Sinha.
\newblock Manifold-theoretic compactifications of configuration spaces.
\newblock {\em Selecta Math. (N.S.)}, 10(3):391--428, 2004.

\bibitem{Sin:OKS}
Dev~P. Sinha.
\newblock Operads and knot spaces.
\newblock {\em J. Amer. Math. Soc.}, 19(2):461--486 (electronic), 2006.

\bibitem{Sta:wha}
Jim Stasheff.
\newblock What is {$\dots$} an operad?
\newblock {\em Notices Amer. Math. Soc.}, 51(6):630--631, 2004.

\bibitem{Sul:inf}
Dennis Sullivan.
\newblock Infinitesimal computations in topology.
\newblock {\em Inst. Hautes \'Etudes Sci. Publ. Math.}, 47:269--331 (1978),
  1977.

\bibitem{Tam:for}
Dmitry~E. Tamarkin.
\newblock Formality of chain operad of little discs.
\newblock {\em Lett. Math. Phys.}, 66(1-2):65--72, 2003.

\bibitem{Wei:EI1}
Michael Weiss.
\newblock Embeddings from the point of view of immersion theory. {I}.
\newblock {\em Geom. Topol.}, 3:67--101 (electronic), 1999.

\end{thebibliography}

\private{\newpage\theendnotes}

\end{document}